\documentclass[11pt]{article}
\usepackage{amsmath,amssymb,amsthm,mathtools,mathrsfs}
\usepackage{bm}
\usepackage[margin=1.05in]{geometry}
\AtBeginDocument{\setlength{\emergencystretch}{2em}}
\usepackage{graphicx}
\usepackage{xcolor}

\usepackage[hidelinks]{hyperref}
\hypersetup{bookmarksdepth=3}
\hypersetup{pdftitle={Airy Path Integral Representation for the Kadomtsev--Petviashvili Equation},
  pdfauthor={Jeremy Quastel and C. Alexander Rodriguez}}
\usepackage{needspace}
\usepackage{booktabs,longtable,enumitem,array}
\newcommand{\E}{\mathbb E}
\newcommand{\Ai}{\operatorname{Ai}}
\newtheorem{theorem}{Theorem}[section]

\newtheorem{lemma}[theorem]{Lemma}
\newtheorem{proposition}[theorem]{Proposition}
\newtheorem{corollary}[theorem]{Corollary}
\theoremstyle{remark}
\newtheorem{remark}[theorem]{Remark}
\theoremstyle{plain}
\title{Airy Path Integral Representation for the Kadomtsev--Petviashvili Equation}
\author{Jeremy Quastel \and C. Alexander Rodriguez}
\date{November 23, 2026}

\providecommand{\R}{\mathbb R}

\providecommand{\sech}{\operatorname{sech}}
\providecommand{\diag}{\operatorname{diag}}

\providecommand{\diff}{\mathop{}\!\mathrm{d}}

\providecommand{\C}{\mathbb C}

\begin{document}
\maketitle
\begin{center}
\begin{minipage}{0.85\textwidth}
\centering\itshape
quin etiam passim nostris in versibus ipsis\\
multa elementa vides multis communia verbis

\vspace{0.35em}
``All through these very lines of mine you see many elements common to many words.''
\end{minipage}

\vspace{0.35em}
{\small --- Lucretius, \emph{De rerum natura} I.823--824, trans. Rouse}
\end{center}
\vspace{0.6em}

\begin{abstract}
We give a representation for solutions of the
Kadomtsev--Petviashvili II (KP-II) equation as path integrals over
the Airy line ensemble. The initial condition enters directly as a
potential, bypassing the usual direct and inverse scattering
constructions. We establish the representation for broad classes of
signed initial data, prove uniqueness in specified solution classes,
and derive a restart identity expressing solutions as Airy path integrals
of their own earlier states. For arbitrary real $L^2$ initial data,
a centered version of the Airy occupation represents the standard
global KP-II flow.

We show how the representation recovers the classical Marchenko
determinant for Korteweg--de Vries (KdV), and identify a Brownian-path
amplitude with the classical KP-II heat-scattering coefficient. We
obtain explicit formulas for soliton solutions and derive KP from a
Toda equation for nonintersecting Brownian bridges. Finally, we use
the restart identity to characterize one-point distributions of the
KPZ fixed point through their singular initial conditions, and derive
sharp lower-tail asymptotics directly from KP.
\end{abstract}
\clearpage
\tableofcontents
\clearpage

\section{Airy line ensemble and the representation}\label{sec:introduction}
We consider the Kadomtsev--Petviashvili II equation on
$\mathbb R^2$, with evolution time $t$, transverse position
$x$, height coordinate $r$, and real initial field $V$.
Our normalization of the Cauchy problem is
\begin{equation}\label{KP}
 \partial_t\phi-\frac12\phi\,\partial_r\phi
 +\frac1{12}\partial_r^3\phi
 +\frac14\partial_x^2\partial_r^{-1}\phi=0,
 \qquad \phi(0,x,r)=V(x,r).
\end{equation}
Here and throughout, $\partial_r^{-1}f(x,r)=-\int_r^\infty f(x,u)\,du$
whenever the integral exists. The sign of the transverse term is that
of KP-II. For background on KP-II, its inverse-scattering theory, and
large-time behavior, see Klein and Saut \cite{KleinSaut}.

The probabilistic object behind the representation is the occupation measure of the rescaled
Airy line ensemble, a sequence of random continuous curves
\[
 \mathcal A^{t,x}_1(y)>\mathcal A^{t,x}_2(y)>\cdots,\qquad y\in\mathbb R.
\]
For each fixed $t>0$ and pairwise distinct spatial positions $y_1,\ldots,y_m\in\mathbb R$, the heights
$\mathcal A^{t,x}_i(y_j)$ form a determinantal point process: for $n\ge1$,
indices $j_1,\ldots,j_n\in\{1,\ldots,m\}$, and any compactly supported
test function $f$ on $\mathbb R^n$,
\begin{equation}\label{eq:airy-correlations}
\begin{aligned}
&\mathbb E\sum_{\substack{i_1,\ldots,i_n\ge1\\
 i_a\ne i_b\ \text{if }j_a=j_b}}
 f\bigl(\mathcal A_{i_1}^{t,x}(y_{j_1}),\ldots,
         \mathcal A_{i_n}^{t,x}(y_{j_n})\bigr)\\
&\qquad=\int_{\mathbb R^n}f(z_1,\ldots,z_n)
 \det\!\left[K^{t,x}(y_{j_a},z_a;y_{j_b},z_b)
 \right]_{a,b=1}^n\,dz_1\cdots dz_n,
\end{aligned}
\end{equation}
where the correlation kernel is
\begin{equation}\label{eq:extended-airy-S}
 K^{t,x}(y_1,z_1;y_2,z_2)
 =
 \Bigl[S_{t,x-y_1}^{*}
 \bigl(\mathbf1_{(-\infty,0]}-\mathbf1_{\{y_1<y_2\}}I\bigr)
 S_{t,y_2-x}\Bigr](z_1,z_2).
\end{equation}
Here the indicator acts by multiplication on $L^2(\mathbb R)$, and
\begin{equation}\label{eq:S}
 S_{t,x}(u,v)=t^{-1/3}e^{\frac{2x^3}{3t^2}+\frac{x(v-u)}{t}}
 \Ai\!\left(t^{-1/3}(v-u)+t^{-4/3}x^2\right)
\end{equation}
is the kernel of $e^{(t/3)\partial^3+x\partial^2}$ for $t>0$,
with $S^*_{t,x}(u,v)=S_{t,x}(v,u)$. At $t=0$ and $x>0$,
$S_{0,x}=e^{x\partial^2}$ is the variance-$2$ Brownian heat
kernel, also denoted by $p_x$:
\begin{equation}\label{eq:heat-kernel}
 S_{0,x}(u,v)=p_x(v-u),\qquad
 p_x(z)=\frac{1}{\sqrt{4\pi x}}\exp\!\left(-\frac{z^2}{4x}\right).
\end{equation}
Removing the parabola gives the stationary Airy line ensemble:
\begin{equation}\label{eq:markov-airy-marginal}
 \{\mathcal A_k^{t,x}(y):k\ge1,\ y\in\mathbb R\}
 \ =\
 \left\{t^{1/3}\bar{\mathcal A}_k(t^{-2/3}y)-\frac{(y-x)^2}{t}:
 k\ge1,\ y\in\mathbb R\right\},
\end{equation}
where $\bar{\mathcal A}$ is stationary in the space variable $y$.

The Airy line ensemble arises as a universal scaling limit in random
growth models, random matrices and systems of nonintersecting paths
(see, e.g., Section~\ref{sec:dyson-brownian-scaling} and
Figure~\ref{fig:db-bridges}). It was introduced by Pr\"ahofer and Spohn
\cite{PS02} as the limit of the multi-layer polynuclear growth model;
see also the earlier soft-edge dynamical correlations of Forrester,
Nagao and Honner \cite{FNH99}, and Okounkov's related work on the
Airy point process \cite{Okounkov}.
The Airy line ensemble is also characterized by its Brownian Gibbs property
\cite{CH,AH}.

For the initial potential $V$, define the occupation functional
\begin{equation}\label{eq:airy-occupation}
 I_V^{t,x}(r)=\sum_{i\ge1}\int_{\mathbb R}
 V\bigl(y,r-\mathcal A_i^{t,x}(y)\bigr)\,dy.
\end{equation}
The KP representation is
\begin{equation}\label{representation}
 \phi(t,x,r)=-2\partial_r^2\log Z_V(t,x,r),
 \qquad Z_V(t,x,r)=\mathbb E\exp\{-I_V^{t,x}(r)\}.
\end{equation}
The initial data $V$ enters directly as a potential, without first
constructing scattering data. The probabilistic representation gives $Z_V>0$.

The stationary Airy point process has one-point density
\[
 \rho_{\rm Ai}(z)=K_{\rm Ai}(z,z)=\Ai'(z)^2-z\Ai(z)^2,
 \quad K_{\rm Ai}=S_{1,0}^{*}\mathbf1_{(-\infty,0]}S_{1,0},
 \quad \rho_{\rm Ai}(z)\sim\frac{\sqrt{-z}}{\pi}\quad(z\to-\infty).
\]
The occupation \eqref{eq:airy-occupation} of a measurable $V\ge0$ has mean
\begin{equation}\label{eq:occ-density}
 \E I_V^{t,x}(r)=\int_{\R^2}V\varrho_{t,x,r},\qquad\varrho_{t,x,r}(y,u)=t^{-1/3}\rho_{\rm Ai}\Bigl(t^{-1/3}\Bigl(r+\frac{(y-x)^2}t-u\Bigr)\Bigr),
\end{equation}
by \eqref{eq:markov-airy-marginal} and the one-point Airy density; for fixed $y$, $\varrho_{t,x,r}(y,u)\sim\pi^{-1}t^{-1/2}\sqrt{u-r-(y-x)^2/t}$ as ${u\to+\infty}$.
Consequently, the ordinary occupation is absolutely finite almost surely whenever
\begin{equation}\label{eq:occ-finite}
 \int_{\mathbb R^2}|V(y,u)|\varrho_{t,x,r}(y,u)\,dy\,du<\infty.
\end{equation}

\begin{remark}[G\'erard's direct formula for Benjamin--Ono]
A striking direct representation of a different kind is available for the
Benjamin--Ono equation through work of P.~G\'erard \cite{GerardBO}; see also
the recent exposition and applications in \cite{GassotGerardMiller}.  On the
line, the solution is recovered explicitly from the initial data through the
resolvent of an operator built from the Hardy-space Lax operator.  Thus, as in
\eqref{representation}, the prescribed initial data enters the
formula directly rather than through the solution of a separate
inverse-scattering problem.  The mechanisms are nevertheless quite
different: G\'erard's formula exploits the Hardy-space and Toeplitz structure
of Benjamin--Ono, whereas the present KP-II formula is built from the Airy
evolution and its path-integral structure.  It would be interesting to
understand whether these direct representations are manifestations of a
common operator-theoretic principle.
\end{remark}

\begin{remark}[Comparison with KP-I]
The representation is specific to KP-II.  Formally, KP-I is obtained by
changing the sign of the transverse term, or equivalently by the complex
rotation $x\mapsto ix$.  At the level of the Lax problem this replaces the
heat operator underlying the KP-II construction by a nonstationary
Schr\"odinger operator.  The corresponding direct scattering problem still
has an explicit oscillatory (Fresnel) expansion in the initial potential, but
the inverse problem no longer collapses to an Airy path integral.  Instead one
must retain the two boundary values of the Schr\"odinger scattering data and
perform the usual triangular, or nonlocal Riemann--Hilbert, factorization.
Fredholm determinant representations of the resulting KP-I tau function are
classical; see, for example, P\"oppe--Sattinger \cite{PS88} and the discussion
of KP-I scattering in \cite{KleinSaut}.  Thus the point of
\eqref{representation} is not simply that KP admits a Fredholm
determinant.  For KP-II the heat structure allows the direct and inverse steps
to collapse into a single Airy path integral in which the initial field
$V(y,r)$ appears unmodified.  An analogous collapse does not appear to occur for KP-I.
\end{remark}

\subsection{Main results}\label{sec:main-results}

\paragraph{Compactly supported data.}
We begin with compactly supported initial data, for which our representation takes its simplest form. In this case we express $Z_V$ as a Fredholm determinant whose kernel combines the Feynman--Kac kernel of the Lax operator $\partial_y-\partial_r^2+V$ of KP-II with the Airy operators \eqref{eq:S}.

\begin{theorem}\label{thm:intro-compact}
Let $V\in C_c^\infty(\R^2)$ be real, with support in $(a,b)\times\R$. Then, for $t>0$ and $x,r\in\R$,
\begin{equation}\label{eq:intro-fredholm}
 Z_V(t,x,r)=\det\bigl(I-\chi_rK_{t,x;V}\chi_r\bigr)_{L^2(\R)},\qquad K_{t,x;V}=S^*_{t,a-x}D_{a,b;V}S_{t,x-b},
\end{equation}
where $\chi_r=\mathbf 1_{[r,\infty)}$, $S_{t,x}$ is the Airy operator \eqref{eq:S}, and $D_{a,b;V}=p_{b-a}-P_{a,b;V}$ is the difference between the heat kernel \eqref{eq:heat-kernel} and the Feynman--Kac kernel of $V$,
\[
 \int_\R P_{a,b;V}(w,w')f(w')\diff w'=\E_{B(a)=w}\Bigl[e^{-\int_a^bV(y,B(y))\diff y}f(B(b))\Bigr],
\]
for a Brownian motion $B$ with diffusion coefficient $2$. Moreover, $\phi_V=-2\partial_r^2\log Z_V$ is a real-analytic solution of the equation in \eqref{KP} on $(0,\infty)\times\R^2$, and $\phi_V(t,\cdot,\cdot)\to V$ in $C^\infty_{\rm loc}(\R^2)$ as $t\downarrow0$.
\end{theorem}

In Section~\ref{sec:direct} we prove Theorem~\ref{thm:intro-compact} directly, by expanding the exponential in $Z_V$, using the correlation functions \eqref{eq:airy-correlations} of the Airy line ensemble, and recognizing the Fredholm series of \eqref{eq:intro-fredholm}.
We note that for nonnegative $V$ it was already known that expectations of multiplicative functionals of nonintersecting paths, and in particular of the Airy line ensemble, are determinantal \cite[Section~4.2]{BCR}.

To see why $\phi_V$ solves KP, we translate by $r$ and write \eqref{eq:intro-fredholm} as a determinant on $L^2(\R_+)$,
\begin{equation}\label{eq:intro-flows}
 Z_V(t,x,r)=\det(I-\mathcal K)_{L^2(\R_+)},\qquad\mathcal K=e^{-\frac t3\partial^3+(a-x)\partial^2+r\partial}\,D_{a,b;V}\,e^{\frac t3\partial^3+(x-b)\partial^2-r\partial}.
\end{equation}
In this form $\mathcal K$ depends on $t$, $x$ and $r$ only through the linear flows generated by $\partial^3$, $\partial^2$ and $\partial$, which imply that $\phi_V$ solves KP; see \cite{Poppe89,QR}.
The initial condition in Theorem~\ref{thm:intro-compact} follows from a law of large numbers for the occupation. For small $t$ the sum \eqref{eq:airy-occupation} has $O(t^{-1/2})$ nonzero terms, each of size $O(t^{1/2})$. Indeed, $\bar{\mathcal A}_i$ is of order $-i^{2/3}$, so by \eqref{eq:markov-airy-marginal} the point $(y,r-\mathcal A_i^{t,x}(y))$ lies above the support of $V$ unless $i=O(t^{-1/2})$ and $|y-x|=O(t^{1/2})$. As $t\downarrow0$ the occupation measure of the curves $r-\mathcal A_i^{t,x}$ converges to the deterministic measure $\frac12(u-r)\diff u$ on the vertical half-line above the point $(x,r)$ (Proposition~\ref{prop:lln} and Figure~\ref{fig:lln}), so that
\[
 I_V^{t,x}(r)\longrightarrow\frac12\int_r^\infty(u-r)V(x,u)\diff u=\frac12\partial_r^{-2}V(x,r).
\]
Consequently,
\[
 \phi_V=-2\partial_r^2\log\E e^{-I_V}\longrightarrow-2\partial_r^2\log e^{-\frac12\partial_r^{-2}V}=V.
\]

\paragraph{Existence.}
To extend Theorem~\ref{thm:intro-compact} from smooth compactly supported data to data decaying sufficiently quickly as $r\to\infty$, we approximate $V$ by functions $V_n\in C_c^\infty(\R^2)$ such that $Z_{V_n}\to Z_V$. To conclude that $\phi_V$ solves KP, one would need the derivatives of $\log Z_{V_n}$ to converge as well, and this does not follow directly from $Z_{V_n}\to Z_V$. We overcome this difficulty in three ways, one for each case of the following theorem.

\begin{theorem}\label{thm:intro-existence}
Let $V\colon\R^2\to\R$ be measurable, and suppose that one of the following conditions holds, with constants $M,C,B<\infty$ and $\delta>0$.
\begin{enumerate}
\item[\textup{(i)}] $V\in L^1(\R^2)$ has compact support, and $V\ge-M$.
\item[\textup{(ii)}] $|V(x,r)|\le Me^{-\delta r_+}$, where $r_+=\max\{r,0\}$.
\item[\textup{(iii)}] $V\in C^5(\R^2)$, $V\ge-M$, and $|\partial_x^j\partial_r^kV(x,r)|\le Ce^{B|x|-\delta r}$ for $j+k\le5$.
\end{enumerate}
Then, on $(0,\infty)\times\R^2$, $0<Z_V<\infty$ and $\phi_V=-2\partial_r^2\log Z_V$ is a classical solution of the equation in \eqref{KP}, and $\phi_V(t,\cdot,\cdot)\to V$ in $\mathcal D'(\R^2)$ as $t\downarrow0$. The solution $\phi_V$ is real analytic in case \textup{(i)} and $C^\infty$ in case \textup{(ii)}.
\end{theorem}

In case~(i) the Fredholm determinants \eqref{eq:intro-fredholm} for $Z_{V_n}$ converge to a Fredholm determinant of the same form for $Z_V$, so \cite[Theorem~1.3]{QR} applies to $Z_V$ directly. In case~(ii) the determinant \eqref{eq:intro-fredholm} is no longer available, and we factor $Z_V$ instead into an exponential prefactor and a regularized determinant on $L^2(\R^2)$. In both cases the derivatives in $t$, $x$ and $r$ fall on the Airy operators \eqref{eq:S}, so $V$ need only be measurable. In case~(iii) we do not use a determinant, and we differentiate under the expectation $Z_V=\E e^{-I_V}$ instead, as in $\partial_rZ_V=-\E[e^{-I_V}I_{\partial_rV}]$. The derivatives now fall on $V$, and case~(iii) therefore assumes bounds on its derivatives up to order five, since the equation involves derivatives of $\log Z_V$ up to that order. Section~\ref{sec:existence-properties} states the three cases separately, in somewhat greater generality, as Theorems~\ref{thm:measurable-data}, \ref{thm:bounded-data} and~\ref{thm:growing-data}.

The well-posedness theory of the Cauchy problem for KP-II has been developed in spaces of Sobolev type. The problem is globally well posed in $L^2(\R^2)$ \cite{Bourgain93,Hadac2008} and, for small data, in the scale-invariant anisotropic Sobolev space \cite{HadacHerrKoch}. Data that do not decay in $x$ have been treated when they are periodic in $x$ and square integrable in $r$, or square integrable perturbations of a KdV soliton \cite{MST2011,Mizumachi15}. The inverse scattering transform has been carried out for smooth localized data \cite{Wickerhauser87,ChenKP}.

The data of Theorem~\ref{thm:intro-existence} need not lie in these spaces, either locally or at infinity. In case~(i) they are integrable but may be unbounded above, so that they need not lie in $L^2$ or in the anisotropic Sobolev spaces above, and yet the solution is real analytic at positive times. In cases~(ii) and~(iii) the data need to decay only as $r\to+\infty$. In every other direction they may be merely bounded and measurable, for example $V(x,r)=-10/((1+e^x)(1+e^r))$, or even grow. The representation \eqref{representation} evolves such data directly. The closest earlier result concerns KdV, which is \eqref{KP} for data independent of $x$. Grudsky and Rybkin \cite{GrudskyRybkin15,GrudskyRybkinKdV} constructed classical solutions of KdV for data that decay in one direction only, through a determinant formula with Hankel operators, and their class contains the data of case~(ii) that do not depend on $x$. For KP-II we know of no earlier existence theorem for bounded data in the plane that decay in one direction only, nor for data that grow.

The path integral $Z_V$ is built on the occupation $I_V^{t,x}(r)$, which is finite almost surely whenever $V$ is integrable against the one-point density of the Airy lines, as in \eqref{eq:occ-finite}. Nonetheless, the representation extends by renormalization to data for which this condition fails. In particular, consider the global solutions $\Phi_t(V)$ of KP-II with real initial data $V\in L^2(\R^2)$ \cite[Theorem~1.6]{Hadac2008}. For such data the occupation $I_V^{t,x}(r)$ may be infinite, but its fluctuation $I_V^{t,x}(r)-\E I_V^{t,x}(r)$ extends by continuity from bounded data with compact support to all of $L^2(\R^2)$. Writing $X_V^{t,x}(r)$ for the extension, we show in Theorem~\ref{thm:L2} that, in the sense of distributions,
\[
 \Phi_t(V)=U_tV-2\partial_r^2\log\E\exp\bigl\{-X_V^{t,x}(r)\bigr\},
\]
where $U_t$ is the linear KP-II group. The two terms come from the mean and from the fluctuation of the occupation. The mean carries the linear evolution, since $2\partial_r^2\E I_V^{t,x}(r)=U_tV(x,r)$ for smooth $V$ with compact support, and the centered expectation, which is a regularized Fredholm determinant, carries the nonlinear interaction.

\paragraph{The restart identity and uniqueness.}
Let $\phi(s)=\phi(s,\cdot,\cdot)$ be a potential at time $s$, not necessarily a solution of KP, and write the equation in \eqref{KP} as $\partial_t\phi=\mathcal N(\phi)$. In Section~\ref{sec:existence-properties} we show that
\begin{equation}\label{eq:intro-generator}
 \partial_tZ_{\phi(s)}(t,x,r)=\frac{d}{d\varepsilon}\Big|_{\varepsilon=0}Z_{\phi(s)+\varepsilon\mathcal N(\phi(s))}(t,x,r)=-\E\Bigl[e^{-I_{\phi(s)}}I_{\mathcal N(\phi(s))}\Bigr],
\end{equation}
where $I_W=I_W^{t,x}(r)$ is the occupation \eqref{eq:airy-occupation} of a potential $W$. Therefore, by the linearity of the occupation in the potential,
\begin{equation}\label{eq:intro-transport}
 (\partial_s-\partial_t)Z_{\phi(s)}(t,x,r)=-\E\Bigl[e^{-I_{\phi(s)}}I_{\partial_s\phi(s)-\mathcal N(\phi(s))}\Bigr].
\end{equation}
When $\phi$ solves KP, the right side of \eqref{eq:intro-transport} vanishes, so $Z_{\phi(s)}(t,x,r)$ depends on $s$ and $t$ only through $s+t$, and the law of large numbers for the occupation identifies it as $\exp\{-\frac12\partial_r^{-2}\phi(s+t,x,r)\}$. The next theorem establishes this formula for solutions that satisfy growth bounds at positive times.

\begin{theorem}\label{thm:intro-restart}
Let $\phi$ be a classical solution of the equation in \eqref{KP} on $(0,T)\times\R^2$. Suppose that for each compact interval $J\subset(0,T)$ there are constants $C$, $B$ and $\delta>0$ such that, for $s\in J$ and $x,r\in\R$,
\[
 \phi(s,x,r)\ge-C,\qquad|\partial_x^j\partial_r^k\phi(s,x,r)|\le Ce^{B|x|-\delta r},\quad j+k\le3.
\]
Then, for $s,t>0$ with $s+t<T$ and $x,r\in\R$,
\begin{equation}\label{eq:intro-restart}
 Z_{\phi(s)}(t,x,r)=\exp\bigl\{-\tfrac12\partial_r^{-2}\phi(s+t,x,r)\bigr\},\qquad\text{and hence}\qquad\phi(s+t)=\phi_{\phi(s)}(t).
\end{equation}
\end{theorem}

We call \eqref{eq:intro-restart} the restart identity: a solution of KP is the path integral of its own earlier states. The bounds of Theorem~\ref{thm:intro-restart} are required only at positive times, and they hold, for instance, for the line solitons of Section~\ref{sec:soliton}. Theorem~\ref{thm:representation} allows, more generally, the growth $e^{B|x|^\beta}$ for any $\beta<2$ and a negative part bounded by $M(1+r_-)^\nu$ with $\nu<1$. For a solution of the form $\phi_V$, the identity $\phi(s+t)=\phi_{\phi(s)}(t)$ would be the semigroup property of the KP flow, if uniqueness were known. But Theorem~\ref{thm:intro-restart} does not mention the formula \eqref{representation} in its hypotheses: the solution $\phi$ is arbitrary, subject only to growth bounds at positive times. In this way, the restart identity becomes our tool for proving uniqueness.

Theorem~\ref{thm:intro-restart} applies to every solution in its class, whether or not the solution was constructed by the formula, and it therefore gives uniqueness once the limit $s\downarrow0$ can be taken in \eqref{eq:intro-restart}. Suppose that $\phi(s)\to V$ as $s\downarrow0$. The right side of \eqref{eq:intro-restart} then tends to $\exp\{-\frac12\partial_r^{-2}\phi(t,x,r)\}$, and if the left side tends to $Z_V(t,x,r)$, then $\phi=\phi_V$. We prove this convergence for classical solutions that satisfy such bounds uniformly down to time $0$ (Theorem~\ref{thm:representation}), and for distributional solutions such that $|\phi(t,x,r)|\le Me^{-\delta r_+}$ and $\phi(t)\to V$ in $\mathcal D'(\R^2)$ as $t\downarrow0$ (Theorem~\ref{thm:bounded-unique}); every such solution is equal to $\phi_V$. In both cases the bounds serve to control the negative part of $\phi(s)$ as $s\downarrow0$. Such control is not automatic, since there is a continuous $V\ge0$ with compact support for which $\phi_V$ is unbounded below as $t\downarrow0$ (Remark~\ref{rem:focusing}).

\paragraph{Classical scattering.}
The classical scattering data of KdV and of KP-II can also be read from the formula. When $V(x,r)=v(r)$ does not depend on $x$, with $v\in C_c^\infty(\R)$, the equation in \eqref{KP} is KdV, and we show in Section~\ref{sec:kdv-marchenko} that $Z_v$ is the Marchenko determinant of $v$, which is built from the reflection coefficient, the bound states and the norming constants of the Schr\"odinger operator $-\partial_r^2+v$. For KP-II itself, the scattering coefficient of the heat operator at a spectral parameter $q+i\kappa$ with $\kappa>0$ is a value of the resolvent kernel of a trace-class operator $\mathsf T_{q,\kappa;V}$ on $L^2(-\kappa,\kappa)$, which is formed from the Fourier transform of the middle factor $D_{a,b;V}$ of \eqref{eq:intro-fredholm} (Section~\ref{sec:classical-heat-scattering}). Since this operator comes from Brownian paths, its determinant can be compared with an expectation over nonintersecting Brownian motions, and the comparison gives
\[
 \det\bigl(I-\mathsf T_{q,\kappa;V}\bigr)\ge\exp\Bigl\{-\frac\kappa\pi\int_{\R^2}V(y,r)\diff y\diff r\Bigr\}>0
\]
for every real $V\in C_c^\infty(\R^2)$. The resolvent therefore exists for potentials of any size and of either sign, and by Theorem~\ref{thm:lt-heat} the direct scattering equation has a unique solution at every such spectral parameter.

\paragraph{Nonintersecting Brownian bridges.}
The Airy line ensemble is the limit, at the edge, of $N$ nonintersecting Brownian bridges as $N\to\infty$ \cite{CH}, and before the limit the path integral satisfies an exact equation. Let $B_1<\dots<B_N$ be nonintersecting Brownian bridges on $[-T,T]$, from height $u$ at $-T$ to height $v$ at $T$, and let
\[
 Z_N^V(u,v)=\E\exp\Bigl\{-\sum_{i=1}^N\int_{-T}^TV\bigl(y,B_i(y)\bigr)\diff y\Bigr\}.
\]
By the Karlin--McGregor formula \cite{KM}, $Z_N^V$ is a ratio of two determinants, and the Desnanot--Jacobi identity for these determinants is the Toda equation
\[
 \partial_u\partial_v\log Z_N^V=\frac N{4T}\Bigl(\frac{Z_{N+1}^VZ_{N-1}^V}{(Z_N^V)^2}-1\Bigr),
\]
in which the number of bridges is the lattice variable (Proposition~\ref{prop:toda}). In the scaling in which the lowest bridges converge to the curves $r-\mathcal A_i^{t,x}$, the Toda equation becomes KP. The formal computation is known from the polynuclear growth model \cite{PNG-Toda}. We make it a proof (Theorem~\ref{thm:toda-kp}), and so obtain a second derivation of the KP equation for $\phi_V$, $V\in C_c^\infty(\R^2)$, which uses neither the determinant \eqref{eq:intro-fredholm} nor the theorem of \cite{QR}.

The bridges also allow the path integral of a soliton to be evaluated in closed form. For the soliton $V(x,r)=-2\kappa^2\sech^2(\kappa r)$ of KdV, for instance,
\[
 \E\exp\Bigl\{2\kappa^2\sum_{i\ge1}\int_\R\sech^2\bigl(\kappa(r-\mathcal A_i^{t,x}(y))\bigr)\diff y\Bigr\}=1+e^{-2\kappa r+2\kappa^3t/3},
\]
which is the tau function of the soliton at time $t$, computed from its profile at time $0$ without the use of the equation. More generally, if $V$ is the initial condition of a line soliton of KP-II whose coefficient matrix is totally nonnegative \cite{ChakravartyKodama,KW}, then, by Theorems~\ref{thm:onestep} and~\ref{thm:general}, $Z_V$ is the Wronskian tau function of the soliton, divided by its dominant term as $r\to+\infty$.

The evaluation rests on a Darboux transformation, which for Brownian paths has a simple meaning: it removes the potential at the price of one more path. In the simplest case $V=-2\partial_r^2\log f$, where $f(x,r)=\sum_jc_je^{p_jr-p_j^2x}$ with $c_j>0$, and the $N$ bridges in the potential $V$ become $N+1$ nonintersecting paths in no potential, of which the lowest and the highest end at free points. The path integral is then the expectation of $f$ at these two endpoints, whose joint law is explicit, and because $f$ is a sum of exponentials, Laplace's method yields its limit as $N\to\infty$. At the highest endpoint only the exponential with the largest $p_j$ survives. At the lowest every exponential survives, and Laplace's method gives each the factor $e^{-p_j^3t/3}$. Their sum is therefore $\sum_jc_je^{p_jr-p_j^2x-p_j^3t/3}$, the tau function of the soliton at time $t$.

\paragraph{The KPZ fixed point.}
Special cases of the formula \eqref{representation} were known from the theory of random growth, with initial data that are not functions. Let $\mathfrak h(t,x;h)$ be the KPZ fixed point \cite{MQR} started from an upper semicontinuous profile $h$ that is bounded above and not identically $-\infty$, and let
\[
 F_h(t,x,r)=\mathbb P\{\mathfrak h(t,x;h)\le r\},\qquad\phi_h=-2\partial_r^2\log F_h.
\]
In \cite{QR} it was proved that $\phi_h$ solves \eqref{KP}. By the variational formula for the KPZ fixed point \cite[Theorem~4.18]{MQR},
\begin{equation}\label{eq:intro-wall}
 F_h(t,x,r)=\mathbb P\bigl\{r-\mathcal A_1^{t,x}(y)\ge h(y)\text{ for all }y\in\R\bigr\},
\end{equation}
and the right side of \eqref{eq:intro-wall} is, formally, $Z_V$ for the potential $V$ that equals $0$ above the graph of $h$ and $+\infty$ below it. Indeed, the occupation of this potential is $0$ if all the curves $r-\mathcal A_i^{t,x}$ stay above the graph, which they do as soon as the lowest curve does, and it is infinite otherwise. The one-point distributions of the KPZ fixed point are in this sense Airy path integrals whose initial data are a wall. The formula of Borodin and Gorin for the stochastic heat equation \cite{BG} and the deformed Airy determinants of Cafasso, Claeys and Ruzza \cite{CCR} are further instances with singular data, concentrated on the line $x=0$; we describe them in Section~\ref{sec:examples-kpz}.

The initial data are not functions, so that none of the previous results apply; it was already observed in \cite{QR} that initial data of this kind lie outside the known well-posedness theories for KP. The restart identity nevertheless applies, and it shows that the KP equation determines the distributions $F_h$. Precisely, $\phi_h$ is the only classical solution of KP that is nonnegative and suitably bounded at positive times, is nonincreasing in $r$ for small times, and has the wall as its initial data in the following sense: as $t\downarrow0$ the solution tends to $0$ locally uniformly above the graph of $h$, while its integrals $\int_{y-\varepsilon}^{y+\varepsilon}\phi(t,x,r)\diff x$ tend to infinity for every point $(y,r)$ below the graph and every $\varepsilon>0$ (Theorem~\ref{thm:kpz-characterization}).

The proof runs the uniqueness argument with the singular data as the limit $s\downarrow0$. For a solution $\phi$ with these properties the potentials $\phi(s)$ converge to the wall, so that $Z_{\phi(s)}(t,x,r)$ converges to the probability in \eqref{eq:intro-wall}, and \eqref{eq:intro-restart} becomes $\exp\{-\frac12\partial_r^{-2}\phi(t,x,r)\}=F_h(t,x,r)$. For $\phi_h$ itself the restart identity reads
\[
 F_h(s+t,x,r)=Z_{\phi_h(s)}(t,x,r),\qquad s,t>0,
\]
which is Corollary~\ref{cor:kpz-restart}. It expresses the law of the KPZ fixed point at time $s+t$ as the path integral, over all the curves of the Airy line ensemble, of the solution of KP at time $s$.

Among the conditions of this characterization, the monotonicity in $r$ is the one that is not apparent for $\phi_h$. It is equivalent to an inequality which we prove for every $h$ (Theorem~\ref{thm:kpz-curvature}):
\[
 \partial_r^3\log F_h(t,x,r)\ge0.
\]
For a narrow wedge, $F_h$ is a GUE Tracy--Widom distribution function, and the inequality reads $-2qq'\ge0$ for the Hastings--McLeod solution $q$ of Painlev\'e~II \cite{HM}. In general it is an inequality for the finite-dimensional distributions of the Airy$_2$ process.

The same inequality lets us extract asymptotics from the KP equation. For the parabolic profiles $h(y)=-ay^2$, $a\ge0$, we show that, for fixed $t>0$ and $x\in\R$,
\[
 \log F_h(t,x,-R)=-\frac{1+(1+at)^{-1/2}}{12t}R^3+o(R^3)\qquad\text{as }R\to\infty,
\]
and the same asymptotics hold for profiles that are locally bounded below and satisfy $h(y)=-ay^2+o(y^2)$ as $|y|\to\infty$ (Theorem~\ref{thm:kpz-quadratic-tail} and Corollary~\ref{cor:kpz-asymptotic-quadratic}). The coefficient equals $1/(6t)$ for flat data, $a=0$, and decreases to $1/(12t)$ as $a\to\infty$; these are the constants in the lower tails of the GOE and GUE Tracy--Widom distributions. This makes rigorous a heuristic argument in \cite[Remark~2.4]{QR}.

\section{Compactly supported initial data}\label{sec:direct}

The next three subsections give the straightforward proof of the representation in the case of compactly supported $V$.

\subsection{From the Airy line ensemble to a Fredholm determinant}

In the case of compact initial data there is a Fredholm determinant representation of \eqref{representation}. We start there.

\begin{theorem}\label{thm:fredholm}
Let $V\in C_c(\R^2)$ with support in $(a,b)\times\R$. Then for $t>0$ and $x,r\in\R$, $0<Z_V(t,x,r)<\infty$ and
\begin{equation}\label{eq:fredholm}
 Z_V(t, x, r)
 =\det\bigl(I-\chi_r K_{t, x; V}\chi_r\bigr)_{L^2(\R)},
\end{equation}
where $\chi_r=\mathbf 1_{[r,\infty)}$, and 
\begin{equation}\label{eq:kernel}
    K_{t, x; V}(u, v) = \iint_{\R^2}S_{t,a-x}(w,u) D_{a, b;V}(w,w') S_{t,x-b}(w',v)\diff w\diff w',
\end{equation}
where $S_{t,x}$ is defined in \eqref{eq:S}, and $D_{a,b;V}$ is the defect 
\begin{equation}\label{eq:defect}
   D_{a, b; V}(w, w') \vcentcolon= p_{b-a}(w'-w) - P_{a,b;V}(w, w'),
\end{equation}
where $p_s$ is the heat kernel \eqref{eq:heat-kernel}, and
$P_{a,b;V}$ is the corresponding Feynman--Kac kernel with potential $V$,
\begin{equation}\label{eq:FK}
 P_{a, b;V}f(w)=\int_{\R}P_{a,b;V}(w,w')f(w')\diff w'=\E_{B(a)=w}\Bigl[e^{-\int_a^bV(y,B(y))\diff y}f(B(b))\Bigr].
\end{equation}
\end{theorem}

\begin{remark}\label{rem:operator-form}
The integral \eqref{eq:kernel} converges absolutely and $\chi_rK_{t,x;V}\chi_r$ is trace class, as shown in Appendix~\ref{app:limit}. We sometimes use the operator notation $K_{t,x;V}=S^*_{t,a-x}D_{a,b;V}S_{t,x-b}$, but this is a bit formal since the outer factors need not be bounded on $L^2(\R)$.

The kernel is independent of the endpoints $a,b$.
To see this, let $a'<a$, with $V=0$ on $(a',a)\times\R$.
The Markov property gives
$P_{a',b;V}=e^{(a-a')\partial^2}P_{a,b;V}$,
and hence
$D_{a',b;V}=e^{(a-a')\partial^2}D_{a,b;V}$.
On the other hand,
$S^*_{t,a'-x}e^{(a-a')\partial^2}=S^*_{t,a-x}$.
Thus the two factors introduced by extending the interval to
the left cancel in \eqref{eq:kernel}, leaving $K_{t,x;V}$
unchanged. Similarly, if $b'>b$ and $V=0$ on
$(b,b')\times\R$, then
$D_{a,b';V}=D_{a,b;V}e^{(b'-b)\partial^2}$ and
$e^{(b'-b)\partial^2}S_{t,x-b'}=S_{t,x-b}$,
so extending the interval to the right also has no effect.
\end{remark}

\begin{figure}[ht]
\centering
\includegraphics{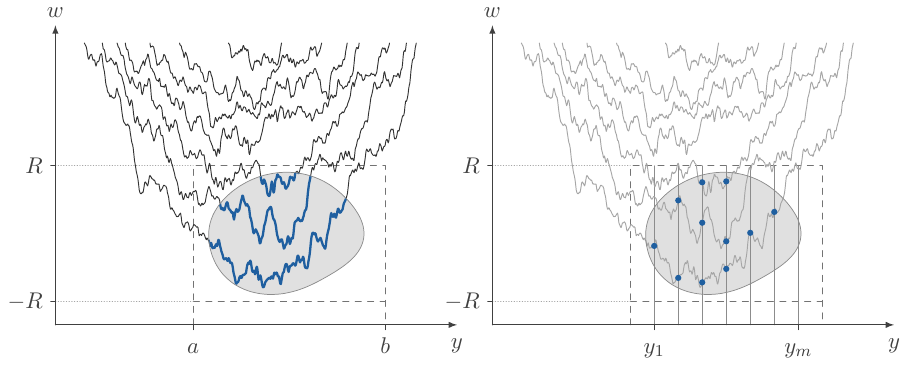}
\caption{The occupation $I_V$ (left) and its discrete version $I_m$ (right).}
\label{fig:occupation-mesh}
\end{figure}

We now prove Theorem~\ref{thm:fredholm}. The starting point is to replace the integral over $y$ in $I_V$ by a Riemann sum (Figure~\ref{fig:occupation-mesh}), so that the Airy correlation functions give a finite-mesh Fredholm determinant. For $m\ge1$, let $\ell=(b-a)/(m+1)$ and $y_j=a+j\ell$ for $j=1,\dots,m$, and let
\[
 Z_m=\E\exp\{-I_m\},\qquad I_m=\ell\sum_{j=1}^m\sum_{i\ge1}V\bigl(y_j,r-\mathcal A_i^{t,x}(y_j)\bigr).
\]
Expanding $e^{-I_m}$ by inclusion--exclusion and taking expectations writes $Z_m$ as a series whose $n$th term is an integral against the $n$-point correlation function of the points $\bigl(y_j,r-\mathcal A_i^{t,x}(y_j)\bigr)$. By \eqref{eq:markov-airy-marginal} and \eqref{eq:airy-correlations}, this correlation function is an $n\times n$ determinant of the extended Airy kernel, and the series can be identified with the Fredholm expansion of the determinant in Lemma~\ref{lem:block}.

\begin{lemma}\label{lem:block}
Let $R>0$ be such that $V$ vanishes outside $(a,b)\times[-R,R]$. Let $g=\operatorname{diag}(g_1,\dots,g_m)$ act on $L^2[-R,R]^m$, where $g_j=1-e^{-\ell V(y_j,\cdot)}$. Let $\widehat K$ be the operator on $L^2[-R,R]^m$ with kernel
\[
 \widehat K_{jk}(w,w')=-\mathbf 1_{j<k}p_{y_k-y_j}(w'-w)+\int_r^\infty S_{t,x-y_j}(w,u)S_{t,y_k-x}(w',u)\diff u.
\]
Then
\[
 Z_m=\det\bigl(I-g\widehat K\bigr)_{L^2[-R,R]^m}.
\]
\end{lemma}

\begin{proof}
Let $G_{ij}=g_j\bigl(r-\mathcal A_i^{t,x}(y_j)\bigr)$ for $i\ge1$ and $1\le j\le m$. Then
\begin{equation}\label{eq:block-expansion}
 \begin{aligned}
 e^{-I_m}&=\prod_{j=1}^m\prod_{i\ge1}e^{-\ell V(y_j,r-\mathcal A_i^{t,x}(y_j))}=\prod_{j=1}^m\prod_{i\ge1}\bigl(1-G_{ij}\bigr)\\
 &=\sum_{n\ge0}\frac{(-1)^n}{n!}\sum_{j_1,\dots,j_n=1}^m\sum_{\substack{i_1,\dots,i_n\ge1\\ i_k\ne i_{k'}\text{ if }j_k=j_{k'}}}\prod_{k=1}^nG_{i_kj_k}.
 \end{aligned}
\end{equation}
The products and sums in \eqref{eq:block-expansion} are finite, since only finitely many $G_{ij}$ are nonzero: $g_j$ vanishes outside $[-R,R]$, and at each $y_j$ only finitely many of the points $r-\mathcal A_i^{t,x}(y_j)$ lie in $[-R,R]$. Taking expectations in \eqref{eq:block-expansion} gives
\[
 \begin{aligned}
 Z_m&=\sum_{n\ge0}\frac{(-1)^n}{n!}\sum_{j_1,\dots,j_n=1}^m\E\sum_{\substack{i_1,\dots,i_n\ge1\\ i_k\ne i_{k'}\text{ if }j_k=j_{k'}}}\prod_{k=1}^nG_{i_kj_k}\\
 &=\sum_{n\ge0}\frac{(-1)^n}{n!}\sum_{j_1,\dots,j_n=1}^m\int_{[-R,R]^n}\prod_{k=1}^ng_{j_k}(w_k)\det\bigl[\widehat K_{j_kj_{k'}}(w_k,w_{k'})\bigr]_{k,k'=1}^n\diff\mathbf w\\
 &=\sum_{n\ge0}\frac{(-1)^n}{n!}\sum_{j_1,\dots,j_n=1}^m\int_{[-R,R]^n}\det\bigl[g_{j_k}(w_k)\widehat K_{j_kj_{k'}}(w_k,w_{k'})\bigr]_{k,k'=1}^n\diff\mathbf w\\
 &=\det\bigl(I-g\widehat K\bigr)_{L^2[-R,R]^m}.
 \end{aligned}
\]
The second equality uses the Airy correlations
\eqref{eq:airy-correlations} with $z_k=r-w_k$, since
$\widehat K_{jk}(w,w')=K^{t,x}(y_j,r-w;y_k,r-w')$ by
\eqref{eq:extended-airy-S}. Lemma~\ref{lem:block-analytic} justifies the
first and last equalities: the expectation can be taken inside the
sums, and $g\widehat K$ is trace class, so its Fredholm series equals
$\det(I-g\widehat K)$.
\end{proof}

For the discretization at the points $y_j$ defined above, let $P_m$ be the Feynman--Kac propagator obtained by replacing the path integral in
\eqref{eq:FK} by $\ell\sum_jV(y_j,B(y_j))$, put
$D_m=e^{(b-a)\partial^2}-P_m$, and define $K_m$ from $D_m$
by the same Airy factors as in \eqref{eq:kernel}.
Recall that Lemma~\ref{lem:block} gives a determinant on $m$ copies of
$L^2[-R,R]$, whereas Theorem~\ref{thm:fredholm} requires one on
$L^2[r,\infty)$. To pass between them, split the extended Airy kernel
into its heat-kernel part and the part factoring through
$L^2[r,\infty)$.

Let $H$ act on $L^2[-R,R]^m$, with blocks
\[
 H_{jk}=
 \begin{cases}
 \mathbf1_{[-R,R]}e^{(y_k-y_j)\partial^2}\mathbf1_{[-R,R]},
 &j<k,\\
 0,&j\ge k.
 \end{cases}
\]
Let $\Phi$ be the column with entries
$\mathbf1_{[-R,R]}S_{t,x-y_j}\chi_r$, and let $\Psi$ be the row
with entries $\chi_rS^*_{t,y_k-x}\mathbf1_{[-R,R]}$.
The kernel in Lemma~\ref{lem:block} then reads
$\widehat K=-H+\Phi\Psi$, and so
\[
 Z_m=\det(I+gH-g\Phi\Psi).
\]
Since $H$ is strictly upper triangular in the mesh-point index,
$(gH)^m=0$,
$\det(I+gH)=1$, and the determinant identity
$\det(I-AB)=\det(I-BA)$ gives
\begin{equation}\label{eq:mesh-reduced}
 Z_m=\det\bigl(I-\Psi(I+gH)^{-1}g\Phi\bigr)_{L^2[r,\infty)}.
\end{equation}

To identify the operator inside the determinant in \eqref{eq:mesh-reduced}, expand the finite inverse
$(I+gH)^{-1}g=\sum_{n=1}^{m}(-gH)^{n-1}g$.
Every term corresponds to a strictly increasing sequence
$j_1<\cdots<j_n$ of mesh points. Inserting the heat
semigroups before and after those times gives the chain
\[
 T_{j_1\cdots j_n}
 =e^{(y_{j_1}-a)\partial^2}g_{j_1}
  e^{(y_{j_2}-y_{j_1})\partial^2}g_{j_2}
  \cdots g_{j_n}e^{(b-y_{j_n})\partial^2}.
\]
These are exactly the terms in the discrete Feynman--Kac
defect. Indeed, by the Markov property,
\[
 P_m=e^{\ell\partial^2}(1-g_1)e^{\ell\partial^2}
       \cdots(1-g_m)e^{\ell\partial^2}.
\]
Expanding this product in the factors $g_j$, the term with
no $g_j$ is the free heat semigroup and cancels in
$D_m=e^{(b-a)\partial^2}-P_m$. Both expansions therefore give
\begin{equation}\label{eq:chains-matched}
 D_m=\sum_{j,k=1}^m e^{(y_j-a)\partial^2}
       \bigl[(I+gH)^{-1}g\bigr]_{jk}
       e^{(b-y_k)\partial^2}.
\end{equation}
Finally, the semigroup property of $S_{t,x}$ absorbs the
heat semigroups in \eqref{eq:chains-matched} into the two outer Airy factors:
\[
 \Psi(I+gH)^{-1}g\Phi=\chi_rK_m\chi_r.
\]
Consequently, the discrete occupation expectation satisfies
\begin{equation}\label{eq:finite-mesh-det}
 Z_m=\det\bigl(I-\chi_rK_m\chi_r\bigr)_{L^2(\R)}.
\end{equation}
The operator estimates justifying these manipulations are in Appendix~\ref{app:mesh}; the limit needed to pass from the discrete expectation to $Z_V$ is stated next.

\begin{proposition}\label{prop:limit}
The expectation $Z_V$ is finite and positive, and $\chi_rK_V\chi_r$ is trace class. As $m\to\infty$,
\[
 Z_m\to Z_V,\qquad\|\chi_rK_m\chi_r-\chi_rK_V\chi_r\|_1\to0,
\]
and hence $\det(I-\chi_rK_m\chi_r)\to\det(I-\chi_rK_V\chi_r)$.
\end{proposition}

Proposition~\ref{prop:limit} is proved in Appendix~\ref{app:limit}.

\begin{proof}[Proof of Theorem~\ref{thm:fredholm}]
By Proposition~\ref{prop:limit}, $0<Z_V<\infty$. Letting $m\to\infty$ in \eqref{eq:finite-mesh-det} gives \eqref{eq:fredholm}.
\end{proof}

\subsection{From the Fredholm determinant to KP-II}\label{sec:fredholm-evolution}

To show that $\phi_V$ solves KP we use the following determinant-flow theorem.

\begin{theorem}[{\cite[Theorem~1.3]{QR}}]\label{thm:qr-flow}
Let $\mathcal K(t,x,r)$, for $t>0$ and $x,r\in\R$, be a trace-class operator on $L^2(\R_+)$ with $\|\mathcal K\|_1$ bounded on compact sets, whose kernel satisfies
\begin{equation}\label{eq:QR-flows}
 \partial_r\mathcal K=(\partial_u+\partial_v)\mathcal K,\qquad\partial_t\mathcal K=-\tfrac13(\partial_u^3+\partial_v^3)\mathcal K,\qquad\partial_x\mathcal K=(\partial_v^2-\partial_u^2)\mathcal K.
\end{equation}
Suppose that $\det(I-\mathcal K)>0$ for all $t>0$ and $x,r\in\R$, that $\mathcal K$ is real analytic in $t$, $x$ and $r$, and that $\|\mathcal K\|_1<1$ for $r$ in some open interval. Then $Q=\partial_r\log\det(I-\mathcal K)$ and $q=\partial_rQ$ satisfy
\begin{equation}\label{eq:QR-KP}
 \partial_tq+\tfrac12\partial_r(q^2)+\tfrac1{12}\partial_r^3q+\tfrac14\partial_x^2Q=0.
\end{equation}
\end{theorem}

Earlier versions in special cases appeared in \cite{ZS74,Poppe84,Poppe89,PS88,HirotaDirect}.

We apply Theorem~\ref{thm:qr-flow} to the kernel of Theorem~\ref{thm:fredholm}. For $V$ as in Theorem~\ref{thm:fredholm}, let
\[
 \mathcal K(t,x,r;u,v)=K_{t,x;V}(r+u,r+v),\qquad u,v>0,
\]
so that $\det(I-\mathcal K)_{L^2(\R_+)}=Z_V(t,x,r)$ by Theorem~\ref{thm:fredholm}.

\begin{proposition}\label{prop:qr-kernel}
\textup{(i)} $\mathcal K$ is trace class, and extends to a trace-norm holomorphic function of $(t,x,r)$ on a complex neighborhood of $(0,\infty)\times\R^2$.

\textup{(ii)} $\mathcal K$ satisfies the flows \eqref{eq:QR-flows}.

\textup{(iii)} For all $j,k,l\ge0$, $\|\partial_t^j\partial_x^k\partial_r^l\mathcal K\|_1\to0$ as $r\to+\infty$, locally uniformly in $(t,x)$.
\end{proposition}

The differential identities \eqref{eq:QR-flows} can be read directly
from the two Airy factors in \eqref{eq:kernel}. Indeed, since
$S_{t,\alpha}=e^{(t/3)\partial^3+\alpha\partial^2}$ and its kernel
$S_{t,\alpha}(w,u)$ depends only on $u-w$,
\[
 \partial_t S_{t,\alpha}(w,u)
 =-\tfrac13\partial_u^3S_{t,\alpha}(w,u),
 \qquad
 \partial_\alpha S_{t,\alpha}(w,u)
 =\partial_u^2S_{t,\alpha}(w,u).
\]
The two Airy factors in \eqref{eq:kernel} have parameters
$(t,a-x)$ and $(t,x-b)$, respectively, while $D_{a,b;V}$
is independent of $t,x,r$. Differentiating under the integral
therefore gives \eqref{eq:QR-flows}, proving
Proposition~\ref{prop:qr-kernel}\textup{(ii)}. The required
derivative bounds are in Appendix~\ref{app:QR-analytic}.

Parts~(i) and~(iii) are proved in Appendix~\ref{app:QR-analytic}.

We can now conclude the KP-II equation for smooth compactly
supported data. Apply Theorem~\ref{thm:qr-flow} to the kernel
$\mathcal K$ of Proposition~\ref{prop:qr-kernel}. Parts~(i) and~(ii)
of the proposition give the trace-class analyticity and the flows,
while part~(iii) gives $\|\mathcal K\|_1<1$ for sufficiently large
$r$. The Fredholm representation in
Theorem~\ref{thm:fredholm} gives
$\det(I-\mathcal K)=Z_V>0$. Since the Fredholm determinant
is holomorphic on trace-class operators, $Z_V$ is real analytic.
Theorem~\ref{thm:qr-flow} yields \eqref{eq:QR-KP} for
$q=\partial_r^2\log Z_V$ and $Q=\partial_r\log Z_V$.
Moreover, $\partial_r Q=q$ and the trace identity
\[
 Q=-\operatorname{tr}\bigl((I-\mathcal K)^{-1}
                        \partial_r\mathcal K\bigr)
\]
together with Proposition~\ref{prop:qr-kernel}(iii) shows that
$Q\to0$ as $r\to+\infty$. Hence $Q=\partial_r^{-1}q$.
Multiplying \eqref{eq:QR-KP} by $-2$ proves that
$\phi_V=-2q$ solves \eqref{KP} for $t>0$.
This completes the KP-II argument.

\subsection{Recovery of the initial data}\label{sec:compact-recovery}

We recover the initial condition in the smooth local topology as $t\downarrow0$. 

\begin{theorem}\label{thm:compact}
Let $V\in C_c^\infty(\R^2;\R)$. The KP-II solution
$\phi_V=-2\partial_r^2\log Z_V$ defined by the Airy representation satisfies,
as $t\downarrow0$,
\begin{equation}\label{eq:compact-trace}
 \phi_V(t,\cdot,\cdot)\longrightarrow V\quad\text{in }C^\infty_{\rm loc}(\R^2).
\end{equation}
\end{theorem}

\begin{figure}[ht]
\centering
\includegraphics{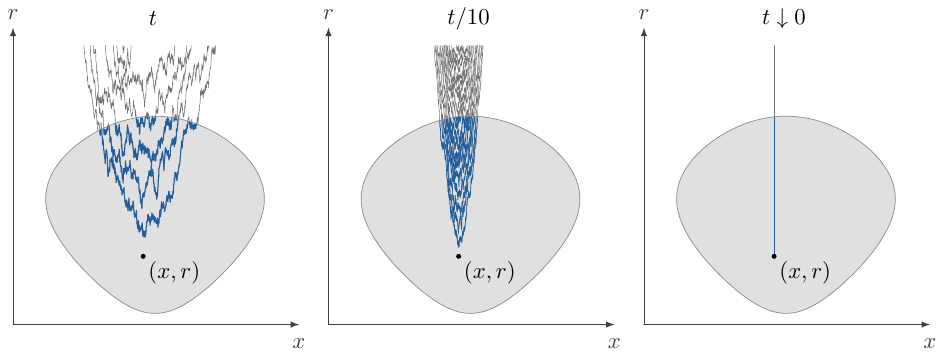}
\caption{The law of large numbers (Proposition~\ref{prop:lln}), at a point $(x,r)$ in the support of $V$ (shaded). Left and middle: the curves $r-\mathcal A_i^{t,x}(\cdot)$ for a finite-size simulation approximating a sample of the Airy line ensemble, at times $t$ and $t/10$. The occupation $I_V^{t,x}(r)$ integrates $V$ along their highlighted parts. Right: in the limit, $I_V^{t,x}(r)$ becomes the integral of $V$ along this half-line against the density $\frac12(u-r)$ at height $u$.}
\label{fig:lln}
\end{figure}

The initial condition in Theorem~\ref{thm:compact} follows from a law of large numbers for the occupation:
$I_V^{t,x}(r)\to\frac12\partial_r^{-2}V(x,r)$ as $t\downarrow0$
(Proposition~\ref{prop:lln} and Figure~\ref{fig:lln}).
Here the right-normalized inverse second derivative is
\[
 \partial_r^{-2}f(x,r):=\int_r^\infty(u-r)f(x,u)\,du,
\]
whenever the integral exists, so that $\partial_r^2\partial_r^{-2}f=f$.
Consequently,
\begin{equation}\label{eq:mV}
 \phi_V=-2\partial_r^2\log\E e^{-I_V}
 \longrightarrow-2\partial_r^2\log e^{-\frac12\partial_r^{-2}V}=V.
\end{equation}

\begin{proposition}[Law of large numbers]\label{prop:lln}
Let $V\in C_c(\R^2)$ be real. As $t\downarrow0$, uniformly for $(x,r)$ in compact subsets of $\R^2$,
\[
 \E I_V^{t,x}(r)\longrightarrow\frac12\partial_r^{-2}V(x,r),\qquad\operatorname{Var}I_V^{t,x}(r)=O(\sqrt t).
\]
In particular, uniformly for $(x,r)$ in compact sets,
\begin{equation}\label{recovery}
 I_V^{t,x}(r)\longrightarrow\frac12\partial_r^{-2}V(x,r)\quad\text{in }L^2.
\end{equation}
\end{proposition}

Here is the short proof of the law of large numbers.  The one-point
Airy intensity gives, after $y=x+v\sqrt t$ and $u=r+v^2-w$,
\begin{equation}\label{eq:lln-mean}
 \E I_V^{t,x}(r)
 =\int_{\R^2}V(x+v\sqrt t,r+v^2-w)\,
 t^{1/6}\rho_{\rm Ai}(t^{-1/3}w)\diff w\diff v.
\end{equation}
Since $\rho_{\rm Ai}(z)\sim\pi^{-1}\sqrt{-z}$ as $z\to-\infty$
and $\rho_{\rm Ai}$ decays rapidly as $z\to+\infty$,
\[
 t^{1/6}\rho_{\rm Ai}(t^{-1/3}w)
 \longrightarrow \pi^{-1}\sqrt{(-w)_+}.
\]
Consequently, by dominated convergence,
\[
 \E I_V^{t,x}(r)\longrightarrow
 \frac1\pi\int_{\R^2}V(x,r+v^2-w)\sqrt{(-w)_+}\diff w\diff v
 =\frac12\int_r^\infty(u-r)V(x,u)\diff u.
\]
For the last equality, put $w=r+v^2-u$ and use
$\pi^{-1}\int_{\R}\sqrt{(u-r-v^2)_+}\diff v
=\tfrac12(u-r)_+$.

For the variance, stationarity of the Airy lines and the same
spatial rescaling give
\[
 I_V^{t,x}(r)=\sqrt t\int_{\R}Y_v\diff v,\qquad
 Y_v=\sum_i f_v\bigl(\bar{\mathcal A}_i(t^{-1/6}v)\bigr),\quad
 f_v(z)=V(x+v\sqrt t,r+v^2-t^{1/3}z).
\]
At each fixed $v$, the Airy points form a determinantal process
with the self-adjoint projection kernel $K_{\rm Ai}$.  The one- and
two-point correlation formulas give
\[
 \operatorname{Var}Y_v
 =\int_{\R}f_v(z)^2\rho_{\rm Ai}(z)\diff z
 -\|K_{\rm Ai}M_{f_v}K_{\rm Ai}\|_{\rm HS}^{\,2}
 \le\int_{\R}f_v(z)^2\rho_{\rm Ai}(z)\diff z,
\]
where $M_{f_v}$ denotes multiplication by $f_v$.  The subtracted term
is nonnegative even for signed $V$; this is the only point where the
determinantal structure is used (see \eqref{eq:occ-variance}).
Since $|V(y,u)|\le Ce^{-u}$, the Airy density Laplace transform
\eqref{app:mu} then gives, locally uniformly in $(x,r)$,
\[
 \operatorname{Var}Y_v\le Ct^{-1/2}e^{-2v^2},
 \qquad 0<t\le1.
\]
Minkowski's inequality now gives
\[
 \sqrt{\operatorname{Var}I_V^{t,x}(r)}
 \le\sqrt t\int_{\R}\sqrt{\operatorname{Var}Y_v}\diff v
 \le C t^{1/4},
\]
which proves the proposition.  In particular, no two-time
Airy correlation estimate is needed.

The domination and locally uniform convergence in the mean calculation,
and the uniform moment bounds behind the variance estimate, are
justified in Appendix~\ref{app:lln}.  The passage from the
occupation law of large numbers to the $C^\infty_{\rm loc}$
initial-data limit in Theorem~\ref{thm:compact} is proved in
Appendix~\ref{app:compact-initial-trace}.

\section{Existence, uniqueness, stability}\label{sec:existence-properties}

\subsection{Extension to decaying initial data}\label{sec:constructs}

We state the three cases of Theorem~\ref{thm:intro-existence} separately, in greater generality.
Here and throughout the paper, a classical solution is one for which $\partial_t\phi$ and the spatial derivatives of $\phi$ through order three are continuous, and the equation holds pointwise.
We write $V_-:=\max\{-V,0\}$ for the negative part of $V$.

\subsubsection{Compact integrable data and data below a height}

We begin with two classes of data for which $Z_V$ is again a Fredholm determinant to which Theorem~\ref{thm:qr-flow} applies.

\begin{theorem}\label{thm:measurable-data}
Let $V\colon\R^2\to\R$ be measurable, and suppose that either
\begin{enumerate}
\item[\textup{(i)}] $V\in L^1(\R^2)$ has compact support and $V_-\le M$, or
\item[\textup{(ii)}] $|V|\le M$ and $V(x,r)=0$ for $r>R$,
\end{enumerate}
for some $M,R<\infty$. Then $\phi_V=-2\partial_r^2\log Z_V$ is a real-analytic solution of the equation in \eqref{KP} on $(0,\infty)\times\R^2$, and $\phi_V(t,\cdot,\cdot)\to V$ in $\mathcal D'(\R^2)$ as $t\downarrow0$.
\end{theorem}

Theorem~\ref{thm:measurable-data} follows from the next two propositions.

\begin{proposition}\label{prop:integrable-det}
Let $V$ be as in Theorem~\ref{thm:measurable-data}(i), with support in $(a,b)\times\R$. Then, for $t>0$ and $x,r\in\R$, $0<Z_V(t,x,r)<\infty$ and \eqref{eq:fredholm} holds, with $K_{t,x;V}$ given by \eqref{eq:kernel}. Moreover, the kernel $\mathcal K_V(t,x,r;u,v)=K_{t,x;V}(r+u,r+v)$ satisfies \textup{(i)--(iii)} of Proposition~\ref{prop:qr-kernel}.
\end{proposition}

We prove Proposition~\ref{prop:integrable-det} in Appendix~\ref{app:extensions}, by approximating $V$ in $L^1(\R^2)$ by functions in $C_c(\R^2)$, to which Theorem~\ref{thm:fredholm} applies. The kernels of the approximations converge in trace norm, because \eqref{eq:kernel} depends on $V$ only through the defect $D_{a,b;V}$, and the defect is Lipschitz in $V$ in the $L^1$ norm, by \eqref{eq:extensions-defect}.

In case~(ii) the support of $V$ need not lie in any $(a,b)\times\R$, so we replace \eqref{eq:kernel} by a kernel that does not involve $a$ and $b$.

\begin{proposition}\label{prop:endpoint-free}
Let $V$ be measurable, with $|V|\le M$ and $V(y,u)=0$ for $u>R$. For $t>0$ and $x,r\in\R$, let $\mathcal K_V$ be the operator on $L^2(\R_+)$ with kernel
\begin{multline}\label{eq:extensions-two-insertions}
\mathcal K_V(u,v)=\int_{\R^2}S_{t,y-x}(w,r+u)V(y,w)S_{t,x-y}(w,r+v)\diff y\diff w\\
-\int_{y<y'}\int_{\R^2}S_{t,y-x}(w,r+u)V(y,w)P_{y,y';V}(w,w')V(y',w')S_{t,x-y'}(w',r+v)\diff w\diff w'\diff y\diff y',
\end{multline}
where $P_{y,y';V}$ is the Feynman--Kac kernel \eqref{eq:FK}. Then $\mathcal K_V$ is trace class, $Z_V(t,x,r)=\det(I-\mathcal K_V)$, and $\mathcal K_V$ satisfies \textup{(i)--(iii)} of Proposition~\ref{prop:qr-kernel}.
\end{proposition}

\begin{proof}[Proof of Proposition~\ref{prop:endpoint-free} for compactly supported $V$]
Let $V$ have compact support in $(a,b)\times\R$. By Proposition~\ref{prop:integrable-det}, it suffices to show that $\mathcal K_V(u,v)=K_{t,x;V}(r+u,r+v)$.

For bounded measurable $V$ and $W$, the Feynman--Kac formula and the Markov property give Duhamel's formula
\begin{equation}\label{eq:duhamel}
 P_{a,b;V}-P_{a,b;W}=-\int_a^bP_{a,y;W}(V_y-W_y)P_{y,b;V}\diff y,
\end{equation}
as a composition of kernels, where $V_y$ is multiplication by $V(y,\cdot)$. Taking $W=0$ in \eqref{eq:duhamel}, and in the same formula with $V$ and $W$ exchanged, gives $P_{a,b;V}=p_{b-a}-\int_a^bp_{y-a}V_yP_{y,b;V}\diff y$ and $P_{y,b;V}=p_{b-y}-\int_y^bP_{y,y';V}V_{y'}p_{b-y'}\diff y'$. Substituting the second identity into the first,
\[
 D_{a,b;V}=\int_a^bp_{y-a}V_yp_{b-y}\diff y-\iint_{a<y<y'<b}p_{y-a}V_yP_{y,y';V}V_{y'}p_{b-y'}\diff y\diff y'.
\]
Insert this into \eqref{eq:kernel}. The heat semigroup identity absorbs the outer heat kernels into the Airy factors, $S_{t,a-x}$ becoming $S_{t,y-x}$ and $S_{t,x-b}$ becoming $S_{t,x-y'}$; Fubini's theorem applies because the heat kernels decay like Gaussians in $w$. After shifting the half-line $[r,\infty)$ to $\R_+$, this gives the kernel \eqref{eq:extensions-two-insertions}. We prove the general case in Appendix~\ref{app:extensions}.
\end{proof}

\begin{proof}[Proof of Theorem~\ref{thm:measurable-data}]
Let $\mathcal K_V$ be as in Proposition~\ref{prop:integrable-det} in case~(i), and as in Proposition~\ref{prop:endpoint-free} in case~(ii). In both cases $Z_V=\det(I-\mathcal K_V)$, and $\mathcal K_V$ satisfies (i)--(iii) of Proposition~\ref{prop:qr-kernel}. Moreover $Z_V>0$, by Proposition~\ref{prop:integrable-det} in case~(i) and by Lemma~\ref{lem:occ-trace}(iii) in case~(ii). So the argument at the end of Section~\ref{sec:fredholm-evolution} applies with $\mathcal K_V$ in place of $\mathcal K$, and shows that $\phi_V$ is real analytic and solves the equation in \eqref{KP} on $(0,\infty)\times\R^2$.

For the initial condition, we show that $\log Z_V(t)\to-\tfrac12\partial_r^{-2}V$ in $L^1_{\rm loc}(\R^2)$ as $t\downarrow0$, so that $\phi_V(t)=-2\partial_r^2\log Z_V(t)\to V$ in $\mathcal D'(\R^2)$. In case~(ii) this is Lemma~\ref{lem:occ-trace}(iii), since $|V(x,r)|\le Me^{R}e^{-r}$. In case~(i) we prove it in Appendix~\ref{app:extensions}.
\end{proof}

\begin{remark}[Scope of the endpoint-free Fredholm kernel]\label{rem:fredholm-scope}
The two-insertion expression \eqref{eq:extensions-two-insertions}
is a trace-class kernel for bounded potentials supported below a
fixed height, by Proposition~\ref{prop:endpoint-free}.
But the ordinary Airy path integral may exist even when the
second integral in that particular kernel is infinite.
For example, choose $0\le h\in C_c^\infty((-2,-1))$, $h\not\equiv0$,
and a smooth function $0\le\chi\le1$ that vanishes for $y\le2$
and equals one for $y\ge3$.  Then
\[
 W(y,w)=\chi(y)e^{-y^2}h(w-y^2)
\]
is a nonnegative Schwartz potential.  At $(t,x,r)=(1,0,0)$,
the occupation $I_W$ has exponential moments of every order,
whereas the two-insertion integral in
\eqref{eq:extensions-two-insertions} diverges for all $u,v\ge0$.
The reason is an uncancelled positive factor
$e^{(y'-y)^3/12}$ when $w$ and $w'$ follow the parabola.
The calculation is in Appendix~\ref{app:parabolic-obstruction}.
This does not exclude other determinant formulas for $W$. Since $0\le W\le\|h\|_\infty e^{-w_+}$, Proposition~\ref{prop:bounded-factor} below expresses $Z_W$ through a regularized determinant on $L^2(\R^2)$. Moreover $W\in L^2(\R^2)$, and in Section~\ref{sec:L2} we show that for every real $V\in L^2(\R^2)$ the path integral with the occupation centered by its mean is a regularized determinant (Theorem~\ref{thm:L2}).
\end{remark}

\subsubsection{Bounded data}\label{sec:determinants}

Next we treat bounded data that decay exponentially as $r\to+\infty$.

\begin{theorem}\label{thm:bounded-data}
Let $V\colon\R^2\to\R$ be measurable, and suppose that for some $M<\infty$ and $\delta>0$,
\begin{equation}\label{eq:bounded-data-bound}
 |V(x,r)|\le Me^{-\delta r_+},\qquad r_+=\max\{r,0\}.
\end{equation}
Then $\phi_V=-2\partial_r^2\log Z_V$ is a $C^\infty$ solution of the equation in \eqref{KP} on $(0,\infty)\times\R^2$, and $\phi_V(t,\cdot,\cdot)\to V$ in $\mathcal D'(\R^2)$ as $t\downarrow0$.
\end{theorem}

\begin{remark}[Step-like data and nonzero backgrounds]
The hypotheses of Theorem~\ref{thm:bounded-data} do not require decay at negative height. For example, $V(r)=-10/(1+e^r)$ tends to $-10$ as $r\to-\infty$ and to $0$ as $r\to+\infty$, and the theorem applies without renormalization. The same is true for genuinely two-dimensional examples such as $V(x,r)=-10/((1+e^x)(1+e^r))$. What is essential here is decay at positive height. A field such as $V(x,r)=-10e^{-x^2}$, independent of $r$, lies outside the ordinary representation: each line has finite occupation, but the sum over infinitely many Airy lines diverges.
\end{remark}

We now prove Theorem~\ref{thm:bounded-data}. For $V$ as in the theorem, the integrals in \eqref{eq:extensions-two-insertions} need not converge (Remark~\ref{rem:fredholm-scope}). The starting point is instead a regularized determinant on $L^2(\R^2)$, which we first state for $V$ with compact support.

\begin{proposition}\label{lem:bounded-det}
Let $V$ be bounded, measurable and compactly supported. For $t>0$ and $x,r\in\R$, let
\begin{equation}\label{eq:bounded-kernel}
 \widehat K^\gamma(y,w;z,w')=t^{-1/3}\int_0^\infty e^{(z-y)\lambda/t^{2/3}}\Ai(\xi_y(w)+\lambda)\Ai(\xi_z(w')+\lambda)\diff\lambda-\mathbf 1_{y<z}p^\gamma(y,w;z,w'),
\end{equation}
where $\xi_y(w)=t^{-1/3}\bigl(r+\tfrac{(y-x)^2}t-w\bigr)$, $p^\gamma(y,w;z,w')=e^{\gamma_z(w')-\gamma_y(w)}p_{z-y}(w'-w)$ and $\gamma_y(w)=\tfrac{2(x-y)^3}{3t^2}+\tfrac{(x-y)(r-w)}t$. Let $\mathsf B_V$ be the operator on $L^2(\R^2)$ with kernel
\[
 \mathsf B_V(y,w;z,w')=|V(y,w)|^{1/2}\,\widehat K^\gamma(y,w;z,w')\operatorname{sgn}V(z,w')\,|V(z,w')|^{1/2}.
\]
Then $\mathsf B_V$ is Hilbert--Schmidt, and
\begin{equation}\label{eq:bounded-det}
 Z_V(t,x,r)=e^{-\E I_V^{t,x}(r)}\det{}_2(I-\mathsf B_V),
\end{equation}
where $\det_2(I-A)=\det((I-A)e^A)$ is the regularized determinant.
\end{proposition}

By \eqref{eq:S}, the factors $e^{\pm\gamma}$ conjugate the kernels $S$, which grow exponentially in the heights, into the Airy kernels of \eqref{eq:bounded-kernel}. This conjugation does not change determinants, and without it the Hilbert--Schmidt bounds in the proof would fail.

\begin{proof}
Let $W=|V|^{1/2}$ and $\widetilde W=\operatorname{sgn}(V)|V|^{1/2}$, and let $\mathsf P:L^2(\R_+)\to L^2(\R^2)$, $\mathsf Q:L^2(\R^2)\to L^2(\R_+)$ and $\mathsf H$ on $L^2(\R^2)$ be the operators with kernels
\begin{equation}\label{eq:bounded-PQ}
 \begin{aligned}
 \mathsf P(y,w;\lambda)&=t^{-1/6}e^{(x-y)\lambda/t^{2/3}}\Ai(\xi_y(w)+\lambda)W(y,w),\\
 \mathsf Q(\lambda;z,w')&=t^{-1/6}e^{(z-x)\lambda/t^{2/3}}\Ai(\xi_z(w')+\lambda)\widetilde W(z,w'),\\
 \mathsf H(y,w;z,w')&=\mathbf 1_{y<z}p^\gamma(y,w;z,w')W(y,w)\widetilde W(z,w').
 \end{aligned}
\end{equation}
By \eqref{eq:bounded-kernel}, $\mathsf B_V=\mathsf P\mathsf Q-\mathsf H$. Conjugation by the unitary map $f\mapsto t^{1/6}f(r+t^{1/3}\cdot)$ from $L^2[r,\infty)$ onto $L^2(\R_+)$ turns the ordered expansion \eqref{eq:ordered-kernel} of Appendix~\ref{app:remarks} into
\[
 \sum_{n\ge1}(-1)^{n+1}\mathsf Q\mathsf H^{n-1}\mathsf P=\mathsf Q(I+\mathsf H)^{-1}\mathsf P.
\]
So, by Theorem~\ref{thm:fredholm} and $\det(I-AB)=\det(I-BA)$,
\[
 Z_V=\det\bigl(I-(I+\mathsf H)^{-1}\mathsf P\mathsf Q\bigr).
\]
Since $I-\mathsf B_V=(I+\mathsf H)\bigl(I-(I+\mathsf H)^{-1}\mathsf P\mathsf Q\bigr)$ and $\det_2(I+\mathsf H)=1$, the product rule for $\det_2$ \cite[Chapter~9]{Simon} gives
\[
 \det{}_2(I-\mathsf B_V)=Z_Ve^{\operatorname{Tr}(\mathsf P\mathsf Q)}=Z_Ve^{\E I_V}.
\]
The operator estimates that justify these steps are in Appendix~\ref{app:bounded}.
\end{proof}

For $V$ as in Theorem~\ref{thm:bounded-data}, the operators $\mathsf P$, $\mathsf Q$ and $\mathsf H$ of \eqref{eq:bounded-PQ} need not be Hilbert--Schmidt, so we cut the plane. Fix $x_0\in\R$, $0<t_0<t_1$ and $r_0\in\R$, and let
\[
 \Theta=(t_0,t_1)\times(x_0-1,x_0+1)\times(r_0,\infty),\qquad\Omega=\{(y,w)\in\R^2:w\le1+\langle y-x_0\rangle\},
\]
where $\langle y\rangle=(1+y^2)^{1/2}$. Let $\mathsf P_\Omega=\mathbf 1_\Omega\mathsf P$, $\mathsf Q_\Omega=\mathsf Q\mathbf 1_\Omega$ and $\mathsf H_\Omega=\mathbf 1_\Omega\mathsf H\mathbf 1_\Omega$, and let $\mathsf E$ be the operator with kernel
\[
 \mathsf E(y,w;z,w')=\bigl(1-\mathbf 1_\Omega(y,w)\mathbf 1_\Omega(z,w')\bigr)\mathsf B_V(y,w;z,w').
\]

\begin{proposition}\label{prop:bounded-factor}
Let $V$ satisfy \eqref{eq:bounded-data-bound}, and let $(t,x,r)\in\Theta$. Then $I+\mathsf H_\Omega$ is invertible; with $\mathsf U=(I+\mathsf H_\Omega)^{-1}$ and $\widetilde{\mathsf B}=\mathsf U(\mathsf P_\Omega\mathsf Q_\Omega+\mathsf E)$, the operator $\mathsf U\mathsf P_\Omega\mathsf Q_\Omega$ is trace class, $\widetilde{\mathsf B}$ is Hilbert--Schmidt, and
\begin{equation}\label{eq:bounded-factor}
 Z_V(t,x,r)=\exp\bigl\{-\operatorname{Tr}(\mathsf U\mathsf P_\Omega\mathsf Q_\Omega)-\E I_{V\mathbf 1_{\Omega^c}}^{t,x}(r)\bigr\}\det{}_2(I-\widetilde{\mathsf B}).
\end{equation}
\end{proposition}

\begin{proof}
Let $V$ have compact support. By \eqref{eq:bounded-kernel},
\begin{equation}\label{eq:bounded-split}
 \mathsf B_V=\mathsf P_\Omega\mathsf Q_\Omega-\mathsf H_\Omega+\mathsf E.
\end{equation}
As in the proof of Proposition~\ref{lem:bounded-det}, the four operators in \eqref{eq:bounded-split} are Hilbert--Schmidt, $I+\mathsf H_\Omega$ is invertible, and $\det_2(I+\mathsf H_\Omega)=1$. By \eqref{eq:bounded-split}, $I-\mathsf B_V=(I+\mathsf H_\Omega)(I-\widetilde{\mathsf B})$, and the product rule for $\det_2$ gives
\[
 \det{}_2(I-\mathsf B_V)=e^{\operatorname{Tr}(\mathsf H_\Omega\widetilde{\mathsf B})}\det{}_2(I-\widetilde{\mathsf B}).
\]
Since $\mathsf H_\Omega$ vanishes off $\Omega\times\Omega$ and $\mathsf E$ vanishes on it, $\mathbf 1_\Omega\widetilde{\mathsf B}\mathbf 1_\Omega=\mathsf U\mathsf P_\Omega\mathsf Q_\Omega$, and with $\mathsf H_\Omega\mathsf U=I-\mathsf U$,
\[
 \operatorname{Tr}(\mathsf H_\Omega\widetilde{\mathsf B})=\operatorname{Tr}(\mathsf H_\Omega\mathsf U\mathsf P_\Omega\mathsf Q_\Omega)=\operatorname{Tr}(\mathsf P_\Omega\mathsf Q_\Omega)-\operatorname{Tr}(\mathsf U\mathsf P_\Omega\mathsf Q_\Omega).
\]
As in the proof of Proposition~\ref{lem:bounded-det}, $\operatorname{Tr}(\mathsf P_\Omega\mathsf Q_\Omega)=\int_\Omega V\varrho=\E I_{V\mathbf 1_\Omega}$. Inserting these into \eqref{eq:bounded-det} gives \eqref{eq:bounded-factor}, since $\E I_V-\E I_{V\mathbf 1_\Omega}=\E I_{V\mathbf 1_{\Omega^c}}$. For general $V$, Lemma~\ref{lem:bounded-operators} gives the invertibility of $I+\mathsf H_\Omega$ and the trace-class and Hilbert--Schmidt properties. Apply \eqref{eq:bounded-factor} to $V_n=V\mathbf 1_{[-n,n]^2}$ and let $n\to\infty$: the right side converges by Lemma~\ref{lem:bounded-operators}(ii), and $Z_{V_n}\to Z_V$ by Lemma~\ref{lem:occ-continuity}.
\end{proof}

Theorem~\ref{thm:bounded-data} follows from the next proposition, which we prove from the factorization \eqref{eq:bounded-factor} in Appendix~\ref{app:bounded}.

\begin{proposition}\label{prop:bounded-stability}
Let $V$ be as in Theorem~\ref{thm:bounded-data}, with $|V(x,r)|\le Me^{-\delta r_+}$.
\begin{enumerate}
\item[\textup{(i)}] $\log Z_V$ is $C^\infty$ on $(0,\infty)\times\R^2$. If $V_n$ are measurable functions with $|V_n(x,r)|\le Me^{-\delta r_+}$ and $V_n\to V$ in $L^1_{\rm loc}(\R^2)$, then $\log Z_{V_n}\to\log Z_V$ in $C^\infty_{\rm loc}((0,\infty)\times\R^2)$.
\item[\textup{(ii)}] Let $0<t_0<t_1$, $r_0\in\R$ and $0<\delta'<\delta$. For every multi-index $\alpha$ in $(t,x,r)$ there is a constant $C$ such that
\begin{equation}\label{eq:bounded-upper}
 |\partial^\alpha\log Z_V(t,x,r)|\le Ce^{-\delta'r},\qquad t_0\le t\le t_1,\ x\in\R,\ r\ge r_0.
\end{equation}
\end{enumerate}
\end{proposition}

\begin{proof}[Proof of Theorem~\ref{thm:bounded-data}]
By Proposition~\ref{prop:bounded-stability}(i), $\phi_V=-2\partial_r^2\log Z_V$ is $C^\infty$ on $(0,\infty)\times\R^2$.

\textit{The KP equation.}
Approximate $V$ by convolving
$V\mathbf1_{[-n,n]^2}$ with a mollifier of radius $1/n$.
Then $V_n\in C_c^\infty$, $|V_n(y,u)|\le Me^\delta
e^{-\delta u_+}$, and $V_n\to V$ in $L^1_{\rm loc}$.
Proposition~\ref{prop:bounded-stability}(i) gives convergence
of $\log Z_{V_n}$ with all spatial derivatives on compact sets,
while \eqref{eq:bounded-upper} supplies the integrable
upper-height bounds. Lemma~\ref{lem:logKP-limit} passes the equation for $\log Z_{V_n}$ to the limit, so
$\phi_V$ solves \eqref{KP}.

\textit{The initial trace.} By Lemma~\ref{lem:occ-trace}(iii), $\log Z_V(t)\to-\bigl(\tfrac12\partial_r^{-2}V\bigr)$ in $L^1_{\rm loc}(\R^2)$ as $t\downarrow0$, so $\phi_V(t)=-2\partial_r^2\log Z_V(t)$ converges to $2\partial_r^2\bigl(\tfrac12\partial_r^{-2}V\bigr)=V$ in $\mathcal D'(\R^2)$.
\end{proof}

\subsubsection{Growing data}

Finally we treat data that may grow in $x$ and as $r\to-\infty$.

\begin{theorem}\label{thm:growing-data}
Let $V\in C^5(\R^2;\R)$, and suppose that its negative part satisfies
\begin{equation}\label{eq:growing-negative-bound}
 V_-(x,r)\le M(1+r_-)^\nu,\qquad r_-=\max\{-r,0\},
\end{equation}
for some $M<\infty$ and $0\le\nu<1$. Suppose that the derivatives $\partial_x^j\partial_r^kV$, $j+k\le5$, obey one of the following, with constants $C_{jk}<\infty$:
\begin{enumerate}
\item[\textup{(i)}] $|\partial_x^j\partial_r^kV(x,r)|\le C_{jk}e^{B(1+|x|)^\beta-\delta r}$, where $B\ge0$, $\delta>0$ and $0\le\beta\le2$;
\item[\textup{(ii)}] for some $m\ge0$, $p>3+m/2$ and $\delta>0$,
\[
 |\partial_x^j\partial_r^kV(x,r)|\le C_{jk}(1+|x|)^m
 \begin{cases}(1+r)^{-p},&r\ge0,\\ e^{\delta r_-},&r<0.\end{cases}
\]
\end{enumerate}
Let $T_*=\delta/B$ in case \textup{(i)} with $\beta=2$ and $B>0$, and $T_*=\infty$ otherwise. Then $\phi_V=-2\partial_r^2\log Z_V$ is a classical solution of the equation in \eqref{KP} on $(0,T_*)\times\R^2$, and $\phi_V(t,\cdot,\cdot)\to V$ in $C^3_{\rm loc}(\R^2)$ as $t\downarrow0$.
\end{theorem}

We prove Theorem~\ref{thm:growing-data} in Appendix~\ref{app:bounded-growing}. For $V$ as in the theorem we do not necessarily have a determinant representation, so we work with the expectation $Z_V=\E e^{-I_V}$ itself. The derivatives of $Z_V$ can be moved onto $V$ by differentiating under the expectation, as in $\partial_rZ_V=-\E[e^{-I_V}I_{\partial_rV}]$, and estimates on the occupation (Appendix~\ref{app:occupation}) bound the resulting expectations.

The local initial-data convergence above is consistent with the instantaneous KP mass constraint: when the height integral $\int_\R\phi(t,x,r)\,dr$ exists with sufficient decay, its second derivative in $x$ vanishes, and suitable decay in $x$ forces the integral to vanish for positive time \cite{MST}. This need not hold for the height integral of the initial field, since the convergence as $t\downarrow0$ is only local.

For the solutions of Theorems~\ref{thm:measurable-data}, \ref{thm:bounded-data} and~\ref{thm:growing-data}, the equation for $\phi_V=-2\partial_r^2\log Z_V$ is equivalent to the Hirota form of KP-II,
\begin{equation}\label{eq:hirota-physical}
 (D_r^4+12D_rD_t+3D_x^2)Z_V\cdot Z_V=0.
\end{equation}

\begin{remark}\label{rem:smooth-data}
If $V\in C^N$ for some $N\ge5$, satisfies \eqref{eq:growing-negative-bound}, and its derivatives through order $N$ satisfy either set of bounds in Theorem~\ref{thm:growing-data}, then $\phi_V(t,\cdot,\cdot)\to V$ in $C^{N-2}_{\rm loc}(\R^2)$ as $t\downarrow0$. The same holds in case~\textup{(i)} of Theorem~\ref{thm:measurable-data} when $V\in C^N_c$ with $N\ge2$.
\end{remark}

\subsection{Solutions are path integrals of their own earlier state}\label{sec:representation}

The restart identity rests on the transport identity \eqref{eq:intro-transport}. For smooth potentials with compact support, it reads as follows.

\begin{proposition}\label{prop:transport}
Let $\psi\in C_c^\infty(\R^3)$ be real, and let
\[
 R_\psi=\partial_s\psi-\frac12\psi\partial_r\psi+\frac1{12}\partial_r^3\psi+\frac14\partial_x^2\partial_r^{-1}\psi
\]
be its residual in KP. Then for $s,x,r\in\R$ and $t>0$,
\begin{equation}\label{eq:transport}
 (\partial_s-\partial_t)Z_{\psi(s)}(t,x,r)=\frac{d}{d\varepsilon}\Big|_{\varepsilon=0}Z_{\psi(s)+\varepsilon R_\psi(s)}(t,x,r)=-\E\Bigl[e^{-I_{\psi(s)}}\,I_{R_\psi(s)}\Bigr],
\end{equation}
where $I_V=I_V^{t,x}(r)$ is the occupation at $(t,x,r)$.
\end{proposition}

\begin{proof}
We write $\psi$ as $\psi(s,y,u)$, with $\partial=\partial_u$. Choose $a<b$ such that all $\psi(s)$ vanish outside $(a,b)\times\R$. Fix $\tau$, $x$ and $r$, and let $K_s=K_{\tau-s,x;\psi(s)}$, so that $Z_{\psi(s)}(\tau-s,x,r)=\det(I-\chi_rK_s\chi_r)$ by Theorem~\ref{thm:fredholm}.

Let $s<\tau$, and write $t=\tau-s$. The left side of \eqref{eq:transport}, at $(s,t)$, is the $s$-derivative of $Z_{\psi(s)}(\tau-s,x,r)=\det(I-\chi_rK_s\chi_r)$. By Jacobi's formula $\partial\det(I-A)=-\det(I-A)\operatorname{tr}\bigl((I-A)^{-1}\partial A\bigr)$,
\begin{equation}\label{eq:representation-jacobi}
 \frac{d}{ds}\det\bigl(I-\chi_rK_s\chi_r\bigr)=-Z_{\psi(s)}(t,x,r)\operatorname{tr}\bigl((I-\chi_rK_s\chi_r)^{-1}\chi_r\partial_sK_s\chi_r\bigr),
\end{equation}
and $I-\chi_rK_s\chi_r$ is invertible because its determinant $Z_{\psi(s)}(t,x,r)$ is positive. We show below that the kernel $K_s$ satisfies
\begin{equation}\label{eq:representation-kernel-identity}
 \frac{d}{ds}K_s=\frac{d}{d\varepsilon}\Big|_{\varepsilon=0}K_{\tau-s,x;\psi(s)+\varepsilon R_\psi(s)}.
\end{equation}
Inserting \eqref{eq:representation-kernel-identity} into \eqref{eq:representation-jacobi}, and applying Jacobi's formula again, now in $\varepsilon$, gives
\begin{equation}\label{eq:representation-variation}
 \begin{aligned}
 \frac{d}{ds}\det\bigl(I-\chi_rK_s\chi_r\bigr)&=\frac{d}{d\varepsilon}\Big|_{\varepsilon=0}\det\bigl(I-\chi_rK_{\tau-s,x;\psi(s)+\varepsilon R_\psi(s)}\chi_r\bigr)\\
 &=\frac{d}{d\varepsilon}\Big|_{\varepsilon=0}Z_{\psi(s)+\varepsilon R_\psi(s)}(\tau-s,x,r).
 \end{aligned}
\end{equation}
The last expression in \eqref{eq:representation-variation} is the derivative with respect to $\varepsilon$ appearing in \eqref{eq:transport}, evaluated at time $t=\tau-s$. Since $Z_{\psi(s)+\varepsilon R_\psi(s)}(\tau-s,x,r)$ is the expectation of $e^{-I_{\psi(s)}-\varepsilon I_{R_\psi(s)}}$, with $I_V=I_V^{\tau-s,x}(r)$, differentiation under the expectation identifies this derivative with $-\E[e^{-I_{\psi(s)}}I_{R_\psi(s)}]$, the final expression in \eqref{eq:transport}.

We prove \eqref{eq:representation-kernel-identity} in the next three steps. The determinant formulas and the derivatives used above are justified in Appendix~\ref{app:representation-transport}. By Remark~\ref{rem:operator-form}, $K_s=S^*_{\tau-s,a-x}D_sS_{\tau-s,x-b}$, where $S_{\tau-s,\alpha}=e^{((\tau-s)/3)\partial^3+\alpha\partial^2}$ has kernel \eqref{eq:S}, $S^*_{\tau-s,\alpha}=e^{-((\tau-s)/3)\partial^3+\alpha\partial^2}$ is its transpose, and $D_s=D_{a,b;\psi(s)}$ is the defect \eqref{eq:defect}.

\textit{The Lax pair.} Let $H_\psi(y)=\partial^2-\psi(s,y,\cdot)$ and
\[
 L=\partial_y+H_\psi,\qquad A=\tfrac13\partial^3-\tfrac12\psi\partial-\tfrac14(\partial\psi)+\tfrac14(\partial_y\partial^{-1}\psi),
\]
where the terms in parentheses act by multiplication. A direct computation gives
\begin{equation}\label{eq:representation-lax}
 [L,A]=\tfrac1{12}\partial^3\psi-\tfrac12\psi\partial\psi+\tfrac14\partial_y^2\partial^{-1}\psi=\partial_sH_\psi+R_\psi,
\end{equation}
where the second equality is the definition of $R_\psi$, since $\partial_sH_\psi=-\partial_s\psi$.

\textit{The defect.} The Feynman--Kac kernels $P_{y,z}=P_{y,z;\psi(s)}$ of \eqref{eq:FK} satisfy $\partial_zP_{y,z}=P_{y,z}H_\psi(z)$ and $\partial_yP_{y,z}=-H_\psi(y)P_{y,z}$. Hence, for $a<y<b$,
\begin{equation}\label{eq:representation-derivative}
 \partial_y\bigl(P_{a,y}A(y)P_{y,b}\bigr)=P_{a,y}\bigl(\partial_yA+[H_\psi,A]\bigr)(y)P_{y,b}=P_{a,y}[L,A](y)P_{y,b}.
\end{equation}
Since $\psi(s)$ vanishes outside $(a,b)\times\R$, $A=A_0\vcentcolon=\frac13\partial^3$ at $y=a$ and $y=b$. So integrating \eqref{eq:representation-derivative} over $a<y<b$ and inserting \eqref{eq:representation-lax} gives
\begin{equation}\label{eq:representation-integrated}
 P_{a,b}A_0-A_0P_{a,b}=\int_a^bP_{a,y}\bigl(\partial_sH_\psi+R_\psi\bigr)(y)P_{y,b}\diff y.
\end{equation}
By Duhamel's formula \eqref{eq:duhamel}, the part of the right side of \eqref{eq:representation-integrated} with $\partial_sH_\psi$ is $\partial_sP_{a,b}$. Solving for $\partial_sP_{a,b}$ gives
\begin{equation}\label{eq:representation-P}
 \partial_sP_{a,b}=[P_{a,b},A_0]-\int_a^bP_{a,y}R_\psi(y)P_{y,b}\diff y.
\end{equation}
In terms of $D_s=e^{(b-a)\partial^2}-P_{a,b}$, the left side of \eqref{eq:representation-P} is $-\partial_sD_s$. The first term on its right side is $[A_0,D_s]$, since $e^{(b-a)\partial^2}$ commutes with $A_0$. Hence
\begin{equation}\label{eq:representation-defect}
 \partial_sD_s=-[A_0,D_s]+\partial_{R_\psi}D_s,\qquad\partial_{R_\psi}D_s=\int_a^bP_{a,y}R_\psi(y)P_{y,b}\diff y.
\end{equation}
Since $\psi(s)$ is smooth with compact support, the kernels $P_{y,z}$ are smooth with Gaussian bounds on their derivatives, so all these kernel integrals converge absolutely.

\textit{Differentiating the Airy factors.} By Duhamel's formula \eqref{eq:duhamel}, $\partial_{R_\psi}D_s$ is the first variation of $D_{a,b;V}$ at $V=\psi(s)$ in the direction $R_\psi(s)$. Since the Airy factors do not depend on $V$, the right side of \eqref{eq:representation-kernel-identity} is therefore $S^*_{\tau-s,a-x}(\partial_{R_\psi}D_s)S_{\tau-s,x-b}$, and \eqref{eq:representation-kernel-identity} becomes
\begin{equation}\label{eq:representation-dressing}
 \partial_sK_s=S^*_{\tau-s,a-x}\bigl(\partial_{R_\psi}D_s\bigr)S_{\tau-s,x-b}.
\end{equation}
By the product rule, the middle factor $D_s$ contributes its derivative \eqref{eq:representation-defect}, $-[A_0,D_s]+\partial_{R_\psi}D_s$, between the Airy factors. The Airy factors contribute $[A_0,D_s]$ between them, since $\partial_sS_{\tau-s,\alpha}=-A_0S_{\tau-s,\alpha}$ and $\partial_sS^*_{\tau-s,\alpha}=A_0S^*_{\tau-s,\alpha}$, and $A_0$ commutes with both. The two commutators cancel, which gives \eqref{eq:representation-dressing}. In the kernel form \eqref{eq:kernel}, moving $A_0$ past the Airy factors is an integration by parts in $w$ and $w'$. There are no boundary terms, since $D_s(w,w')$ is smooth and, with its derivatives, decays like a Gaussian in $|w|+|w'|$, while $S_{\tau-s,\alpha}(w,u)$ grows at most exponentially in $w$ (proof of Lemma~\ref{lem:S}).
\end{proof}

For a solution $\phi$ of KP, $R_\phi=0$. Proposition~\ref{prop:transport} does not apply to $\phi$ directly, since solutions of KP do not in general have compact support.

\begin{theorem}\label{thm:representation}
Let $T>0$, and let $\phi$ be a classical solution of the equation in \eqref{KP} on $(0,T)\times\R^2$. Suppose that for each compact interval $J\subset(0,T)$ there are constants $C_{jk}<\infty$, $B\ge0$, $0\le\beta<2$, $\delta>0$, $M<\infty$ and $0\le\nu<1$ such that, for $s\in J$ and $x,r\in\R$,
\begin{equation}\label{eq:representation-bounds}
 |\partial_x^j\partial_r^k\phi(s,x,r)|\le C_{jk}e^{B(1+|x|)^\beta-\delta r},\quad j+k\le3,\qquad
 \phi_-(s,x,r)\le M(1+r_-)^\nu.
\end{equation}
\begin{enumerate}
\item[\textup{(i)}] For $s,t>0$ with $s+t<T$ and $x,r\in\R$,
\[
 Z_{\phi(s)}(t,x,r)=\exp\{-\tfrac12\partial_r^{-2}\phi(s+t,x,r)\}.
\]
In particular, $\phi(s+t,x,r)=-2\partial_r^2\log Z_{\phi(s)}(t,x,r)$.
\item[\textup{(ii)}] Suppose moreover that for each $T'<T$ there are constants $C,B,M<\infty$, $0\le\beta<2$, $\delta>0$ and $0\le\nu<1$ such that, for $0<s<T'$ and $x,r\in\R$,
\[
 |\phi(s,x,r)|\le Ce^{B(1+|x|)^\beta-\delta r},\qquad\phi_-(s,x,r)\le M(1+r_-)^\nu,
\]
and that $\phi(s)\to V$ locally uniformly as $s\downarrow0$. Then $Z_V(t,x,r)=\exp\{-\tfrac12\partial_r^{-2}\phi(t,x,r)\}$, and so $\phi(t)=-2\partial_r^2\log Z_V(t)$, for $0<t<T$ and $x,r\in\R$. In particular, two such solutions with the same $V$ coincide.
\end{enumerate}
\end{theorem}

Theorem~\ref{thm:representation} follows from the next proposition, which extends Proposition~\ref{prop:transport}, in integrated form, to potentials without compact support. We prove it in Appendix~\ref{app:transport-noncompact}.

\begin{proposition}\label{prop:transport-noncompact}
Let $\psi\colon(s_-,s_+)\times\R^2\to\R$ have continuous $\partial_s\psi$ and spatial derivatives through order three, satisfying \eqref{eq:representation-bounds} on compact time intervals, with $|\partial_s\psi|\le C\bigl(e^{B(1+|x|)^\beta-\delta r}+e^{2B(1+|x|)^\beta-2\delta r}\bigr)$ there, and let $R_\psi$ be its residual in KP, as in Proposition~\ref{prop:transport}. Then, for $s_-<s_1<s_2<s_+$, $\tau>s_2$ and $x,r\in\R$,
\begin{equation}\label{eq:kpz-transport}
 Z_{\psi(s_2)}(\tau-s_2,x,r)-Z_{\psi(s_1)}(\tau-s_1,x,r)
 =-\int_{s_1}^{s_2}\E\Bigl[
 e^{-I_{\psi(s)}^{\tau-s,x}(r)}
 I_{R_\psi(s)}^{\tau-s,x}(r)\Bigr]\diff s.
\end{equation}
If $\tau<s_+$, the same holds for $s_2=\tau$, with $Z_{\psi(\tau)}(0,x,r)=\exp\{-\frac12\partial_r^{-2}\psi(\tau,x,r)\}$.
\end{proposition}

\begin{proof}[Proof of Theorem~\ref{thm:representation}(i)]
Since $\phi$ solves KP, $R_\phi=0$ and $\partial_s\phi=\frac12\phi\partial_r\phi-\frac1{12}\partial_r^3\phi-\frac14\partial_x^2\partial_r^{-1}\phi$, which satisfies the bound of Proposition~\ref{prop:transport-noncompact} by \eqref{eq:representation-bounds}. So Proposition~\ref{prop:transport-noncompact} applies to $\psi=\phi$ on $(0,T)$. With $s_1=s$ and $s_2=\tau=s+t$, the right side of \eqref{eq:kpz-transport} vanishes, and \eqref{eq:kpz-transport} reads $Z_{\phi(s)}(t,x,r)=\exp\{-\frac12\partial_r^{-2}\phi(s+t,x,r)\}$. Applying $-2\partial_r^2\log$ gives the second identity.
\end{proof}

\begin{proof}[Proof of Theorem~\ref{thm:representation}(ii)]
Fix $0<\tau<T'<T$ and $x,r\in\R$. By part~(i), $Z_{\phi(s)}(\tau-s,x,r)=\exp\{-\tfrac12\partial_r^{-2}\phi(\tau,x,r)\}$ for $0<s<\tau$, which does not depend on $s$. So if
\begin{equation}\label{eq:zero-limit}
 Z_{\phi(s)}(\tau-s,x,r)\longrightarrow Z_V(\tau,x,r)\qquad\text{as }s\downarrow0,
\end{equation}
then $Z_V(\tau,x,r)=\exp\{-\tfrac12\partial_r^{-2}\phi(\tau,x,r)\}$ and $\phi(\tau)=2\partial_r^2\bigl(\tfrac12\partial_r^{-2}(\phi(\tau))\bigr)=-2\partial_r^2\log Z_V(\tau)$. A locally uniform limit keeps pointwise bounds, so $V$ is continuous and satisfies the two bounds of the hypothesis on $(0,T')$. Lemma~\ref{lem:occ-continuity}, with $V_n=\phi(s_n)$ and $t_n=\tau-s_n$ for $s_n\downarrow0$, gives \eqref{eq:zero-limit}.
\end{proof}

Finally, we state a converse of Theorem~\ref{thm:representation}(i).

\begin{proposition}\label{prop:kpz-converse}
Let $T>0$, and let $\phi$ be a continuous function on $(0,T)\times\R^2$ whose spatial derivatives of order at most three are continuous and satisfy the bounds \eqref{eq:representation-bounds} on each compact $J\subset(0,T)$. If $Z_{\phi(s)}(t,x,r)=\exp\{-\frac12\partial_r^{-2}\phi(s+t,x,r)\}$ for all $s,t>0$ with $s+t<T$ and all $x,r\in\R$, then $\phi$ is a classical solution of the equation in \eqref{KP} on $(0,T)\times\R^2$.
\end{proposition}

We prove Proposition~\ref{prop:kpz-converse} in Appendix~\ref{app:kpz-converse-proof}.

\subsection{Uniqueness}\label{sec:initial-data}

Theorem~\ref{thm:representation}(ii) gives uniqueness for classical solutions satisfying bounds down to time $0$. A similar argument applies to bounded weak solutions.

For distributional solutions, the nonlinear term $\phi\partial_r\phi$ in \eqref{KP} is understood as $\frac12\partial_r(\phi^2)$.

\begin{theorem}\label{thm:bounded-unique}
Let $T>0$, and let $\phi$ be a real distributional solution of the equation in \eqref{KP} on $(0,T)\times\R^2$ such that for some $M<\infty$ and $\delta>0$,
\[
 |\phi(t,x,r)|\le Me^{-\delta r_+}\qquad\text{for almost every }(t,x,r)\in(0,T)\times\R^2.
\]
Suppose that $\phi(t)\to V$ in $\mathcal D'(\R^2)$ as $t\downarrow0$, for $t$ outside a set of measure zero. Then $\phi=\phi_V$ almost everywhere, where $\phi_V$ is the solution of Theorem~\ref{thm:bounded-data}; moreover, $V$ is represented by a measurable function satisfying $|V|\le Me^{-\delta r_+}$.
\end{theorem}

We prove Theorem~\ref{thm:bounded-unique} in Appendix~\ref{app:zero}. The proof first shows that $\phi$ is smooth at positive times and satisfies the bounds of Theorem~\ref{thm:representation} on compact time intervals (Lemma~\ref{lem:representation-smoothing}), so that part~(i) of that theorem applies.

\begin{remark}\label{rem:smooth-compact-unique}
For $V\in C_c^\infty$, Proposition~\ref{prop:L2-weighted} and Theorem~\ref{thm:representation}(ii) imply that $\phi_V$ satisfies, for every $T>0$, a bound
\[
 |\phi_V(t,x,r)|\le M_Te^{-\delta_T r_+},\qquad 0<t<T.
\]
Theorem~\ref{thm:bounded-unique} therefore gives uniqueness among distributional solutions with initial data $V$ satisfying such a bound on every finite time interval.
\end{remark}

\begin{remark}\label{rem:focusing}
The bounds down to $t=0$ in Theorems~\ref{thm:representation}(ii) and~\ref{thm:bounded-unique} can fail even for $\phi_V$ with continuous compactly supported $V$. Indeed, there is such a $V\colon\R^2\to[0,1]$ with $\phi_V(t_n,x_n,r_n)\le-n$ for some $t_n\downarrow0$ and $(x_n,r_n)\to(0,0)$ (Appendix~\ref{app:focusing}). Neither theorem applies to this example. Uniqueness without bounds uniform down to $t=0$ remains open.
\end{remark}

\begin{remark}\label{rem:bounded-restart}
We do not know whether the solutions of Theorems~\ref{thm:bounded-data} and~\ref{thm:growing-data} satisfy the hypotheses of Theorem~\ref{thm:representation}, since our bounds on $\phi_V(s)$ hold only above each fixed height. For the solutions of Theorem~\ref{thm:bounded-data}, one inequality holds: for $s,t>0$ and $x,r\in\R$,
\[
 0<Z_{\phi_V(s)}(t,x,r)\le Z_V(s+t,x,r),
\]
with equality in the limit $t\downarrow0$ with $s+t$ fixed. The proof of this inequality is in Appendix~\ref{app:zero}. Whether equality holds for all $s$ and $t$ is open.
\end{remark}

\subsection{Comparison and stability in the initial data}
\label{sec:stability}

A useful consequence of the probabilistic representation is a comparison
principle for $-\log Z_V$.
Whenever the expectations are finite and positive,
\begin{equation}\label{eq:stability-order}
 V_1\le V_2\quad\Longrightarrow\quad
 -\log Z_{V_1}\le-\log Z_{V_2}.
\end{equation}
Let $V_\theta=V+\theta W$, and suppose that $V_\theta$ satisfies the
hypotheses of Theorem~\ref{thm:growing-data} for $\theta$ in an open
interval. On that interval, and for $0<t<T_*$, differentiation under the expectation
(Lemmas~\ref{lem:occ-moments} and~\ref{lem:occ-exp}) gives
\begin{align}
 \frac{d}{d\theta}(-\log Z_{V_\theta})
 &=\E_{V_\theta}I_W,\label{eq:stability-first-variation}\\
 \frac{d^2}{d\theta^2}(-\log Z_{V_\theta})
 &=-\operatorname{Var}_{V_\theta}(I_W)\le0,
 \label{eq:stability-second-variation}
\end{align}
where $d\mathbb P_V=e^{-I_V}d\mathbb P/Z_V$.
Thus $W\mapsto\E_{V_\theta}I_W$ is linear and nonnegative for
$W\ge0$.

The one-section Airy intensity gives an explicit deterministic
positive operator. With $\rho_{\rm Ai}(z)=K_{\rm Ai}(z,z)$, define
\begin{equation}\label{eq:stability-intensity}
 (\mathcal T_tV)(x,r)=\mathbb{E} I_V^{t,x}(r)=\int_{\mathbb R^2}V(y,v)\,t^{-1/3}\rho_{\rm Ai}\!\left(\frac{r+(y-x)^2/t-v}{t^{1/3}}\right)\,dy\,dv.
\end{equation}
The formula is valid for nonnegative measurable $V$, and by
absolute integrability for signed $V$.

\begin{theorem}
\label{thm:positive-stability}
Let $V_1,V_2\in C_c^\infty(\mathbb R^2)$ be nonnegative. Then for
$t>0$ and all $x,r$,
\begin{equation}\label{eq:stability-pointwise}
 |\log Z_{V_1}(t,x,r)-\log Z_{V_2}(t,x,r)|
\le \exp\!\left\{\max\bigl((\mathcal T_tV_1)(x,r),
                                    (\mathcal T_tV_2)(x,r)\bigr)\right\}
 (\mathcal T_t|V_1-V_2|)(x,r).
\end{equation}
Consequently, if $\Omega\subset\mathbb R^2$ and
$\mathcal C\subset(0,\infty)\times\mathbb R^2$ are compact, and
$V_1,V_2$ are supported in $\Omega$ with
$\|V_1\|_{L^1}+\|V_2\|_{L^1}\le M$, then
\begin{equation}\label{eq:stability-L1}
 \|\log Z_{V_1}-\log Z_{V_2}\|_{C(\mathcal C)}
 \le C_{\mathcal C,\Omega,M}\|V_1-V_2\|_{L^1(\mathbb R^2)}.
\end{equation}
In particular, for every $t>0$ and
$\eta\in C_c^\infty(\R^2)$ with
$\{t\}\times\operatorname{supp}\eta\subset\mathcal C$,
\begin{equation}\label{eq:stability-distribution}
 |\langle\phi_{V_1}(t)-\phi_{V_2}(t),\eta\rangle|
 \le 2C_{\mathcal C,\Omega,M}\|V_1-V_2\|_{L^1}
          \|\partial_r^2\eta\|_{L^1}.
\end{equation}
\end{theorem}

\begin{proof}
Write $V_\theta=(1-\theta)V_1+\theta V_2$. By
\eqref{eq:stability-first-variation} and $V_\theta\ge0$,
\[
 \left|\frac{d}{d\theta}(-\log Z_{V_\theta})\right|
 \le\frac{\mathbb{E} I_{|V_2-V_1|}}{\mathbb{E} e^{-I_{V_\theta}}}
 \le e^{\mathbb{E} I_{V_\theta}}\mathbb{E} I_{|V_2-V_1|}.
\]
The first inequality uses $e^{-I_{V_\theta}}\le1$; the second is
Jensen's inequality. Since $\mathbb{E} I_V=\mathcal T_tV$, integration
in $\theta$ gives \eqref{eq:stability-pointwise}.
For compact $\Omega,\mathcal C$, the kernel in \eqref{eq:stability-intensity} is
bounded uniformly for $(t,x,r)\in\mathcal C$ and $(y,v)\in\Omega$. Thus
$\mathcal T_tV\le C_{\mathcal C,\Omega}\|V\|_{L^1}$ for nonnegative $V$
supported in $\Omega$, proving \eqref{eq:stability-L1}.
Finally $\phi_{V_1}-\phi_{V_2}
=-2\partial_r^2(\log Z_{V_1}-\log Z_{V_2})$
in distributions, giving \eqref{eq:stability-distribution}.
\end{proof}

\section{A second proof via nonintersecting Brownian bridges}\label{sec:dyson-toda-kp}

Section~\ref{sec:direct} proved that $\phi_V$ solves KP by writing $Z_V$ as a Fredholm determinant and applying the determinant flows (Theorems~\ref{thm:fredholm} and~\ref{thm:qr-flow}). Here we give a second proof, for $V\in C_c^\infty$, through $N$ nonintersecting Brownian bridges whose edge converges to the Airy lines. The bridge analogue $Z_N^V$ of the expectation $Z_V$ satisfies an exact Toda equation (Proposition~\ref{prop:toda}). The formal passage from a Toda equation of this kind to KP is known from the polynuclear growth model \cite{PNG-Toda}; we turn it into a proof (Theorem~\ref{thm:toda-kp}). Together with the law of large numbers (Proposition~\ref{prop:lln}), which identifies the initial data, this shows that $\phi_V$ solves KP-II with initial data $V$, without using Theorems~\ref{thm:fredholm} and~\ref{thm:qr-flow}.

\subsection{The bridge limit}\label{sec:dyson-brownian-scaling}

\begin{figure}[ht]
\centering
\includegraphics{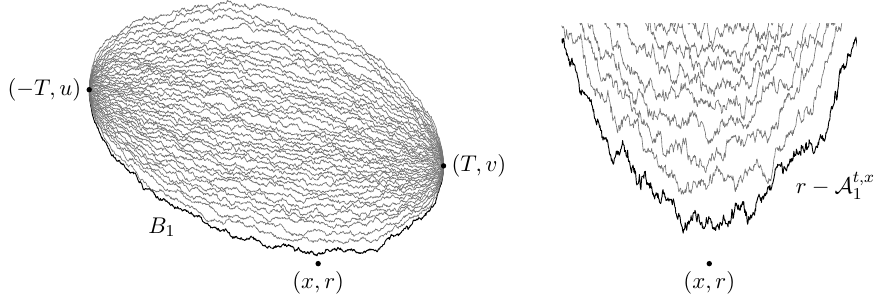}
\caption{The bridges $B_1<\dots<B_N$ for $N=40$, and the curves $r-\mathcal A_i^{t,x}$ to which the lowest bridges converge.}
\label{fig:db-bridges}
\end{figure}

Fix $t>0$ and $x,r\in\R$. For $N\ge1$ let $B_1<\dots<B_N$ be nonintersecting Brownian bridges with diffusion coefficient $2$ on $[-T,T]$, from height $u$ at $-T$ to height $v$ at $T$, where
\begin{equation}\label{eq:db-endpoints}
 T=(t^2N)^{1/3},\qquad u=2\sqrt{NT}+r+\frac{x^2}t+\frac{2xT}t,\qquad v=2\sqrt{NT}+r+\frac{x^2}t-\frac{2xT}t.
\end{equation}
As $N\to\infty$ the lowest bridges converge to the curves $r-\mathcal A_i^{t,x}$ of the occupation \eqref{eq:airy-occupation} (Figure~\ref{fig:db-bridges}): for each $k$,
\begin{equation}\label{eq:db-target}
 \bigl(B_1,\dots,B_k\bigr)\Longrightarrow\bigl(r-\mathcal A_1^{t,x},\dots,r-\mathcal A_k^{t,x}\bigr)
\end{equation}
in law, for the topology of uniform convergence on compact sets of $y$.

The limit \eqref{eq:db-target} follows from the classical limit for nonintersecting Brownian bridges $\widetilde B_1>\dots>\widetilde B_N$ with diffusion coefficient $2$ on $[-T,T]$ that start and end at $0$. At each $y$ the heights $\widetilde B_i(y)$ are, up to scaling, the eigenvalues of an $N\times N$ GUE matrix, as for Dyson's Brownian motion, and for large $N$ the highest bridge is close to $2\sqrt{NT}\sqrt{1-y^2/T^2}$, which is approximately $2\sqrt{NT}-y^2/t$ near $y=0$. The top bridges converge to the Airy lines of \eqref{eq:markov-airy-marginal}: $\widetilde B_i-2\sqrt{NT}\Rightarrow\mathcal A_i^{t,0}$ \cite[Theorem~3.1]{CH}. Since $B_i(y)$ has the law of $2\sqrt{NT}+r+x^2/t-2xy/t-\widetilde B_i(y)$, this gives \eqref{eq:db-target}, by \eqref{eq:markov-airy-marginal}; see Appendix~\ref{app:edge-scaled}.

For a potential $V$, the analogue of the expectation $Z_V(t,x,r)$ for the bridges is
\[
 Z_N^V(u,v)=\E\exp\Bigl\{-\sum_{i=1}^N\int_{-T}^TV\bigl(y,B_i(y)\bigr)\diff y\Bigr\}.
\]

\begin{proposition}[Bridge limit]\label{prop:bridge-limit}
Let $V\in C_c^\infty(\R^2)$ be real. Then $Z_V$ is finite, positive and smooth on $(0,\infty)\times\R^2$. Moreover, for all $j,k\ge0$,
\[
 \partial_x^j\partial_r^k\log Z_N^V(u,v)\longrightarrow\partial_x^j\partial_r^k\log Z_V(t,x,r)\qquad(N\to\infty),
\]
uniformly for $(t,x,r)$ in compact subsets of $(0,\infty)\times\R^2$, where $T$, $u$ and $v$ are given by \eqref{eq:db-endpoints}.
\end{proposition}

The Fredholm argument in Section~\ref{sec:fredholm-evolution} already gives the first assertion. We prove it again in Appendix~\ref{app:bridge-limit}, together with the convergence, so that this section does not use Theorems~\ref{thm:fredholm} and~\ref{thm:qr-flow}.

\subsection{The Toda equation and its KP limit}\label{sec:toda-kp-limit}

We will show that $Z_N^V(u,v)$ is a ratio of two determinants, and that it satisfies a Toda equation in $u$, $v$ and $N$. Let $P_{-T,T;V}(u,v)$ be the Feynman--Kac kernel \eqref{eq:FK}, from height $u$ at time $-T$ to height $v$ at time $T$; for $V=0$ it is $p_{2T}(v-u)$. Put $\tau_0^V=1$, $Z_0^V=1$ and, for $N\ge1$,
\[
 \tau_N^V(u,v)=\det\bigl[\partial_u^i\partial_v^jP_{-T,T;V}(u,v)\bigr]_{i,j=0}^{N-1}.
\]

\begin{proposition}[Toda equation]\label{prop:toda}
Let $V\in C_c(\R^2)$ vanish for $|y|\ge Y$, and let $T>Y$ and $u,v\in\R$. Then for $N\ge1$
\begin{equation}\label{eq:db-ratio}
 Z_N^V(u,v)=\frac{\tau_N^V(u,v)}{\tau_N^0(u,v)}>0,
\end{equation}
and
\begin{equation}\label{eq:db-toda}
 \partial_u\partial_v\log Z_N^V=\frac N{4T}\Bigl(\frac{Z_{N+1}^VZ_{N-1}^V}{(Z_N^V)^2}-1\Bigr).
\end{equation}
\end{proposition}

We prove Proposition~\ref{prop:toda} in Appendix~\ref{app:toda-kp-proofs} (see also O'Connell and Warren \cite[Section~3]{OW}). The ratio \eqref{eq:db-ratio} is the Karlin--McGregor formula \cite{KM} for nonintersecting paths, in the limit where their starting points coalesce at $u$ and their final points at $v$. The Toda equation \eqref{eq:db-toda} then follows from the Desnanot--Jacobi identity for the matrix of $\tau_{N+1}^V$, which reads
\begin{equation}\label{eq:toda-desnanot}
 \tau_N^V\,\partial_u\partial_v\tau_N^V-(\partial_u\tau_N^V)(\partial_v\tau_N^V)=\tau_{N+1}^V\tau_{N-1}^V.
\end{equation}
Indeed, the $u$-derivative of each row of the matrix of $\tau_N^V$ is the next row, so $\partial_u\tau_N^V$ is the determinant with only its last row differentiated; similarly, $\partial_v\tau_N^V$ has only its last column differentiated. Each of the four determinants on the left of \eqref{eq:toda-desnanot} is therefore obtained from the matrix of $\tau_{N+1}^V$ by deleting one of its last two rows and one of its last two columns. Dividing \eqref{eq:toda-desnanot} by $(\tau_N^V)^2$ and subtracting the corresponding identity for $V=0$ gives \eqref{eq:db-toda}.

\begin{theorem}\label{thm:toda-kp}
Let $V\in C_c^\infty(\R^2)$ be real, and let $F=\log Z_V$. Then
\begin{equation}\label{eq:toda-kp-log}
 \partial_t\partial_rF+\frac14\partial_x^2F+\frac1{12}\partial_r^4F+\frac12\bigl(\partial_r^2F\bigr)^2=0\qquad\text{on }(0,\infty)\times\R^2,
\end{equation}
and $\phi_V=-2\partial_r^2F$ solves the equation in \eqref{KP} for $t>0$, with $\partial_r^{-1}\phi_V=-2\partial_rF$.
\end{theorem}

We prove Theorem~\ref{thm:toda-kp} in Appendix~\ref{app:weak-toda-limit} by letting $N\to\infty$ in the Toda equation \eqref{eq:db-toda}. To take this limit, we change variables to
\[
 s=\frac N{T^3},\qquad X=\frac{u-v}{4T},\qquad R=\frac{u+v}2-2\sqrt{NT}.
\]
By \eqref{eq:db-endpoints}, $s=1/t^2$, $X=x/t$ and $R=r+x^2/t$. For convenience, let $h=1/T$. We write $f_h(s,X,R)=\log Z_N^V(u,v)$ and $f(s,X,R)=\log Z_V(t,x,r)$. Since $4\partial_u\partial_v=\partial_R^2-\frac{h^2}4\partial_X^2$ and $N/T=s/h^2$, the Toda equation \eqref{eq:db-toda} becomes
\begin{equation}\label{eq:db-scaled-toda}
 \partial_R^2f_h-\frac{h^2}4\partial_X^2f_h=\frac s{h^2}\bigl(e^{\delta_hf_h}-1\bigr),
\end{equation}
where $\delta_h=\mathsf S_h^++\mathsf S_h^--2$ and
\[
 \mathsf S_h^\pm g(s,X,R)=g\Bigl(s\pm h^3,\,X,\,R+\frac{2\sqrt s-2\sqrt{s\pm h^3}}{h^2}\Bigr).
\]
The shift in $R$ keeps $u$ and $v$ fixed, so $\mathsf S_h^\pm f_h=\log Z_{N\pm1}^V(u,v)$, and $e^{\delta_hf_h}$ is the ratio $Z_{N+1}^VZ_{N-1}^V/(Z_N^V)^2$ in \eqref{eq:db-toda}.

If $f_h$ were smooth in $s$, we could expand $\delta_hf_h$ in powers of $h$ in \eqref{eq:db-scaled-toda}: the terms of order $1$ cancel, and the terms of order $h^2$ give \eqref{eq:toda-kp-log} in the variables $(s,X,R)$. But Proposition~\ref{prop:bridge-limit} controls only the derivatives of $f_h$ in $X$ and $R$, not how $f_h$ changes with $N$. So we multiply by a test function $\psi\in C_c^\infty\bigl((0,\infty)\times\R^2\bigr)$, sum over $N$ and integrate over $X$ and $R$. Summation by parts in $N$ and a translation of $R$ then move $\delta_h$ from $f_h$ onto $\psi$:
\[
 \sum_{N\ge1}\int_{\R^2}(\delta_hf_h)\,\psi\diff X\diff R=\sum_{N\ge1}\int_{\R^2}f_h\,(\delta_h\psi)\diff X\diff R,\qquad s=Nh^3.
\]
On the left, \eqref{eq:db-scaled-toda} gives $\delta_hf_h$ in terms of derivatives of $f_h$ in $X$ and $R$; on the right, the expansion in powers of $h$ is made for $\psi$, which is smooth. The terms of order $h^2$ on the two sides are equal, and those of order $h^4$ give \eqref{eq:toda-kp-log} in the limit $h\to0$ (Proposition~\ref{prop:weak-toda}).

\begin{remark}[Other discrete and multilayer realizations]
The Toda structure is not special to the Brownian-bridge realization used
above. It can also be obtained from multiline polynuclear growth: the
multiline PNG observables satisfy an exact Toda relation, and the
formal $1{:}2{:}3$ KPZ scaling leads from this discrete identity to the
corresponding continuum relation; see \cite{PNG-Toda}. Related Toda structures occur in other integrable
particle and cellular-automaton models. For example, Kamioka and Takagaki
\cite{KamiokaTakagaki} give nonintersecting-path formulas for the discrete
Toda molecule and, after ultradiscretization, for the ultradiscrete Toda
molecule, with the resulting dynamics compared explicitly with the
box--ball system. We do not pursue these alternative realizations here.
\end{remark}

\section{Renormalized $L^2$ initial data and scattering}\label{sec:renormalized-extensions}\label{sec:renormalized-solutions}\label{sec:classical-scattering}

\subsection{Renormalized $L^2$ initial data}\label{sec:L2}

For real $V\in L^2(\R^2)$, KP has a global solution with initial data $V$, unique in a Bourgain-type solution space \cite[Theorem~1.6]{Hadac2008}. These solutions need not satisfy the bounds of Theorem~\ref{thm:representation}, and for $V\in L^2$ the occupation $I_V^{t,x}(r)$ need not be finite. However, after centering the occupation, the path integral still represents these solutions, and the centered path integral is a regularized Fredholm determinant.

Let $U_t$ be the linear KP-II group, the solution operator of \eqref{KP} without its nonlinear term, and let
\[
 \widetilde V_{t,x,r}(s,z)=t^{2/3}V\bigl(x+t^{2/3}s,\,r+t^{1/3}(s^2-z)\bigr).
\]
Define the operator $\widetilde{\mathcal K}_V$ on $L^2(\R^2)$ by its kernel
\[
 \widetilde{\mathcal K}_V(Y,Y')=|\widetilde V_{t,x,r}(Y)|^{1/2}K_{\rm ext}(Y,Y')\operatorname{sgn}\widetilde V_{t,x,r}(Y')\,|\widetilde V_{t,x,r}(Y')|^{1/2},
\]
where, for $Y=(s,z)$ and $Y'=(s',z')$, the stationary extended Airy kernel is
\[
 K_{\rm ext}(Y,Y')=
 \begin{cases}
 \displaystyle\int_0^\infty e^{-(s-s')\lambda}\Ai(z+\lambda)\Ai(z'+\lambda)\diff\lambda,&s\ge s',\\[6pt]
 \displaystyle-\int_0^\infty e^{-(s'-s)\lambda}\Ai(z-\lambda)\Ai(z'-\lambda)\diff\lambda,&s<s'.
 \end{cases}
\]
It is related to the kernel \eqref{eq:extended-airy-S} by
\[
 K^{1,0}(s,z-s^2;s',z'-s'^2)
 =e^{H(s,z)-H(s',z')}K_{\rm ext}(Y,Y'),\qquad
 H(s,z)=\frac{s^3}{3}-sz.
\]
The gauge factors cancel in correlation determinants.
Recall that for a Hilbert--Schmidt operator $A$, $\det_2(I-A)$ denotes the regularized Fredholm determinant; for trace-class $A$ it equals $\det((I-A)e^A)$.

For real $V\in L^2(\R^2)$, let $X_V^{t,x}(r)$ be the $L^2(\mathbb P)$ limit of
\[
 I_{V_n}^{t,x}(r)-\E I_{V_n}^{t,x}(r),
\]
where $V_n$ are real bounded compactly supported functions converging to $V$ in $L^2(\R^2)$. Proposition~\ref{prop:L2-centered} shows that this limit exists and is independent of the approximation.

\begin{theorem}\label{thm:L2}
Let $V\in L^2(\R^2)$ be real, and let $\Phi(V)\in C([0,\infty);L^2(\R^2))$ be the standard global mild solution of KP-II with initial data $V$, defined through the Fourier group \eqref{eq:L2-group} and its Duhamel formulation \cite[Theorem~1.6]{Hadac2008}. Then
\[
 \Phi_t(V)=U_tV-2\partial_r^2\log\tau_V^{\rm ren}(t)\qquad\text{in }\mathcal D'\bigl((0,\infty)\times\R^2\bigr),
\]
where
\[
 \tau_V^{\rm ren}(t,x,r)=\E e^{-X_V^{t,x}(r)}=\det{}_2\bigl(I-\widetilde{\mathcal K}_V\bigr).
\]
\end{theorem}

Theorem~\ref{thm:L2} follows from the three results below, which we prove in Appendix~\ref{app:L2}. The first establishes the stated $L^2$ extension of the centered occupation and its determinant identity.

\begin{proposition}\label{prop:L2-centered}
For real $V_1,V_2\in L^2(\R^2)$ and $t>0$, $x,r\in\R$, the operators $\widetilde{\mathcal K}_{V_1}$ and $\widetilde{\mathcal K}_{V_2}$ are Hilbert--Schmidt, and
\begin{equation}\label{eq:L2-HS}
 \|\widetilde{\mathcal K}_{V_1}\|_{\rm HS}^2\le Ct^{1/3}\|V_1\|_2^2,\qquad
 \|\widetilde{\mathcal K}_{V_1}-\widetilde{\mathcal K}_{V_2}\|_{\rm HS}^2
 \le Ct^{1/3}(\|V_1\|_2+\|V_2\|_2)\|V_1-V_2\|_2.
\end{equation}
The map $V\mapsto X_V^{t,x}(r)$ extends uniquely from real bounded compactly supported functions to a continuous linear map from real $L^2(\R^2)$ to $L^2(\mathbb P)$. For every real $V\in L^2(\R^2)$ and $p\in\R$,
\begin{equation}\label{eq:L2-det2}
 \E|X_V^{t,x}(r)|^2\le Ct^{1/3}\|V\|_2^2,\qquad
 \E e^{pX_V^{t,x}(r)}=\det{}_2(I+p\widetilde{\mathcal K}_V)
 \le e^{Cp^2t^{1/3}\|V\|_2^2}.
\end{equation}
\end{proposition}

The change of variables from the stationary lines to the parabolic
lines then gives a continuous renormalized expectation.

\begin{lemma}\label{lem:L2-tau}
Let $V$ be real and in $L^2(\R^2)$, and let $T>0$.
\begin{enumerate}
\item[\textup{(i)}] For $t>0$ and $x,r\in\R$, $\tau_V^{\rm ren}(t,x,r)=\det_2(I-\widetilde{\mathcal K}_V)$ and
\begin{equation}\label{eq:L2-tau-bounds}
 1\le\tau_V^{\rm ren}(t,x,r)\le e^{Ct^{1/3}\|V\|_2^2}.
\end{equation}
Moreover $\tau_V^{\rm ren}$ is continuous on $(0,\infty)\times\R^2$.
\item[\textup{(ii)}] If $V_n\to V$ in $L^2(\R^2)$, then $\log\tau_{V_n}^{\rm ren}\to\log\tau_V^{\rm ren}$ uniformly on $(0,T]\times\R^2$.
\item[\textup{(iii)}] If $V\in C_c^\infty(\R^2;\R)$, then
\begin{equation}\label{eq:L2-factor}
 Z_V(t,x,r)=e^{-\E I_V^{t,x}(r)}\tau_V^{\rm ren}(t,x,r).
\end{equation}
\end{enumerate}
\end{lemma}

Finally, the mean of the uncentered occupation supplies exactly
the linear KP-II evolution.

\begin{lemma}\label{lem:L2-mean}
For $V\in C_c^\infty(\R^2;\R)$, $t>0$ and $x,r\in\R$,
\[
 2\partial_r^2\E I_V^{t,x}(r)=U_tV(x,r).
\]
\end{lemma}

By Lemmas~\ref{lem:L2-tau}(iii) and~\ref{lem:L2-mean}, the mean occupation accounts for the linear KP-II evolution, and the logarithm of the centered expectation for the nonlinear correction. On the Fredholm side, subtracting the first cumulant removes the linear term in the determinant expansion, leaving $\det_2$, which is defined for Hilbert--Schmidt operators even when the ordinary Fredholm determinant is not.

\subsection{KdV: the one-dimensional check}
\label{sec:kdv-marchenko}

We next compare the Airy representation with classical scattering.
We first check the KdV reduction, where the Airy
determinant coincides with the classical Marchenko determinant.
Let $v\in C_c^\infty(\R;\R)$ and $V(y,r)=v(r)$. Stationarity makes
$Z_v(t,x,r)$ independent of $x$. Write $\phi(t,r)$ for the corresponding
classical KdV solution, obtained by suppressing the transverse term in
\eqref{KP}. Let $-\kappa_j^2$ be the negative eigenvalues of
$-\partial_r^2+v$. Denote by $f_+(r,k)$ the right Jost solution of
$(-\partial_r^2+v)f=k^2f$, normalized by
$f_+(r,k)\sim e^{ikr}$ as $r\to+\infty$. The right norming constants are
\[
 c_j=\left(\int_{\R}f_+(r,i\kappa_j)^2\diff r\right)^{-1}.
\]
For real $k\ne0$, let $f_-(r,k)\sim e^{-ikr}$ as $r\to-\infty$,
and write $f_-(r,k)=a(k)e^{-ikr}+b(k)e^{ikr}$ to the right of
the support of $v$. The right reflection coefficient is
$\mathcal R(k)=b(k)/a(k)$. With these conventions the Marchenko input is
\begin{equation}\label{eq:airy-marchenko-classical-F}
 F_t^{\rm M}(s)=\sum_jc_je^{-\kappa_js+2\kappa_j^3t/3}
 +\frac1{2\pi}\int_{\R}\mathcal R(k)e^{iks+2ik^3t/3}\,dk .
\end{equation}
The corresponding Hankel operator on $L^2([r,\infty))$ is
$(H_{F_t^{\rm M},r}f)(u)=\int_r^\infty F_t^{\rm M}(u+w)f(w)\,dw$.
The classical Marchenko reconstruction formula \cite[(5.3), Theorem~6.2]{GrudskyRybkinKdV} gives
\begin{equation}\label{eq:am-classical-tau}
 \phi(t,r)=-2\partial_r^2\log\det(I+H_{F_t^{\rm M},r}),
 \qquad \det(I+H_{F_t^{\rm M},r})>0.
\end{equation}

For every real $v\in C_c^\infty(\R)$, the Airy representation agrees
with the classical KdV Marchenko determinant:
\begin{equation}\label{eq:am-complete-identification}
 Z_v(t,r)=\det(I+H_{F_t^{\rm M},r})
 =\exp\{-\tfrac12\partial_r^{-2}\phi(t,r)\}.
\end{equation}
Indeed, the classical KdV solution for compact smooth data satisfies
the one-sided exponential decay and derivative bounds of
Theorem~\ref{thm:representation}(ii) on bounded time intervals,
by one-sided weighted KdV energy estimates. For $\lambda>0$ and
$g=e^{\lambda r}\phi$, the KdV equation gives
\[
 \frac{d}{dt}\|g\|_2^2=-\frac\lambda2\|g_r\|_2^2
 +\frac{\lambda^3}{6}\|g\|_2^2-\frac{2\lambda}{3}\int_\R e^{2\lambda r}\phi^3\,dr.
\]
Since the last term is at most $\frac{2\lambda}{3}\|\phi\|_\infty\|g\|_2^2$,
this bounds $\|g(t)\|_2$ on bounded time intervals.
The same argument for the spatial derivatives of $\phi$, with
smooth truncated weights and the classical bounds for smooth solutions of KdV, gives
exponentially weighted Sobolev estimates on bounded time
intervals, and by Sobolev embedding the derivative bounds
required by Theorem~\ref{thm:representation}(ii); the negative part of $\phi$ is
bounded, and $\phi(t)\to v$ locally uniformly as $t\downarrow0$.
Theorem~\ref{thm:representation}(ii) therefore applies, and the two representations
have the same second logarithmic derivative. They are equal because their logarithms and
first derivatives vanish as $r\to+\infty$.

\subsection{Classical heat scattering and the Airy representation}
\label{sec:classical-heat-scattering}\label{sec:kp-soft-marchenko}

The inverse-scattering method for KP-II begins with a heat operator,
rather than the Schr\"odinger operator of KdV.  The essential
difference is that the normalized heat eigenfunction is not
holomorphic in its complex spectral parameter.  We explain how
this leads directly to the inverse $\bar\partial$ problem.  For
the classical construction and rigorous treatments of its
scattering equations, see
\cite{ABF,Wickerhauser87,GrinevichKP,ChenKP}.

\textit{Direct heat scattering.}
Let $V\in C_c^\infty(\R^2;\R)$ be the initial potential in
our KP-II equation.  The spatial Lax equation is
\begin{equation}\label{eq:comparison-heat-lax}
 (\partial_y-\partial_r^2+V(y,r))\Psi_z(y,r)=0.
\end{equation}
For $V=0$, \eqref{eq:comparison-heat-lax} has the solutions $e^{zr+z^2y}$, $z\in\mathbb C$.
Accordingly, for general $V$ write
$\Psi_z(y,r)=e^{zr+z^2y}\chi_z(y,r)$ and choose the
normalized solution $\chi_z\to1$ at spatial infinity.
Substitution into \eqref{eq:comparison-heat-lax} yields
\begin{equation}\label{eq:comparison-direct-equation}
 (\partial_y-\partial_r^2-2z\partial_r)\chi_z=-V\chi_z,
 \qquad \chi_z=1-G_z*(V\chi_z),
\end{equation}
where $G_z$ is the Green function of
$L_z=\partial_y-\partial_r^2-2z\partial_r$.

\textit{The $\bar\partial$ derivative of the Green function.}
The Green function $G_z$ has the Fourier multiplier
$(\xi^2+i\eta-2iz\xi)^{-1}$, where $\xi,\eta$ are
dual to $r,y$.  Put $z=q+i\kappa$ and first suppose
$\kappa>0$.  Inverting in $\eta$ by residues gives
\begin{equation}\label{eq:comparison-green-band}
 G_z(r,y)=\frac1{2\pi}\int_\R
 e^{i\xi(r+2qy)-\xi(\xi+2\kappa)y}
 \bigl(\mathbf1_{\{y>0\}}
       -\mathbf1_{\{-2\kappa<\xi<0\}}\bigr)\,d\xi.
\end{equation}

In the integrand of \eqref{eq:comparison-green-band},
write
$F=e^{i\xi(r+2qy)-\xi(\xi+2\kappa)y}$.
Since $\partial_qF=2i\xi yF$ and
$\partial_\kappa F=-2\xi yF$, we have
$\bar\partial_zF=0$, where
$\bar\partial_z=\frac12(\partial_q+i\partial_\kappa)$.
Thus the entire $\bar\partial_z$ derivative comes
from the band $-2\kappa<\xi<0$ of the second indicator, whose endpoint moves with $\kappa$.  Indeed,
$\partial_\kappa\mathbf1_{\{-2\kappa<\xi<0\}}
=2\delta(\xi+2\kappa)$, so the Fourier integral gives
$\bar\partial_zG_z=-i(2\pi)^{-1}
e^{-2i\kappa r-4iq\kappa y}$.
The calculation for $\kappa<0$ is the same, with
the opposite sign.  Setting
$\zeta(z)=(-2\kappa,-4q\kappa)$, we obtain
\begin{equation}\label{eq:comparison-dbar-green}
 \bar\partial_zG_z(r,y)=
 -\frac{i\,\operatorname{sgn}(\kappa)}{2\pi}
 e^{i(r,y)\cdot\zeta(z)},\qquad \kappa\ne0.
\end{equation}
Since $G_z$ is a convolution kernel, its derivative at two
spatial points $X,X'\in\R^2$ is
\begin{equation}\label{eq:comparison-dbar-green-rank-one}
 \bar\partial_zG_z(X-X')
 =-\frac{i\,\operatorname{sgn}(\kappa)}{2\pi}
 e^{iX\cdot\zeta(z)}e^{-iX'\cdot\zeta(z)}.
\end{equation}
Thus $\bar\partial_zG_z$ is a rank-one operator on functions
supported in a bounded set.

Differentiating \eqref{eq:comparison-direct-equation} with
respect to $\bar z$ and using \eqref{eq:comparison-dbar-green-rank-one} gives
\begin{equation}\label{eq:comparison-dbar-chi-equation}
 (I+G_zV)\bar\partial_z\chi_z
 =-(\bar\partial_zG_z)(V\chi_z)
 =\mathsf{a}(z)e^{i(r,y)\cdot\zeta(z)},
\end{equation}
where
\begin{equation}\label{eq:comparison-direct-a}
 \mathsf{a}(z)=
 \frac{i\,\operatorname{sgn}(\kappa)}{2\pi}
 \int_{\R^2}V(y,r)\chi_z(y,r)
 e^{-i(r,y)\cdot\zeta(z)}\,dr\,dy
\end{equation}
is the scattering coefficient in our normalization.
In first approximation $\chi_z=1$, so
$\mathsf{a}$ is a Fourier transform of $V$;
the exact transform is nonlinear because $\chi_z$
depends on $V$.

\textit{The spectral $\bar\partial$ equation.}
By \eqref{eq:comparison-dbar-chi-equation}, $\bar\partial_z\chi_z$ is $\mathsf{a}(z)$ times the function
$(I+G_zV)^{-1}e^{i(r,y)\cdot\zeta(z)}$,
which can be identified using a symmetry of the Green
function.  The Green function satisfies
$G_{\bar z}=\overline{G_z}$; moreover, the same
Fourier integral, after a shift of its frequency
variable by $2\kappa$, gives
\[
 G_z(r,y)=e^{i(r,y)\cdot\zeta(z)}G_{\bar z}(r,y).
\]
Since $V$ is real, conjugating \eqref{eq:comparison-direct-equation} and using
the uniqueness of its solution gives $\chi_{\bar z}=\overline{\chi_z}$.
Multiplying the equation for $\chi_{\bar z}$
by $e^{i(r,y)\cdot\zeta(z)}$ then gives
\[
 (I+G_zV)
 \bigl(e^{i(r,y)\cdot\zeta(z)}\overline{\chi_z(y,r)}\bigr)
 =e^{i(r,y)\cdot\zeta(z)}.
\]
Comparing this with \eqref{eq:comparison-dbar-chi-equation}, and using the invertibility
of $I+G_zV$, we obtain the spectral $\bar\partial$ equation
\[
 \bar\partial_z\chi_z(y,r)=
 \mathsf{a}(z)e^{i(r,y)\cdot\zeta(z)}
 \overline{\chi_z(y,r)}.
\]

\textit{The inverse $\bar\partial$ problem.}
The second equation of the KP-II Lax pair makes
the evolution of the scattering data explicit.
In the normalization \eqref{KP}, the free spectral phase is
$\Theta(z;t,x,r)=zr+z^2x-\frac t3z^3$.
Let
$P(z;t,x,r)=
e^{\overline{\Theta(z;t,x,r)}-\Theta(z;t,x,r)}$.
At $t=0$ this is
$e^{i(r,x)\cdot\zeta(z)}$, exactly the factor
in the direct spectral equation.

The inverse construction reverses that equation:
given $\mathsf{a}(z)$, find $m(z;t,x,r)\to1$
as $|z|\to\infty$ satisfying
\begin{equation}\label{eq:comparison-inverse-dbar}
 \bar\partial_zm(z;t,x,r)
 =\mathsf{a}(z)P(z;t,x,r)\,
 \overline{m(z;t,x,r)}.
\end{equation}
By the Cauchy--Green formula this is the integral
equation $m=1+\mathcal C(\mathsf{a}P\overline m)$, where
$\mathcal Cf(z)=\pi^{-1}\int_{\mathbb C}
f(w)/(z-w)\,dA(w)$.
Thus the inverse step is a real-linear equation over
the whole complex spectral plane, rather than the
contour problem familiar from KdV.

When $m$ admits an expansion
\[
 m(z)=1+\frac{m^{(1)}}z+O(z^{-2})
\]
as $z\to\infty$ in a specified direction, the normalized spatial
Lax equation is
\[
 (\partial_x-\partial_r^2-2z\partial_r+\phi)m=0.
\]
Substituting the expansion and comparing terms of order $z^0$ gives
\[
 -2\partial_r m^{(1)}+\phi=0,
\]
hence $\phi=2\partial_r m^{(1)}$.
The coefficient $m^{(1)}$ depends on the choice of primitive.
Equivalently, when the integrals and derivatives
are justified,
\begin{equation}\label{eq:comparison-inverse-reconstruction}
 \phi(t,x,r)=\frac2\pi\int_{\mathbb C}\partial_r\!\left(\mathsf{a}(z)P(z;t,x,r)\overline{m(z;t,x,r)}\right)dA(z).
\end{equation}
The derivative is kept inside the integral because
the undifferentiated coefficient integral need not
be absolutely convergent near the real spectral axis.
The full construction requires estimates for the
direct and inverse integral equations and for the
reconstruction; see \cite{Wickerhauser87,GrinevichKP,ChenKP}.
For example, an expansion justified for compactly
supported \emph{spectral} data cannot simply be assumed
from compact support of the initial potential.

The classical method therefore consists of
determining $\mathsf{a}$ from $V$ by
\eqref{eq:comparison-direct-a}, multiplying it by the explicit
spectral phase $P$ defined above, solving the inverse
$\bar\partial$ equation \eqref{eq:comparison-inverse-dbar},
and recovering $\phi$ by \eqref{eq:comparison-inverse-reconstruction}.
The Airy representation avoids these two successive
spectral integral equations: it expresses the
solution directly in terms of $V$.

\textit{The Brownian transform and the scattering coefficient.}
Our construction starts with the original potential rather than
the nonlinear scattering coefficient.  Let
$V\in C_c^\infty(\R^2;\R)$ and choose $a<b$ strictly outside its
support in the transverse variable $y$.  Recall the
defect $D_{a,b;V}$ from \eqref{eq:defect}.  Its two-frequency transform is
\begin{equation}\label{eq:tf-data}
 T_V(k,\ell)=\frac{e^{-ak^2+b\ell^2}}{2\pi}
 \int_{\R^2}e^{iku-i\ell v}
 D_{a,b;V}(u,v)\diff u\diff v,
 \qquad k,\ell\in\mathbb C.
\end{equation}
The ordered Duhamel expansion \eqref{eq:tf-series} proves that
$T_V$ is entire in both variables and independent of $a,b$.
The dependence on the original potential is explicit at each order,
in contrast with the spectral coefficient $\mathsf{a}(z)$ of
\eqref{eq:comparison-direct-a}.

The Fourier representation of the Airy kernels $S_{t,x}$
expresses the Fredholm kernel \eqref{eq:kernel} in terms of $T_V$.
Choose
$\alpha,\beta>0$ such that $t\alpha>x-a$ and $t\beta>b-x$,
and set $\Gamma_-=\R-i\alpha$ and $\Gamma_+=\R+i\beta$, both oriented
from left to right. Then
\[
 \begin{aligned}
 K_{t,x;V}(u,v)=\frac1{2\pi}\int_{\Gamma_-}\int_{\Gamma_+}
 &e^{-itk^3/3+xk^2-iku}\\
 &{}\times e^{it\ell^3/3-x\ell^2+i\ell v}
 T_V(k,\ell)\diff\ell\diff k.
 \end{aligned}
\]
The contour shifts give Gaussian decay in the Fourier integrals,
justifying the substitution of \eqref{eq:tf-data} into \eqref{eq:kernel}.
Thus the Airy determinant of Theorem~\ref{thm:fredholm}
is recovered directly from $T_V$, without solving a
heat-scattering eigenfunction equation or an inverse
$\bar\partial$ problem.

To compare with direct scattering, write $z=q+i\kappa$,
where $q\in\R$ and $\kappa>0$, and shear the potential:
\[
 W_q(y,u)=V(y,u-2qy).
\]
Restrict $T_{W_q}$ to the frequency band $(-\kappa,\kappa)$:
\[
 (\mathsf T_{q,\kappa;V}f)(k)
 =\int_{-\kappa}^{\kappa}T_{W_q}(k,\ell)f(\ell)\,d\ell,
 \qquad f\in L^2(-\kappa,\kappa).
\]
Because the kernel is entire, $\mathsf T_{q,\kappa;V}$ is
trace class.  The nonlinear amplitude is the endpoint value of
its resolvent kernel,
\begin{equation}\label{eq:tf-amp}
 \mathcal A_q(\kappa;V)
 =2\pi\bigl[\mathsf T_{q,\kappa;V}
       (I-\mathsf T_{q,\kappa;V})^{-1}\bigr](\kappa,-\kappa).
\end{equation}
For every real $V\in C_c^\infty(\R^2)$ the determinant obeys
\begin{equation}\label{eq:comparison-finite-band-positive}
 \det(I-\mathsf T_{q,\kappa;V})
 \ge\exp\left\{-\frac{\kappa}{\pi}
              \int_{\R^2}V(y,r)\diff y\diff r\right\}>0.
\end{equation}
Thus the inverse in \eqref{eq:tf-amp} exists without a sign or
smallness assumption.  This quantitative result, proved in
Lemma~\ref{lem:finite-band-invertible}, is stronger than mere
Fredholm invertibility.  The proof approximates the frequency
band by the lowest Dirichlet eigenfunctions on a large interval,
interprets the finite-dimensional determinant as an exponential
moment for stationary nonintersecting Brownian motions, and passes
Jensen's inequality to the whole-line determinant using
trace-norm convergence.

Here is how this finite-band resolvent enters the heat equation.
The shear $u=r+2qy$ changes the normalized direct-scattering
equation \eqref{eq:comparison-direct-equation} to
\[
 (\partial_y-\partial_u^2-2i\kappa\partial_u)
 \chi^q_{i\kappa}=-W_q\chi^q_{i\kappa},
 \qquad
 \chi^q_{i\kappa}(y,u)
       =\chi_{q+i\kappa}(y,u-2qy).
\]
For this operator the Green function is
\begin{equation}\label{eq:tf-Green}
 G_{i\kappa}(u,s)
 =e^{-i\kappa u}\bigg[
   \mathbf1_{\{s>0\}}e^{\kappa^2s}p_s(u)
   -\frac1{2\pi}\int_{-\kappa}^{\kappa}
          e^{imu+(\kappa^2-m^2)s}\,dm\bigg].
\end{equation}
The first term vanishes for $s\le0$, and the second is supported
on a finite Fourier band.  Writing
$h=W_q\chi^q_{i\kappa}$, one obtains
\begin{equation}\label{eq:lt-heat-eqn}
 h=W_q(1-G_{i\kappa}h),
\end{equation}
where $G_{i\kappa}$ acts by convolution in the differences
of the two variables.

\begin{theorem}
\label{thm:lt-heat}
Let $V\in C_c^\infty(\R^2)$ be real,
$q\in\R$ and $\kappa>0$, and let $E$ be a compact rectangle containing the support of $W_q$.
Then equation \eqref{eq:lt-heat-eqn}, interpreted on $L^2(E)$, has a
unique solution $h$, and
\begin{equation}\label{eq:lt-heat-amplitude}
 \int_{\R^2}e^{2i\kappa u}h(y,u)\,dy\,du
       =\mathcal A_q(\kappa;V).
\end{equation}
Undoing the shear gives the normalized eigenfunction
$\chi_{q+i\kappa}$ of \eqref{eq:comparison-direct-equation}.
\end{theorem}

The proof is in Appendix~\ref{app:heat-scattering}.
Undoing the shear in the direct coefficient
\eqref{eq:comparison-direct-a} now yields
\begin{equation}\label{eq:comparison-chen-identification}
 \mathsf{a}(q+i\kappa)=\frac{i}{2\pi}\,
       \mathcal A_q(\kappa;V),\qquad \kappa>0.
\end{equation}

\textit{Comparison of the two reconstructions.}
For real smooth compactly supported $V$, the Airy representation
gives the KP-II solution and the preceding calculation identifies
its Brownian resolvent amplitude with the classical direct
heat-scattering coefficient.  The two reconstructions then give the same
solution in the function classes where the inverse
scattering problem and the KP-II initial-value problem are known to be
well posed.  Formally comparing
$\phi=2\partial_r m^{(1)}$ with
$\phi=-2\partial_r^2\log Z_V$
suggests
\begin{equation}\label{eq:comparison-m1-tau}
 m^{(1)}=-\partial_r\log Z_V,
\end{equation}
with the integration constant fixed by normalization at
$r=+\infty$.  This is an identity when the coefficient expansion,
reconstruction formula and normalization are justified in a common
class. There is also a formal determinant interpretation of the inverse
$\bar\partial$ step, which is worth distinguishing from the theorem
above.  For real scattering data, the reflection relation
$m(\bar z)=\overline{m(z)}$ suggests, after conjugating by the
exponential phase, a spectral Cauchy operator of the schematic form
\[
 L_{\mathsf{a}}(z,w)=
 \frac{\mathsf{a}(\bar w)e^{\Theta(z;t,x,r)-\Theta(\bar w;t,x,r)}}
      {\pi(z-\bar w)}.
\]
Its $r$ derivative has the separated kernel
\[
 \partial_rL_{\mathsf{a}}(z,w)=
 \frac{\mathsf{a}(\bar w)e^{\Theta(z;t,x,r)-\Theta(\bar w;t,x,r)}}{\pi},
\]
and is therefore formally rank one.  Formally, the rank-one derivative suggests a Fredholm determinant relation
\begin{equation}\label{eq:chen-formal-determinant}
 m^{(1)}=-\partial_r\log\det(I-L_{\mathsf{a}}).
\end{equation}
This identity is only formal: we have not established that
$L_{\mathsf{a}}$ is trace class, nor that its Fredholm series
converges.

The essential distinction is therefore not the mere occurrence
of a Fredholm determinant, which also arises in inverse scattering.
It is that in our formula the initial potential enters directly
under an Airy occupation expectation, without first computing
$\mathsf{a}(z)$ and solving \eqref{eq:comparison-inverse-dbar}.
Extensions of the Airy formula beyond compact data, such as
Theorem~\ref{thm:L2}, do not on their own extend the pointwise
scattering identification \eqref{eq:comparison-chen-identification}
to all of $L^2$.

\section{Direct evaluation for line solitons}\label{sec:soliton}\label{app:gram}

\subsection{Direct evaluation of the Airy path integral}\label{subsec:goal}

In this section we ask whether, for soliton initial data, $Z_V$ can be
evaluated directly from the Brownian-path structure of the Airy ensemble,
without using the KP equation.
For example, for the KdV single-soliton potential
\[
 V(x,r)=-2\kappa^2\sech^2(\kappa r),\qquad \kappa>0,
\]
we shall derive the exact identity
\begin{equation}\label{eq:one-example}
 \displaystyle
 \E_{\mathcal A}\exp\left\{2\kappa^2\sum_{i\ge1}\int_\R
 \sech^2\bigl(\kappa(r-\mathcal A_i^{t,x}(y))\bigr)\diff y\right\}
 =1+e^{-2\kappa r+2\kappa^3t/3}.
\end{equation}
Thus we derive the tau function of the KdV soliton at time $t$ directly
from the initial potential $V$, using the probabilistic representation.

\begin{figure}[!htbp]
\centering
\includegraphics{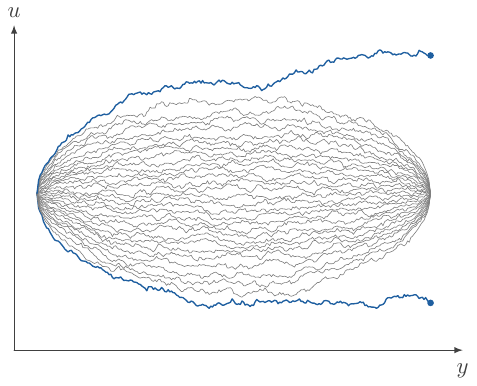}
\caption{$N+1$ nonintersecting Brownian paths with two free endpoints.}
\label{fig:s8-endpoints}
\end{figure}

The idea of the derivation is the following. As in
Section~\ref{sec:dyson-toda-kp}, we approximate the Airy lines by $N$
nonintersecting Brownian bridges, so that $Z_V$ is the limit as $N\to\infty$
of an expectation over $N$ nonintersecting Brownian paths in the potential
$V$. The soliton potential is $V=-2\partial_r^2\log f$, where $f$ is the tau
function of the soliton at time $0$, and a Darboux transformation converts
the $N$ paths in the potential $V$ into $N+1$ nonintersecting Brownian paths
with no potential, weighted by $f$ at their endpoints. All but two of these
endpoints are fixed; the lowest and the highest are free
(Figure~\ref{fig:s8-endpoints}). As $N\to\infty$, the expectation of $f$ at
the lowest free endpoint converges to the tau function at time $t$, and the
highest free endpoint gives its normalization. With $n$ Darboux
transformations, which add $n$ paths and leave the $n$ lowest and the $n$
highest endpoints free, the same idea gives the Wronskian tau functions of
general line solitons.

We begin with one Darboux transformation. Let $p_1<\dots<p_m$ and
$c_1,\dots,c_m>0$, and let
\begin{equation}\label{eq:seed-potential}
 f(x,r)=\sum_{j=1}^m c_je^{p_jr-p_j^2x},\qquad
 V_f(x,r)=-2\partial_r^2\log f(x,r).
\end{equation}
The function $f$ is positive and solves the backward heat equation
$\partial_xf+\partial_r^2f=0$, and $V_f$ is the initial value of the
regular KP line soliton $-2\partial_r^2\log\tau(t,x,r)$, where
\begin{equation}\label{eq:theta}
 \tau(t,x,r)=\sum_{j=1}^m c_je^{\theta(p_j)},\qquad
 \theta(p)=pr-p^2x-\frac{p^3t}{3}.
\end{equation}
Since $f$ is a positive exponential
sum, $V_f\le0$ and is bounded, and it decays exponentially as
$r\to+\infty$ for each fixed $x$.

\begin{theorem}\label{thm:onestep}
Let $V_f$ and $\tau$ be as in \eqref{eq:seed-potential} and
\eqref{eq:theta}. For every $t>0$ and $x,r\in\R$,
\begin{equation}\label{eq:seed-identity}
 \displaystyle
 Z_{V_f}(t,x,r)=
 \frac{\tau(t,x,r)}{c_m e^{\theta(p_m)}}
 =\sum_{j=1}^m\frac{c_j}{c_m}
 e^{\theta(p_j)-\theta(p_m)}.
\end{equation}
\end{theorem}

The denominator $c_me^{\theta(p_m)}$ in \eqref{eq:seed-identity} is the
dominant term of $\tau$ as $r\to+\infty$, so $Z_{V_f}(t,x,r)\to1$.

\begin{proof}[Proof of Theorem~\ref{thm:onestep}]

Let $T=(t^2N)^{1/3}$, as in Section~\ref{sec:dyson-toda-kp}, and let $X_1<\dots<X_N$ be nonintersecting Brownian bridges with diffusion coefficient $2$ on the interval $[a,b]$, all starting and ending at the height $R$, where
\begin{equation}\label{eq:s8-edge}
 a=x-T,\qquad b=x+T,\qquad R=2\sqrt{NT}+r.
\end{equation}
For a potential $V(y,u)$, let
\[
 Z_N(a,b;R;V)=\E\exp\Bigl\{-\sum_{i=1}^N\int_a^bV\bigl(y,X_i(y)\bigr)\diff y\Bigr\}.
\]
In Appendix~\ref{app:s8-bridge} we prove that, for bounded continuous $V$ satisfying $|V(y,u)|\le Ce^{c|y|-\delta u}$ with $C,c,\delta>0$, and in particular for $V=V_f$,
\begin{equation}\label{eq:bridge-limit}
 Z_N(a,b;R;V)\longrightarrow Z_V(t,x,r)\qquad(N\to\infty).
\end{equation}

Next we compute $Z_N(a,b;R;V_f)$ by exploiting the probabilistic
structure of a Darboux transformation.
Write $P_{a,b;Q}(u,v)$ for the variance-$2$ Feynman--Kac
kernel with potential $Q$, and $p_s(v-u)$ for the free heat kernel
of \eqref{eq:heat-kernel}.
Let $\psi>0$ satisfy $(\partial_y+\partial_u^2-Q)\psi=0$, and let
\begin{equation}\label{eq:darboux-potential}
 Q^\psi=Q-2\partial_u^2\log\psi.
\end{equation}
Then the Darboux intertwining identity is
\begin{equation}\label{eq:darboux}
 (\partial_y+\partial_u^2-Q^\psi)
 (\partial_u-\partial_u\log\psi)
 = (\partial_u-\partial_u\log\psi)
 (\partial_y+\partial_u^2-Q).
\end{equation}
The one-particle identity holds for a general base potential.

\begin{lemma}[One path over a base potential]\label{lem:onepath-Q}
Let $a<b$, and let $Q,\psi$ be smooth on $[a,b]\times\R$,
with $\psi>0$ satisfying $(\partial_y+\partial_u^2-Q)\psi=0$.
Assume $Q,Q^\psi$ are bounded and
$\psi(y,u)\le Ce^{C|u|}$. Define the tilted transition density
\begin{equation}\label{eq:tilted-kernel}
 p^{\psi,Q}_{y,b}(u,w)
 =\frac{P_{y,b;Q}(u,w)\psi(b,w)}{\psi(y,u)},
 \qquad a\le y<b.
\end{equation}
Then $\int_\R p^{\psi,Q}_{y,b}(u,w)\diff w=1$ and
\begin{equation}\label{eq:Q-onepath}
 P_{y,b;Q^\psi}(u,v)
 =-\frac{\psi(y,u)}{\psi(b,v)}\partial_u
 \int_{-\infty}^{v}p^{\psi,Q}_{y,b}(u,w)\diff w.
\end{equation}
\end{lemma}
The proof is in Appendix~\ref{app:onepath}.

To see the simplest case, recall from
\eqref{eq:seed-potential} that
$f(x,r)=\sum_j c_j e^{p_jr-p_j^2x}$.
Each exponential satisfies
$(\partial_x+\partial_r^2)e^{p_jr-p_j^2x}=0$,
so $f$ solves the backward heat equation and is strictly positive.
Thus we may take $Q=0$ and $\psi=f$ in the Darboux identity:
the original operator is the free heat operator, while the transformed
potential is
\[
 Q^f=0-2\partial_r^2\log f=V_f.
\]
This is precisely the potential introduced in
\eqref{eq:seed-potential}; the Darboux transformation converts the free
heat equation into the heat equation with potential $V_f$.

For $Q=0$ and $\psi=f$, the tilted kernel becomes
$p^f_{a,b}(u,w)=p_{b-a}(w-u)f(b,w)/f(a,u)$.
Lemma~\ref{lem:onepath-Q} then gives
\begin{equation}\label{eq:one-path}
 P_{a,b;V_f}(u,v)
 =-\frac{f(a,u)}{f(b,v)}\partial_u
 \int_{-\infty}^{v}p^f_{a,b}(u,w)\diff w.
\end{equation}
Integrating these derivatives over interlaced starting
points produces differences at the interval endpoints.

For $u=(u_1<\cdots<u_{N+1})$ and $w=(w_1<\cdots<w_N)$ write
$u\prec w$ if $u_1<w_1<u_2<\cdots<w_N<u_{N+1}$.
Applying Lemma~\ref{lem:onepath-Q}
to each row of a Karlin--McGregor determinant gives the \emph{one-more-path identity}
\begin{equation}\label{eq:one-more-path}
 \int_{u\prec w}\det[P_{a,b;Q^\psi}(w_i,v_j)]_{i,j\le N}
 \prod_{i=1}^N\frac{\diff w_i}{\psi(a,w_i)}
 =\frac{\displaystyle\int_{z\prec v}
 \det[P_{a,b;Q}(u_i,z_j)]_{i,j\le N+1}
 \prod_{j=1}^{N+1}\psi(b,z_j)\diff z_j}
 {\displaystyle\prod_{i=1}^{N+1}\psi(a,u_i)
 \prod_{j=1}^N\psi(b,v_j)}.
\end{equation}
The determinants in \eqref{eq:one-more-path} have a concrete
probabilistic meaning. By the Karlin--McGregor formula~\cite{KM}, the
determinant on the left side is the Feynman--Kac kernel, with potential
$Q^\psi$, of $N$ nonintersecting Brownian paths from $w$ to $v$, and
the determinant on the right side is the Feynman--Kac kernel, with
potential $Q$, of $N+1$ nonintersecting Brownian paths from $u$ to $z$.
We prove \eqref{eq:one-path} and \eqref{eq:one-more-path} in
Appendix~\ref{app:onepath}.

We now let the points $u$ and $v$ in \eqref{eq:one-more-path} coalesce
at $R$, where the bridges of $Z_N(a,b;R;V_f)$ start and end. Take
$Q=0$ and $\psi=f$ in \eqref{eq:one-more-path}, and divide by the same
identity with $f=1$, for which $V_f=0$. As $u_1,\dots,u_{N+1}\to R$,
the interlacing $u\prec w$ forces $w_1,\dots,w_N\to R$. As
$v_1,\dots,v_N\to R$, the interlacing $z\prec v$ forces
$z_2,\dots,z_N\to R$, while the lowest and the highest endpoints remain
free, at $z_1=R-L$ and $z_{N+1}=R+U$ with $L,U>0$. After this division,
the left side tends to $Z_N(a,b;R;V_f)/f(a,R)^N$, as in the proof of
Proposition~\ref{prop:toda}, and the right side tends to
$\E[f(b,R-L)f(b,R+U)]/\bigl(f(a,R)^{N+1}f(b,R)\bigr)$, where the
expectation is over the limiting law of the two free endpoints.
Therefore
\begin{equation}\label{eq:two-ends}
 \displaystyle
 Z_N(a,b;R;V_f)=
 \frac{\E[f(b,R-L)f(b,R+U)]}{f(a,R)f(b,R)}.
\end{equation}
The limits are justified in Appendix~\ref{app:interlacing}.
The pair $(L,U)$ has the explicit density
\begin{equation}\label{eq:LU-density}
 \mathbb P(L\in\diff\ell,U\in\diff u)
 \ \propto\ (\ell+u)(\ell u)^{N-1}
 \exp\left(-\frac{\ell^2+u^2}{8T}\right)
 \diff\ell\diff u,\quad \ell,u>0.
\end{equation}
Indeed, we show in Appendix~\ref{app:interlacing} that $-L$ and $U$ are the two nonzero eigenvalues of
$\left(\begin{smallmatrix}0&\xi\\\xi^*&G\end{smallmatrix}\right)$,
where $\xi$ is an $N$-component complex Gaussian vector and $G$ an
independent real Gaussian of variance $4T$. The variables $U-L$ and
$LU$ are then independent, Gaussian and Gamma, respectively.

We now evaluate the two free endpoints as $N\to\infty$. Let $\mu_N$ be the probability
measure on $(0,\infty)$ with density proportional to
$\eta^{N-1}e^{-\eta^2/(8T)}$, and let $\nu_0$ be its mean.
The joint density \eqref{eq:LU-density} is nearly a product: it is the
density of two independent variables with law $\mu_N$, multiplied by the
factor $(\ell+u)/(2\nu_0)$. This factor tends to
$1$ in the exponentially weighted expectations below, since both
endpoints are concentrated near $2\sqrt{NT}$ on a scale small compared with $\sqrt{NT}$.
The main calculation therefore reduces to exponential moments of $\mu_N$.

The elementary Laplace estimate (proved in Appendix~\ref{app:soliton-proof}) is
\begin{equation}\label{eq:scalar-laplace}
 \log\E_{\mu_N}e^{p\eta}
 =2p\sqrt{NT}+Tp^2+\frac{t}{3}p^3+o(1),
 \qquad
 \nu_p=2\sqrt{NT}+2Tp+tp^2+o(1).
\end{equation}
Here $
\nu_p$ is the mean after exponential tilting by $e^{p\eta}$.
The estimate is uniform for $p$ in compact intervals.
Apply these exponential-moment estimates to the two endpoint weights
$f(b,R-\eta)$ and $f(b,R+\eta)$.
At the \emph{lower} endpoint,
\begin{equation}\label{eq:lower-limit}
 \E_{\mu_N}f(b,R-\eta)
 =\sum_jc_je^{p_jR-p_j^2b}
   \E_{\mu_N}e^{-p_j\eta}
 \longrightarrow\sum_jc_j
   e^{p_jr-p_j^2x-p_j^3t/3}
 =\tau(t,x,r),
\end{equation}
since $R-2\sqrt{NT}=r$ and $b-T=x$. All phases survive.
At the \emph{upper} endpoint, by contrast, $R+\eta$ is large and positive,
so only the largest $p_m$ survives; \eqref{eq:scalar-laplace} gives
\begin{equation}\label{eq:upper-limit}
 \frac{\E_{\mu_N}f(b,R+\eta)}{f(a,R)f(b,R)}
 \longrightarrow\frac{1}{c_m e^{\theta(p_m)}}.
\end{equation}
The corresponding exponentially tilted means are $2\sqrt{NT}+O(T)$,
and $T/\sqrt{NT}\to0$, which justifies discarding the correction factor
$(L+U)/(2\nu_0)$ in the limit.
Equations \eqref{eq:two-ends}--\eqref{eq:upper-limit} and the bridge limit
\eqref{eq:bridge-limit} prove Theorem~\ref{thm:onestep}.
\end{proof}

Three orders of endpoint asymptotics have distinct roles in one Brownian
endpoint law. Its \emph{linear} contribution locates the Airy edge $2\sqrt{NT}$;
its \emph{quadratic} term compensates the heat-equation phase $-p^2b$;
and its \emph{cubic} term produces the KP dispersion $-p^3t/3$.
The full tau function and the normalizing dominant phase come from
opposite free endpoints.

\subsection{Examples: KdV and KP solitons}\label{subsec:examples}\label{sec:s8-examples}

\paragraph{KdV single soliton.}
With $p_1=-\kappa,p_2=\kappa$ and $c_1=c_2=1/2$,
\[
 f(x,r)=e^{-\kappa^2x}\cosh(\kappa r),\qquad
 V_f(x,r)=-2\kappa^2\sech^2(\kappa r).
\]
Theorem~\ref{thm:onestep} gives precisely \eqref{eq:one-example}.
The constant term comes from the largest phase $p_2=\kappa$; the moving
exponential comes from the phase $p_1=-\kappa$.

\paragraph{KP resonant Y soliton.}
For $p_1=-\kappa,p_2=0,p_3=\kappa$ and all $c_j=1$, put
\[
 E_1=e^{-\kappa(r-\kappa x-\kappa^2t/3)},\qquad
 E_2=e^{-2\kappa(r-\kappa^2t/3)}.
\]
Then
\begin{equation}\label{eq:Y}
 Z_{V_f}=1+E_1+E_2,\qquad
 -2\partial_r^2\log Z_{V_f}
 =-\frac{2\kappa^2(E_1+4E_2+E_1E_2)}{(1+E_1+E_2)^2}.
\end{equation}
The dominant pairs among $1,E_1,E_2$ determine the three branches of the
resonant Y. In the Brownian picture the three terms are not three
independent paths: they are the three available exponential phases at the
\emph{same} lower free endpoint.

\paragraph{A KP soliton of arbitrary slope.}
For
\[
 V(x,r)=-2\kappa^2\sech^2\bigl(\kappa(r-\mu x-r_0)\bigr),
\]
take the phases $p_1=\mu/2-\kappa$, $p_2=\mu/2+\kappa$ and
$c_1/c_2=e^{2\kappa r_0}$. The same identity yields
\begin{equation}\label{eq:slope}
 Z_V(t,x,r)=1+\exp\left\{-2\kappa(r-r_0-\mu x)
 +\left(\frac{2\kappa^3}{3}+\frac{\kappa\mu^2}{2}\right)t\right\}.
\end{equation}
This includes the KdV example when $\mu=r_0=0$.

\subsection{General KP solitons: $n$ Darboux steps}\label{subsec:general}\label{sec:s8-several}

A general regular KP line soliton is a Wronskian rather than a single
sum of exponentials. Let
\begin{equation}\label{eq:Wronskians}
 f_i(x,r)=\sum_{j=1}^m c_{ij}e^{p_jr-p_j^2x},\qquad
 W_k(x,r)=\det[\partial_r^{j-1}f_i(x,r)]_{i,j=1}^k,
 \quad W_0=1.
\end{equation}
Suppose for the moment that $W_k>0$ for all $k\le n$. The ratios
\begin{equation}\label{eq:psi-chain}
 \psi_k=\frac{W_k}{W_{k-1}},\qquad
 V_k=-2\partial_r^2\log W_k,
\end{equation}
satisfy
\[
 (\partial_x+\partial_r^2-V_{k-1})\psi_k=0,\qquad
 V_k=V_{k-1}-2\partial_r^2\log\psi_k.
\]
As in \eqref{eq:two-ends}, where one Darboux transformation removed $V_f$, we may therefore remove $V_n$ by $n$ successive Darboux transformations.
At the level of Brownian paths, this changes $N$ paths in the potential
$V_n$ into $N+n$ free paths, after $n$ intertwinings.
When the original endpoints coalesce, $N-n$ of the new endpoints are
pinned at $R$; $n$ are free below $R$ and $n$ above.

To state Proposition~\ref{prop:n-endpoints} below, we introduce the
\emph{divided Wronskian}. For distinct $z_1<\cdots<z_n$, let
\begin{equation}\label{eq:divided-W}
 \mathcal W_x(z_1,\ldots,z_n)
 =C_n\frac{\det[f_i(x,z_j)]_{i,j=1}^n}
 {\Delta(z)},\qquad
 \Delta(z)=\prod_{i<j}(z_j-z_i),\qquad
 C_n=\prod_{j=0}^{n-1}j!.
\end{equation}
The divided Wronskian is continuous when the points coalesce, with
$\mathcal W_x(R,\ldots,R)=W_n(x,R)$. For $n=2$ the intermediate interlacing
integral can already be seen directly:
\begin{equation}\label{eq:two-W}
 \frac{\det[f_i(x,z_j)]_{i,j=1}^{2}}{z_2-z_1}
 =\frac{f_1(x,z_1)f_1(x,z_2)}{z_2-z_1}
 \int_{z_1}^{z_2}\frac{W_2(x,w)}{f_1(x,w)^2}\diff w.
\end{equation}
This identifies the weighted interlacing volume with the divided
Wronskian. Its $n$-level generalization is proved in
Appendix~\ref{app:soliton-proof}.

\begin{proposition}[The $2n$ free endpoints]\label{prop:n-endpoints}
Let $N\ge n$, and let $\xi$ be an $N\times n$ complex
Gaussian matrix with independent entries satisfying
$\E|\xi_{ij}|^2=2(b-a)$. Let $G$ be an independent $n\times n$ GUE
matrix with diagonal variance $2(b-a)$ and
$\E|G_{ij}|^2=2(b-a)$ above the diagonal. The Hermitian matrix
\begin{equation}\label{eq:border}
 M=\begin{pmatrix}0_{N\times N}&\xi\\\xi^*&G\end{pmatrix}
\end{equation}
has $N-n$ zero eigenvalues, $n$ negative eigenvalues
$\lambda^-_1<\cdots<\lambda^-_n$ and $n$ positive eigenvalues
$\lambda^+_1<\cdots<\lambda^+_n$, almost surely. Then
\begin{equation}\label{eq:n-endpoints}
 \displaystyle
 Z_N(a,b;R;V_n)=
 \frac{\E[\mathcal W_b(R+\lambda^-)\mathcal W_b(R+\lambda^+)]}
 {W_n(a,R)W_n(b,R)}.
\end{equation}
\end{proposition}

We prove Proposition~\ref{prop:n-endpoints} in Appendix~\ref{app:interlacing}.

For a general regular KP line soliton, we take the coefficient matrix
$c=(c_{ij})$ in \eqref{eq:Wronskians} to be an $n\times m$ matrix
of rank $n$ with all maximal minors nonnegative. For an $n$-subset
$I\subset\{1,\ldots,m\}$, put
\begin{equation}\label{eq:alpha}
 c_I=\det(c_{ij})_{i\le n,j\in I},\quad
 \alpha_I=c_I\Delta(p_I),\quad
 \theta_I=\sum_{j\in I}\theta(p_j).
\end{equation}
Then the classical tau function is
\[
 \tau_n(t,x,r)=\sum_{|I|=n}\alpha_I e^{\theta_I(t,x,r)}.
\]
At $t=0$ it is the Wronskian of \eqref{eq:Wronskians}: $\tau_n(0,x,r)=W_n(x,r)$.
Among the sets $I$ with $\alpha_I>0$, there is a unique basis $I_*$
maximizing $\sum_{j\in I}p_j$ (the greedy basis of the associated
matroid). The term $\alpha_{I_*}e^{\theta_{I_*}}$ dominates the tau function as $r\to+\infty$.

\begin{theorem}[Direct Airy evaluation for totally nonnegative KP solitons]
\label{thm:general}
Let $c$ be an $n\times m$ matrix of rank $n$ with all maximal minors
nonnegative, and let $V=-2\partial_r^2\log W_n$. Then, for $t>0$ and $x,r\in\R$,
\begin{equation}\label{eq:general-tau}
 \displaystyle
 Z_V(t,x,r)=
 \sum_{|I|=n}\frac{\alpha_I}{\alpha_{I_*}}
 e^{\theta_I-\theta_{I_*}}
 =\frac{\tau_n(t,x,r)}{\alpha_{I_*}e^{\theta_{I_*}}}.
\end{equation}
\end{theorem}

\begin{proof}[Brownian argument]
A change of row basis makes every intermediate Wronskian positive,
so Proposition~\ref{prop:n-endpoints} applies. This positivity
construction is developed in Appendix~\ref{app:positive}.
Expanding the divided Wronskians by Cauchy--Binet, the lower
free endpoints retain all bases $I$ and produce $\tau_n$,
while the upper endpoints select the unique dominant basis
$I_*$ and supply its normalization. The endpoint groups are
not independent at finite $N$: the joint Laplace analysis in
Appendix~\ref{app:asymptotics} proves that their coupling
cancels asymptotically. These facts, together with the
noncompact-potential bridge limit \eqref{eq:bridge-limit},
give \eqref{eq:general-tau}. The complete interlacing
construction and endpoint evaluation are in
Appendix~\ref{app:soliton-proof}.
\end{proof}

We remark that Theorems~\ref{thm:onestep} and~\ref{thm:general} can also
be proved using Theorem~\ref{thm:representation}. The classical soliton tau function
gives an explicit KP solution, so the representation theorem identifies
its normalized tau function with $Z_V$. The point of this section, however,
is to derive these identities directly from
nonintersecting Brownian bridges, explaining both the soliton
coefficients and the normalization.

\paragraph{KdV $n$-solitons.}
Take $0<\kappa_1<\cdots<\kappa_n$ and $\rho_j>0$, and choose
\[
 f_j(x,r)=e^{-\kappa_j^2x}
 \bigl(e^{\kappa_jr}+(-1)^{j-1}\rho_je^{-\kappa_jr}\bigr).
\]
The nonzero Grassmannian coefficients correspond to choosing one of
$\pm\kappa_j$ in each row. If $S$ indexes the chosen negative phases,
Theorem~\ref{thm:general} yields
\begin{equation}\label{eq:KdV-n}
 Z_V(t,x,r)=\sum_{S\subset\{1,\ldots,n\}}\gamma_S
 \prod_{i\in S}e^{-2\kappa_ir+2\kappa_i^3t/3},\quad
 \gamma_S=\prod_{i\in S}\rho_i
 \prod_{i\in S}\prod_{j
\notin S}
 \frac{\kappa_i+\kappa_j}{|\kappa_i-\kappa_j|}.
\end{equation}
\subsection{Airy wanderers}\label{subsec:wanderers}

Given the free endpoints $z_1=R-L$ and $z_{N+1}=R+U$, the $N+1$ free paths
in the proof of Theorem~\ref{thm:onestep} consist of $N-1$ bridges from
$R$ to $R$ and two wanderers, paths that end at $z_1$ and $z_{N+1}$
instead of at $R$. Near the edge of the bridges, the scaling limits of
bridges with wanderers are the Airy processes with wanderers of Adler,
Ferrari and van Moerbeke~\cite{AFvM}, whose line ensembles are constructed
in~\cite{CH,Dimitrov}. With the weight
$f(b,z_1)f(b,z_{N+1})$ of \eqref{eq:two-ends}, the free paths should converge
near the lower edge of the bridges to an Airy line ensemble with one
wanderer, and the lowest free path should converge to the wanderer. Our
computation does not use this limit, so we do not prove it and
only identify, formally, where the wanderer leaves the lower edge of the
bridges.

Near $y=x$ this edge is the parabola $y\mapsto r+(y-x)^2/t$
(Section~\ref{sec:dyson-toda-kp}), and when $f$ has a single phase $p$ the
lowest free path ends, up to an error of order $\sqrt T$, on its tangent
$\ell_p(y)=r+2p(y-x)-p^2t$ at $y=x+pt$. Indeed,
$f(b,z_1)$ is then a multiple of $e^{pz_1}$, and under this weight, by
\eqref{eq:LU-density} and \eqref{eq:scalar-laplace}, the endpoint $z_1=R-L$ has
mean $r+2Tp-tp^2+o(1)=\ell_p(b)+o(1)$ and standard deviation of order
$\sqrt T$. For general $f$ the weight $f(b,z_1)$ is a sum of $m$ such
exponentials, and the lowest free path ends on one of the tangents
$\ell_{p_1},\dots,\ell_{p_m}$, whose heights at $y=b$ differ by order $T$.
By \eqref{eq:LU-density}--\eqref{eq:upper-limit}, the probability that
it ends on $\ell_{p_j}$ tends to $c_je^{\theta(p_j)}/\tau(t,x,r)$.

Formally, in the limit the
wanderer therefore leaves the parabola along $\ell_{p_j}$, at $y=x+p_jt$,
with probability $c_je^{\theta(p_j)}/\tau(t,x,r)$. For the solitons of
Theorem~\ref{thm:general} the $n$ lowest free paths are wanderers, and the
estimates of Appendix~\ref{app:asymptotics} show in the same way that in
the limit they leave the parabola at the points $x+p_it$, $i\in I$, where
the random set $I$ has probability $\alpha_Ie^{\theta_I}/\tau_n(t,x,r)$.

\section{The one-point laws of the KPZ fixed point}\label{sec:kpz-fixed-point}

The KPZ fixed point gives solutions of KP whose initial data are not even functions: formally, the initial condition is an escarpment, $0$ above the graph of a profile $h$ and $+\infty$ below it. The KPZ fixed point is a Markov process on upper semicontinuous profiles, the scaling limit of the KPZ universality class, with transition probabilities given by Fredholm determinants \cite{MQR}. Let $\mathfrak h(t,x;h)$ be the KPZ fixed point started from $h\colon\R\to[-\infty,\infty)$, and let
\[
 F_h(t,x,r)=\mathbb P\{\mathfrak h(t,x;h)\le r\},\qquad\phi_h=-2\partial_r^2\log F_h.
\]
The field $\phi_h$ solves \eqref{KP} for $t>0$, with the escarpment as initial condition \cite[Theorem~1.1 and (1.25)]{QR}. Throughout this section, $h$ is upper semicontinuous, bounded above and not identically $-\infty$.

By the variational formula \cite[Theorem~4.18 and Remark~4.19]{MQR}, and the stationarity and reflection symmetry of the Airy$_2$ process $\bar{\mathcal A}_1$,
\begin{equation}\label{eq:kpz-wall}
 F_h(t,x,r)=\mathbb P\bigl\{r-\mathcal A_1^{t,x}(y)\ge h(y)\text{ for all }y\in\R\bigr\}.
\end{equation}
Formally, the right side is $Z_V$ with $V$ the escarpment; the occupation $I_V^{t,x}(r)$ is then $0$ if every curve $r-\mathcal A_i^{t,x}$ stays above the graph of $h$, and $+\infty$ otherwise, so $Z_V$ is the probability that all the curves stay above the graph. Since the curves are ordered, this is the probability that $r-\mathcal A_1^{t,x}$ does. However, the escarpment takes the value $+\infty$, so it is outside the theory of Section~\ref{sec:direct}, and $\phi_h$ is unbounded as $t\downarrow0$, so Theorems~\ref{thm:representation}(ii) and~\ref{thm:bounded-unique} do not apply.

We show that $\phi_h$, and hence the one-point law $F_h$ of the KPZ fixed point, is nonetheless characterized by KP and its behaviour as $t\downarrow0$.

\begin{theorem}\label{thm:kpz-characterization}
The function $\phi_h$ is the only classical solution of the equation in \eqref{KP} on $(0,\infty)\times\R^2$ that satisfies the following four admissibility conditions.
\begin{enumerate}
\item[\textup{$(\aleph 1)$}] For each compact $J\subset(0,\infty)$ there are $\delta>0$ and, for each $\lambda>0$, a constant $C_\lambda<\infty$ such that $0\le\phi(t,x,r)\le C_\lambda e^{\lambda r_--\delta r_+}$ for $t\in J$ and $x,r\in\R$.
\item[\textup{$(\aleph 2)$}] $\phi(t)\to0$ locally uniformly on $\mathrm{epi}_S(h):=\{(y,r):r>h(y)\}$ as $t\downarrow0$.
\item[\textup{$(\aleph 3)$}] $\int_{y-\varepsilon}^{y+\varepsilon}\phi(t,x,r)\diff x\to\infty$ as $t\downarrow0$, for every $(y,r)\in\mathrm{hypo}_S(h):=\{(y,r):r<h(y)\}$ and every $\varepsilon>0$.
\item[\textup{$(\aleph 4)$}] There is $t_*>0$ such that $\phi(t,x,\cdot)$ is nonincreasing for $0<t<t_*$ and $x\in\R$.
\end{enumerate}
\end{theorem}

Conditions $(\aleph 2)$ and $(\aleph 3)$ encode the initial condition, and $(\aleph 1)$ lets us use the restart identity of Theorem~\ref{thm:representation}(i). In Section~\ref{subsec:kpz-uniqueness} we show that, for every classical solution $\phi$ of KP that satisfies $(\aleph 1)$,
\begin{equation}\label{eq:kpz-chain}
 \exp\{-\tfrac12\partial_r^{-2}\phi(s+t,x,r)\}=Z_{\phi(s)}(t,x,r),\qquad s,t>0,\quad x,r\in\R.
\end{equation}
Now let $s\downarrow0$. The left side tends to $\exp\{-\tfrac12\partial_r^{-2}\phi(t,x,r)\}$. On the right, the occupation $I_{\phi(s)}^{t,x}(r)$ tends to zero when the lowest curve $r-\mathcal A_1^{t,x}$ stays above the graph of $h$, by $(\aleph 2)$, and to infinity when the curve crosses the graph, by $(\aleph 3)$ and $(\aleph 4)$. So the right side tends to the probability that the curve stays above the graph (Lemma~\ref{lem:kpz-wall}), which is $F_h(t,x,r)$ by \eqref{eq:kpz-wall}. Hence $\phi(t)=-2\partial_r^2\log F_h(t)=\phi_h(t)$.

In Section~\ref{subsec:kpz-membership} we show that $\phi_h$ satisfies $(\aleph 1)$. So \eqref{eq:kpz-chain} holds for $\phi_h$, and its left side is $F_h(s+t,x,r)$:

\begin{corollary}\label{cor:kpz-restart}
For $s,t>0$ and $x,r\in\R$, $F_h(s+t,x,r)=Z_{\phi_h(s)}(t,x,r)$.
\end{corollary}

\subsection{Uniqueness}\label{subsec:kpz-uniqueness}

Let $\phi$ be a classical solution of KP that satisfies $(\aleph 1)$. By Lemma~\ref{lem:kpz-smoothing}, $\phi$ is smooth and satisfies the bounds \eqref{eq:representation-bounds} of Theorem~\ref{thm:representation}, so Theorem~\ref{thm:representation}(i) gives \eqref{eq:kpz-chain}. It remains to show that $Z_{\phi(s)}(t,x,r)\to F_h(t,x,r)$ as $s\downarrow0$. This follows from the next lemma, which we prove in Appendix~\ref{app:kpz-wall}.

\begin{lemma}\label{lem:kpz-wall}
Let $t>0$ and $0<s_0<t_*$, and let $\phi\colon(0,s_0]\times\R^2\to[0,\infty)$ be continuous and satisfy $(\aleph 2)$, $(\aleph 3)$ and $(\aleph 4)$. Suppose that for some $C<\infty$ and $\delta>0$,
\begin{equation}\label{eq:kpz-uniform-tail}
 1-Z_{\phi(s)}(t,x,r)\le Ce^{-\delta r},\qquad0<s\le s_0,\quad x,r\in\R.
\end{equation}
Then $Z_{\phi(s)}(t,x,r)\to F_h(t,x,r)$ as $s\downarrow0$, for $x,r\in\R$.
\end{lemma}

\begin{proof}[Proof of uniqueness in Theorem~\ref{thm:kpz-characterization}]
Let $\phi$ satisfy the four admissibility conditions, fix $t>0$ and $0<s_0<t_*$, and write $F(\tau)=\exp\{-\frac12\partial_r^{-2}\phi(\tau)\}$. By the restart identity \eqref{eq:kpz-chain} and $(\aleph 1)$ on $[t,t+s_0]$,
\[
 1-Z_{\phi(s)}(t,x',r')=1-F(s+t,x',r')\le\tfrac12\partial_r^{-2}\phi(s+t,x',r')\le Ce^{-\delta r'}
\]
for $0<s\le s_0$ and $r'\ge0$, and the left side is at most $1$ for $r'<0$. So Lemma~\ref{lem:kpz-wall} applies, and $Z_{\phi(s)}(t,x,r)\to F_h(t,x,r)$. Since also $Z_{\phi(s)}(t,x,r)=F(s+t,x,r)\to F(t,x,r)$, $F(t)=F_h(t)$, and $\phi(t)=-2\partial_r^2\log F(t)=\phi_h(t)$.
\end{proof}

\subsection{Verification of the admissibility conditions}\label{subsec:kpz-membership}

We now show that $\phi_h$ satisfies the four admissibility conditions of Theorem~\ref{thm:kpz-characterization}. The main difficulty is condition $(\aleph 4)$, the monotonicity of $\phi_h$ in $r$. This follows from the following theorem, proved in Appendix~\ref{app:kpz-curvature}.

\begin{theorem}\label{thm:kpz-curvature}
For $t>0$ and $x,r\in\R$, $\partial_r^3\log F_h(t,x,r)\ge0$. Equivalently, $\phi_h(t,x,\cdot)$ is nonincreasing.
\end{theorem}

\begin{proof}[Proof that $\phi_h$ satisfies \textup{$(\aleph 1)$}, \textup{$(\aleph 2)$}, \textup{$(\aleph 3)$} and \textup{$(\aleph 4)$}]
By \cite[Theorems~1.1 and~1.4]{QR}, $F_h$ is positive and real analytic on $(0,\infty)\times\R^2$, and $\phi_h$ solves the equation in \eqref{KP}. Write $\ell=\log F_h$ and $\bar h=\sup h$. Two comparisons with Tracy--Widom distribution functions follow from \eqref{eq:kpz-wall}: for every $y$ with $h(y)>-\infty$,
\begin{equation}\label{eq:kpz-TW}
 F_{\rm GOE}\bigl(4^{1/3}t^{-1/3}(r-\bar h)\bigr)\le F_h(t,x,r)\le F_{\rm GUE}\bigl(t^{-1/3}\bigl(r-h(y)+(x-y)^2/t\bigr)\bigr).
\end{equation}
The lower bound uses Johansson's identity $\mathbb P\{\sup_y(\bar{\mathcal A}_1(y)-y^2)\le m\}=F_{\rm GOE}(4^{1/3}m)$ \cite{Johansson03}, \cite[Theorem~1]{CQR2013}. For the estimates below we need only the nonsharp two-sided cubic bounds
$cS^3\le-\log F_{\rm GOE/GUE}(-S)\le CS^3$ for large $S$,
and a stretched-exponential GOE right-tail bound; these follow, for
example, from the finite-matrix estimates of \cite{LedouxRider2010}
by soft-edge convergence. In particular $\ell\to0$ as $r\to\infty$.

By Theorem~\ref{thm:kpz-curvature}, $\ell_{rr}$ is nondecreasing in $r$. Since $\ell\le0$ and $\ell\to0$ as $r\to\infty$, it follows that $\ell_{rr}\le0$, $\ell_r\ge0$ and $\ell_r\to0$. This gives $(\aleph 4)$ and $\phi_h\ge0$. For $d>0$,
\begin{equation}\label{eq:kpz-differences}
 \phi_h(t,x,r)\le\frac2{d^2}\bigl(-\ell(t,x,r-2d)\bigr),\qquad
 \phi_h(t,x,b)\ge\frac2{d^2}\bigl(-\ell(t,x,b)+2\ell(t,x,b+d)\bigr).
\end{equation}
Both follow from the monotonicity of $\ell_{rr}$, since a second difference of $\ell$ with step $d$ is the integral of $\ell_{rr}$ over a square of side $d$; the first also uses $\ell(r-d)\le\ell(r)\le0$, and the second $\ell(b+2d)\le0$.

\textit{$(\aleph 1)$.} Take $d=1$ in the first inequality of \eqref{eq:kpz-differences} and use the lower bound in \eqref{eq:kpz-TW}. For $t$ in a compact $J\subset(0,\infty)$, the Tracy--Widom tails give $\phi_h\le C(1+r_-)^3e^{-\delta r_+}$, for every $\delta>0$, which implies $(\aleph 1)$.

\textit{$(\aleph 2)$.} Let $(x_0,r_0)\in\mathrm{epi}_S(h)$. Since $h$ is upper semicontinuous, there are $\varepsilon,\rho>0$ with $h\le r_0-4\rho$ on $[x_0-2\varepsilon,x_0+2\varepsilon]$. Let $|x-x_0|\le\varepsilon$ and $|r-r_0|\le2\rho$. The lowest curve $r-\mathcal A_1^{t,x}$ falls below $h$ somewhere on $[x_0-2\varepsilon,x_0+2\varepsilon]$ only if $\sup_y\mathcal A_1^{t,x}(y)>2\rho$, which has probability $1-F_{\rm GOE}(2\cdot4^{1/3}t^{-1/3}\rho)$. For $y\ge x_+=x_0+2\varepsilon$ we have $(y-x)^2\ge(y-x_+)^2+\varepsilon^2$, so $\mathcal A_1^{t,x}(y)\le\mathcal A_1^{t,x_+}(y)-\varepsilon^2/t$, and the curve falls below $h\le\bar h$ there with probability at most $1-F_{\rm GOE}(4^{1/3}t^{-1/3}(r-\bar h+\varepsilon^2/t))$; similarly for $y\le x_0-2\varepsilon$. These probabilities tend to $0$ uniformly, so $\ell\to0$ uniformly on $|x-x_0|\le\varepsilon$, $|r-r_0|\le2\rho$. The first inequality of \eqref{eq:kpz-differences}, with $d=\rho/2$, gives $\phi_h\to0$ uniformly on $|x-x_0|\le\varepsilon$, $|r-r_0|\le\rho$.

\textit{$(\aleph 3)$.} Let $b<h(y)$, choose $c_0>0$ with $c_0^2<(h(y)-b)/2$, and choose $d$ with $b+d>\bar h+1$. For $|x-y|\le c_0\sqrt t$, the upper bound in \eqref{eq:kpz-TW} and the GUE tail give $-\ell(t,x,b)\ge c_1/t$ for small $t$, while the lower bound gives $\ell(t,x,b+d)\to0$ uniformly in $x$. By the second inequality of \eqref{eq:kpz-differences}, $\phi_h(t,x,b)\ge c_2/t$ for $|x-y|\le c_0\sqrt t$ and small $t$, so $\int_{y-\varepsilon}^{y+\varepsilon}\phi_h(t,x,b)\diff x\ge2c_0c_2t^{-1/2}\to\infty$.
\end{proof}

\begin{proof}[Proof of Corollary~\ref{cor:kpz-restart}]
By the proof above, $\phi_h$ satisfies $(\aleph 1)$, so the equality in \eqref{eq:kpz-chain} holds for it. Since $\ell$ and $\ell_r$ tend to $0$ as $r\to\infty$, integrating $\ell_{rr}=-\frac12\phi_h$ twice gives $\exp\{-\frac12\partial_r^{-2}\phi_h\}=F_h$.
\end{proof}

If $h$ is finite only at finitely many spatial points $y_1,\dots,y_m$, then \eqref{eq:kpz-wall} reads
\begin{equation}\label{eq:kpz-finite-sources}
 F_h(t,x,r)=\mathbb P\Bigl\{\bar{\mathcal A}_1(t^{-2/3}y_j)\le t^{-1/3}\Bigl(r-h(y_j)+\frac{(y_j-x)^2}t\Bigr),\ 1\le j\le m\Bigr\},
\end{equation}
a distribution function of the Airy$_2$ process at $m$ points, taken along a line in its arguments. Theorem~\ref{thm:kpz-curvature} says that its logarithm has a nonnegative third derivative along the line. For $m=1$, $F_h$ is a GUE Tracy--Widom distribution function, and $(\log F_{\rm GUE})'''=-2qq'\ge0$, since the Hastings--McLeod solution $q$ of Painlev\'e II is positive and decreasing \cite{HM}.

Theorem~\ref{thm:kpz-curvature} is a limit of the same inequality for $N$ independent Brownian bridges conditioned to be ordered at $n$ given times. In the lemma, $X_{ij}$ is the value of the $i$th bridge at the $j$th time.

\begin{lemma}\label{lem:kpz-gaussian}
Let $\Gamma$ be a symmetric positive definite $n\times n$ matrix with nonpositive off-diagonal entries and nonnegative row sums. Let $X_1,\dots,X_N$ be independent Gaussian vectors with means $\mu_1,\dots,\mu_N\in\R^n$ and the same covariance matrix $\Gamma^{-1}$, and let $c_1,\dots,c_n\in(-\infty,\infty]$. Then
\[
 P(\theta)=\mathbb P\bigl\{X_{1j}\le\theta+c_j\text{ for all }j\bigm|X_{1j}>X_{2j}>\dots>X_{Nj}\text{ for all }j\bigr\}
\]
satisfies $(\log P)'''\ge0$.
\end{lemma}

Both results are proved in Appendix~\ref{app:kpz-curvature}.

\subsection{Lower tails from KP and quadratic initial data}\label{subsec:kpz-lower-tails}

For parabolic and asymptotically quadratic initial profiles, the KP equation,
together with curvature and coarse probabilistic tail bounds, determines the
leading cubic lower tail of the one-point distribution of the KPZ fixed point.
The calculation uses the curvature inequality of
Theorem~\ref{thm:kpz-curvature} and nonsharp cubic bounds for the GOE and GUE
laws, but \emph{not} their known sharp left-tail constants. The argument is a PDE analogue of the method used by
Hastings and McLeod \cite{HM} for Painlev\'e II. Parabolic initial data
reduce KP exactly to a one-dimensional evolution equation. At the flat-data
endpoint, positivity, monotonicity, and coarse tail bounds select the
distinguished far-left scaling, without using the known sharp Tracy--Widom
asymptotics. A weighted stability estimate then extends the identified limit
to general parabolic initial data. We give the argument here and the technical
estimates in Appendix~\ref{app:kpz-quadratic-tail}.

For $a\ge0$ let $h_a(y)=-ay^2$, and write $F_a=F_{h_a}$ and
$\phi_a=-2\partial_r^2\log F_a$.

\begin{theorem}[Parabolic lower tails from KP]\label{thm:kpz-quadratic-tail}
For every $a\ge0$, as $R\to\infty$,
\begin{equation}\label{eq:kpz-lower-tail-quadratic}
-\frac1{R^3}\log F_a(t,\sqrt R\,X,Rz)
\longrightarrow
\frac{1+(1+at)^{-1/2}}{12t}
\left(-z-\frac{aX^2}{1+at}\right)_+^3,
\end{equation}
and
\begin{equation}\label{eq:kpz-lower-tail-phi}
\frac1R\phi_a(t,\sqrt R\,X,Rz)
\longrightarrow
\frac{1+(1+at)^{-1/2}}t
\left(-z-\frac{aX^2}{1+at}\right)_+.
\end{equation}
Both convergences are locally uniform in $(t,X,z)$ for $t>0$,
and locally uniform in $a\ge0$.
\end{theorem}

The theorem concerns the left tail: in the scaling of the theorem,
this is the region $z+aX^2/(1+at)<0$; outside this region the
$R^3$-scaled logarithm has limit zero.

For the main calculation, let $\bar{\mathcal A}_1$ be the stationary
Airy$_2$ process, and set
\[
 G_b(s)=\mathbb P\left\{\sup_y
       \bigl(\bar{\mathcal A}_1(y)-by^2\bigr)\le s\right\},
 \qquad g_b=-2\partial_s^2\log G_b,\quad b\ge1.
\]
Completing the square in the variational formula gives the exact identity
\begin{equation}\label{eq:kpz-parabolic-exact}
 F_a(t,x,r)=G_{1+at}\left(t^{-1/3}
            \left(r+\frac{ax^2}{1+at}\right)\right).
\end{equation}
Thus the dependence on $x$ is known before taking any limit.
Substituting into KP-II gives
\begin{equation}\label{eq:kpz-parabolic-kdv}
 (b-1)\partial_b g_b-\frac s3\partial_sg_b
 -\frac12g_b\partial_sg_b+\frac1{12}\partial_s^3g_b
 +\left(\frac{b-1}{2b}-\frac23\right)g_b=0.
\end{equation}
This holds for $b>1$; at $b=1$ the stationary equation follows
directly from flat-data KP and scaling.

Writing $v_R(b,z)=R^{-1}g_b(Rz)$, the dispersive term in
\eqref{eq:kpz-parabolic-kdv} is multiplied by $R^{-3}$. The candidate limit
\[
 v_*(b,z)=K(b)(-z)_+,\qquad K(b)=1+b^{-1/2},
\]
solves the resulting first-order equation. Its coefficient satisfies
\[
 (b-1)K'=\left(1-\frac{b-1}{2b}\right)K-\frac12K^2.
\]
The point is to prove that the dispersive solutions $v_R$ select $v_*$.
The proof in Appendix~\ref{app:kpz-quadratic-tail} has three steps.
First, the stationary equation at $b=1$, the curvature inequality and
nonsharp cubic estimates select $v_*(1,z)=2(-z)_+$, without using the GOE
cubic constant. Second, positive association of the Airy$_2$ process
controls the entrance uniformly as $b\downarrow1$. Third, a weighted
$L^2$ estimate for the reduced KdV equation propagates the entrance
condition. Smoothing the corner of $v_*$ is essential for this last step;
its residual vanishes in the energy estimate. Integrating the resulting
limit of $g_b$ twice proves \eqref{eq:kpz-lower-tail-quadratic}.
The curvature inequality and monotone second differences then give
\eqref{eq:kpz-lower-tail-phi}.

The parabolic result extends to initial profiles which have only
quadratic \emph{asymptotics} at infinity.

\begin{corollary}[Asymptotically quadratic initial profiles]
\label{cor:kpz-asymptotic-quadratic}
Let $h:\R\to\R$ be upper semicontinuous and locally bounded below,
with $\sup_\R h<\infty$, and suppose that, for some $a\ge0$,
\[
 h(y)=-ay^2+o(y^2)\qquad (|y|\to\infty).
\]
Then the two limits \eqref{eq:kpz-lower-tail-quadratic} and
\eqref{eq:kpz-lower-tail-phi} hold with $F_h,\phi_h$ in place of
$F_a,\phi_a$. They are locally uniform in $(t,X,z)$ for $t>0$.
In particular the leading coefficient is unchanged by asymmetric
or nonsmooth perturbations of size $o(y^2)$.
\end{corollary}

\begin{proof}
Suppose first that $a>0$. For every $0<\delta<a$ there is
$C_\delta<\infty$ such that
\[
 -(a+\delta)y^2-C_\delta
 \le h(y)\le -(a-\delta)y^2+C_\delta,\qquad y\in\R.
\]
The inequalities hold outside a compact set by the hypothesis at
infinity, and on the remaining compact set by local boundedness.
Monotonicity of the KPZ fixed point in its initial data, remembering
that distribution functions have the reverse order, gives
\begin{equation}\label{eq:kpz-quadratic-comparison}
 F_{a-\delta}(t,x,r-C_\delta)
 \le F_h(t,x,r)
 \le F_{a+\delta}(t,x,r+C_\delta).
\end{equation}
Apply Theorem~\ref{thm:kpz-quadratic-tail} with
$x=\sqrt R\,X$ and $r=Rz$. The shifts $C_\delta/R$ disappear,
and the limiting expression is continuous in $a,t,X,z$.
Letting $\delta\downarrow0$ proves the first limit.

If $a=0$, use the upper comparison $h\le M=\sup h$ and the
lower comparison $h(y)\ge-\delta y^2-C_\delta$:
\[
 F_{\rm flat}(t,x,r-M)
 \le F_h(t,x,r)
 \le F_\delta(t,x,r+C_\delta).
\]
Letting $R\to\infty$, then $\delta\downarrow0$, gives the flat
coefficient. Uniformity on compact sets follows from the uniform
parabolic limit and the same comparison. Finally,
$\partial_r^3\log F_h\ge0$ by Theorem~\ref{thm:kpz-curvature};
monotone second-difference bounds applied to the locally uniform
cubic limit give the asserted limit of $\phi_h$, including at the
edge where the limiting second derivative vanishes.
\end{proof}

For the flat profile $a=0$ the coefficient in
\eqref{eq:kpz-lower-tail-quadratic} is $1/(6t)$ and the limiting field
is $2(-z)_+/t$. As $a\to\infty$, the displayed coefficients and
limiting fields approach the narrow-wedge expressions
$1/(12t)$ and $t^{-1}(-z-X^2/t)_+$, respectively.
This interpolation is an identity between the explicit limiting
profiles; it does not require an interchange of the limits
$a\to\infty$ and $R\to\infty$.

Finally, with $u_R(t,X,z)=R^{-1}\phi_h(t,\sqrt R\,X,Rz)$,
the KP-II equation becomes
\begin{equation}\label{eq:kpz-dispersionless-scaled}
 \partial_tu_R-\tfrac12u_R\partial_zu_R
 +\tfrac14\partial_X^2\partial_z^{-1}u_R
 +\frac1{12R^3}\partial_z^3u_R=0.
\end{equation}
The field on the right of \eqref{eq:kpz-lower-tail-phi}
is a distributional solution of the dispersionless KP-II
equation obtained by dropping the last term. In this sense
the theorem identifies the limit selected by the original KP
solution, rather than merely exhibiting a formal solution.

\subsection{Approximate solutions and the restart identity}\label{subsec:kpz-stability}

The proof of uniqueness has a quantitative form. If a function satisfies $(\aleph 2)$, $(\aleph 3)$ and $(\aleph 4)$ and solves KP up to a small residual, then the distribution function it defines is close to $F_h$. As in Proposition~\ref{prop:transport}, $R_\phi$ denotes the residual of $\phi$ in KP.

\begin{theorem}\label{thm:kpz-stability}
Let $T>0$, and let $\phi\colon(0,T)\times\R^2\to[0,\infty)$ have continuous derivatives $\partial_t\phi$ and $\partial_x^j\partial_r^k\phi$, $j+k\le3$, and satisfy $(\aleph 2)$, $(\aleph 3)$ and $(\aleph 4)$. Suppose that for some $\delta>0$ and each compact $J\subset(0,T)$,
\[
 |\partial_t\phi|+\sum_{j+k\le3}|\partial_x^j\partial_r^k\phi|\le C_Je^{-\delta r}\quad\text{on }J\times\R^2,
\]
and that $|R_\phi(t,x,r)|\le\nu(t)e^{-\delta r}$ with $\nu\in L^1(0,T)$. Then, for $0<t<T$ and $x,r\in\R$,
\begin{equation}\label{eq:kpz-stability}
 \bigl|\exp\{-\tfrac12\partial_r^{-2}\phi(t,x,r)\}-F_h(t,x,r)\bigr|\le\frac{e^{\delta^3t/12-\delta r}}{2\delta^2}\int_0^t\nu(s)\diff s.
\end{equation}
\end{theorem}

The proof of uniqueness used the equation only through the equality in \eqref{eq:kpz-chain}. So the one-point laws of the KPZ fixed point are also determined by their restart identity, Corollary~\ref{cor:kpz-restart}, and their behaviour as $t\downarrow0$:

\begin{corollary}\label{cor:kpz-restart-characterization}
The function $\phi_h$ is the unique continuous function on $(0,\infty)\times\R^2$ satisfying conditions $(\aleph 1)$--$(\aleph 4)$ of Theorem~\ref{thm:kpz-characterization} and the following restart identity:
\[
 Z_{\phi(s)}(t,x,r)=\exp\{-\tfrac12\partial_r^{-2}\phi(s+t,x,r)\},\qquad s,t>0,\quad x,r\in\R.
\]
\end{corollary}

We prove Theorem~\ref{thm:kpz-stability} in Appendix~\ref{app:kpz-transport}.

\subsection{The escarpment does not determine the solution}\label{subsec:kpz-thinning}

For a narrow wedge, the limits of $F_h$ above and below the graph are shared by a family of solutions of KP that satisfy $(\aleph 1)$ and $(\aleph 2)$ but not $(\aleph 3)$ and $(\aleph 4)$. Let $h$ be the narrow wedge at $y_0$ of height $h_0$, so that $h(y_0)=h_0$ and $h=-\infty$ elsewhere, and for $0<p\le1$ let
\begin{equation}\label{eq:kpz-thinned}
 F^{(p)}(t,x,r)=\det(I-pK_{\rm Ai})_{L^2(\sigma,\infty)},\qquad\sigma=t^{-1/3}\Bigl(r-h_0+\frac{(x-y_0)^2}t\Bigr),
\end{equation}
and $\phi^{(p)}=-2\partial_r^2\log F^{(p)}$. We call these the thinned laws. Then $F^{(1)}=F_h$, and for $p<1$, $F^{(p)}(t,x,\cdot)$ is the distribution function of the largest point of the Airy point process after each point is kept independently with probability $p$, in the scaling of the narrow wedge. In the path integral, $F^{(p)}(t,x,r)$ is the mean of $(1-p)$ to the power of the number of curves $r-\mathcal A_i^{t,x}$ that pass below the tip $(y_0,h_0)$, since a curve does so exactly when $\bar{\mathcal A}_i(t^{-2/3}y_0)>\sigma$. Reading the occupation of a potential $\delta_{y_0}(x)g(r)$ as $\sum_ig(r-\mathcal A_i^{t,x}(y_0))$, as in Section~\ref{sec:examples-kpz}, this says that $F^{(p)}=Z_{V_p}$ for
\[
 V_p(x,r)=-\log(1-p)\,\delta_{y_0}(x)\mathbf 1_{\{r<h_0\}},
\]
which has the form \eqref{CCRinitial}, moved to the point $y_0$.

\begin{proposition}\label{prop:kpz-thinned}
Let $0<p<1$. Then $\phi^{(p)}$ is a classical solution of the equation in \eqref{KP} on $(0,\infty)\times\R^2$ that satisfies $(\aleph 1)$ and $(\aleph 2)$, and $\phi^{(p)}(t,x,\cdot)$ is not monotone for any $t>0$ and $x\in\R$. As $t\downarrow0$, $F^{(p)}(t)\to1$ locally uniformly on $\mathrm{epi}_S(h)$, $F^{(p)}(t,y_0,r)\to0$ for $r<h_0$, $\phi^{(p)}(t)\to V_p$ in $\mathcal D'(\R^2)$, and
\[
 \int_{y_0-\varepsilon}^{y_0+\varepsilon}\phi^{(p)}(t,x,b)\diff x\longrightarrow-\log(1-p),\qquad b<h_0,\quad\varepsilon>0.
\]
\end{proposition}

So $F^{(p)}$ has the same limits as $F_h$ above and below the graph, and their limits differ only at the tip itself, where $F^{(p)}(t,y_0,h_0)=\det(I-pK_{\rm Ai})_{L^2(0,\infty)}$ for every $t$. But $\phi^{(p)}$ tends to the finite potential $V_p$, and $(\aleph 3)$ fails. The thinned laws fail $(\aleph 4)$ as well, so they leave open whether $(\aleph 3)$ or $(\aleph 4)$ could be dropped while the other is kept. We prove Proposition~\ref{prop:kpz-thinned} in Appendix~\ref{app:kpz-thinning}.

\subsection{Other singular initial conditions: Borodin--Gorin and Cafasso--Claeys--Ruzza}\label{sec:examples-kpz}

\label{sec:kpz}
The examples below have potentials concentrated at one transverse position.
Unlike the infinite potential below a graph used for the KPZ fixed
point, these potentials assign finite height-dependent weights to the
Airy lines at a single point. Consequently the characterization theorem
of this section does not directly apply. Extending the restart-based uniqueness argument to such singular
initial potentials would require an additional limiting argument. On
the other hand, the formulas themselves are known identities; we record
them here as precursor examples of the representation.

Consider the multiplicative stochastic heat equation
\[
 \partial_t\mathcal Z=\tfrac14\partial_x^2\mathcal Z
 +\xi\mathcal Z,\qquad \mathcal Z(0,x)=\delta_0(x),
\]
where $\xi$ is space--time white noise with
$\E[\xi(t,x)\xi(s,y)]=\delta(t-s)\delta(x-y)$. As shown in \cite{QR},
\[
-2\partial_r^2\log\mathbb E_{\rm KPZ}
\exp\{-\tfrac12e^{-r+t/12}\mathcal Z(t,x)\}
\]
satisfies our KP-II equation \eqref{KP}.
The corresponding singular initial potential is
\begin{equation}\label{BGentrance}
 V(x,r)=\delta_0(x)\log(1+e^{-r}).
\end{equation}
Our representation formula gives
\[
\mathbb E_{\rm KPZ}
 \exp\{-\tfrac12e^{-r+t/12}\mathcal Z(t,x)\}
 =\mathbb E_{\mathcal A}
 \exp\left\{-\sum_{i\ge1}
 \log(1+e^{-r+\mathcal A_i^{t,x}(0)})\right\}.
\]
This is the Borodin--Gorin formula \cite{BG}.

\label{sec:CCR}
The deformed Airy kernel determinants of Cafasso--Claeys--Ruzza
\cite{CCR} are a closely related generalization, with the Fermi factor replaced by a general
multiplicative weight.  They correspond to the singular initial potential
\begin{equation}\label{CCRinitial}
 V_\sigma(x,r)=-\delta_0(x)\log(1-\sigma(-r)),
\end{equation}
where $\sigma:\mathbb R\to[0,1]$ satisfies the assumptions of
\cite{CCR}.  The corresponding Airy functional is
\[
 F_\sigma(t,x,r)=\mathbb E_{\mathcal A}\prod_{i\ge1}
 \bigl(1-\sigma(\mathcal A_i^{t,x}(0)-r)\bigr).
\]
Since $\mathcal A_i^{t,x}(0)=t^{1/3}\bar{\mathcal A}_i(0)-x^2/t$,
it depends on $x,r$ only through $s=r+x^2/t$.  Writing
$F_\sigma(t,s)=F_\sigma(t,x,r)$ and
$\Phi_\sigma(t,s)=-2\partial_s^2\log F_\sigma(t,s)$,
the KP equation reduces to cylindrical KdV.

In the notation of \cite{CCR}, their deformed Airy determinant is
\begin{equation}\label{CCRidentity}
 Q_\sigma(\xi,\tau)
 =F_\sigma\left(\tau^{-2},-\frac{\xi}{\tau}\right)
 =Z_{V_\sigma}\left(\tau^{-2},0,-\frac{\xi}{\tau}\right).
\end{equation}
Their field
\begin{equation}\label{CCRu}
 u_\sigma(\xi,\tau)
 =\partial_\xi^2\log Q_\sigma(\xi,\tau)+\frac{\xi}{2\tau}
 =-\frac1{2\tau^2}\Phi_\sigma\left(\tau^{-2},-\frac{\xi}{\tau}\right)
 +\frac{\xi}{2\tau}
\end{equation}
satisfies
\begin{equation}\label{CCRKdV}
 \partial_\tau u_\sigma+2u_\sigma\partial_\xi u_\sigma
 +\frac16\partial_\xi^3u_\sigma=0.
\end{equation}
This follows from the standard change of variables from cylindrical
KdV to KdV.  Thus the solutions of \cite{CCR} arise from the same
one-slice Airy construction as Borodin--Gorin; the parabolic Airy
scaling accounts for both the change of variables and the linear
correction $\xi/(2\tau)$.

\appendix
\numberwithin{equation}{section}

\section{Fredholm identities and analytic estimates}\label{app:fredholm}

\subsection{Lemma~\ref{lem:block} and the finite-mesh determinant identity}\label{app:mesh}

Let $R>0$ be such that $V$ vanishes outside
$(a,b)\times[-R,R]$, and let $M=\sup|V|$.
The proofs of Lemma~\ref{lem:block} and \eqref{eq:finite-mesh-det}
use the following bound and two further auxiliary lemmas.

\begin{lemma}\label{lem:S}
Let $t>0$, $r\in\R$, and $L>0$. There is a constant $C$
such that
\[
 \|S_{t,x}(u,\cdot)\|_{L^2[r,\infty)}\le Ce^{C|u|}
 \qquad(u\in\R,\ |x|\le L).
\]
\end{lemma}

\begin{proof}
By \cite[\S9.7]{DLMF}, $|\Ai(\xi)|\le C_0e^{-\frac23\xi_+^{3/2}}$
for real $\xi$. For $|x|\le L$, the exponential prefactor
in \eqref{eq:S} is at most $e^{2L^3/3t^2+L|v-u|/t}$,
and the Airy argument is at least $t^{-1/3}(v-u)$, so
\[
 |S_{t,x}(u,v)|\le C_1\exp\Bigl\{
 \frac{L|v-u|}{t}-\frac23t^{-1/2}(v-u)_+^{3/2}
 \Bigr\},\qquad C_1=t^{-1/3}C_0e^{2L^3/3t^2}.
\]
In $\|S_{t,x}(u,\cdot)\|^2_{L^2[r,\infty)}$, the values
$v\ge u$ contribute at most a constant depending on $t$ and $L$,
and the values $r\le v<u$ at most
$C_1^2(u-r)_+e^{2L(u-r)_+/t}$.
Hence $\|S_{t,x}(u,\cdot)\|_{L^2[r,\infty)}
\le Ce^{C|u|}$, with $C$ depending on $t$, $r$ and $L$.
\end{proof}

The next lemma bounds the exponential moments of the occupation and its Riemann sums.

\begin{lemma}\label{lem:counts}
For $c,y\in\R$, let $N_c(y)$ be the number of $i$ with $\bar{\mathcal A}_i(y)\ge c$.
\begin{enumerate}
\item[\textup{(i)}] $\E N_c(y)$ is finite and does not depend on $y$, and $\E e^{\theta N_c(y)}\le\exp\{(e^\theta-1)\E N_c(0)\}$ for $\theta\ge0$.
\item[\textup{(ii)}] Let $c=t^{-1/3}(r-R)$. For $p\ge0$,
\[
 \E e^{pI_{|V|}}\le\exp\bigl\{(e^{pM(b-a)}-1)\E N_c(0)\bigr\},
\]
and the same bound holds for the Riemann sum $I_m$ with $|V|$ in place of $V$.
\end{enumerate}
\end{lemma}

\begin{proof}
(i) Let $K_c$ be the operator on $L^2[c,\infty)$ with kernel $K_{\rm Ai}(z,z')$. For $\theta\ge0$,
\[
 \E e^{\theta N_c(y)}=\det\bigl(I+(e^\theta-1)K_c\bigr)\le e^{(e^\theta-1)\operatorname{tr}K_c},\qquad
 \operatorname{tr}K_c=\int_c^\infty\rho_{\rm Ai}(z)\diff z=\E N_c(y).
\]
The equality follows from \eqref{eq:airy-correlations}, by expanding $e^{\theta N_c(y)}=\prod_i\bigl(1+(e^\theta-1)\mathbf 1_{[c,\infty)}(\bar{\mathcal A}_i(y))\bigr)$. The inequality holds because $K_c$ is positive by \eqref{eq:extended-airy-S}, and $\det(I+A)\le e^{\operatorname{tr}A}$ for positive trace class $A$. The trace does not depend on $y$, and it is finite since $\rho_{\rm Ai}(z)$ decays like $e^{-\frac43z^{3/2}}$ as $z\to+\infty$.

(ii) We first show that
\begin{equation}\label{eq:counts-pathwise}
 I_{|V|}\le M\int_a^bN_c(t^{-2/3}y)\diff y.
\end{equation}
By \eqref{eq:markov-airy-marginal} we may take $r-\mathcal A_i^{t,x}(y)=r+(y-x)^2/t-t^{1/3}\bar{\mathcal A}_i(t^{-2/3}y)$. If $\bar{\mathcal A}_i(t^{-2/3}y)<c$, then $t^{1/3}\bar{\mathcal A}_i(t^{-2/3}y)<r-R$, so $r-\mathcal A_i^{t,x}(y)>R$ and $V(y,r-\mathcal A_i^{t,x}(y))=0$. Hence at each $y$, at most $N_c(t^{-2/3}y)$ terms of $\sum_i|V(y,r-\mathcal A_i^{t,x}(y))|$ are nonzero, each is at most $M$, and all of them vanish unless $a<y<b$. Therefore, by \eqref{eq:counts-pathwise}, Jensen's inequality and (i),
\[
\begin{aligned}
 \E e^{pI_{|V|}}&\le\E\exp\Bigl\{pM\int_a^bN_c(t^{-2/3}y)\diff y\Bigr\}\\
 &\le\frac1{b-a}\int_a^b\E e^{pM(b-a)N_c(t^{-2/3}y)}\diff y\le\exp\bigl\{(e^{pM(b-a)}-1)\E N_c(0)\bigr\}.
\end{aligned}
\]
For the Riemann sum $I_m$, at $y=y_j$ the sum over $i$ is likewise at most $MN_c(t^{-2/3}y_j)$, and Jensen's inequality gives
\[
 \E\exp\Bigl\{pM\ell\sum_{j=1}^mN_c(t^{-2/3}y_j)\Bigr\}\le\frac1m\sum_{j=1}^m\E e^{pMm\ell N_c(t^{-2/3}y_j)}\le\exp\bigl\{(e^{pM(b-a)}-1)\E N_c(0)\bigr\},
\]
where the last inequality uses $m\ell\le b-a$ and (i).
\end{proof}

We also use two analytic facts in the proof of Lemma~\ref{lem:block}.

\begin{lemma}\label{lem:block-analytic}
Let $g$, $\widehat K$ and $G_{ij}$ be as in Lemma~\ref{lem:block} and its proof.
\begin{enumerate}
\item[\textup{(i)}] $g\widehat K$ is a trace class operator on $L^2[-R,R]^m$, and $\det(I-g\widehat K)$ is given by its Fredholm series.
\item[\textup{(ii)}] $\E\prod_{j=1}^m\prod_{i\ge1}\bigl(1+|G_{ij}|\bigr)<\infty$.
\end{enumerate}
\end{lemma}

\begin{proof}
(i) In the operator notation of Remark~\ref{rem:operator-form},
\begin{equation}\label{eq:Khat-blocks}
 \widehat K_{jk}=\mathbf 1_{[-R,R]}\Bigl(-\mathbf 1_{j<k}\,e^{(y_k-y_j)\partial^2}+S_{t,x-y_j}\,\chi_r\,S^*_{t,y_k-x}\Bigr)\mathbf 1_{[-R,R]}.
\end{equation}
By Lemma~\ref{lem:S}, $\mathbf 1_{[-R,R]}S_{t,x-y_j}\chi_r$ and $\chi_rS^*_{t,y_k-x}\mathbf 1_{[-R,R]}$ are Hilbert--Schmidt, so the second term above is trace class. For $j<k$, the first term is $-\bigl(\mathbf 1_{[-R,R]}e^{s\partial^2/2}\bigr)\bigl(e^{s\partial^2/2}\mathbf 1_{[-R,R]}\bigr)$ with $s=y_k-y_j>0$, again a product of Hilbert--Schmidt operators. Thus $\widehat K$, and hence $g\widehat K$, is trace class. Since its kernel is continuous, the determinant is given by the usual Fredholm series.

(ii) We have
\[
 \prod_{j=1}^m\prod_{i\ge1}\bigl(1+|G_{ij}|\bigr)\le\exp\Bigl\{\ell\sum_{j=1}^m\sum_{i\ge1}\bigl|V\bigl(y_j,r-\mathcal A_i^{t,x}(y_j)\bigr)\bigr|\Bigr\},
\]
since $G_{ij}=1-e^{-\ell V(y_j,r-\mathcal A_i^{t,x}(y_j))}$ and $1+|1-e^{-z}|\le e^{|z|}$. The right side is $e^{I_m}$ with $|V|$ in place of $V$, which has finite expectation by Lemma~\ref{lem:counts}(ii) with $p=1$.
\end{proof}

\begin{proof}[Analytic details for \eqref{eq:finite-mesh-det}]
The operators $\Phi$ and $\Psi$ introduced in the main text are
Hilbert--Schmidt by Lemma~\ref{lem:S}. Each nonzero block of $H$
is a heat semigroup compressed to $[-R,R]$, hence is trace class
by factorization into two Hilbert--Schmidt operators. The
block-triangular identities and the Fredholm determinant
manipulations in \eqref{eq:mesh-reduced} are therefore valid.

It remains to check that the formal heat chains used in
\eqref{eq:chains-matched} define the asserted kernels. Each $g_j$
is supported in $[-R,R]$, and the heat factors on either side
of a chain have strictly positive time. Its kernel is consequently
bounded by $C_m e^{-c(w^2+w'^2)}$ for some $C_m,c>0$.
This Gaussian decay dominates the at-most-exponential
growth in $w,w'$ of the outer Airy kernels $S_{t,a-x}$
and $S_{t,x-b}$, so their kernel integrals converge
absolutely. The finite sum defining $D_m$ has the same
bound. Finally, $\Psi(I+gH)^{-1}g\Phi$ is trace class as a
product of Hilbert--Schmidt and bounded operators, and the
identification with $\chi_rK_m\chi_r$ holds in trace class.
\end{proof}

\subsection{Proof of Proposition~\ref{prop:limit}}\label{app:limit}

We first show that the defects $D_m$ converge to $D_V$, with a Gaussian bound uniform in $m$.

\begin{lemma}\label{lem:gauss}
There are constants $C,c>0$ such that, for all $m\ge1$ and $w,w'\in\R$,
\[
 |D_m(w,w')|+|D_V(w,w')|\le Ce^{-c(w^2+w'^2)},
\]
and $D_m(w,w')\to D_V(w,w')$ as $m\to\infty$.
\end{lemma}

\begin{proof}
Let $\E^{w,w'}$ denote expectation over a Brownian bridge $B$ with diffusion coefficient $2$ from $B(a)=w$ to $B(b)=w'$. Conditioning on $B(b)$ in \eqref{eq:FK} gives
\[
 D_V(w,w')=p_{b-a}(w'-w)\E^{w,w'}\Bigl[1-e^{-\int_a^bV(y,B(y))\diff y}\Bigr],
\]
and the same holds for $D_m$ with $\ell\sum_jV(y_j,B(y_j))$ in place of the integral. We prove the Gaussian bound for $D_V$:
\[
 \begin{aligned}
 |D_V(w,w')|&\le e^{M(b-a)}p_{b-a}(w'-w)\E^{w,w'}\int_a^b|V(y,B(y))|\diff y\\
 &=e^{M(b-a)}\int_a^b\int_\R p_{y-a}(\eta-w)|V(y,\eta)|p_{b-y}(w'-\eta)\diff\eta\diff y\le Ce^{-c(w^2+w'^2)}.
 \end{aligned}
\]
The first inequality holds because $|1-e^{-z}|\le|z|e^{z_-}$ for real $z$, and $|V|\le M$. The last inequality holds because on the support of $V$ the times $y-a$ and $b-y$ lie in $[\delta,b-a]$ for some $\delta>0$ and $|\eta|\le R$, so $p_{y-a}(\eta-w)\le Ce^{-cw^2}$ and $p_{b-y}(w'-\eta)\le Ce^{-cw'^2}$ there. The same argument, with $\int_a^b\diff y$ replaced by $\ell\sum_j$ at $y=y_j$, proves the bound for $D_m$, since $m\ell<b-a$.

Finally, $D_m(w,w')\to D_V(w,w')$ by dominated convergence, since along every continuous path the Riemann sums $\ell\sum_jV(y_j,B(y_j))$ converge to $\int_a^bV(y,B(y))\diff y$ and are bounded by $M(b-a)$.
\end{proof}

\begin{proof}[Proof of Proposition~\ref{prop:limit}]
\noindent\textit{Step~1: The kernels.} First, the integral \eqref{eq:kernel} converges absolutely for all $u,v$, by Lemma~\ref{lem:gauss} and the bound $|S_{t,\alpha}(w,u)|\le C_{\alpha,u}e^{|\alpha||w|/t}$ from \eqref{eq:S}. Next we show that $\chi_rK_V\chi_r$ is trace class and that $\|\chi_rK_m\chi_r-\chi_rK_V\chi_r\|_1\to0$. For $w\in\R$ let $\phi_w=\chi_rS_{t,a-x}(w,\cdot)$ and $\psi_w=\chi_rS_{t,x-b}(w,\cdot)$. By Lemma~\ref{lem:S}, $\phi_w$ and $\psi_w$ lie in $L^2(\R)$ with norms at most $Ce^{C|w|}$, and by dominated convergence they depend continuously on $w$ in $L^2(\R)$. The operator with kernel $\phi_w(u)\psi_{w'}(v)$ has rank one and trace norm $\|\phi_w\|\|\psi_{w'}\|$, and by Lemma~\ref{lem:gauss},
\[
 \iint|D_V(w,w')|\|\phi_w\|\|\psi_{w'}\|\diff w\diff w'\le C\iint e^{-c(w^2+w'^2)+C(|w|+|w'|)}\diff w\diff w'<\infty.
\]
Hence the integral \eqref{eq:kernel}, with $\chi_r$ on both sides, converges in trace norm, and $\chi_rK_V\chi_r$ is trace class. The same holds with $D_m$ in place of $D_V$, and
\[
 \|\chi_r(K_m-K_V)\chi_r\|_1\le\iint|D_m(w,w')-D_V(w,w')|\|\phi_w\|\|\psi_{w'}\|\diff w\diff w'.
\]
The right side tends to zero by dominated convergence and Lemma~\ref{lem:gauss}. The Fredholm determinant is continuous in trace norm, so $\det(I-\chi_rK_m\chi_r)\to\det(I-\chi_rK_V\chi_r)$.

\noindent\textit{Step~2: The expectations.} By Lemma~\ref{lem:counts}(ii) with $p=2$, $\E e^{2I_{|V|}}<\infty$, and the same bound holds for $I_m$ with $|V|$ in place of $V$, for all $m$. In particular, the series defining $I_V$ converges absolutely almost surely, so $Z_V>0$. Since $|I_V|\le I_{|V|}$, and likewise for $I_m$, the variables $e^{-I_V}$ and $e^{-I_m}$ are bounded in $L^2$ uniformly in $m$. Hence $Z_V<\infty$, and the $e^{-I_m}$ are uniformly integrable, so $Z_m\to Z_V$ once $I_m\to I_V$ in probability. We show that $I_m\to I_V$ in $L^1$.

Fix $n\ge1$, and split $I_m$ and $I_V$ into the contributions of the lines $i\le n$ and $i>n$. For each $i\le n$, the Riemann sums $\ell\sum_jV(y_j,r-\mathcal A_i^{t,x}(y_j))$ converge almost surely to $\int_a^bV(y,r-\mathcal A_i^{t,x}(y))\diff y$, since the Airy lines and $V$ are continuous. These sums and integrals are at most $M(b-a)$ in absolute value, so by bounded convergence the contributions of the lines $i\le n$ converge in $L^1$. By the proof of Lemma~\ref{lem:counts}(ii), at each $y$ only the top $N_c(t^{-2/3}y)$ lines have $V(y,r-\mathcal A_i^{t,x}(y))\ne0$, so at most $(N_c(t^{-2/3}y)-n)_+$ of the lines $i>n$ do. By stationarity, the contributions of the lines $i>n$ to $I_m$ and to $I_V$ therefore have expected absolute value at most $M(b-a)\E(N_c(0)-n)_+$. Hence
\[
 \limsup_{m\to\infty}\E|I_m-I_V|\le2M(b-a)\E(N_c(0)-n)_+,
\]
which tends to zero as $n\to\infty$ by dominated convergence, since $\E N_c(0)<\infty$ by Lemma~\ref{lem:counts}(i).
\end{proof}

\subsection{The ordered Feynman--Kac expansion}\label{app:remarks}

Expanding the exponential in \eqref{eq:FK} in powers of $V$ gives $K_{t,x;V}=\sum_{n\ge1}(-1)^{n+1}K^{[n]}_{t,x;V}$, where
\begin{equation}\label{eq:ordered-kernel}
 \begin{aligned}
 K^{[n]}_{t,x;V}(u,v)=\int_{a<y_1<\cdots<y_n<b}&\int_{\R^n}S_{t,y_1-x}(w_1,u)S_{t,x-y_n}(w_n,v)\\
 &\times\prod_{j=1}^nV(y_j,w_j)\prod_{j=1}^{n-1}p_{y_{j+1}-y_j}(w_{j+1}-w_j)\diff\mathbf w\diff\mathbf y.
 \end{aligned}
\end{equation}
If $|V|\le M$ and $V$ vanishes outside $(a,b)\times[-R,R]$, then $\|\chi_rK^{[n]}_{t,x;V}\chi_r\|_1\le2CR(M(b-a))^n/n!$ with $C$ independent of $n$ and $M$, so $\sum_{n\ge1}(-1)^{n+1}\chi_rK^{[n]}_{t,x;V}\chi_r$ converges to $\chi_rK_{t,x;V}\chi_r$ in trace norm.

\begin{proof}[Proof of the expansion]
Expanding the exponential in \eqref{eq:FK} in powers of $V$ and using the Markov property writes $D_{a,b;V}$ as a series whose $n$th term is at most $Ce^{-c(w^2+w'^2)}(M(b-a))^n/n!$ in absolute value, as in the proof of Lemma~\ref{lem:gauss}. So the series may be substituted into \eqref{eq:kernel}, and integrating its $n$th term over $w$ and $w'$ by the heat semigroup identity gives $K^{[n]}_{t,x;V}$. For the trace-norm bound, Lemma~\ref{lem:S} bounds the $L^2[r,\infty)$ norms of $S_{t,y_1-x}(w_1,\cdot)$ and $S_{t,x-y_n}(w_n,\cdot)$ in \eqref{eq:ordered-kernel} by a constant $C_1$ where $V(y_1,w_1)$ and $V(y_n,w_n)$ are nonzero. So, as in Step~1 of the proof of Proposition~\ref{prop:limit}, $\|\chi_rK^{[n]}_{t,x;V}\chi_r\|_1$ is at most $C_1^2$ times the integral of $\prod_j|V(y_j,w_j)|\prod_jp_{y_{j+1}-y_j}(w_{j+1}-w_j)$, which is at most $2R(M(b-a))^n/n!$.
\end{proof}

\subsection{Proof of Proposition~\ref{prop:qr-kernel}}\label{app:QR-analytic}

We first extend the bound of Lemma~\ref{lem:S} to complex parameters, using the principal branch of $t^{1/3}$.

\begin{lemma}\label{lem:S-complex}
Let $t_0>0$ and $x_0,r_0\in\R$. Then there is a complex neighborhood $U$ of $(t_0,x_0,r_0)$ such that, for each $j\ge0$,
\[
 |\partial_u^jS_{t,\alpha}(w,r+u)|\le C_je^{C_j|w|^{3/2}-u^{3/2}/C_j}
\]
for all $(t,x,r)\in U$, $\alpha\in\{a-x,x-b\}$, $w\in\R$ and $u\ge0$.
\end{lemma}

\begin{proof}
By \eqref{eq:S}, $\partial_u^jS_{t,\alpha}(w,r+u)$ is a linear combination of the terms
\[
 (\alpha/t)^{j-k}t^{-(k+1)/3}e^{\frac{2\alpha^3}{3t^2}}e^{\frac{\alpha(r+u-w)}t}\Ai^{(k)}(\xi),\qquad\xi=t^{-1/3}(r+u-w)+t^{-4/3}\alpha^2,
\]
for $0\le k\le j$. Take $U$ bounded, with $|t|\ge t_0/2$ and $|\arg t|\le\pi/4$ on $U$. Then on $U$ the product $(\alpha/t)^{j-k}t^{-(k+1)/3}e^{2\alpha^3/(3t^2)}$ is bounded, $|e^{\alpha(r+u-w)/t}|\le Ce^{C(u+|w|)}$, and $|\arg t^{-1/3}|\le\pi/12$. By \cite[\S9.7]{DLMF}, $|\Ai^{(k)}(\xi)|\le C(1+|\xi|)^Ce^{\frac23|\xi|^{3/2}}$ for all $\xi\in\C$, and $|\Ai^{(k)}(\xi)|\le C(1+|\xi|)^Ce^{-|\xi|^{3/2}/C}$ for $|\arg\xi|\le\pi/6$. If $u\ge2|w|+K$ for a large constant $K$, then $|\arg\xi|\le\pi/6$ and $|\xi|\ge u/C$, so each term is at most $Ce^{Cu-u^{3/2}/C}$. If $u<2|w|+K$, then $|\xi|\le C(1+|w|)$, so each term is at most $Ce^{C(1+|w|)^{3/2}}$. Summing over $k$, both cases give the bound of the lemma.
\end{proof}

\begin{proof}[Proof of Proposition~\ref{prop:qr-kernel}(i) and (iii)]
Fix $(t_0,x_0,r_0)$ with $t_0>0$, and let $U$ be as in Lemma~\ref{lem:S-complex}. For $w\in\R$ and $(t,x,r)\in U$ let $\phi_w(u)=S_{t,a-x}(w,r+u)$ and $\psi_w(v)=S_{t,x-b}(w,r+v)$ for $u,v>0$. By Lemma~\ref{lem:S-complex}, $\phi_w$ and $\psi_w$ lie in $L^2(\R_+)$ with norms at most $Ce^{C|w|^{3/2}}$. Since they are holomorphic in $(t,x,r)$ pointwise, they are holomorphic as $L^2(\R_+)$-valued functions. For $f,g\in L^2(\R_+)$, write $f\otimes g$ for the operator with kernel $f(u)g(v)$; it has rank one and $\|f\otimes g\|_1=\|f\|\|g\|$.

(i) By Lemma~\ref{lem:gauss},
\[
 \iint|D_{a,b;V}(w,w')|\|\phi_w\|\|\psi_{w'}\|\diff w\diff w'\le C\iint e^{-c(w^2+w'^2)+C(|w|^{3/2}+|w'|^{3/2})}\diff w\diff w'<\infty,
\]
uniformly on $U$. So, as in Step~1 of the proof of Proposition~\ref{prop:limit}, the integral $\iint D_{a,b;V}(w,w')\phi_w\otimes\psi_{w'}\diff w\diff w'$ converges in trace norm, and for real $(t,x,r)$ its kernel is $\mathcal K$. Its integrand is holomorphic on $U$, since $\phi_w$ and $\psi_{w'}$ are, so the integral is holomorphic on $U$.

(iii) Let $(t,x)$ range over a compact set $\Lambda\subset(0,\infty)\times\R$, and let $r\ge r_0$ be real. By (ii), $\partial_t^j\partial_x^k\partial_r^l\mathcal K$ is a finite linear combination of integrals $\iint D_{a,b;V}(w,w')\phi_w^{(p)}\otimes\psi_{w'}^{(q)}\diff w\diff w'$, each of trace norm at most $\iint|D_{a,b;V}(w,w')|\|\phi_w^{(p)}\|\|\psi_{w'}^{(q)}\|\diff w\diff w'$. The norm $\|\phi_w^{(p)}\|=\|\partial_u^pS_{t,a-x}(w,\cdot)\|_{L^2[r,\infty)}$ decreases in $r$ and tends to zero as $r\to\infty$, uniformly in $(t,x)\in\Lambda$ by Dini's theorem. By Lemma~\ref{lem:S-complex} on a finite cover of $\Lambda\times\{r_0\}$, this norm is at most $Ce^{C|w|^{3/2}}$ for $r\ge r_0$. The same holds for $\psi_{w'}^{(q)}$, so (iii) follows by dominated convergence, as in (i).
\end{proof}

\section{Airy occupation: intensity, moments, and concentration}
\label{app:section2}

\subsection{Airy intensity and occupation bounds}\label{app:occupation}

We begin by writing the occupation $I_V^{t,x}(r)$ in terms of the lines $\bar{\mathcal A}$ of \eqref{eq:markov-airy-marginal}. Since $\bar{\mathcal A}$ is stationary, we may take the lines in \eqref{eq:airy-occupation} to be $\mathcal A_i^{t,x}(y)=t^{1/3}\bar{\mathcal A}_i(t^{-2/3}(y-x))-(y-x)^2/t$. Substituting $y=x+t^{2/3}s$ gives
\begin{equation}\label{eq:occ-scaled}
 I_V^{t,x}(r)=t^{2/3}\sum_{i\ge1}\int_\R V\bigl(x+t^{2/3}s,r+t^{1/3}(s^2-\bar{\mathcal A}_i(s))\bigr)\diff s.
\end{equation}
With $s=t^{-1/6}v$ and $f_v(z)=V(x+v\sqrt t,r+v^2-t^{1/3}z)$, this becomes
\begin{equation}\label{eq:occ-sections}
 I_V^{t,x}(r)=\sqrt t\int_\R Y_v\diff v,\qquad Y_v=\sum_{i\ge1}f_v\bigl(\bar{\mathcal A}_i(t^{-1/6}v)\bigr).
\end{equation}
Minkowski's integral inequality, applied to $Y_v$ and to $Y_v-\E Y_v$, gives, for $1\le q<\infty$,
\begin{equation}\label{eq:occ-minkowski}
 \|I_V^{t,x}(r)\|_q\le\sqrt t\int_\R\|Y_v\|_q\diff v,\qquad\|I_V^{t,x}(r)-\E I_V^{t,x}(r)\|_2\le\sqrt t\int_\R(\operatorname{Var}Y_v)^{1/2}\diff v.
\end{equation}
By stationarity, $Y_v$ has the law of the linear statistic $\sum_if_v(\bar{\mathcal A}_i(0))$ of the Airy point process $\{\bar{\mathcal A}_i(0)\}_{i\ge1}$, so the moment bounds for the occupation below reduce to bounds on linear statistics of this point process.

\subsubsection*{Intensity and determinantal statistics}

The intensity $\rho_{\rm Ai}(z)=K_{\rm Ai}(z,z)=\int_0^\infty\Ai(z+v)^2\diff v$ of the Airy point process has Laplace transform
\begin{equation}\label{app:mu}
 \mu(q)=\int_\R e^{qz}\rho_{\rm Ai}(z)\diff z=\frac{e^{q^3/12}}{2\sqrt{\pi}q^{3/2}},\qquad q>0,
\end{equation}
by Tonelli's theorem and $\int_\R e^{qz}\Ai(z)^2\diff z=e^{q^3/12}/(2\sqrt{\pi q})$ \cite[Lemma~2.6]{Okounkov}. For $0<t\le1$,
\begin{equation}\label{eq:airy-density-scaled-bound}
 t^{1/6}\rho_{\rm Ai}(t^{-1/3}q)\le C(1+\sqrt{(-q)_+})e^{-q_+},\qquad q\in\R,
\end{equation}
since, by the Airy asymptotics \cite[\S9.7]{DLMF}, there is a $c_0>0$ with
\begin{equation}\label{eq:airy-density-bound}
 \rho_{\rm Ai}(z)\le C\bigl(1+\sqrt{(-z)_+}\bigr)e^{-c_0z_+^{3/2}},\qquad z\in\R.
\end{equation}

The next lemma bounds the moments, exponential moments and variance of linear statistics $\sum_if(\bar{\mathcal A}_i(0))$ by those of a Poisson process with intensity $\rho_{\rm Ai}$. Let $\Pi_q$ be the set of partitions of $\{1,\dots,q\}$.

\begin{lemma}\label{lem:one-position}
For measurable $f\ge0$, integers $q\ge1$ and $\theta\ge0$,
\begin{equation}\label{eq:occ-one-position}
 \E\Bigl(\sum_if(\bar{\mathcal A}_i(0))\Bigr)^q\le\sum_{\pi\in\Pi_q}\prod_{J\in\pi}\int_\R f^{|J|}\rho_{\rm Ai}\diff z,\qquad
 \E e^{\theta\sum_if(\bar{\mathcal A}_i(0))}\le\exp\Bigl\{\int_\R(e^{\theta f}-1)\rho_{\rm Ai}\diff z\Bigr\}.
\end{equation}
For real measurable $f$ with $\E\bigl(\sum_i|f(\bar{\mathcal A}_i(0))|\bigr)^2<\infty$,
\begin{equation}\label{eq:occ-variance}
 \operatorname{Var}\sum_if(\bar{\mathcal A}_i(0))\le\int_\R f^2\rho_{\rm Ai}\diff z.
\end{equation}
\end{lemma}

\begin{proof}
By truncation we may take $f$ bounded with compact support. Grouping the indices $i_1,\dots,i_q$ in $\E\bigl(\sum_if(\bar{\mathcal A}_i(0))\bigr)^q$ by which of them are equal, and applying \eqref{eq:airy-correlations} at one position, gives
\[
 \E\Bigl(\sum_if(\bar{\mathcal A}_i(0))\Bigr)^q=\sum_{\pi\in\Pi_q}\int\prod_{J\in\pi}f(z_J)^{|J|}\det\bigl[K_{\rm Ai}(z_J,z_{J'})\bigr]_{J,J'\in\pi}\prod_{J\in\pi}\diff z_J.
\]
Since $K_{\rm Ai}$ is positive, Hadamard's inequality bounds the determinant by $\prod_{J\in\pi}\rho_{\rm Ai}(z_J)$, which gives the first bound in \eqref{eq:occ-one-position}. The second follows as in the proof of Lemma~\ref{lem:counts}(i), with $g=e^{\theta f}-1$ in place of $(e^\theta-1)\mathbf 1_{[c,\infty)}$. Finally, \eqref{eq:airy-correlations} and $K_{\rm Ai}^2=K_{\rm Ai}$ give \eqref{eq:occ-variance}, since $\operatorname{Var}\sum_if(\bar{\mathcal A}_i(0))=\int f^2\rho_{\rm Ai}-\|K_{\rm Ai}M_fK_{\rm Ai}\|_{\rm HS}^2$, with $M_f$ multiplication by $f$.
\end{proof}

\subsubsection*{Uniform moments and exponential integrability}

Fix $0<T<\infty$. The occupation estimates below assume that $V$ satisfies one of the following bounds, with a constant $A\ge0$:
\begin{equation}\label{eq:occ-exponential-bound}
 |V(y,u)|\le Ae^{B(1+|y|)^\beta-\delta u},
\end{equation}
where $B\ge0$, $\delta>0$, and either $0\le\beta<2$, or $\beta=2$ and $BT<\delta$; or
\begin{equation}\label{eq:occ-algebraic-bound}
 |V(y,u)|\le A(1+|y|)^m\psi(u),\qquad\psi(u)=\begin{cases}(1+u)^{-p},&u\ge0,\\ e^{-\delta u},&u<0,\end{cases}
\end{equation}
where $m\ge0$, $p>2+m/2$ and $\delta>0$.

\begin{lemma}\label{lem:occ-moments}
Let $V$ satisfy \eqref{eq:occ-exponential-bound} or \eqref{eq:occ-algebraic-bound}, let $1\le q<\infty$, and put $\omega(r)=e^{-\delta r}$, respectively $\omega(r)=(1+r_+)^{2+m/2-p}$. Let $R>0$ and $r_0\in\R$.
\begin{enumerate}
\item[\textup{(i)}] For $0<t\le T$, $|x|\le R$ and $r\ge r_0$,
\begin{equation}\label{eq:occ-moment-bounds}
 \|I_V^{t,x}(r)\|_q\le CA\omega(r),\qquad\|I_V^{t,x}(r)-\E I_V^{t,x}(r)\|_2\le CAt^{1/4}\omega(r).
\end{equation}
Under \eqref{eq:occ-exponential-bound}, $C$ does not depend on $r_0$, so \eqref{eq:occ-moment-bounds} holds for all $r$. If $B=0$, respectively $m=0$, $C$ does not depend on $R$, so \eqref{eq:occ-moment-bounds} holds for all $x$.
\item[\textup{(ii)}] Let $L\ge2R$ and $N\ge1$. If $V$ vanishes on $[-L,L]\times[-N,N]$, then for $0<t\le T$, $|x|\le R$ and $|r|\le R$,
\[
 \|I_V^{t,x}(r)\|_q\le CA\bigl(e^{-cL^2}+e^{-cN}\bigr),\quad\text{respectively}\quad\|I_V^{t,x}(r)\|_q\le CA\bigl(L^{m+4-2p}+N^{-c}\bigr).
\]
\end{enumerate}
The constants $C,c>0$ depend on $V$ only through the parameters in \eqref{eq:occ-exponential-bound} or \eqref{eq:occ-algebraic-bound}, as appropriate.
\end{lemma}

\begin{proof}
Since $L^q$ norms increase with $q$, we may take $q\ge2$ to be an integer, and by linearity we may take $A=1$. Let $f_v$ and $Y_v$ be as in \eqref{eq:occ-sections}.

\textit{Step~1: Reduction to one position.} Suppose that for some function $w\ge0$,
\begin{equation}\label{eq:occ-section-condition}
 \int_\R|f_v|^j\rho_{\rm Ai}\diff z\le Ct^{-1/2}w(v)^j,\qquad1\le j\le q,\ v\in\R.
\end{equation}
Then
\begin{equation}\label{eq:occ-step-one}
 \|I_V^{t,x}(r)\|_q\le C\int_\R w(v)\diff v,\qquad\|I_V^{t,x}(r)-\E I_V^{t,x}(r)\|_2\le Ct^{1/4}\int_\R w(v)\diff v.
\end{equation}
Indeed, since $t\le T$, Lemma~\ref{lem:one-position} and \eqref{eq:occ-section-condition} give
\[
 \E\Bigl(\sum_i|f_v(\bar{\mathcal A}_i(0))|\Bigr)^q\le\sum_{\pi\in\Pi_q}\prod_{J\in\pi}Ct^{-1/2}w(v)^{|J|}\le Ct^{-q/2}w(v)^q,
\]
so $\sqrt t\|Y_v\|_q\le Cw(v)$. By \eqref{eq:occ-variance} and \eqref{eq:occ-section-condition} with $j=2$, $\sqrt t(\operatorname{Var}Y_v)^{1/2}\le Ct^{1/4}w(v)$. Inserting these two bounds in \eqref{eq:occ-minkowski} gives \eqref{eq:occ-step-one}.

\textit{Step~2: Part \textup{(i)} under \eqref{eq:occ-exponential-bound}.} \eqref{eq:occ-section-condition} holds with $w(v)=e^{B(1+|x+v\sqrt t|)^\beta-\delta(r+v^2)}$, since $|f_v(z)|\le w(v)e^{\delta t^{1/3}z}$ and $\mu(j\delta t^{1/3})\le Ct^{-1/2}$ by \eqref{app:mu}. For $|x|\le R$ and $0<t\le T$, moreover, $w(v)\le Ce^{-\delta r-cv^2}$. This is clear if $\beta<2$, with $C$ independent of $R$ if $B=0$. If $\beta=2$, it follows by choosing $\kappa>0$ with $B(1+\kappa)T<\delta$ and using
\begin{equation}\label{eq:occ-beta-two}
 (1+|x+v\sqrt t|)^2\le(1+\kappa^{-1})(1+|x|)^2+(1+\kappa)tv^2.
\end{equation}
So \eqref{eq:occ-step-one} gives \eqref{eq:occ-moment-bounds} with $\omega(r)=e^{-\delta r}$, for all $r$.

\textit{Step~3: Part \textup{(i)} under \eqref{eq:occ-algebraic-bound}.} We may take $r_0\le0$. For $r_0\le r<0$, \eqref{eq:occ-scaled} gives $I_V^{t,x}(r)=I_{V'}^{t,x}(r-r_0)$ with $V'(y,u)=V(y,u+r_0)$. Since $\psi(u+r_0)\le C\psi(u)$, $V'$ satisfies \eqref{eq:occ-algebraic-bound} with a larger constant. So it suffices to take $r\ge0$. We show that for $|x|\le R$,
\begin{equation}\label{eq:occ-section-integral}
 \int_\R|f_v(z)|^j\rho_{\rm Ai}(z)\diff z\le C(1+|v|)^{jm}t^{-1/2}(1+r+v^2)^{3/2-jp},\qquad j\ge1.
\end{equation}
For $|x|\le R$, $|f_v(z)|\le C(1+|v|)^m\psi(u)$ with $u=r+v^2-t^{1/3}z$. For $z\le0$, $\psi(u)=(1+r+v^2+t^{1/3}|z|)^{-p}$ and $\rho_{\rm Ai}(z)\le C(1+\sqrt{|z|})$ by \eqref{eq:airy-density-bound}, and a change of variables bounds this part of the integral by the right side of \eqref{eq:occ-section-integral}. For $z\ge0$ with $t^{1/3}z\le(r+v^2)/2$, we have $u\ge(r+v^2)/2$, so $\psi(u)\le C(1+r+v^2)^{-p}$, and this part is at most $C(1+|v|)^{jm}(1+r+v^2)^{-jp}$. For the remaining $z$, $u\ge-t^{1/3}z$ gives $\psi(u)\le e^{\delta t^{1/3}z}$, and $z\ge(r+v^2)/(2T^{1/3})$. Since $j\delta t^{1/3}z\le\frac{c_0}2z^{3/2}+C$, \eqref{eq:airy-density-bound} bounds this part by $C(1+|v|)^{jm}(1+r+v^2)^{-jp}$ as well. As $t^{-1/2}(1+r+v^2)^{3/2}\ge T^{-1/2}$, this proves \eqref{eq:occ-section-integral}.

Since $(1+r+v^2)^{3/2-jp}\le\bigl((1+r+v^2)^{3/2-p}\bigr)^j$ for $j\ge1$, \eqref{eq:occ-section-integral} gives \eqref{eq:occ-section-condition} with $w(v)=C(1+|v|)^m(1+r+v^2)^{3/2-p}$. The substitution $v=s\sqrt{1+r}$ bounds $\int w(v)\diff v$ by $C(1+r)^{2+m/2-p}$, where the integral in $s$ converges because $p>2+m/2$. So \eqref{eq:occ-step-one} gives \eqref{eq:occ-moment-bounds}, with $C$ independent of $R$ if $m=0$.

\textit{Step~4: Part \textup{(ii)}.} Write $V=V_1+V_2$, with $V_1=V\mathbf 1_{\{|y|>L\}}$ and $V_2=V\mathbf 1_{\{|y|\le L,\,|u|>N\}}$.

The function $f_v$ of $V_1$ vanishes unless $|x+v\sqrt t|>L$, and for $|x|\le R\le L/2$ this requires $|v|>L/(2\sqrt T)$. So \eqref{eq:occ-section-condition} holds for $V_1$ with $w(v)\mathbf 1_{\{|v|>L/(2\sqrt T)\}}$ in place of $w(v)$, where $w$ is the weight of Step~2, respectively of Step~3 after the translation with $r_0=-R$. This weight is at most $Ce^{-cv^2}$, respectively $C(1+|v|)^{m+3-2p}$, so \eqref{eq:occ-step-one} gives $\|I_{V_1}^{t,x}(r)\|_q\le Ce^{-cL^2}$, respectively $CL^{m+4-2p}$.

For $V_2$ under \eqref{eq:occ-exponential-bound}, choose $\delta_1<\delta$, with $BT<\delta_1$ if $\beta=2$. Then $e^{-\delta u}\mathbf 1_{\{u>N\}}\le e^{-(\delta-\delta_1)N}e^{-\delta_1u}$ and $e^{-\delta u}\mathbf 1_{\{u<-N\}}\le e^{-N}e^{-(\delta+1)u}$. So part (i), with $\delta_1$, respectively $\delta+1$, in place of $\delta$, gives $\|I_{V_2}^{t,x}(r)\|_q\le Ce^{-cN}$. Under \eqref{eq:occ-algebraic-bound}, choose $2+m/2<p'<p$. Then $(1+u)^{-p}\mathbf 1_{\{u>N\}}\le(1+N)^{p'-p}(1+u)^{-p'}$. So part (i), with $p'$ in place of $p$, and with $\delta+1$ in place of $\delta$ for $u<-N$, gives $\|I_{V_2}^{t,x}(r)\|_q\le CN^{-c}$.
\end{proof}

\begin{lemma}\label{lem:occ-exp}
Let $V\ge0$ satisfy \eqref{eq:occ-exponential-bound} or \eqref{eq:occ-algebraic-bound}, and suppose that $V(y,u)\le M(1+u_-)^\nu$ for some $M<\infty$ and $0\le\nu<1$. Let $\theta>0$, $R>0$ and $r_0\in\R$. Then $\E e^{\theta I_V^{t,x}(r)}\le C$ for $0<t\le T$, $|x|\le R$ and $r\ge r_0$, where $C$ depends on $V$ only through $A$, $M$, $\nu$ and the parameters in the corresponding bound. If $B=0$, respectively $m=0$, this holds for all $x$ with the same $C$.
\end{lemma}

\begin{proof}
Under \eqref{eq:occ-algebraic-bound} we may take $r_0=0$, by the translation in Step~3 of the proof of Lemma~\ref{lem:occ-moments}, which preserves the bound $V\le M(1+u_-)^\nu$ with a larger $M$. Let $f_v$ and $Y_v$ be as in \eqref{eq:occ-sections}; $f_v\ge0$ since $V\ge0$.

Choose $\gamma$ with $\frac12<\gamma<\frac32-\nu$, and under \eqref{eq:occ-algebraic-bound} also $\gamma<p-\frac32-\frac m2$; this is possible because $\nu<1$ and $p>2+\frac m2$. Let $g(v)=c_\gamma(1+v^2)^{-\gamma}$, with $c_\gamma$ chosen so that $g$ is a probability density. Put $\theta_v=\theta\sqrt t/g(v)$, so that $\theta_v\le C\sqrt t(1+v^2)^\gamma$. By \eqref{eq:occ-sections}, $\theta I_V^{t,x}(r)=\int g(v)\theta_vY_v\diff v$. So Jensen's inequality for $g(v)\diff v$, stationarity and the exponential bound in Lemma~\ref{lem:one-position} give
\[
 \E e^{\theta I_V^{t,x}(r)}\le\int_\R g(v)\E e^{\theta_vY_v}\diff v\le\sup_{v\in\R}e^{Q(v)},\qquad Q(v)=\int_\R\bigl(e^{\theta_vf_v(z)}-1\bigr)\rho_{\rm Ai}(z)\diff z.
\]
It remains to bound $Q(v)$ uniformly in $v\in\R$, $0<t\le T$, $|x|\le R$ and $r\ge r_0$. By Steps~2 and~3 of the proof of Lemma~\ref{lem:occ-moments},
\begin{equation}\label{eq:occ-first-moment}
 \int_\R f_v\rho_{\rm Ai}\diff z\le Ct^{-1/2}e^{-cv^2},\quad\text{respectively}\quad\int_\R f_v\rho_{\rm Ai}\diff z\le Ct^{-1/2}(1+|v|)^{m+3-2p}.
\end{equation}

\textit{Large $|v|$, small $z$.} Fix $\varepsilon\in(0,1)$, with $B(1+\kappa)T<\delta(1-\varepsilon)$ for some $\kappa>0$ if $\beta=2$. For $t^{1/3}z\le\varepsilon v^2$, the height $r+v^2-t^{1/3}z$ is at least $r_0+(1-\varepsilon)v^2$. Thus \eqref{eq:occ-exponential-bound} gives $f_v(z)\le Ce^{-cv^2}$, using \eqref{eq:occ-beta-two} if $\beta=2$, while \eqref{eq:occ-algebraic-bound} gives $f_v(z)\le C(1+|v|)^{m-2p}$. Since $\theta_v\le C(1+v^2)^\gamma$, and $2\gamma+m-2p<0$ under \eqref{eq:occ-algebraic-bound}, there is a $v_0$ such that $\theta_vf_v\le1$ on this set for $|v|\ge v_0$. There $e^{\theta_vf_v}-1\le2\theta_vf_v$, so by \eqref{eq:occ-first-moment} the contribution of this set to $Q(v)$ is at most $C(1+v^2)^\gamma e^{-cv^2}$, respectively $C(1+|v|)^{2\gamma+m+3-2p}$. Both are bounded, the second because $\gamma<p-\frac32-\frac m2$.

\textit{Large $|v|$, large $z$.} For $|v|\ge v_0$ and $t^{1/3}z\ge\varepsilon v^2$, we have $t^{1/3}z\ge1$ after enlarging $v_0$. The height is at least $r_0-t^{1/3}z$, so the bound $V\le M(1+u_-)^\nu$ gives $f_v(z)\le C(t^{1/3}z)^\nu$, and
\[
 \frac{\theta_vf_v(z)}{z^{3/2}}\le C|v|^{2\gamma}t^{1/2+\nu/3}\Bigl(\frac{\varepsilon v^2}{t^{1/3}}\Bigr)^{\nu-3/2}=Ct|v|^{2\gamma+2\nu-3}.
\]
Since $2\gamma+2\nu-3<0$, after enlarging $v_0$ we have $\theta_vf_v\le\frac{c_0}2z^{3/2}$ on this set. By \eqref{eq:airy-density-bound}, its contribution to $Q(v)$ is at most $C\int_0^\infty e^{-c_0z^{3/2}/2}\diff z$.

\textit{Bounded $|v|$.} For $|v|\le v_0$, $\theta_v\le C\sqrt t$. For $z\le0$, the height is at least $r_0$, so $f_v\le M(1+|r_0|)^\nu$ and $e^{\theta_vf_v}-1\le C\theta_vf_v$. By \eqref{eq:occ-first-moment}, this part is therefore at most $C\sqrt t\int f_v\rho_{\rm Ai}\le C$. For $z\ge0$, the height is at least $r_0-t^{1/3}z$, so $\theta_vf_v(z)\le C(1+z^\nu)$. Since $e^{C(1+z^\nu)}$ is integrable against the bound \eqref{eq:airy-density-bound}, this part is bounded too.

So $Q$ is bounded, and if $B=0$, respectively $m=0$, none of these bounds depends on $x$.
\end{proof}

\textit{Differences of exponentials.} For potentials $V$ and $V'$ whose occupations are realized on one probability space, and for $1\le q<\infty$,
\begin{equation}\label{eq:occ-Z-difference}
 \bigl\|e^{-I_V}-e^{-I_{V'}}\bigr\|_q\le\|I_V-I_{V'}\|_{2q}\bigl(\|e^{I_{V_-}}\|_{2q}+\|e^{I_{V'_-}}\|_{2q}\bigr),
\end{equation}
and the same holds with a constant $a$ in place of $I_{V'}$ and $e^{a_-}$ in place of $e^{I_{V'_-}}$. This follows from $|e^{-a}-e^{-b}|\le|a-b|(e^{a_-}+e^{b_-})$, which is the mean value theorem, together with $I_V\ge-I_{V_-}$ and the Cauchy--Schwarz inequality.

\subsubsection*{Continuity, differentiation, and traces}

\begin{lemma}\label{lem:diff-under-E}
Let $V\in C_c^\infty(\R^2;\R)$ and $t>0$. In the realization \eqref{eq:occ-scaled}, almost surely, $I_V^{t,x}(r)$ is a smooth function of $(x,r)$, with $\partial_x^j\partial_r^kI_V^{t,x}(r)=I_{\partial_x^j\partial_r^kV}^{t,x}(r)$. Moreover, $Z_V(t,\cdot,\cdot)$ is smooth, with $\partial_x^j\partial_r^kZ_V=\E\partial_x^j\partial_r^ke^{-I_V}$.
\end{lemma}

\begin{proof}
Fix $(x_0,r_0)\in\R^2$. We show that, almost surely, the series \eqref{eq:occ-scaled} for $I_V$ and its termwise derivatives converge uniformly for $|x-x_0|,|r-r_0|\le1$. Then $I_V$ is smooth there, with $\partial_x^j\partial_r^kI_V=I_{\partial_x^j\partial_r^kV}$, and covering $\R^2$ by countably many such squares gives the first assertion. For $G\in C_c(\R^2)$, the function $\widetilde G(y,u)=\sup\{|G(y',u')|:|y'-y|,|u'-u|\le1\}$ is continuous with compact support, and
\begin{equation}\label{eq:diff-majorant}
 t^{2/3}\sum_{i\ge1}\int_\R\sup_{|x-x_0|,|r-r_0|\le1}\bigl|G\bigl(x+t^{2/3}s,r+t^{1/3}(s^2-\bar{\mathcal A}_i(s))\bigr)\bigr|\diff s\le I_{\widetilde G}^{t,x_0}(r_0).
\end{equation}
By Lemma~\ref{lem:occ-moments}(i), the right side has finite moments of all orders, so it is finite almost surely. Applied to $G=\partial_x^j\partial_r^kV$ for all $j,k$, \eqref{eq:diff-majorant} gives the uniform convergence.

By the chain rule, $\partial_x^j\partial_r^ke^{-I_V}$ is $e^{-I_V}$ times a polynomial in the occupations $I_{\partial_x^a\partial_r^bV}$, $a\le j$, $b\le k$. For $|x-x_0|,|r-r_0|\le1$, \eqref{eq:diff-majorant} gives $|I_{\partial_x^a\partial_r^bV}|\le I_{\widetilde G}^{t,x_0}(r_0)$ with $G=\partial_x^a\partial_r^bV$, and $e^{-I_V}\le e^{I_{V_-}}\le e^{I_{\widetilde G}^{t,x_0}(r_0)}$ with $G=V_-$. By Lemma~\ref{lem:occ-moments}(i), Lemma~\ref{lem:occ-exp} and H\"older's inequality, the supremum of $|\partial_x^j\partial_r^ke^{-I_V}|$ over the square therefore has finite expectation. Dominated convergence then shows that $Z_V=\E e^{-I_V}$ is smooth there, with $\partial_x^j\partial_r^kZ_V=\E\partial_x^j\partial_r^ke^{-I_V}$.
\end{proof}

\begin{lemma}\label{lem:occ-continuity}
Let $V_n$ and $V$ be measurable and satisfy \eqref{eq:occ-exponential-bound} with common constants, let $V_n\to V$ in $L^1_{\rm loc}(\R^2)$, and let $t_n\to t$ in $(0,T]$. Then, in the realization \eqref{eq:occ-scaled}, $I_{V_n}^{t_n,x}(r)\to I_V^{t,x}(r)$ in $L^q$ for $x,r\in\R$ and $1\le q<\infty$. If moreover $(V_n)_-,V_-\le M(1+u_-)^\nu$ with $0\le\nu<1$, then $Z_{V_n}(t_n,x,r)\to Z_V(t,x,r)$.
\end{lemma}

\begin{proof}
Fix $x$ and $r$. For a function $G$ on $\R^2$ and $\tau>0$, let $G^\tau(s,z)=\tau^{2/3}G\bigl(x+\tau^{2/3}s,r+\tau^{1/3}(s^2-z)\bigr)$, so that $I_G^{\tau,x}(r)=\sum_i\int_\R G^\tau(s,\bar{\mathcal A}_i(s))\diff s$ by \eqref{eq:occ-scaled}. Then
\[
 |I_{V_n}^{t_n,x}(r)-I_V^{t,x}(r)|\le\Delta_n\vcentcolon=\sum_{i\ge1}\int_\R\bigl|V_n^{t_n}-V^t\bigr|\bigl(s,\bar{\mathcal A}_i(s)\bigr)\diff s.
\]
We show that $\E\Delta_n\to0$. Since $\Delta_n\le I_{|V_n|}^{t_n,x}(r)+I_{|V|}^{t,x}(r)$, Lemma~\ref{lem:occ-moments}(i) bounds $\|\Delta_n\|_{2q}$ uniformly in $n$, so H\"older's inequality then gives $\|\Delta_n\|_q\to0$.

By stationarity, $\E\Delta_n=\iint_{\R^2}|V_n^{t_n}-V^t|(s,z)\rho_{\rm Ai}(z)\diff s\diff z$. For $\tau\in\{t,t_1,t_2,\dots\}$, a compact subset of $(0,T]$, \eqref{eq:occ-exponential-bound} gives $|V_n^\tau|,|V^\tau|\le Ce^{-cs^2+Cz_+-cz_-}$, as in Step~2 of the proof of Lemma~\ref{lem:occ-moments}. This majorant is integrable against $\rho_{\rm Ai}(z)\diff s\diff z$ by \eqref{eq:airy-density-bound}, so we may restrict the integral to a compact set $D$. Since the map $(s,z)\mapsto(x+\tau^{2/3}s,r+\tau^{1/3}(s^2-z))$ has Jacobian $\tau$ and $\rho_{\rm Ai}$ is bounded on $D$, $\iint_D|G^\tau|\rho_{\rm Ai}\le C\int_{D'}|G|$ for these $\tau$ and every $G$, with $C$ and a compact set $D'$ independent of $\tau$ and $G$. Split $V_n^{t_n}-V^t=(V_n-V)^{t_n}+(V^{t_n}-V^t)$. The integral of $|(V_n-V)^{t_n}|\rho_{\rm Ai}$ over $D$ is at most $C\int_{D'}|V_n-V|$, which tends to zero. That of $|V^{t_n}-V^t|\rho_{\rm Ai}$ tends to zero for continuous $V$, since then $V^{t_n}\to V^t$ uniformly on $D$. It therefore tends to zero for general $V$ as well, by approximating $V$ in $L^1(D')$ by continuous functions and using the same bound.

The convergence of $Z$ follows from \eqref{eq:occ-Z-difference} with $q=1$ and Lemma~\ref{lem:occ-exp}, applied to $(V_n)_-$ and $V_-$.
\end{proof}

We write $\E I_V(t)$ for the function $(x,r)\mapsto\E I_V^{t,x}(r)$.

\begin{lemma}\label{lem:occ-trace}
\begin{enumerate}
\item[\textup{(i)}] Let $V\in L^1(\R^2)$ vanish outside $[-R,R]^2$. Then $\int_{[-R,R]^2}\E I_{|V|}^{t,x}(r)\diff x\diff r\le C\|V\|_1$ for $0<t\le1$, with $C$ depending only on $R$, and $\E I_V(t)\to\frac12\partial_r^{-2}V$ in $L^1_{\rm loc}(\R^2)$ as $t\downarrow0$.
\item[\textup{(ii)}] Let $V_n$ and $V$ be continuous and satisfy \eqref{eq:occ-exponential-bound} or \eqref{eq:occ-algebraic-bound} with common constants, let $V_n\to V$ locally uniformly, and let $t_n\downarrow0$. Then, for $1\le q<\infty$ and locally uniformly in $(x,r)$, $I_{V_n}^{t_n,x}(r)\to\frac12\partial_r^{-2}V(x,r)$ in $L^q$. If moreover $(V_n)_-\le M(1+u_-)^\nu$ with $0\le\nu<1$, then also $e^{-I_{V_n}^{t_n,x}(r)}\to e^{-\frac12\partial_r^{-2}V(x,r)}$ in $L^q$, locally uniformly in $(x,r)$, and in particular $Z_{V_n}(t_n,x,r)\to e^{-\frac12\partial_r^{-2}V(x,r)}$.
\item[\textup{(iii)}] Let $V$ be measurable and satisfy \eqref{eq:occ-exponential-bound} or \eqref{eq:occ-algebraic-bound}, and let $V_-\le M(1+u_-)^\nu$ with $0\le\nu<1$. Then $0<Z_V<\infty$ on $(0,T]\times\R^2$. As $t\downarrow0$, $\log Z_V(t)+\E I_V(t)\to0$ locally uniformly, and $\E I_V(t)\to\frac12\partial_r^{-2}V$ and $\log Z_V(t)\to-\frac12\partial_r^{-2}V$ in $L^1_{\rm loc}(\R^2)$. Hence $-2\partial_r^2\log Z_V(t)\to V$ in $\mathcal D'(\R^2)$.
\end{enumerate}
\end{lemma}

\begin{proof}
(i) For $|y|,|u|\le R$ and $0<t\le1$, substitute $x=y+v\sqrt t$ in the integral of $\varrho_{t,x,r}(y,u)$ over $[-R,R]^2$, extend the $v$-integration to $\R$, and apply \eqref{eq:airy-density-scaled-bound}. The resulting integral is bounded by a constant depending only on $R$. Integrating against $|V(y,u)|$ and using \eqref{eq:occ-density} gives the bound. Since $V$ also vanishes outside larger squares, it suffices to prove the limit in $L^1([-R,R]^2)$. Let $G\in C_c(\R^2)$ vanish outside $[-R,R]^2$. By Proposition~\ref{prop:lln}, $\E I_G(t)\to\frac12\partial_r^{-2}G$ uniformly on $[-R,R]^2$, and $\|\partial_r^{-2}G\|_{L^1([-R,R]^2)}\le C\|G\|_1$. Approximating $V$ in $L^1$ by such $G$ and using the bound proves the limit.

(ii) Fix $R>0$ and $L\ge2R$. Let $\Theta\in C_c(\R^2)$ take values in $[0,1]$, equal $1$ on $[-L,L]^2$, and vanish outside $[-2L,2L]^2$. Write $I_G$ for $I_G^{t_n,x}(r)$, and $\partial_r^{-2}G$ for its value at $(x,r)$. Then
\[
 I_{V_n}-\tfrac12\partial_r^{-2}V=\bigl(I_{V\Theta}-\tfrac12\partial_r^{-2}(V\Theta)\bigr)+I_{(V_n-V)\Theta}+I_{V_n(1-\Theta)}-\tfrac12\partial_r^{-2}\bigl(V(1-\Theta)\bigr).
\]
The following bounds are uniform for $|x|,|r|\le R$. Since $V\Theta\in C_c(\R^2)$, the first term tends to zero in $L^2$ by Proposition~\ref{prop:lln}, and hence in every $L^q$, by Lemma~\ref{lem:occ-moments}(i) and H\"older's inequality. Let $\epsilon_n$ be the supremum of $|V_n-V|$ on $[-2L,2L]^2$, which tends to zero. Then $|(V_n-V)\Theta|\le\epsilon_ne^{2L}e^{-u}$, so Lemma~\ref{lem:occ-moments}(i) bounds the second term in $L^q$ by $C\epsilon_ne^{2L}$. The function $V_n(1-\Theta)$ vanishes on $[-L,L]^2$ and satisfies the same bound as $V_n$. So by Lemma~\ref{lem:occ-moments}(ii), the third term is bounded in $L^q$, uniformly in $n$, by a quantity that tends to zero as $L\to\infty$. The last term also tends to zero as $L\to\infty$, since for $|x|\le R$ it involves only heights $u>L$, where the assumed decay bound is integrable against $(u-r)\diff u$. Letting $n\to\infty$ and then $L\to\infty$ proves the first assertion. The second follows from \eqref{eq:occ-Z-difference} with the constant $\frac12\partial_r^{-2}V(x,r)$, which is bounded for $|x|,|r|\le R$, and Lemma~\ref{lem:occ-exp}, applied to $(V_n)_-$.

(iii) By Lemma~\ref{lem:occ-moments}(i), $I_{|V|}^{t,x}(r)$ has finite moments, so the series for $I_V^{t,x}(r)$ converges absolutely almost surely and $Z_V>0$. By Lemma~\ref{lem:occ-exp}, applied to $V_-$, $Z_V\le\E e^{I_{V_-}}<\infty$. Fix $R>0$; the bounds below are uniform for $|x|,|r|\le R$ and $0<t\le T$. By \eqref{eq:occ-Z-difference} with $q=1$ and the constant $\E I_V$, Lemma~\ref{lem:occ-moments}(i) and Lemma~\ref{lem:occ-exp},
\[
 |Z_V-e^{-\E I_V}|\le\|I_V-\E I_V\|_2\bigl(\|e^{I_{V_-}}\|_2+e^{\E I_{V_-}}\bigr)\le Ct^{1/4}.
\]
By Jensen's inequality, $Z_V\ge e^{-\E I_V}\ge e^{-\E I_{|V|}}$, which is bounded below by Lemma~\ref{lem:occ-moments}(i). Hence $\log Z_V(t)+\E I_V(t)\to0$ uniformly.

It remains to show that $\E I_V(t)\to\frac12\partial_r^{-2}V$ in $L^1([-R,R]^2)$, which also gives the convergence of $\log Z_V(t)$. Take $L\ge2R$ and $\Theta$ as in (ii). Since $V\Theta$ is integrable and vanishes outside $[-2L,2L]^2$, (i) gives $\E I_{V\Theta}(t)\to\frac12\partial_r^{-2}(V\Theta)$ in $L^1([-R,R]^2)$. By Lemma~\ref{lem:occ-moments}(ii), $|\E I_{V(1-\Theta)}(t)|$ is bounded by a quantity that tends to zero as $L\to\infty$, uniformly in $t$, and $\frac12\partial_r^{-2}(V(1-\Theta))\to0$ uniformly, as in (ii).

For the last assertion, let $\varphi\in C_c^\infty(\R^2)$. By \eqref{eq:occ-exponential-bound} or \eqref{eq:occ-algebraic-bound}, $(x,r)\mapsto\int_r^\infty(u-r)|V(x,u)|\diff u$ is locally bounded. Fubini's theorem and $\int_{-\infty}^u(u-r)\partial_r^2\varphi(x,r)\diff r=\varphi(x,u)$ therefore give $\int_{\R^2}(\partial_r^{-2}V)\partial_r^2\varphi=\int_{\R^2}V\varphi$. With the convergence of $\log Z_V(t)$ in $L^1_{\rm loc}(\R^2)$,
\[
 -2\int_{\R^2}\log Z_V(t)\partial_r^2\varphi\longrightarrow\int_{\R^2}(\partial_r^{-2}V)\partial_r^2\varphi=\int_{\R^2}V\varphi.
\]
\end{proof}

\subsection{Small-time concentration of the occupation}\label{app:lln}

\begin{proof}[Proof of Proposition~\ref{prop:lln}]
The mean calculation, including the limit and the elementary integral
over $v$, is given in Section~\ref{sec:compact-recovery}.
Here we justify the dominated convergence and its local uniformity,
using the representation \eqref{eq:lln-mean}.

To apply dominated convergence, let $R>0$ be such that $V=0$ outside $[-R,R]^2$. For $(x,r)\in[-R,R]^2$ and $0<t\le1$, the integrand in \eqref{eq:lln-mean} is at most
\[
 C\|V\|_\infty e^{-w_+}\mathbf 1_{\{w\ge-2R,\ v^2\le w+2R\}},
\]
which is integrable in $(w,v)$ and independent of $t$, $x$ and $r$. Indeed, $V(x+v\sqrt t,r+v^2-w)$ vanishes unless $|r+v^2-w|\le R$, which for $|r|\le R$ forces $w\ge-2R$ and $v^2\le w+2R$; and for $w\ge-2R$, \eqref{eq:airy-density-scaled-bound} gives $t^{1/6}\rho_{\rm Ai}(t^{-1/3}w)\le Ce^{-w_+}$.

To prove the uniformity, let $t_n\downarrow0$, and let $(x_n,r_n)\in[-R,R]^2$ converge to $(x,r)$. By the continuity of $V$ and dominated convergence, $\E I_V^{t_n,x_n}(r_n)\to\frac12\partial_r^{-2}V(x,r)$. Convergence along every such sequence, together with the continuity of $\frac12\partial_r^{-2}V$, gives uniform convergence for $(x,r)\in[-R,R]^2$. Since $R$ can be enlarged, the convergence is uniform on every compact set.

\textit{The variance.} Since $V=0$ outside $[-R,R]^2$, the potential $V$ satisfies $|V(y,u)|\le\|V\|_\infty e^{R-u}$, which is \eqref{eq:occ-exponential-bound} with $A=\|V\|_\infty e^R$, $B=0$ and $\delta=1$. Lemma~\ref{lem:occ-moments}(i), with $T=1$, therefore gives $\|I_V^{t,x}(r)-\E I_V^{t,x}(r)\|_2\le Ct^{1/4}e^{-r}$ for $0<t\le1$ and all $x,r\in\R$. Hence $\operatorname{Var}I_V^{t,x}(r)=O(\sqrt t)$, uniformly for $(x,r)$ in compact sets. By the triangle inequality in $L^2$, the limit of the mean and the bound on the variance give \eqref{recovery}.
\end{proof}

\section{Initial-value, extension, and restart arguments}
\label{app:initial-value-methods}

Section~\ref{sec:fredholm-evolution} established the KP-II equation
for smooth compactly supported data. Here we prove the smooth
initial-data limit using the occupation estimates, and then give
the technical arguments for less localized potentials and for
the restart identities.

\subsection{Compact data: recovery of the initial condition}\label{app:compact-initial-trace}

\begin{proof}[Proof of Theorem~\ref{thm:compact}]
We show that $Z_V\to e^{-\frac12\partial_r^{-2}V}$ in $C^\infty_{\rm loc}(\R^2)$ as $t\downarrow0$. Since $e^{-\frac12\partial_r^{-2}V}$ is positive, this convergence and the chain \eqref{eq:mV} give \eqref{eq:compact-trace}. The convergence follows from the law of large numbers for the occupations of $V$ and of its derivatives. Write $\partial^\eta=\partial_x^{\eta_1}\partial_r^{\eta_2}$ for a multi-index $\eta$. By Lemma~\ref{lem:diff-under-E} and the chain rule, $\partial^\eta Z_V$ is the expectation of $e^{-I_V}$ times a polynomial in the occupations $I_{\partial^\gamma V}$, $\gamma\le\eta$. As $t\downarrow0$, $e^{-I_V}\to e^{-\frac12\partial_r^{-2}V}$ and $I_{\partial^\gamma V}\to\frac12\partial_r^{-2}\partial^\gamma V$ in every $L^q$, $q<\infty$, uniformly for $(x,r)$ in compact sets. These limits follow from Lemma~\ref{lem:occ-trace}(ii): as in the proof of Proposition~\ref{prop:lln}, $V$ and each $\partial^\gamma V$ satisfy \eqref{eq:occ-exponential-bound} with $B=0$, and $V_-$ is bounded. H\"older's inequality then shows that $\partial^\eta Z_V$ converges, uniformly on compact sets, to $e^{-\frac12\partial_r^{-2}V}$ times the same polynomial in the functions $\frac12\partial_r^{-2}\partial^\gamma V$. Since $\partial^\gamma\partial_r^{-2}V=\partial_r^{-2}\partial^\gamma V$, the chain rule shows that this limit is $\partial^\eta e^{-\frac12\partial_r^{-2}V}$.
\end{proof}

\subsection{Extensions beyond compact smooth data}\label{app:extensions}

We extend the compact construction to integrable and one-sided
data by combining Fredholm estimates with
the preceding occupation bounds.
In this subsection and the next, $\partial^\eta=\partial_x^{\eta_1}\partial_r^{\eta_2}$ for a multi-index $\eta=(\eta_1,\eta_2)$.

\subsubsection{Compact integrable data}

\begin{proof}[Proof of Proposition~\ref{prop:integrable-det}]
Choose $\varepsilon>0$ and $R$ such that $V$ vanishes outside $\mathcal C=[a+\varepsilon,b-\varepsilon]\times[-R,R]$, and put $\ell=b-a$. The proof follows Section~\ref{sec:direct} for continuous data, with a Gaussian bound on the defect $D_{a,b;V}$ of \eqref{eq:defect} that uses only $\|V\|_1$ and $M$.

\textit{The defect.} The proof of Lemma~\ref{lem:gauss} gives, for $W$ with the same properties as $V$,
\begin{equation}\label{eq:extensions-defect}
 \begin{aligned}
 |D_{a,b;V}(w,w')|&\le Ce^{M\ell}\|V\|_1e^{-c(w^2+w'^2)},\\
 |D_{a,b;V}(w,w')-D_{a,b;W}(w,w')|&\le Ce^{M\ell}\|V-W\|_1e^{-c(w^2+w'^2)}.
 \end{aligned}
\end{equation}
That proof uses $|V|\le M$ only through $\int_a^bV(y,B(y))\diff y\ge-M\ell$, which holds here because $V_-\le M$. The same chain shows that this integral converges absolutely for almost every bridge. Since both $\int_a^bV(y,B(y))\diff y$ and $\int_a^bW(y,B(y))\diff y$ are at least $-M\ell$, the same chain with $|V-W|$ in place of $|V|$ gives the second bound.

\textit{Finiteness.} Since $\varrho_{t,x,r}$ is bounded on $\mathcal C$, \eqref{eq:occ-density} gives $\E I_{|V|}^{t,x}(r)<\infty$, so the occupation $I_V$ converges absolutely almost surely, and $Z_V>0$. Since $e^{-I_V}\le e^{MI_{\mathbf 1_{\mathcal C}}}$, Lemma~\ref{lem:occ-exp} gives $Z_V<\infty$.

\textit{The determinant.} Truncating $V_+$ at level $n$ and mollifying with a nonnegative kernel gives continuous $V_n$ with $(V_n)_-\le M$, supports in a fixed rectangle $\mathcal C'\subset(a,b)\times\R$, and $V_n\to V$ in $L^1$. By Theorem~\ref{thm:fredholm}, $Z_{V_n}=\det(I-\chi_rK_{t,x;V_n}\chi_r)$. As $n\to\infty$, the right side converges to $\det(I-\chi_rK_{t,x;V}\chi_r)$: the computation in Step~1 of the proof of Proposition~\ref{prop:limit} (Appendix~\ref{app:limit}), combined with \eqref{eq:extensions-defect}, gives $\|\chi_r(K_{t,x;V_n}-K_{t,x;V})\chi_r\|_1\le C\|V_n-V\|_1$. The left side converges to $Z_V$ by dominated convergence: $I_{V_n}\to I_V$ in probability, because $\E|I_{V_n}-I_V|\le\sup_{\mathcal C'}\varrho_{t,x,r}\|V_n-V\|_1\to0$ by \eqref{eq:occ-density}, and $e^{-I_{V_n}}\le e^{MI_{\mathbf 1_{\mathcal C'}}}$, which is integrable by Lemma~\ref{lem:occ-exp}. Therefore $Z_V(t,x,r)=\det(I-\chi_rK_{t,x;V}\chi_r)$.

\textit{The kernel.} The kernel $\mathcal K_V$ satisfies (i)--(iii) of Proposition~\ref{prop:qr-kernel}: its proof uses $V$ only through the Gaussian bound of Lemma~\ref{lem:gauss}, which \eqref{eq:extensions-defect} replaces.
\end{proof}

\begin{proof}[Proof of the initial condition in Theorem~\ref{thm:measurable-data}(i)]
Let $\widetilde{\mathcal C}\subset\R^2$ be compact. We show that $\log Z_V(t)\to-\bigl(\tfrac12\partial_r^{-2}V\bigr)$ in $L^1(\widetilde{\mathcal C})$ as $t\downarrow0$. Then, as at the end of the proof of Lemma~\ref{lem:occ-trace}, $\phi_V(t)=-2\partial_r^2\log Z_V(t)\to V$ in $\mathcal D'(\R^2)$.

By Lemma~\ref{lem:occ-trace}(i), $\E I_V(t)\to \bigl(\tfrac12\partial_r^{-2}V\bigr)$ in $L^1(\widetilde{\mathcal C})$, so it suffices to show that $\E I_V+\log Z_V\to0$ in $L^1(\widetilde{\mathcal C})$. For $k\ge1$,
\[
 0\le\E I_V+\log Z_V\le\E I_{V-\min\{V,k\}}+\bigl(\E I_{\min\{V,k\}}+\log Z_{\min\{V,k\}}\bigr).
\]
Here, the first inequality follows from Jensen's inequality. The function $\min\{V,k\}$ is bounded, vanishes outside $\mathcal C$, and has negative part $V_-\le M$, so by Lemma~\ref{lem:occ-trace}(iii) the last term tends to $0$ uniformly on $\widetilde{\mathcal C}$ as $t\downarrow0$. By Lemma~\ref{lem:occ-trace}(i), the $L^1(\widetilde{\mathcal C})$ norm of $\E I_{V-\min\{V,k\}}(t)$ is at most $C\|V-\min\{V,k\}\|_1$ for $0<t\le1$, which tends to $0$ as $k\to\infty$.
\end{proof}

\subsubsection{Data below a height}

\begin{lemma}\label{lem:extensions-rows}
Let $t_0>0$, $r_0,R\in\R$ and $\lambda>0$. There are a complex neighborhood $U$ of $t_0$, an $\epsilon_1>0$, and constants $C,c>0$ such that
\[
 |S_{t,\alpha}(w,r+v)|\le Ce^{-\lambda\operatorname{Re}(r+v-w)}e^{C(\operatorname{Re}\alpha)_+-c(\operatorname{Re}\alpha)_-^3}
\]
for $t\in U$, complex $\alpha$ with $|\operatorname{Im}\alpha|\le\epsilon_1$, complex $r$ with $\operatorname{Re}r\ge r_0$ and $|\operatorname{Im}r|\le1$, real $w\le R$, and $v\ge0$. After $U$, $\epsilon_1$ and the strip in $r$ are shrunk, the same bound holds for every derivative of fixed order in $t$, $\alpha$, $r$ and $v$.
\end{lemma}

\begin{proof}
Put $q=r+v-w$, so that $\operatorname{Re}q\ge r_0-R$ and $|\operatorname{Im}q|\le1$, and recall from \eqref{eq:S} that
\[
 S_{t,\alpha}(w,r+v)=t^{-1/3}\exp\Bigl\{\frac{2\alpha^3}{3t^2}+\frac{\alpha q}t\Bigr\}\Ai(\xi),\qquad\xi=t^{-4/3}\alpha^2+t^{-1/3}q.
\]
We take $U$ and $\epsilon_1$ so small that $\arg t$ is small, and $\arg(\pm\alpha)$ is small when $\pm\operatorname{Re}\alpha$ is large. We use $|\Ai(\xi)|\le Ce^{-\frac23\operatorname{Re}\xi^{3/2}}$ in a sector $|\arg\xi|\le\pi/3$ \cite[\S9.7]{DLMF}, with the asymptotic $\Ai(\xi)\sim(2\sqrt\pi)^{-1}\xi^{-1/4}e^{-\frac23\xi^{3/2}}$ there. Fix a large $Y_0$.

\textit{The case $\operatorname{Re}\alpha\ge Y_0$.} Write $\xi=t^{-4/3}\alpha^2(1+z)$ with $z=qt/\alpha^2$. The exponent in $S$ combines with the Airy exponent into
\[
 \frac{2\alpha^3}{3t^2}+\frac{\alpha q}t-\frac23t^{-2}\alpha^3(1+z)^{3/2}=\frac{\alpha^3}{t^2}H(z),\qquad H(z)=\frac23+z-\frac23(1+z)^{3/2}=-\frac{z^2}2\int_0^1\frac{1-s}{\sqrt{1+sz}}\diff s,
\]
so that $\frac{\alpha^3}{t^2}H(z)=-\frac{q^2}{2\alpha}\int_0^1(1-s)(1+sqt/\alpha^2)^{-1/2}\diff s$. If $\operatorname{Re}q$ is bounded, then so is $|q|$, and $\frac{\alpha^3}{t^2}H(z)$ and $e^{\lambda\operatorname{Re}q}$ are bounded. If $\operatorname{Re}q$ is large, all the arguments involved are small, and the real part of the integral is at least $c(1+|z|)^{-1/2}$, so
\[
 \operatorname{Re}\frac{\alpha^3}{t^2}H(z)\le-\frac{c(\operatorname{Re}q)^2}{\operatorname{Re}\alpha+\sqrt{\operatorname{Re}q}}\le C_\lambda\operatorname{Re}\alpha-\lambda\operatorname{Re}q.
\]
For the last inequality, bound the denominator by $2\operatorname{Re}\alpha$ if $\operatorname{Re}q\le(\operatorname{Re}\alpha)^2$, and by $2\sqrt{\operatorname{Re}q}$ otherwise. The prefactor $t^{-1/3}\xi^{-1/4}$ is bounded.

\textit{The case $\operatorname{Re}\alpha\le-Y_0$.} Write $\alpha=-A$. Then $\xi=t^{-4/3}A^2+t^{-1/3}q$ has a small argument and $|\xi|\ge c(|A|^2+\operatorname{Re}q)-C$, so $|\Ai(\xi)|\le Ce^{-c|A|^3-c(\operatorname{Re}q)_+^{3/2}}$. In the prefactor, $\operatorname{Re}(-2A^3/3t^2)\le-c|A|^3$ and $\operatorname{Re}(-Aq/t)\le C|A|$. Together, $|S_{t,\alpha}|\le Ce^{-c|\operatorname{Re}\alpha|^3-\lambda\operatorname{Re}q}$.

\textit{The case $|\operatorname{Re}\alpha|\le Y_0$.} If $\operatorname{Re}q$ is large, $\xi$ has a small argument and modulus at least $c\operatorname{Re}q$, so $|S_{t,\alpha}|\le Ce^{C\operatorname{Re}q-c(\operatorname{Re}q)^{3/2}}\le C_\lambda e^{-\lambda\operatorname{Re}q}$. If $\operatorname{Re}q$ is bounded, $S_{t,\alpha}$ is bounded.

The three cases prove the bound. The bound for derivatives follows from Cauchy's formula.
\end{proof}

\begin{proof}[Proof of Proposition~\ref{prop:endpoint-free}]
Let $V$ be measurable with $|V|\le M$ and $V(y,u)=0$ for $u>R$. Fix $t_0>0$, $x_0\in\R$ and $r_0\in\R$, and let $U$, $\epsilon_1$ and the strip in $r$ be as in Lemma~\ref{lem:extensions-rows}, for a fixed $\lambda>0$. Let $\Omega$ be the set of complex $(t,x,r)$ with $t\in U$, $|x-x_0|<\epsilon_1$ and $r$ in the strip; for real $y$, the lemma applies to $\alpha=\pm(y-x)$ on $\Omega$.

\textit{Trace bounds.} By Lemma~\ref{lem:extensions-rows} with $\alpha=\pm(y-x)$, the functions $v\mapsto S_{t,y-x}(w,r+v)$ and $v\mapsto S_{t,x-y}(w,r+v)$ have $L^2(\R_+)$ norms at most $Ce^{-\lambda\operatorname{Re}(r-w)}E_1(y)$ and $Ce^{-\lambda\operatorname{Re}(r-w)}E_2(y)$, where
\[
 E_1(y)=e^{Cy_+-cy_-^3},\qquad E_2(y)=e^{Cy_--cy_+^3}.
\]
Each integrand in \eqref{eq:extensions-two-insertions} is a rank-one kernel, with trace norm the product of these two norms. The Feynman--Kac formula gives $0\le P_{y,y';V}(w,w')\le e^{M(y'-y)}p_{y'-y}(w'-w)$, and
\[
 \iint_{w,w'\le R}e^{\lambda(w+w')}p_{y'-y}(w'-w)\diff w\diff w'\le\lambda^{-1}e^{2\lambda R}.
\]
Hence the first term of \eqref{eq:extensions-two-insertions} has trace norm at most $Ce^{-2\lambda\operatorname{Re}r}\int_\R E_1E_2\diff y$, and the second at most
\[
 Ce^{-2\lambda\operatorname{Re}r}\iint_{y<y'}E_1(y)E_2(y')e^{M(y'-y)}\diff y\diff y'.
\]
Both integrals are finite, by the cubic decay of $E_1$ at $-\infty$ and of $E_2$ at $+\infty$; in the second, $y<y'$, so $y<0$ or $y'>0$, and $|y'|\le|y|$ if $y'<0$, $|y|\le y'$ if $y>0$. The same bounds hold for every fixed derivative of the integrands in $(t,x,r)$, by the last statement of Lemma~\ref{lem:extensions-rows}. As in the proof of Proposition~\ref{prop:qr-kernel}(i), $\mathcal K_V$ is therefore trace-norm holomorphic on $\Omega$, with
\begin{equation}\label{eq:extensions-below-trace}
 \|\mathcal K_V\|_1\le Ce^{-2\lambda\operatorname{Re}r}.
\end{equation}
The constants depend on $V$ only through $M$ and $R$.

\textit{The determinant.} Let $V_n=V\mathbf 1_{[-n,n]^2}$. Since $V_n$ is bounded with compact support, Proposition~\ref{prop:endpoint-free} for compactly supported $V$ (Section~\ref{sec:constructs}) gives $Z_{V_n}=\det(I-\mathcal K_{V_n})$. As $n\to\infty$, $\mathcal K_{V_n}\to\mathcal K_V$ in trace norm at each point of $\Omega$: for fixed $y<y'$ and $w,w'$, every continuous bridge path from $(y,w)$ to $(y',w')$ stays in $[-n,n]^2$ once $n$ is large, so $P_{y,y';V_n}(w,w')\to P_{y,y';V}(w,w')$ by dominated convergence, and the trace majorants above do not depend on $n$. Since $\|\mathcal K_{V_n}\|_1$ is bounded on $\Omega$ uniformly in $n$, Cauchy's formula bounds the derivatives of $\mathcal K_{V_n}$ uniformly in $n$ on compact subsets of $\Omega$, so the convergence is uniform there. Also $Z_{V_n}\to Z_V$, by Lemma~\ref{lem:occ-continuity}, since $|V_n|,|V|\le Me^{R}e^{-u}$. Therefore $Z_V=\det(I-\mathcal K_V)$ at real points.

\textit{The flows.} The kernels $\mathcal K_{V_n}$ satisfy the flows \eqref{eq:QR-flows}, by Proposition~\ref{prop:integrable-det}, and their derivatives converge locally uniformly on $\Omega$, so $\mathcal K_V$ satisfies the flows as well. By \eqref{eq:extensions-below-trace} and Cauchy's formula, every derivative of $\mathcal K_V$ in $(t,x,r)$ tends to zero in trace norm as $r\to+\infty$.
\end{proof}

\subsubsection{The example of Remark~\ref{rem:fredholm-scope}}
\label{app:parabolic-obstruction}

\begin{proof}[Proof of the claims in Remark~\ref{rem:fredholm-scope}]
Let $h$, $\chi$ and $W(y,w)=\chi(y)e^{-y^2}h(w-y^2)$ be as in the remark, and let $(t,x,r)=(1,0,0)$. Since $W(y,w)$ vanishes unless $w<y^2-1$, we have $0\le W(y,w)\le\|h\|_\infty e^{-w}$, so Lemma~\ref{lem:occ-exp} gives $\E e^{\theta I_W}<\infty$ for every $\theta>0$.

We show that, for all $u,v\ge0$, the first integral in \eqref{eq:extensions-two-insertions} is finite and the second is $+\infty$. Substitute $w=y^2+z$ and $w'=y'^2+z'$. By \eqref{eq:S},
\begin{equation}\label{eq:parabolic-S}
 S_{1,y}(w,u)=e^{-y^3/3+y(u-z)}\Ai(u-z),\qquad S_{1,-y'}(w',v)=e^{y'^3/3-y'(v-z')}\Ai(v-z').
\end{equation}
In the product $S_{1,y}(w,u)S_{1,-y}(w,v)$ the cubic terms cancel, so the first integral is
\[
 \int_\R\chi(y)e^{-y^2+y(u-v)}\diff y\int_\R h(z)\Ai(u-z)\Ai(v-z)\diff z,
\]
which is finite. In the second integral, $h(z)h(z')$ vanishes unless $z,z'\in(-2,-1)$, and then $\Ai(u-z)$ and $\Ai(v-z')$ are positive, so the integrand is nonnegative. By \eqref{eq:FK}, $P_W(y,y';w,w')\ge e^{-\sqrt\pi\|h\|_\infty}p_{y'-y}(w'-w)$, because $\int_y^{y'}W(s,B(s))\diff s\le\sqrt\pi\|h\|_\infty$ for every path $B$. By \eqref{eq:parabolic-S} and \eqref{eq:heat-kernel}, with $d=y'-y$,
\begin{equation}\label{eq:parabolic-cubic}
 \begin{aligned}
 &S_{1,y}(w,u)p_d(w'-w)S_{1,-y'}(w',v)\\
 &\quad=\frac{\Ai(u-z)\Ai(v-z')}{\sqrt{4\pi d}}\exp\Bigl\{yu-y'v+\frac{d^3}{12}+\frac d2(z+z')-\frac{(z'-z)^2}{4d}\Bigr\}.
 \end{aligned}
\end{equation}
Let $J\subset(-2,-1)$ be a closed interval on which $h$ is bounded below by a positive constant. If $3\le y\le4$, $y'\ge5$ and $z,z'\in J$, then $d\ge y'-4$, and the lower bound for $P_W$ and \eqref{eq:parabolic-cubic} show that the integrand is at least
\[
 \frac c{\sqrt{y'}}\exp\Bigl\{\frac{(y'-4)^3}{12}-y'^2-(v+2)y'\Bigr\},
\]
where $c>0$ does not depend on $y'$. The integral of this function over $y'\ge5$ is infinite, so the second integral is $+\infty$.

The failure persists under compact approximation. Let $W_n\in C_c^\infty(\R^2)$ with $0\le W_n\uparrow W$. Then $Z_{W_n}\to Z_W$ by monotone convergence, and $Z_{W_n}=\det(I-\mathcal K_{W_n})$ by Theorem~\ref{thm:fredholm}, with $\mathcal K_{W_n}$ given by \eqref{eq:extensions-two-insertions}. But $\langle f,\mathcal K_{W_n}f\rangle\to-\infty$ for $f=\mathbf 1_{[0,1]}$, so the operators $\mathcal K_{W_n}$ do not converge, even weakly. Indeed, the first integral for $W_n$ is at most the first integral for $W$, which is bounded for $0\le u,v\le1$. The integrand of the second integral for $W_n$ is nonnegative and converges pointwise to the integrand for $W$, since $P_{W_n}\to P_W$ by \eqref{eq:FK}. So, by Fatou's lemma, the second integral for $W_n$, integrated over $0\le u,v\le1$, tends to $+\infty$.
\end{proof}

\subsection{Bounded and growing data}\label{app:bounded-growing}
This subsection proves Proposition~\ref{prop:bounded-stability} and Theorem~\ref{thm:growing-data}.
For a function $F$ of $(t,x,r)$, let
\begin{equation}\label{eq:extensions-logKP}
 G[F]=\tfrac1{12}F_{rrrr}+\tfrac12F_{rr}^2+\tfrac14F_{xx}.
\end{equation}
For $V$ as in Theorem~\ref{thm:measurable-data}(i), \eqref{eq:QR-KP} and
Proposition~\ref{prop:qr-kernel}(iii) give
\begin{equation}\label{eq:extensions-logKP-identity}
 \partial_t\partial_r\log Z_V=-G[\log Z_V].
\end{equation}

\begin{lemma}\label{lem:logKP-limit}
Suppose smooth $F_n$ on $(0,T)\times\R^2$ satisfy
$\partial_t\partial_rF_n=-G[F_n]$ and converge to $F$
with all $(x,r)$ derivatives through order five, locally uniformly.
Suppose that for every compact $J\subset(0,T)$,
$R<\infty$, and $r_0\in\R$, $F$ and these derivatives
are bounded on $J\times[-R,R]\times[r_0,\infty)$ by
$C\omega(r)$, with $\omega\in L^1([r_0,\infty))$ and
$\omega(r)\to0$ as $r\to\infty$.
Then $F_t,F_{tr},F_{trr}$ exist and are continuous, and
$\phi=-2F_{rr}$ solves the equation in \eqref{KP} on $(0,T)\times\R^2$, with
$\partial_r^{-1}\phi=-2F_r$.
\end{lemma}
\begin{proof}
Integrating \eqref{eq:extensions-logKP-identity} for $F_n$ from $s$ to $t$
and passing to the limit gives
$F_r(t)-F_r(s)=-\int_s^t G[F](v)\diff v$.
These bounds, together with
$F(\infty)=F_r(\infty)=0$, allow us to integrate in $r$:
\[
 F(t,x,r)-F(s,x,r)
   =\int_s^t\int_r^\infty G[F](v,x,u)\diff u\diff v.
\]
Thus $F_t=\int_r^\infty G[F](t,x,u)\diff u$,
$F_{tr}=-G[F]$, and $F_{trr}=-\partial_rG[F]$.
Since $F_r(\infty)=0$, we also have
$\partial_r^{-1}\phi=-2F_r$. Substitution in
\eqref{KP} completes the proof.
\end{proof}

\subsubsection{Bounded data}\label{app:bounded}

\begin{proof}[The operator estimates in the proof of Proposition~\ref{lem:bounded-det}]
Let $V$ vanish outside $(a,b)\times(-R,R)$, and put $M=\sup|V|$.

\textit{Hilbert--Schmidt.} For $y,z\in[a,b]$ and $w,w'\in[-R,R]$, the factors $e^{\pm\gamma}$ and the Airy integral in \eqref{eq:bounded-kernel} are bounded, and $\int_{-R}^R\int_{-R}^Rp_s(w'-w)^2\diff w\diff w'\le2R(8\pi s)^{-1/2}$. So
\begin{equation}\label{eq:bounded-local-HS}
 \int_{[-R,R]^2}|\widehat K^\gamma(y,w;z,w')|^2\diff w\diff w'\le C\bigl(1+|z-y|^{-1/2}\bigr),\qquad y,z\in[a,b].
\end{equation}
Since $|z-y|^{-1/2}$ is integrable, $\mathsf B_V$ and $\mathsf H$ are Hilbert--Schmidt. So are $\mathsf P$ and $\mathsf Q$, since $\Ai$ decays faster than any exponential.

\textit{The series.} Suppose first that $V$ is continuous, so that Theorem~\ref{thm:fredholm} and the expansion of Appendix~\ref{app:remarks} apply. By \eqref{eq:S}, with $u=r+t^{1/3}\lambda$, the kernels of $\mathsf P$ and $\mathsf Q$ are $t^{1/6}e^{-\gamma_y(w)}S_{t,x-y}(w,u)W(y,w)$ and $t^{1/6}e^{\gamma_z(w')}S_{t,z-x}(w',u)\widetilde W(z,w')$. So in the kernel of $\mathsf Q\mathsf H^{n-1}\mathsf P$ the factors $e^{\pm\gamma}$ telescope, and the products $\widetilde WW$ give back $V$ at each point. After conjugation by $f\mapsto t^{1/6}f(r+t^{1/3}\cdot)$, this is the kernel of $\chi_rK^{[n]}_{t,x;V}\chi_r$ in \eqref{eq:ordered-kernel}. In the kernel of $\mathsf H^n$, the heat kernels compose, and the intermediate positions are ordered, so they range over a simplex of volume $(z-y)^{n-1}/(n-1)!$. Schur's test therefore gives $\|\mathsf H^n\|\le C(M(b-a))^n/n!$. So the Neumann series for $(I+\mathsf H)^{-1}$ converges, and $(I+\mathsf H)^{-1}\mathsf P\mathsf Q$ is trace class. By the same bound, $\mathsf H$ has spectral radius $0$, so it has no nonzero eigenvalues, and $\det_2(I+\mathsf H)=1$.

\textit{The trace.} On the diagonal $(z,w')=(y,w)$, the exponentials in $\mathsf P$ and $\mathsf Q$ cancel, and the integral over $\lambda$ is $t^{-1/3}\rho_{\rm Ai}(\xi_y(w))=\varrho(y,w)$, the density in \eqref{eq:occ-density}. So $\operatorname{Tr}(\mathsf P\mathsf Q)=\int_{\R^2}V\varrho=\E I_V$.

\textit{Measurable $V$.} Take continuous $V_n$ vanishing outside $(a,b)\times(-R,R)$, with $|V_n|\le M$ and $V_n\to V$ almost everywhere. By \eqref{eq:bounded-local-HS} and dominated convergence, $\mathsf B_{V_n}\to\mathsf B_V$ in Hilbert--Schmidt norm, and $\det_2$ is continuous in that norm. Lemma~\ref{lem:occ-continuity} gives $\E I_{V_n}\to\E I_V$ and $Z_{V_n}\to Z_V$.
\end{proof}

We prove Proposition~\ref{prop:bounded-stability} from Proposition~\ref{prop:bounded-factor}, with $\Theta$, $\Omega$ and the operators defined before it in Section~\ref{sec:determinants}.

We write $\mathfrak S_1$ and $\mathfrak S_2$ for the trace class and the Hilbert--Schmidt class, with norms $\|\cdot\|_1$ and $\|\cdot\|_2$, and $\|\cdot\|$ for the operator norm.

\begin{lemma}\label{lem:bounded-operators}
Let $V$ satisfy \eqref{eq:bounded-data-bound}.
\begin{enumerate}
\item[\textup{(i)}] $I+\mathsf H_\Omega$ is invertible for every $(t,x,r)\in\Theta$, and
\[
 (I+\mathsf H_\Omega)^{-1}\mathsf P_\Omega\mathsf Q_\Omega\in\mathfrak S_1,\qquad\widetilde{\mathsf B}\in\mathfrak S_2,\qquad\E I^{t,x}_{V\mathbf 1_{\Omega^c}}(r)\in\R
\]
are $C^\infty$ functions of $(t,x,r)\in\Theta$. Each of their derivatives is bounded on $\Theta$ by a constant that depends on $V$ only through $M$ and $\delta$, and not on $x_0$.
\item[\textup{(ii)}] Let $V_n$ satisfy \eqref{eq:bounded-data-bound} with the same $M$ and $\delta$, and let $V_n\to V$ in $L^1_{\rm loc}(\R^2)$. Then every derivative of each function in \textup{(i)} for $V_n$ converges to the same derivative for $V$, uniformly on compact subsets of $\Theta$.
\end{enumerate}
\end{lemma}

We prove Lemma~\ref{lem:bounded-operators} at the end of this section.

\begin{proof}[Proof of Proposition~\ref{prop:bounded-stability}]
\textit{Scalar bounds.} We first show that for every $T>0$ and $r_0\in\R$ there are $c>0$ and $C$, depending only on $M$, $\delta$, $T$ and $r_0$, such that
\begin{equation}\label{eq:bounded-scalar}
 Z_V\ge c,\qquad|Z_V-1|+|\log Z_V|\le Ce^{-\delta r},\qquad0<t\le T,\ x\in\R,\ r\ge r_0.
\end{equation}
By Lemmas~\ref{lem:occ-moments}(i) and~\ref{lem:occ-exp}, applied to $|V|$ with $A=M$ and $B=0$, $\|I_{|V|}^{t,x}(r)\|_2\le Ce^{-\delta r}$ and $\E e^{2I_{|V|}^{t,x}(r)}\le C$ for all $x\in\R$. Jensen's inequality gives $Z_V\ge e^{-\E I_{|V|}}\ge c$. By \eqref{eq:occ-Z-difference} with $q=1$ and the constant $0$, $|Z_V-1|\le\|I_V\|_2(\|e^{I_{V_-}}\|_2+1)\le Ce^{-\delta r}$, and, since $Z_V\ge c$, $|\log Z_V|\le|Z_V-1|/\min\{c,1\}$. In particular, $I_V^{t,x}(r)$ converges absolutely almost surely, and $0<Z_V<\infty$.

\textit{Smoothness and stability.} The function $A\mapsto\det_2(I-A)$ is entire on $\mathfrak S_2$, with $|\det_2(I-A)|\le e^{\|A\|_2^2/2}$ \cite[Chapter~9]{Simon}, so by Cauchy's estimates its derivatives are bounded on bounded sets. By \eqref{eq:bounded-factor}, Lemma~\ref{lem:bounded-operators}(i) and the chain rule, $Z_V$ is therefore $C^\infty$ on $\Theta$, with derivatives bounded as in Lemma~\ref{lem:bounded-operators}(i). If $V_n$ are as in Proposition~\ref{prop:bounded-stability}(i), Lemma~\ref{lem:bounded-operators}(ii) gives $\partial^\alpha Z_{V_n}\to\partial^\alpha Z_V$ uniformly on compact subsets of $\Theta$. Since $Z_V\ge c>0$ by \eqref{eq:bounded-scalar}, with $c$ depending on $V$ only through $M$ and $\delta$, the chain rule gives the same two statements for $\log Z_V$. As the sets $\Theta$ cover $(0,\infty)\times\R^2$, this proves Proposition~\ref{prop:bounded-stability}(i).

\textit{Bounds at large heights.} We interpolate between the decay of $\log Z_V$ and the boundedness of its derivatives. Let $k=|\alpha|$, and choose an integer $N>k$ with $\delta(1-k/N)\ge\delta'$. Fix $t_0\le t\le t_1$, $x\in\R$ and $r\ge r_0$, and let $\mathcal R$ be the cube of side $\min\{1,t_0\}$ centered at $(t,x,r)$. On $\mathcal R$, \eqref{eq:bounded-scalar} gives $|\log Z_V|\le Ce^{-\delta r}$, and the previous step, applied with $x_0=x$ and a set $\Theta$ containing $\mathcal R$, gives $|\partial^\beta\log Z_V|\le C$ for $|\beta|\le N$. Neither constant depends on $x$. The Gagliardo--Nirenberg inequality on $\mathcal R$,
\[
 \sup_{\mathcal R}|\partial^\alpha f|\le C\Bigl(\sup_{\mathcal R}|f|+\bigl(\sup_{\mathcal R}|f|\bigr)^{1-k/N}\max_{|\beta|=N}\bigl(\sup_{\mathcal R}|\partial^\beta f|\bigr)^{k/N}\Bigr),
\]
applied to $f=\log Z_V$, gives $|\partial^\alpha\log Z_V(t,x,r)|\le Ce^{-\delta(1-k/N)r}\le Ce^{-\delta'r}$, which is \eqref{eq:bounded-upper}.
\end{proof}

\begin{proof}[Proof of Lemma~\ref{lem:bounded-operators}]
\textit{Step~1: $\mathsf P_\Omega$ and $\mathsf Q_\Omega$.} It is elementary to show, from the decay of $\Ai$ \cite[\S9.7]{DLMF}, that, for $(t,x,r)\in\Theta$, $(y,w)\in\Omega$ and $\lambda\ge0$,
\[
 \bigl|e^{(x-y)\lambda/t^{2/3}}\Ai(\xi_y(w)+\lambda)\bigr|\le Ce^{-c(|y-x|^3+\lambda^{3/2}+w_-^{3/2})},
\]
since $(x-y)\lambda/t^{2/3}\le\frac{2}{3\sqrt3}\bigl((y-x)^2/t^{4/3}+\lambda\bigr)^{3/2}$. On $\Omega$, $|W|^2\le M$ and $w\le1+\langle y-x_0\rangle$, so the kernel of $\mathsf P_\Omega$ is square integrable over $\Omega\times\R_+$, uniformly on $\Theta$. The same holds for $\mathsf Q_\Omega$. Derivatives in $(t,x,r)$ multiply these kernels by polynomials in $y-x$, $\lambda$ and $\xi_y(w)$, or replace $\Ai$ by some $\Ai^{(k)}$, which obeys the same bound up to a polynomial factor. The exponential decay in the bound above absorbs these polynomials, so the derivatives of the kernels are square integrable as well.

\textit{Step~2: $\mathsf H_\Omega$ is bounded and compact.} Write $\mathsf H_\Omega=W\mathsf H_0\widetilde W$, where
\[
 \mathsf H_0(y,w;z,w')=\mathbf 1_\Omega(y,w)\mathbf 1_{y<z}p^\gamma(y,w;z,w')\mathbf 1_\Omega(z,w').
\]
Completing the square in the exponent of $p^\gamma$, and using $w,w'\le3+|y-x|+s$ on $\Omega$, one finds that for $(t,x,r)\in\Theta$ and $(y,w),(z,w')\in\Omega$ with $s=z-y>0$,
\begin{equation}\label{eq:bounded-heat}
 |p^\gamma(y,w;z,w')|\le Ce^{Cs-cs((y-x)^2+s^2+w_-+r_+)}p_s\Bigl(w'-w+\frac{(2x-y-z)s}t\Bigr).
\end{equation}
Integrating in $w'$ and then in $s$ gives
\begin{equation}\label{eq:bounded-schur}
 \int|\mathsf H_0(y,w;z,w')|\diff z\diff w'\le\frac C{1+(y-x)^2+w_-},
\end{equation}
and the integral over $(y,w)$ is bounded as well, so $\mathsf H_0$ is bounded by Schur's test. It is also compact: with $B_L=\{(y,w)\in\Omega:|y-x_0|\le L,\ w\ge-L\}$, the operator $\mathbf 1_{B_L}\mathsf H_0$ is Hilbert--Schmidt, and by \eqref{eq:bounded-schur} and Schur's test, $\|\mathbf 1_{\Omega\setminus B_L}\mathsf H_0\|\le CL^{-1/2}$. Since $|W|,|\widetilde W|\le M^{1/2}$, $\mathsf H_\Omega$ is also bounded and compact. Derivatives in $(t,x,r)$ multiply $p^\gamma$ by polynomials in the argument of $p_s$ and in derivatives of the exponent, and the decay above absorbs these, so every $\partial^\alpha\mathsf H_0$ and $\partial^\alpha\mathsf H_\Omega$ is also bounded and compact.

\textit{Step~3: the inverse.} Since $|V|\le M$ and $p^\gamma>0$, the argument for $\mathsf H^n$ in the proof of Proposition~\ref{lem:bounded-det}, with the indicators of $\Omega$ at the intermediate points bounded by $1$, gives
\[
 |\mathsf H_\Omega^n(y,w;z,w')|\le\frac{M^ns^{n-1}}{(n-1)!}\mathsf H_0(y,w;z,w').
\]
By \eqref{eq:bounded-heat}, as in Step~2,
\[
 \|\mathsf H_\Omega^n\|\le\frac{CM^n}{(n-1)!}\int_0^\infty s^{n-1}e^{Cs-cs^3}\diff s\le C_0C_1^n\frac{\Gamma(n/3)}{\Gamma(n)}.
\]
So the Neumann series for $(I+\mathsf H_\Omega)^{-1}$ converges in operator norm, uniformly on $\Theta$, and by Step~2 the derivatives of $(I+\mathsf H_\Omega)^{-1}$ are bounded.

\textit{Step~4: $\mathsf E$ is Hilbert--Schmidt.} The kernel of $\mathsf E$ is $W\widehat K^\gamma\widetilde W$ off $\Omega\times\Omega$. For $y<z$, the integral in \eqref{eq:bounded-kernel} taken over all $\lambda\in\R$ equals $p^\gamma$, by the full-line Airy semigroup identity, so
\begin{equation}\label{eq:bounded-kernel-minus}
 \widehat K^\gamma(y,w;z,w')=-t^{-1/3}\int_0^\infty e^{-(z-y)\lambda/t^{2/3}}\Ai(\xi_y(w)-\lambda)\Ai(\xi_z(w')-\lambda)\diff\lambda,\qquad y<z.
\end{equation}
Fix $0<s_0<\delta t_0/16$. When $z-y\ge s_0$ we bound $\widehat K^\gamma$ using \eqref{eq:bounded-kernel-minus}, and when $z-y<s_0$ using \eqref{eq:bounded-kernel}, where $e^{(z-y)\lambda/t^{2/3}}\le e^{s_0\lambda/t^{2/3}}$ in its integral. For these integrals, it is elementary to show that, for $(t,x,r)\in\Theta$ and $(y,w)\in\R^2$,
\begin{equation}\label{eq:bounded-airy-integrals}
 |V(y,w)|\int_0^\infty\Bigl(e^{s_0\lambda/t^{2/3}}\Ai(\xi_y(w)+\lambda)^2+e^{-s_0\lambda/t^{2/3}}\Ai(\xi_y(w)-\lambda)^2\Bigr)\diff\lambda\le Ce^{-c((y-x)^2+|w|)},
\end{equation}
since $|V(y,w)|\le Me^{-\delta w_+}$ and $s_0\xi_y(w)_-/t^{2/3}\le\delta w_+/16+C$. By \eqref{eq:bounded-airy-integrals} and the Cauchy--Schwarz inequality in $\lambda$, the kernel of $\mathsf E$ without the heat term of \eqref{eq:bounded-kernel} for $z-y<s_0$ is square integrable.

It remains to bound the heat term in \eqref{eq:bounded-kernel} for $0<s=z-y<s_0$. Since $\mathsf E$ vanishes on $\Omega\times\Omega$, one of the two points lies outside $\Omega$. Completing the square as in Step~2, but for arbitrary heights, gives \eqref{eq:bounded-heat} with an extra factor $e^{s(w_++w'_+)/(2t)}\le e^{\delta(w_++w'_+)/32}$, so by \eqref{eq:bounded-data-bound}
\[
 |V(y,w)V(z,w')|p^\gamma(y,w;z,w')^2\le Cs^{-1/2}e^{-\delta(w_++w'_+)/2}p_{s/2}\Bigl(w'-w+\frac{(2x-y-z)s}t\Bigr).
\]
Since $e^{-\delta w/2}$ is integrable over the complement of $\Omega$, the right side is integrable over pairs of points with $0<z-y<s_0$ and one point outside $\Omega$.

Derivatives in $(t,x,r)$ multiply the integrands and $p^\gamma$ by polynomials, as in Steps~1 and~2, and the decay above absorbs these. So every derivative of the kernel of $\mathsf E$ is square integrable, uniformly on $\Theta$.

\textit{Step~5: the mean above the cut, and $\widetilde{\mathsf B}$.} Since $\varrho(y,w)=t^{-1/3}\int_0^\infty\Ai(\xi_y(w)+\lambda)^2\diff\lambda$, \eqref{eq:bounded-airy-integrals} gives $|V\varrho|\le Ce^{-c((y-x)^2+|w|)}$. Derivatives in $(t,x,r)$ produce $\varrho$ or products $\Ai^{(k)}(\xi_y(w))\Ai^{(k')}(\xi_y(w))$, times polynomials in $\xi_y(w)$ and $y-x$, and the same decay absorbs these. So, by \eqref{eq:occ-density}, $\E I_{V\mathbf 1_{\Omega^c}}=\int_{\Omega^c}V\varrho$ is $C^\infty$ on $\Theta$ with bounded derivatives. By Steps 1--4 and the Leibniz rule, so are $\mathsf U\mathsf P_\Omega\mathsf Q_\Omega$, into $\mathfrak S_1$, and $\widetilde{\mathsf B}$, into $\mathfrak S_2$.

\textit{Step~6: dependence on $V$.} A subscript $n$ marks the objects built from $V_n$. The bounds of Steps 1--5 depend on $V$ only through $M$ and $\delta$, so they hold uniformly in $n$. In particular, the derivatives of each order are Lipschitz on $\Theta$, uniformly in $n$, so it suffices to prove convergence at each point of $\Theta$. Since every subsequence of $(V_n)$ has a further subsequence that converges almost everywhere, we may also assume that $V_n\to V$ almost everywhere. Then $W_n\to W$ and $\widetilde W_n\to\widetilde W$ almost everywhere, and $|W_n|,|\widetilde W_n|\le M^{1/2}$.

For $\mathsf P_\Omega$, $\mathsf Q_\Omega$ and $\mathsf E$, the kernels of their derivatives converge almost everywhere, and Steps 1 and 4 give square-integrable majorants. So the derivatives converge in $\mathfrak S_2$, by dominated convergence. The same argument, with the integrable majorant of Step~5, applies to the mean.

For $\mathsf H_\Omega$, multiplication by $W_n-W$ or by $\widetilde W_n-\widetilde W$ tends to zero strongly, by dominated convergence, and is bounded by $2M^{1/2}$. Since $\partial^\alpha\mathsf H_0$ is compact,
\[
 \partial^\alpha\mathsf H_{\Omega,n}-\partial^\alpha\mathsf H_\Omega=(W_n-W)(\partial^\alpha\mathsf H_0)\widetilde W_n+W(\partial^\alpha\mathsf H_0)(\widetilde W_n-\widetilde W)\longrightarrow0
\]
in operator norm. Then $\mathsf U_n-\mathsf U=\mathsf U_n(\mathsf H_\Omega-\mathsf H_{\Omega,n})\mathsf U\to0$, since $\|\mathsf U_n\|$ is bounded uniformly in $n$ by Step~3. The derivatives of $\mathsf U$ are sums of products of $\mathsf U$ and derivatives of $\mathsf H_\Omega$, so they converge as well. Finally, by the Leibniz rule, the derivatives of $\mathsf U\mathsf P_\Omega\mathsf Q_\Omega$ and $\widetilde{\mathsf B}$ converge in $\mathfrak S_1$ and $\mathfrak S_2$, respectively.
\end{proof}

\subsubsection{Growing data}

\begin{proof}[Proof of Theorem~\ref{thm:growing-data}]
Fix $T<T_*$ and a closed rectangle $\mathcal C=[-R,R]\times[r_0,r_1]$. Since $T<T_*$, each $\partial^\eta V$, $|\eta|\le5$, satisfies \eqref{eq:occ-exponential-bound} or \eqref{eq:occ-algebraic-bound}, so Lemma~\ref{lem:occ-moments} applies to it, and Lemma~\ref{lem:occ-exp} applies to $W=V_-$. By Lemma~\ref{lem:occ-trace}(iii), $0<Z_V<\infty$, and $Z_V\ge e^{-\E I_{|V|}}$ by Jensen's inequality, so $Z_V$ is bounded below by a positive constant on $(0,T]\times[-R,R]\times[r_0,\infty)$.

The proof approximates $V$ by the cutoffs $V_n(y,u)=\chi(y/n)\chi(u/n)V(y,u)$, where $\chi\in C_c^\infty(\R)$, $0\le\chi\le1$, $\chi=1$ on $[-1,1]$ and $\chi=0$ outside $[-2,2]$. By Theorem~\ref{thm:measurable-data}(i), the compact-data functions
$F_n=\log Z_{V_n}$ are smooth and satisfy
\eqref{eq:extensions-logKP-identity}. We verify convergence
of their spatial derivatives and the upper-height bounds required
by Lemma~\ref{lem:logKP-limit}.

\textit{Step~1: cutoff errors.} By the product rule, and since the derivatives of $\chi(y/n)\chi(u/n)$ are bounded uniformly in $n$, each $|\partial^\eta(V_n-V)|$, $|\eta|\le5$, satisfies the corresponding growth and decay bound in Theorem~\ref{thm:growing-data}, with a constant independent of $n$, and vanishes on $[-n,n]^2$. So Lemma~\ref{lem:occ-moments}(ii), with $L=N=n$, gives, for every $q<\infty$, as $n\to\infty$,
\begin{equation}\label{eq:extensions-cutoff}
 \sup\bigl\{\|I^{t,x}_{\partial^\eta(V_n-V)}(r)\|_q:0<t\le T,\ (x,r)\in\mathcal C\bigr\}\longrightarrow0,\qquad|\eta|\le5.
\end{equation}
Moreover $(V_n)_-\le V_-$ and $|V_n|\le|V|$, so the exponential moments of $I_{(V_n)_-}$ and the lower bound $Z_{V_n}\ge e^{-\E I_{|V|}}$ hold uniformly in $n$.

\textit{Step~2: differentiation under the expectation.} By the chain rule, $\partial^\eta e^{-I}=e^{-I}\mathcal P_\eta$, where $\mathcal P_\eta$ is a polynomial in the derivatives $\partial^\zeta I$, $0<\zeta\le\eta$. Let $F\in C_c^N$ for some $N$, and $t>0$. The proof of Lemma~\ref{lem:diff-under-E}, applied to derivatives of order at most $N$, shows that in the realization \eqref{eq:occ-scaled} $\partial^\zeta I_F=I_{\partial^\zeta F}$ for $|\zeta|\le N$, and that, for $|\eta|\le N$,
\begin{equation}\label{eq:extensions-bell}
 \partial^\eta Z_F=\E\bigl[e^{-I_F}\mathcal P_\eta\bigl((I_{\partial^\zeta F})_\zeta\bigr)\bigr].
\end{equation}

\textit{Step~3: positive time.} Next we show that, for $|\eta|\le5$, $\partial^\eta Z_{V_n}$ converges, uniformly on $(0,T]\times\mathcal C$, to the right side of \eqref{eq:extensions-bell} with $F=V$. Indeed, \eqref{eq:extensions-bell} holds for $F=V_n\in C_c^5$. Since $I_{V_n},I_V\ge-I_{V_-}$, the mean value theorem gives $|e^{-I_{V_n}}-e^{-I_V}|\le|I_{V_n-V}|e^{I_{V_-}}$, and $\mathcal P_\eta((I_{\partial^\zeta V_n})_\zeta)-\mathcal P_\eta((I_{\partial^\zeta V})_\zeta)$ is a sum of products, each containing a factor $I_{\partial^\zeta(V_n-V)}$. The convergence then follows from \eqref{eq:extensions-cutoff}, the moment bounds of Lemmas~\ref{lem:occ-moments}(i) and~\ref{lem:occ-exp}, and H\"older's inequality. Since the $\partial^\eta Z_{V_n}$ are continuous, the spatial derivatives of $Z_V$ through order five exist, are continuous in $(t,x,r)$, and are given by \eqref{eq:extensions-bell} with $F=V$. By the uniform lower bound on $Z_{V_n}$ and $Z_V$ and the chain rule, $\partial^\eta\log Z_{V_n}\to\partial^\eta\log Z_V$ uniformly on $(0,T]\times\mathcal C$, for $|\eta|\le5$.

\textit{Step~4: large height.} With $\omega$ as in Lemma~\ref{lem:occ-moments}(i), for $|x|\le R$, $r\ge r_0$ and $0<t\le T$,
\begin{equation}\label{eq:extensions-log-tail}
 |\partial^\eta \log Z_V(t,x,r)|\le C\omega(r),\qquad|\eta|\le5.
\end{equation}
For $|\eta|>0$, $\mathcal P_\eta$ has no constant term, so \eqref{eq:extensions-bell} with $F=V$, H\"older's inequality, \eqref{eq:occ-moment-bounds} and Lemma~\ref{lem:occ-exp} bound $\partial^\eta Z_V$ by $C\omega(r)$. For $\eta=0$, \eqref{eq:occ-Z-difference} with $a=0$ and the same bounds give $|Z_V-1|\le C\omega(r)$. Since $Z_V$ is bounded below on this set, the chain rule gives \eqref{eq:extensions-log-tail}.

\textit{Step~5: the initial trace.} Next we show that, as $t\downarrow0$, $Z_V(t)\to e^{-\frac12\partial_r^{-2}V}$ in $C^5(\mathcal C)$. Indeed, by Lemma~\ref{lem:occ-trace}(ii), $I_{\partial^\zeta V}^{t,x}(r)\to\tfrac12\partial_r^{-2}\partial^\zeta V(x,r)$ for $|\zeta|\le5$, and $e^{-I_V^{t,x}(r)}\to e^{-\frac12\partial_r^{-2}V(x,r)}$, in every $L^q$, $q<\infty$, uniformly for $(x,r)\in\mathcal C$. So, by \eqref{eq:extensions-bell} with $F=V$ and H\"older's inequality, $\partial^\eta Z_V(t)$ converges, uniformly on $\mathcal C$, to $e^{-\frac12\partial_r^{-2}V}\mathcal P_\eta\bigl((\tfrac12\partial_r^{-2}\partial^\zeta V)_\zeta\bigr)$, for $|\eta|\le5$. This limit is $\partial^\eta e^{-\frac12\partial_r^{-2}V}$, because $\partial_r^{-2}\partial^\zeta V=\partial^\zeta\partial_r^{-2}V$, by differentiation under the integral defining $\partial_r^{-2}V$. Since $Z_V$ is bounded below, the chain rule gives $\log Z_V(t)\to-\tfrac12\partial_r^{-2}V$ in $C^5(\mathcal C)$, and hence $\phi_V(t)=-2\partial_r^2\log Z_V(t)\to V$ in $C^3(\mathcal C)$.

\textit{Step~6: time derivatives and KP.}
Step~3 gives convergence of $F_n=\log Z_{V_n}$ and all
their spatial derivatives through order five, locally
uniformly at positive times. Step~4 supplies the weight
$\omega$ in \eqref{eq:extensions-log-tail}.
It is integrable at $+\infty$: in the case of algebraic decay,
this follows from $p>3+m/2$.
Lemma~\ref{lem:logKP-limit} therefore gives the required
time derivatives and the KP equation for
$\phi_V=-2\partial_r^2\log Z_V$.
\end{proof}

\begin{proof}[Proof of Remark~\ref{rem:smooth-data}]
For the first claim, the proof of Theorem~\ref{thm:growing-data} with $N$ in place of five gives, in Step~5, $\log Z_V(t)\to-\tfrac12\partial_r^{-2}V$ in $C^N_{\rm loc}$, and hence $\phi_V(t)\to V$ in $C^{N-2}_{\rm loc}$. For the second, let $V\in C_c^N$ with $N\ge2$. Then $V_-$ is bounded, and every $\partial^\eta V$ with $|\eta|\le N$ satisfies the bound in Theorem~\ref{thm:growing-data}(i) with $B=0$. Steps~2 and~5 of the proof of that theorem, with $N$ in place of five and applied to $V$ directly, use no derivative of $V$ of order above $N$, and give the same conclusion.
\end{proof}

\subsection{Transport and restart of Airy expectations}\label{app:representation}

\subsubsection{The transport identity}\label{app:representation-transport}

\textit{The analytic steps in the proof of Proposition~\ref{prop:transport}.} Since $R_\psi(s)$ contains $\frac14\partial_y^2\partial^{-1}\psi(s)$, which is independent of $u$ below the support of $\psi(s)$, the potential $\psi(s)+\varepsilon R_\psi(s)$ need not have compact support. It is bounded, however, and vanishes above the support of $\psi(s)$. So, by Proposition~\ref{prop:endpoint-free}, the determinant formula for $Z_{\psi(s)+\varepsilon R_\psi(s)}$ holds, with the kernel \eqref{eq:kernel}: the computation in the proof of that proposition extends, with Lemma~\ref{lem:extensions-rows} for Fubini's theorem. As in the proof of that proposition, $\chi_rK_{\tau-s,x;\psi(s)+\varepsilon R_\psi(s)}\chi_r$ and $\chi_rK_{t,x;\psi(s)}\chi_r$ are continuously differentiable in trace norm, in $\varepsilon$ and in $(s,t)$ respectively, with Duhamel's formula \eqref{eq:duhamel} for the derivatives of the Feynman--Kac kernels. Since $|\psi(s)|+|R_\psi(s)|\le Ce^{-u}$, Lemmas~\ref{lem:occ-moments}(i) and~\ref{lem:occ-exp} give $\E[|I_{R_\psi(s)}|e^{I_{|\psi(s)|}+I_{|R_\psi(s)|}}]<\infty$, which justifies differentiating under the expectation.

Proposition~\ref{prop:transport} assumes $\psi\in C_c^\infty(\R^3)$. The proof of Proposition~\ref{prop:transport-noncompact} in Appendix~\ref{app:representation-classical} needs it for a potential multiplied by cutoff functions, which has only three continuous spatial derivatives. The next lemma gives Proposition~\ref{prop:transport}, in integrated form, at that regularity.

\begin{lemma}\label{lem:representation-transport}
Let $\tau\in\R$, and let $\psi(s,y,u)$ be defined for $s$ in an open interval containing $[s_1,s_2]$, where $s_2<\tau$. Suppose that all $\psi(s)$ vanish outside a common compact subset of $\R^2$, and that $\partial_s\psi$ and the spatial derivatives of $\psi$ through order three exist and are continuous. Then, for $x,r\in\R$,
\[
 Z_{\psi(s_2)}(\tau-s_2,x,r)-Z_{\psi(s_1)}(\tau-s_1,x,r)=-\int_{s_1}^{s_2}\E\Bigl[e^{-I_{\psi(s)}^{\tau-s,x}(r)}I_{R_\psi(s)}^{\tau-s,x}(r)\Bigr]\diff s.
\]
\end{lemma}

\begin{proof}
Mollify $\psi$ in $(s,y,u)$ at scale $1/n$, and call the result $\psi_n$. The functions $\psi_n$ are smooth and vanish outside a common compact set, and $\psi_n$, $\partial_s\psi_n$ and the spatial derivatives of $\psi_n$ through order three converge uniformly on $[s_1,s_2]\times\R^2$ to those of $\psi$. Proposition~\ref{prop:transport} also needs compact support in $s$. Multiplying $\psi_n$ by a smooth compactly supported function of $s$ that equals $1$ near $[s_1,s_2]$ provides it, and changes neither side of \eqref{eq:transport} for $s$ near $[s_1,s_2]$. So Proposition~\ref{prop:transport} gives the derivative of $s\mapsto Z_{\psi_n(s)}(\tau-s,x,r)$ on $[s_1,s_2]$ as $-\E[e^{-I_{\psi_n(s)}^{\tau-s,x}(r)}I_{R_{\psi_n}(s)}^{\tau-s,x}(r)]$, and integrating gives the lemma for $\psi_n$.

It remains to let $n\to\infty$. Since $\partial_y^2\partial^{-1}\psi_n=\partial^{-1}\partial_y^2\psi_n$, and $\partial^{-1}$ is bounded in the supremum norm on functions supported in a fixed compact set, $R_{\psi_n}\to R_\psi$ uniformly on $[s_1,s_2]\times\R^2$. The functions $\psi_n(s)$ and $R_{\psi_n}(s)$, and their limits, vanish for large $|y|$ and above a common height, since $\partial^{-1}$ integrates from $+\infty$. So, uniformly in $s\in[s_1,s_2]$, $|\psi_n(s)|+|R_{\psi_n}(s)|\le Ce^{-u}$, and $|\psi_n(s)-\psi(s)|$ and $|R_{\psi_n}(s)-R_\psi(s)|$ are at most a multiple of $e^{-u}$ that tends to zero as $n\to\infty$. The times $\tau-s$ in the path integral lie in $[\tau-s_2,\tau-s_1]$. So Lemma~\ref{lem:occ-moments}(i) gives $\|I_{\psi_n(s)-\psi(s)}\|_2+\|I_{R_{\psi_n}(s)-R_\psi(s)}\|_2\to0$ uniformly in $s$, and Lemma~\ref{lem:occ-exp} bounds the exponential moments of $I_{|\psi_n(s)|}$ uniformly in $n$ and $s$. With $|e^{-a}-e^{-b}|\le|a-b|(e^{a_-}+e^{b_-})$, both sides of the lemma for $\psi_n$ converge to those for $\psi$, and the integrand converges boundedly in $s$.
\end{proof}

\subsubsection{Transport for noncompact potentials}\label{app:transport-noncompact}\label{app:representation-classical}

\begin{proof}[Proof of Proposition~\ref{prop:transport-noncompact}]
Fix $\tau$ and $x,r\in\R$, and let $s_-<s_1<s_2<s_+$ with $s_2\le\tau$. Proposition~\ref{prop:transport} requires compact support. We therefore apply its integrated form, Lemma~\ref{lem:representation-transport}, to a cutoff $\psi_{L,N}$ of $\psi$ that equals $\psi$ on a box $[-L,L]\times[-N,N]$, and let $L,N\to\infty$. The error vanishes in the limit because $R_{\psi_{L,N}}-R_\psi$ and $\psi_{L,N}-\psi$ both vanish on the box, apart from a small nonlocal term in the first, and the occupation of a potential that vanishes on the box is small.

Fix a compact interval $J\subset(s_-,s_+)$ whose interior contains $[s_1,s_2]$. On $J$ the bounds \eqref{eq:representation-bounds} hold with constants $C_{jk}$, $B$, $\beta$, $\delta$, $M$ and $\nu$. As in Appendix~\ref{app:representation-transport}, we write the variables of $\psi(s)$ as $(y,u)$, and $\partial=\partial_u$. With
\[
 \bar\omega(y,u)=e^{B(1+|y|)^\beta-\delta u},
\]
the first bound in \eqref{eq:representation-bounds} reads $|\partial_y^j\partial^k\psi(s)|\le C_{jk}\bar\omega$, and the bound on $\partial_s\psi$ reads $|\partial_s\psi(s)|\le C(\bar\omega+\bar\omega^2)$. By Lemma~\ref{lem:occ-trace}(iii), $Z_{\psi(s)}(\tau-s,x,r)$ is finite and positive for $s\in J$, $s<\tau$. From now on, all occupations are taken at $(\tau-s,x,r)$, and the constants $C$ and $c$ do not depend on $L$, $N$ or $s$. We first prove \eqref{eq:kpz-transport} for $s_2<\tau$.

\textit{Cutting $\psi$ off.} Let $\eta,\chi\in C^\infty(\R)$ take values in $[0,1]$, equal $1$ on $[-1,1]$ and vanish outside $[-2,2]$. For $N\ge1$ and $L\ge\max\{1,2|x|,2|r|\}$ put $\eta_L(y)=\eta(y/L)$, $\chi_N(u)=\chi(u/N)$ and
\[
 \psi_{L,N}=\eta_L\chi_N\psi.
\]
Then $\psi_{L,N}=\psi$ on the box $[-L,L]\times[-N,N]$, and $\psi_{L,N}$ vanishes outside $[-2L,2L]\times[-2N,2N]$. It satisfies the hypotheses of Lemma~\ref{lem:representation-transport} on the interior of $J$, so
\begin{equation}\label{eq:restart-cutoff}
 Z_{\psi_{L,N}(s_2)}(\tau-s_2,x,r)-Z_{\psi_{L,N}(s_1)}(\tau-s_1,x,r)=-\int_{s_1}^{s_2}\E\Bigl[e^{-I_{\psi_{L,N}(s)}}I_{R_{\psi_{L,N}}(s)}\Bigr]\diff s.
\end{equation}
We show that replacing $\psi_{L,N}$ by $\psi$ changes the integrand on the right, and each term on the left, by at most $C(e^{-cN}+e^{-cL^2})$.

\textit{The residual of $\psi_{L,N}$.} The bounds on $\partial_y\psi$ and $\partial_y^2\psi$ let $\partial_y$ pass under $\partial^{-1}$, and the product rule gives
\begin{equation}\label{eq:restart-residual}
 \begin{aligned}
 R_{\psi_{L,N}}={}&\eta_L\chi_NR_\psi\\
 &+\tfrac14\eta_L[\partial^{-1},\chi_N]\partial_y^2\psi\\
 &+\tfrac12\eta_L\chi_N(1-\chi_N)\psi\,\partial\psi-\tfrac12\eta_L\chi_N\chi_N'\psi^2+\tfrac1{12}\eta_L\bigl(\chi_N'''\psi+3\chi_N''\partial\psi+3\chi_N'\partial^2\psi\bigr)\\
 &+\tfrac12\eta_L(1-\eta_L)\chi_N\psi\,\partial(\chi_N\psi)+\tfrac12\eta_L'\partial^{-1}(\chi_N\partial_y\psi)+\tfrac14\eta_L''\partial^{-1}(\chi_N\psi),
 \end{aligned}
\end{equation}
where $([\partial^{-1},\chi_N]f)(y,u)=\int_u^\infty(\chi_N(u)-\chi_N(v))f(y,v)\diff v$. At points $(y,u)\in[-L,L]\times[-N,N]$, the first line equals $R_\psi$, and all terms in the last two lines vanish, since the cutoff functions are identically one in this rectangle. The only remaining contribution to $R_{\psi_{L,N}}-R_\psi$ is the commutator term in the second line, which need not vanish because $\partial^{-1}$ is nonlocal. For $|u|\le N$ it involves only heights above $N$, since $\chi_N=1$ on $[u,N]$, and so
\[
 \Bigl|\tfrac14\eta_L(y)\int_N^\infty(1-\chi_N(v))\partial_y^2\psi(y,v)\diff v\Bigr|\le C\int_N^\infty\bar\omega(y,v)\diff v\le Ce^{-\delta N/2}e^{B(1+|y|)^\beta-\delta u/2}.
\]
Off the box, the last three lines are at most $C(\bar\omega+\bar\omega^2)$: the bounds \eqref{eq:representation-bounds} bound $\psi$, $\partial\psi$, $\partial^2\psi$, $\partial_y\psi$ and $\partial_y^2\psi$ by $C\bar\omega$, the derivatives of $\eta_L$ and $\chi_N$ are bounded, and $\partial^{-1}$ and the commutator take a function bounded by $\bar\omega$ to one bounded by $\bar\omega/\delta$. By the same argument, with the bounds on $\partial^3\psi$ and $\partial_s\psi$, $|R_\psi|\le C(\bar\omega+\bar\omega^2)$, so $R_{\psi_{L,N}}-R_\psi$ is at most $C(\bar\omega+\bar\omega^2)$ off the box. Lemma~\ref{lem:occ-moments}(i), applied to the second line on the box, and Lemma~\ref{lem:occ-moments}(ii), applied to the restrictions of $\bar\omega$ and $\bar\omega^2$ to the complement of the box, give, for $s\in J$ with $s<\tau$,
\begin{equation}\label{eq:representation-residual-small}
 \bigl\|I_{R_{\psi_{L,N}}(s)}-I_{R_\psi(s)}\bigr\|_2\le C\bigl(e^{-cN}+e^{-cL^2}\bigr).
\end{equation}

\textit{Removing the cutoff.} The difference $\psi(s)-\psi_{L,N}(s)=(1-\eta_L\chi_N)\psi(s)$ vanishes on the box and is at most $C_{00}\bar\omega$, so Lemma~\ref{lem:occ-moments}(ii) gives $\|I_{\psi(s)-\psi_{L,N}(s)}\|_4\le C(e^{-cN}+e^{-cL^2})$. Since $\eta_L$ and $\chi_N$ take values in $[0,1]$, $\psi_{L,N}(s)_-\le\psi(s)_-\le M(1+u_-)^\nu$ and $|\psi_{L,N}(s)|\le|\psi(s)|\le C_{00}\bar\omega$, so Lemma~\ref{lem:occ-exp} bounds $\E\exp\{4I_{\psi(s)_-}\}$ and $\E\exp\{4I_{\psi_{L,N}(s)_-}\}$ uniformly in $L$, $N$ and $s\in J$ with $s<\tau$. The bound \eqref{eq:occ-Z-difference} with $q=2$ then gives
\begin{equation}\label{eq:restart-exp-difference}
 \bigl\|e^{-I_{\psi(s)}}-e^{-I_{\psi_{L,N}(s)}}\bigr\|_2\le C\bigl(e^{-cN}+e^{-cL^2}\bigr),\qquad s\in J,\ s<\tau.
\end{equation}
In particular, each term on the left of \eqref{eq:restart-cutoff} changes by at most $C(e^{-cN}+e^{-cL^2})$ when $\psi_{L,N}$ is replaced by $\psi$.

\textit{The integrand.} The integrand of \eqref{eq:restart-cutoff} minus the integrand of \eqref{eq:kpz-transport} is
\[
 \E\Bigl[e^{-I_{\psi_{L,N}(s)}}\bigl(I_{R_{\psi_{L,N}}(s)}-I_{R_\psi(s)}\bigr)\Bigr]+\E\Bigl[\bigl(e^{-I_{\psi_{L,N}(s)}}-e^{-I_{\psi(s)}}\bigr)I_{R_\psi(s)}\Bigr].
\]
By the Cauchy--Schwarz inequality, the first term is at most $C(e^{-cN}+e^{-cL^2})$, by \eqref{eq:representation-residual-small} and the exponential moments of $I_{\psi_{L,N}(s)_-}$. So is the second, by \eqref{eq:restart-exp-difference} and the bound $\|I_{R_\psi(s)}\|_2\le C$ of Lemma~\ref{lem:occ-moments}(i). Both bounds are uniform in $s\in[s_1,s_2]$, so letting $L,N\to\infty$ in \eqref{eq:restart-cutoff} proves \eqref{eq:kpz-transport}.

\textit{The endpoint.} Now let $s_2=\tau<s_+$, so that the interior of $J$ contains $[s_1,\tau]$. The bounds above hold for all $s\in J$ with $s<\tau$. So \eqref{eq:kpz-transport} holds with every $s'\in(s_1,\tau)$ in place of $s_2$, and its integrand is bounded on $[s_1,\tau)$, by Lemmas~\ref{lem:occ-moments}(i) and~\ref{lem:occ-exp}. As $s\uparrow\tau$, $\psi(s)\to\psi(\tau)$ locally uniformly, since $\psi$ is continuous. The potentials $\psi(s)$, $s\in J$, share the bounds $|\psi(s)|\le C_{00}\bar\omega$ and $\psi(s)_-\le M(1+u_-)^\nu$. So Lemma~\ref{lem:occ-trace}(ii), with $V_n=\psi(s_n)$, $V=\psi(\tau)$ and $t_n=\tau-s_n$ for any sequence $s_n\uparrow\tau$, gives $Z_{\psi(s)}(\tau-s,x,r)\to\exp\{-\frac12\partial_r^{-2}\psi(\tau,x,r)\}$. Letting $s'\uparrow\tau$ proves \eqref{eq:kpz-transport} for $s_2=\tau$.
\end{proof}

\subsubsection{The converse of the restart identity}\label{app:kpz-converse-proof}

\begin{proof}[Proof of Proposition~\ref{prop:kpz-converse}]
Write the equation in \eqref{KP} as $\partial_t\phi=\mathcal N(\phi)$, with
\[
 \mathcal N(V)=\tfrac12V\partial_rV-\tfrac1{12}\partial_r^3V-\tfrac14\partial_x^2\partial_r^{-1}V.
\]
Fix $s\in(0,T)$ and let $V=\phi(s)$. The function $\psi(\sigma)=V$, $\sigma\in\R$, satisfies the hypotheses of the transport identity \eqref{eq:kpz-transport}, with residual $R_\psi=-\mathcal N(V)$. Taking $s_2=\tau$ in \eqref{eq:kpz-transport} and writing $t=\tau-s_1$,
\[
 Z_V(t,x,r)=\exp\{-\tfrac12\partial_r^{-2}V(x,r)\}-\int_0^t\E\Bigl[e^{-I_V^{t',x}(r)}I_{\mathcal N(V)}^{t',x}(r)\Bigr]\diff t'.
\]
As $t'\downarrow0$, Lemma~\ref{lem:occ-trace}(ii) gives $e^{-I_V^{t',x}(r)}\to\exp\{-\frac12\partial_r^{-2}V(x,r)\}$ and, applied separately to the quadratic and the linear terms of $\mathcal N(V)$, each of which satisfies a bound of the form \eqref{eq:occ-exponential-bound}, $I_{\mathcal N(V)}^{t',x}(r)\to\frac12\partial_r^{-2}\mathcal N(V)(x,r)$ in $L^2$, uniformly for $(x,r)$ in compact sets. Hence, uniformly on compact sets,
\[
 \log Z_{\phi(s)}(t,x,r)=-\tfrac12\partial_r^{-2}\phi(s,x,r)-\tfrac t2\partial_r^{-2}\mathcal N(\phi(s))(x,r)+o(t)\qquad(t\downarrow0).
\]
By hypothesis, the left side is $-\frac12\partial_r^{-2}\phi(s+t,x,r)$. For $\eta\in C_c^\infty(\R^2)$, integration by parts gives $\int\frac12\partial_r^{-2}g\,\partial_r^2\eta=\frac12\int g\eta$, so pairing with $\partial_r^2\eta$ gives
\[
 \lim_{t\downarrow0}\frac1t\int_{\R^2}\bigl(\phi(s+t)-\phi(s)\bigr)\eta=\int_{\R^2}\mathcal N(\phi(s))\eta.
\]
The right side is continuous in $s$, by the continuity of the derivatives of $\phi$ and the bounds \eqref{eq:representation-bounds}. A continuous function whose right derivative is continuous is continuously differentiable, so for $0<s_1<s_2<T$,
\[
 \phi(s_2)-\phi(s_1)=\int_{s_1}^{s_2}\mathcal N(\phi(\sigma))\diff\sigma,
\]
pointwise, since both sides are continuous. So $\partial_s\phi=\mathcal N(\phi)$, and $\phi$ is a classical solution.
\end{proof}

\subsection{Recovering the initial time}\label{app:zero}

Let $\phi$ be as in Theorem~\ref{thm:bounded-unique}, replaced by its smooth representative as in the proof below. As in the proof of Theorem~\ref{thm:representation}(ii) in Section~\ref{sec:representation}, Theorem~\ref{thm:representation}(i) applies to $\phi$, and it remains to prove \eqref{eq:zero-limit} for $0<\tau<T$ and $x,r\in\R$.

\begin{proof}[Proof of Theorem~\ref{thm:bounded-unique}]
For a real test function $\psi$ and almost every $s\in(0,T)$,
\[
 \Bigl|\int\phi(s)\psi\Bigr|\le\int Me^{-\delta r_+}|\psi|.
\]
This inequality passes to $V$ as $s\downarrow0$ outside the null set. So $\psi\mapsto\langle V,\psi\rangle$ extends to a linear functional of norm at most one on $L^1(Me^{-\delta r_+}\diff x\diff r)$, and $V=\vartheta Me^{-\delta r_+}$ for a measurable $\vartheta$ with $|\vartheta|\le1$; in particular $|V|\le Me^{-\delta r_+}$.

By Lemma~\ref{lem:representation-smoothing} below, with $J'=(0,T)$, $C=M$ and $B=\beta=q=0$, $\phi$ agrees almost everywhere with a smooth function, which solves \eqref{KP} classically and satisfies on each compact interval $J\subset(0,T)$ the bounds \eqref{eq:representation-bounds} with $\beta=1$ and some $A<\infty$ and $\lambda>0$ in place of $B$ and $\delta$; the cap $\phi_-\le M$ holds everywhere by continuity. We replace $\phi$ by this function, so Theorem~\ref{thm:representation}(i) applies to it. Fix $0<\tau<T'<T$. By the last assertion of Lemma~\ref{lem:representation-smoothing}, with $J'=(0,T')$, $\phi(s)\to V$ in $L^2_{\rm loc}(\R^2)$ as $s\downarrow0$. Since $e^{-\delta u_+}\le e^{-\delta u}$, the potentials $\phi(s)$, $0<s<T'$, and $V$ are bounded by $Me^{-\delta u}$, and their negative parts by $M$. Lemma~\ref{lem:occ-continuity}, with $V_n=\phi(s_n)$ and $t_n=\tau-s_n$ for $s_n\downarrow0$, gives \eqref{eq:zero-limit}. As in the proof of Theorem~\ref{thm:representation}(ii), $\phi(\tau)=-2\partial_r^2\log Z_V(\tau)=\phi_V(\tau)$ for $0<\tau<T$.
\end{proof}

The proof of Theorem~\ref{thm:bounded-unique} used the following lemma: a weak solution with bounds on its values alone is smooth at positive times, with derivative bounds of the form \eqref{eq:representation-bounds}.

\begin{lemma}\label{lem:representation-smoothing}
Let $\phi$ be a distributional solution of \eqref{KP} on $(0,T)\times\R^2$, and let $J'\subset(0,T)$ be an interval such that, almost everywhere on $J'\times\R^2$,
\[
 |\phi(t,x,r)|\le Ce^{B(1+|x|)^\beta}(1+r_-)^qe^{-\delta r_+},
\]
with $C,B<\infty$, $0\le\beta<2$, $q\ge0$ and $\delta>0$. Then $\phi$ agrees almost everywhere with a function that is $C^\infty$ on $(\operatorname{int}J')\times\R^2$. Let $\kappa=\max\{1,\beta\}$. For every compact interval $J\subset\operatorname{int}J'$ and every $m\ge0$ there are $C_m,A<\infty$ and $\lambda>0$ such that
\[
 |\partial_x^j\partial_r^k\phi(t,x,r)|\le C_me^{A(1+|x|)^\kappa-\lambda r},\qquad t\in J,\quad j+k\le m.
\]
If $J'=(0,T')$, and $\phi(t)\to V$ in $\mathcal D'(\R^2)$ as $t\downarrow0$ outside a null set, for a measurable $V$ with the same bound, then the smooth function satisfies $\phi(t)\to V$ in $L^2_{\rm loc}(\R^2)$ as $t\downarrow0$.
\end{lemma}

\begin{proof}
We multiply \eqref{KP} by a weight that grows exponentially in $r$, which makes its linear part dissipative.

\textit{The weighted equation.} Let $\langle x\rangle=(1+x^2)^{1/2}$. For $a>0$ and $0<\lambda<\delta$, let $W_{a,\lambda}(x,r)=e^{\lambda r-a\langle x\rangle^\kappa}$ and $w=W_{a,\lambda}\phi$. Multiplying \eqref{KP} by $W_{a,\lambda}$ gives
\begin{equation}\label{eq:representation-weighted}
 \begin{aligned}
 w_t+\Lambda w&=g,\qquad\Lambda=\tfrac1{12}Q^3+\tfrac14Q^{-1}\partial_x^2,\\
 g&=\tfrac14Q\bigl((W_{a/2,\lambda/2}\phi)^2\bigr)+\tfrac12Q^{-1}\partial_x(hw)-\tfrac14Q^{-1}\bigl((h^2+h')w\bigr),
 \end{aligned}
\end{equation}
where $Q=\partial_r-\lambda$, $Q^{-1}f(x,r)=-\int_r^\infty e^{\lambda(r-v)}f(x,v)\diff v$ and $h=\partial_x\log W_{a,\lambda}$: conjugation by $W_{a,\lambda}$ turns $\partial_r$, $\partial_r^{-1}$ and $\partial_x$ into $Q$, $Q^{-1}$ and $\partial_x-h$.

We may assume that $B\ge0$. Since $(1+|x|)^\beta\le2\langle x\rangle^\kappa$, the hypothesis gives $|W_{a,\lambda}\phi|\le Ce^{(2B-a)\langle x\rangle^\kappa+\lambda r}(1+r_-)^qe^{-\delta r_+}$ almost everywhere on $J'\times\R^2$. So if $a>4B$, then $W_{a,\lambda}\phi$ and $W_{a/2,\lambda/2}\phi$ are bounded in $L^2(\R^2)\cap L^\infty(\R^2)$, uniformly for $t\in J'$, and the three functions $(W_{a/2,\lambda/2}\phi)^2$, $hw$ and $(h^2+h')w$ in $g$ are bounded in $L^2$, since $h$ and $h'$ grow at most polynomially in $x$.

\textit{The linear part.} In Fourier variables $(k,\xi)$ for $(x,r)$, the operator $\Lambda$ is multiplication by
\[
 \tfrac1{12}(i\xi-\lambda)^3-\frac{k^2}{4(i\xi-\lambda)},\qquad\text{whose real part is}\quad\frac\lambda4\Xi-\frac{\lambda^3}{12},\qquad\Xi=\xi^2+\frac{k^2}{\lambda^2+\xi^2}.
\]
So $e^{-t\Lambda}$ is a strongly continuous semigroup on every Sobolev space $H^\sigma(\R^2)$, and it is elementary to show that, for $0<t\le T$, $\alpha\ge0$ and $\sigma\in\R$,
\begin{equation}\label{eq:representation-semigroup}
 \begin{aligned}
 \|e^{-t\Lambda}f\|_{H^{\sigma+\alpha}}+\|e^{-t\Lambda}Qf\|_{H^{\sigma+\alpha}}+\|e^{-t\Lambda}Q^{-1}\partial_xf\|_{H^{\sigma+\alpha}}&+\|e^{-t\Lambda}Q^{-1}f\|_{H^{\sigma+\alpha}}\\
 &\le C(1+t^{-\alpha-1/2})\|f\|_{H^\sigma},
 \end{aligned}
\end{equation}
since $1+\xi^2+k^2\le C(1+\Xi)^2$ and the symbols of $Q$, $Q^{-1}\partial_x$ and $Q^{-1}$ are at most $C(1+\Xi)^{1/2}$ in modulus. Let $a>4B$. By \eqref{eq:representation-weighted} and the bounds above, $w_t$ is bounded in $H^{-3}$ on $J'$. So, after modification on a null set of times, $w$ is continuous in $H^{-3}$, and for $t_0\le t$ in $J'$,
\begin{equation}\label{eq:representation-mild}
 w(t)=e^{-(t-t_0)\Lambda}w(t_0)+\int_{t_0}^te^{-(t-v)\Lambda}g(v)\diff v.
\end{equation}

\textit{Regularity.} We show by induction on $n\ge0$ that, for every $a>2^{n+2}B$, every $0<\lambda<\delta$ and every compact interval $J\subset\operatorname{int}J'$,
\begin{equation}\label{eq:representation-bootstrap}
 W_{a,\lambda}\phi\in C\bigl(J;H^{n/4}\bigr).
\end{equation}
In \eqref{eq:representation-mild}, take $t_0\in\operatorname{int}J'$ below $J$. By \eqref{eq:representation-semigroup}, the first term is continuous in every $H^\sigma$ for $t\in J$, since $t-t_0$ is bounded below. For $n=0$, \eqref{eq:representation-semigroup} with $\sigma=\alpha=0$ and the bounds of the first step bound the integrand by $C(1+(t-v)^{-1/2})$ in $L^2$. This bound is integrable in $v$, so the integral is continuous in $L^2$. Suppose that \eqref{eq:representation-bootstrap} holds for $n$, and let $a>2^{n+3}B$. Then the three functions in $g$ are bounded in $H^{n/4}$ on compact subintervals of $\operatorname{int}J'$: the square by the inequality $\|f^2\|_{H^\sigma}\le C\|f\|_\infty\|f\|_{H^\sigma}$, $\sigma\ge0$, applied to $f=W_{a/2,\lambda/2}\phi$, and $hw$ and $(h^2+h')w$ because they are $W_{a/2,\lambda}\phi$ multiplied by $he^{-a\langle x\rangle^\kappa/2}$ and $(h^2+h')e^{-a\langle x\rangle^\kappa/2}$, which are smooth with bounded derivatives. So \eqref{eq:representation-semigroup} with $\sigma=n/4$ and $\alpha=\frac14$ bounds the integrand by $C(1+(t-v)^{-3/4})$ in $H^{(n+1)/4}$, and \eqref{eq:representation-bootstrap} follows for $n+1$.

Given $m\ge0$, let $n>4(m+1)$ and $a>2^{n+2}B$. By \eqref{eq:representation-bootstrap} and Sobolev embedding in two dimensions, the derivatives of $w$ in $(x,r)$ of order at most $m$ are bounded and continuous on $J\times\R^2$. Since $\phi=e^{a\langle x\rangle^\kappa-\lambda r}w$, and each derivative of $e^{a\langle x\rangle^\kappa}$ is at most a polynomial in $x$ times $e^{a\langle x\rangle^\kappa}$, the bound of the lemma holds for any $A>a$. Finally, $H^\sigma$ is an algebra for $\sigma>1$, so \eqref{eq:representation-weighted} and \eqref{eq:representation-bootstrap} show that $\partial_tw$ is continuous in every $H^\sigma$ when $a$ is large enough, and so is every $\partial_t^jw$, by differentiating \eqref{eq:representation-weighted} in $t$. Hence $\phi$ is $C^\infty$ on $(\operatorname{int}J')\times\R^2$.

\textit{The limit at $t=0$.} Let $J'=(0,T')$ and $V$ be as in the last assertion, fix $a>4B$ and $0<\lambda<\delta$, and write $W=W_{a,\lambda}$. Then $WV\in L^2$, and the bounds of the first step hold uniformly for $0<t<T'$. For every test function $\eta$, $\int w(t)\eta=\langle\phi(t),W\eta\rangle\to\langle V,W\eta\rangle$ as $t\downarrow0$ outside the null set, and $\|w(t)\|_2$ is bounded, so $w(t)\to WV$ weakly in $L^2$ along these times. Letting $t_0\downarrow0$ along them in \eqref{eq:representation-mild} gives, for $0<t<T'$,
\[
 w(t)=e^{-t\Lambda}WV+\int_0^te^{-(t-v)\Lambda}g(v)\diff v,
\]
since the adjoint semigroup is strongly continuous on $L^2$ too, and the integrand is bounded as in the case $n=0$. The integral is at most $C(\sqrt t+t)$ in $L^2$, and $e^{-t\Lambda}WV\to WV$ in $L^2$ as $t\downarrow0$. So $w(t)\to WV$ in $L^2$, and, since $W$ is bounded below on compact sets, $\phi(t)\to V$ in $L^2_{\rm loc}(\R^2)$.
\end{proof}

\begin{proof}[Proof of the claims in Remark~\ref{rem:bounded-restart}]
Let $V$ be measurable with $|V|\le Me^{-\delta r_+}$, and write $\phi=\phi_V$. Choose $V_n\in C_c^\infty(\R^2;\R)$ converging to $V$ in $L^1_{\rm loc}$ with $|V_n|\le M'e^{-\delta r_+}$ for one $M'\ge M$ (mollified truncations of $V$, for instance), and write $\phi_n=\phi_{V_n}$. We use three facts about these fields.

First, by Proposition~\ref{prop:bounded-stability}(i), $\log Z_{V_n}\to\log Z_V$ in $C^\infty_{\rm loc}((0,\infty)\times\R^2)$, so $\phi_n\to\phi$ locally uniformly at positive times. Proposition~\ref{prop:bounded-stability}(ii), applied to $V_n$ and $V$ with $M'$ in place of $M$ and to the derivative $\partial_r^2$, gives for every $0<d<\delta$, $0<s_0<s_1$ and $h\in\R$ a constant $C$, not depending on $n$, such that
\begin{equation}\label{eq:zero-strip}
 |\phi_n(s,y,u)|+|\phi(s,y,u)|\le Ce^{-du},\qquad s_0\le s\le s_1,\quad y\in\R,\quad u\ge h.
\end{equation}

Second, the same bound for $\log Z_V$ and $\partial_r\log Z_V$ lets us insert $\phi=-2\partial_r^2\log Z_V$ into $\bigl(\tfrac12\partial_r^{-2}(\phi(t))\bigr)$ and integrate by parts twice, which gives
\begin{equation}\label{eq:zero-normalization}
 Z_V(t,x,r)=\exp\{-\tfrac12\partial_r^{-2}\phi(t,x,r)\},\qquad t>0,\quad x,r\in\R.
\end{equation}

Third, by Proposition~\ref{prop:L2-weighted} and the compact-smooth case of Theorem~\ref{thm:L2}, $\phi_n$ agrees with the standard global KP-II solution with initial data $V_n$ and satisfies the hypotheses of Theorem~\ref{thm:representation} for every $T$. Also, \eqref{eq:zero-normalization} holds for $V_n$ by the argument just given. So part~(i) of that theorem gives
\begin{equation}\label{eq:zero-compact-restart}
 Z_{\phi_n(s)}(\tau,x,r)=\exp\{-\tfrac12\partial_r^{-2}\phi_n(s+\tau,x,r)\}=Z_{V_n}(s+\tau,x,r),\qquad s,\tau>0,\quad x,r\in\R.
\end{equation}

\textit{The one-sided restart.} Fix $s,\tau>0$, $x,r\in\R$ and $0<d<\delta$, and use the realization \eqref{eq:occ-scaled}, in which the line $i$ visits the height $u_i(z)=r+\tau^{1/3}(z^2-\bar{\mathcal A}_i(z))$ at the position $x+\tau^{2/3}z$. The lowest visited height is $r-\tau^{1/3}\mathcal M$, with $\mathcal M=\sup_{z\in\R}(\bar{\mathcal A}_1(z)-z^2)$ one random variable for all $\tau$. The top line $\bar{\mathcal A}_1$ is the Airy$_2$ process, and $\mathcal M<\infty$ almost surely by the parabolic-maximum theorem \cite[Theorem~1]{CQR2013}. By Lemma~\ref{lem:occ-moments}(i), the occupation $\tau^{2/3}\sum_i\int_\R e^{-du_i(z)}\diff z$ of the potential $e^{-du}$ is finite almost surely. On the event where it and $\mathcal M$ are finite, choose $h\le r-\tau^{1/3}\mathcal M$; then every visited height is at least $h$, and \eqref{eq:zero-strip} with $s_0<s<s_1$ gives
\[
 |\phi_n(s,x+\tau^{2/3}z,u_i(z))|+|\phi(s,x+\tau^{2/3}z,u_i(z))|\le C_he^{-du_i(z)}
\]
for all $i$, $z$ and $n$. Since $\sum_i\int_\R e^{-du_i(z)}\diff z<\infty$ and $\phi_n(s)\to\phi(s)$ pointwise, dominated convergence gives $I_{\phi_n(s)}^{\tau,x}(r)\to I_{\phi(s)}^{\tau,x}(r)$ almost surely, and the sum and integral defining $I_{\phi(s)}^{\tau,x}(r)$ converge absolutely. By \eqref{eq:zero-compact-restart} and Proposition~\ref{prop:bounded-stability}(i), $\E e^{-I_{\phi_n(s)}^{\tau,x}(r)}=Z_{V_n}(s+\tau,x,r)\to Z_V(s+\tau,x,r)$, so Fatou's lemma, \eqref{eq:zero-normalization} and the finiteness of $I_{\phi(s)}^{\tau,x}(r)$ give
\begin{equation}\label{eq:zero-one-sided}
 0<Z_{\phi(s)}(\tau,x,r)\le Z_V(s+\tau,x,r)=\exp\{-\tfrac12\partial_r^{-2}\phi(s+\tau,x,r)\}.
\end{equation}

\textit{The terminal limit.} Fix $t>0$ and $\tau_k\downarrow0$ with $\tau_k\le t/2$. Take a smooth $\chi$ with $\chi(u)=0$ for $u\le r-2$ and $\chi(u)=1$ for $u\ge r-1$, and put $F_k(y,u)=\chi(u)\phi(t-\tau_k,y,u)$. By \eqref{eq:zero-strip} with $[s_0,s_1]=[t/2,t]$ and $h=r-2$, the functions $F_k$ and $\chi\phi(t)$ are bounded by $Ce^{-du}$; they are continuous, and $F_k\to\chi\phi(t)$ locally uniformly. Lemma~\ref{lem:occ-trace}(ii) gives $I_{F_k}^{\tau_k,x}(r)\to \bigl(\tfrac12\partial_r^{-2}(\chi\phi(t))\bigr)(x,r)=\bigl(\tfrac12\partial_r^{-2}(\phi(t))\bigr)(x,r)$ in $L^2$, the equality because $\chi=1$ on $[r,\infty)$. On the event $\{\tau_k^{1/3}\mathcal M\le1\}$, whose probability tends to one, every visited height is at least $r-1$, so $I_{\phi(t-\tau_k)}^{\tau_k,x}(r)=I_{F_k}^{\tau_k,x}(r)$. Hence $I_{\phi(t-\tau_k)}^{\tau_k,x}(r)\to \bigl(\tfrac12\partial_r^{-2}(\phi(t))\bigr)(x,r)$ in probability, and almost surely along a subsequence, along which Fatou's lemma and \eqref{eq:zero-normalization} give $\liminf Z_{\phi(t-\tau_k)}(\tau_k,x,r)\ge Z_V(t,x,r)$. Since $(\tau_k)$ is arbitrary, \eqref{eq:zero-one-sided} gives $Z_{\phi(s)}(t-s,x,r)\to Z_V(t,x,r)$ as $s\uparrow t$.
\end{proof}

\section{Centered Airy occupation and \texorpdfstring{$L^2$}{L2} renormalization}
\label{app:centered-occupation-methods}

Centering removes the divergent first cumulant of the Airy occupation for $L^2$ data. We construct the centered statistic, identify its
regularized Fredholm determinant, and match the resulting formula
with the classical KP-II flow. The final subsection develops
a separate small-time focusing example.

\subsection{Centered \texorpdfstring{$L^2$}{L2} occupation and regularized determinants}\label{app:L2}

Theorem~\ref{thm:L2} has two parts: for every real $F\in L^2(\R^2)$ the centered occupation of the stationary lines is defined and sub-Gaussian (Proposition~\ref{prop:L2-centered}), so that $\tau_V^{\rm ren}$ is defined for real $V\in L^2(\R^2)$ (Lemma~\ref{lem:L2-tau}); and the formula gives the standard $L^2$ flow of KP-II. We prove the second for $V\in C_c^\infty$ by Theorem~\ref{thm:representation}(ii), using that $2\partial_r^2\E I_V$ is the linear evolution of $V$ (Lemma~\ref{lem:L2-mean}) and that the standard solution decays as $r\to+\infty$ (Proposition~\ref{prop:L2-weighted}), and extend it to $L^2$ by approximation.

The Hilbert--Schmidt bounds of Proposition~\ref{prop:L2-centered} rest on the following estimate for the square of the extended kernel.

\begin{lemma}\label{lem:L2-square}
For $\sigma\in\{-1,1\}$ and a bounded function $f$ on $\R$ with compact support, let
\[
 M_f^\sigma(a,b)=\int_\R f(z)\Ai(z+\sigma a)\Ai(z+\sigma b)\diff z.
\]
Then
\begin{equation}\label{eq:L2-spectral}
 \int_0^\infty\int_0^\infty\frac{|M_f^\sigma(a,b)|^2}{a+b}\diff a\diff b\le C\|f\|_2^2.
\end{equation}
Consequently the operator $\mathcal Q$ on $L^2(\R^2)$ with kernel $|K_{\rm ext}(Y,Y')|^2$ is bounded.
\end{lemma}

We use the Airy bounds \cite[\S9.7]{DLMF}
\begin{equation}\label{eq:L2-airy-bounds}
 |\Ai(z)|\le C(1+|z|)^{-1/4},\qquad |\Ai(z)|\le Ce^{-\frac23z^{3/2}}\quad(z\ge0),
\end{equation}
and the fact that, for $\sigma=\pm1$, the map $f\mapsto\int\Ai(z+\sigma a)f(z)\diff z$ is unitary from $L^2(\diff z)$ to $L^2(\diff a)$: it is a composition of reflections with convolution by $\Ai$, whose Fourier transform $e^{ik^3/3}$ has modulus one.

\begin{proof}[Proof of Lemma~\ref{lem:L2-square}]
\textit{The region $a+b\ge1$.} There $(a+b)^{-1}\le2(1+a)^{-1}$. Extending the $b$-integral to $\R$ and applying the unitarity of the Airy transform to $z\mapsto f(z)\Ai(z+\sigma a)$ gives $\int_0^\infty|M_f^\sigma(a,b)|^2\diff b\le\int|f(z)|^2\Ai(z+\sigma a)^2\diff z$. So this region contributes at most
\[
 2\int_\R|f(z)|^2\int_0^\infty\frac{\Ai(z+\sigma a)^2}{1+a}\diff a\diff z.
\]
The inner integral is bounded uniformly in $z$ by H\"older's inequality, since $\Ai^2\in L^3(\R)$ by \eqref{eq:L2-airy-bounds} and $(1+a)^{-1}\in L^{3/2}(0,\infty)$.

\textit{The region $a+b<1$.} Cover it by the shells $\varepsilon/2<a+b\le\varepsilon$, $\varepsilon=2^{-j}$, $j\ge0$, and split $f=f_0+f_1$ with $f_0=f\mathbf 1_{\{z\ge-4\varepsilon^{-2}\}}$. For $0\le a,b\le\varepsilon$, \eqref{eq:L2-airy-bounds} gives $\Ai(z+\sigma a)^2\Ai(z+\sigma b)^2\le C(1+|z|)^{-1}$ for $z\le0$ and $\le Ce^{-z}$ for $z\ge0$, hence
\[
 \int_{-4\varepsilon^{-2}}^\infty\Ai(z+\sigma a)^2\Ai(z+\sigma b)^2\diff z\le C(1+|\log\varepsilon|).
\]
By the Cauchy--Schwarz inequality, $|M_{f_0}^\sigma(a,b)|^2\le C(1+|\log\varepsilon|)\|f\|_2^2$ on the shell. The shell has area at most $\varepsilon^2$ and $(a+b)^{-1}\le2/\varepsilon$ on it, so $f_0$ contributes at most $C\varepsilon(1+|\log\varepsilon|)\|f\|_2^2$. On the support of $f_1$ and for $0\le a\le\varepsilon$, $|z+\sigma a|\ge3\varepsilon^{-2}$, so $\Ai(z+\sigma a)^2\le C\varepsilon$ by \eqref{eq:L2-airy-bounds}. The shell lies in $\{0\le a\le\varepsilon,\ b\ge0\}$, and extending the $b$-integral to $\R$ and using unitarity as before bounds the contribution of $f_1$ by
\[
 \frac2\varepsilon\int_0^\varepsilon\int_{z<-4\varepsilon^{-2}}|f(z)|^2\Ai(z+\sigma a)^2\diff z\diff a\le C\varepsilon\|f\|_2^2.
\]
Since $|M_f^\sigma|^2\le2|M_{f_0}^\sigma|^2+2|M_{f_1}^\sigma|^2$, the shell contributes at most $C\varepsilon(1+|\log\varepsilon|)\|f\|_2^2$, and summing over $\varepsilon=2^{-j}$ proves \eqref{eq:L2-spectral}.

\textit{The operator $\mathcal Q$.} Let $f,g\ge0$ be bounded with compact support in $\R^2$, and write $f(s)=f(s,\cdot)$. Let $K_A$ be the stationary extended Airy kernel defined in Section~\ref{sec:renormalized-solutions}, with its spectral integrals restricted to $[0,A]$. For $s>s'$, squaring the first branch gives
\[
 \begin{aligned}
 &\iint_{s>s'}f(s,z)|K_A(s,z;s',z')|^2g(s',z')\diff(s,z)\diff(s',z')\\
 &\qquad=\int_0^A\int_0^A\iint_{s>s'}e^{-(s-s')(a+b)}M^{+1}_{f(s)}(a,b)M^{+1}_{g(s')}(a,b)\diff s\diff s'\diff a\diff b,
 \end{aligned}
\]
where all integrals converge absolutely because $\Ai$ is bounded and $A$ is finite. Convolution in $s$ with $\mathbf 1_{h>0}e^{-(a+b)h}$ has norm at most $(a+b)^{-1}$ on $L^2(\R)$, so the right side is at most
\[
 \int_0^\infty\int_0^\infty\frac{\|M^{+1}_{f(\cdot)}(a,b)\|_{L^2_s}\|M^{+1}_{g(\cdot)}(a,b)\|_{L^2_s}}{a+b}\diff a\diff b\le C\|f\|_2\|g\|_2,
\]
by the Cauchy--Schwarz inequality in $(a,b)$ and \eqref{eq:L2-spectral} applied at each $s$. The branch $s<s'$ gives the same bound with $\sigma=-1$. As $A\to\infty$, $K_A\to K_{\rm ext}$ off the null set $\{s=s'\}$, and Fatou's lemma gives $\langle f,\mathcal Qg\rangle\le C\|f\|_2\|g\|_2$. Monotone approximation extends this to nonnegative $f,g\in L^2(\R^2)$, and $|\langle f,\mathcal Qg\rangle|\le\langle|f|,\mathcal Q|g|\rangle$ to all of $L^2(\R^2)$.
\end{proof}

For real bounded compactly supported $F$ on $\R^2$, let $I_F=\sum_{i\ge1}\int_\R F\bigl(s,\bar{\mathcal A}_i(s)\bigr)\diff s$ and $X_F=I_F-\E I_F$. With $V(y,u)=F(y,y^2-u)$, \eqref{eq:occ-scaled} at $(t,x,r)=(1,0,0)$ gives $I_F=I_V^{1,0}(0)$, which converges absolutely almost surely by Lemma~\ref{lem:occ-moments}(i) applied to $|V|$. For real $F\in L^2(\R^2)$ let $\widetilde{\mathcal K}_F$ be the operator with kernel
\begin{equation}\label{eq:L2-BF}
 \widetilde{\mathcal K}_F(Y,Y')=|F(Y)|^{1/2}K_{\rm ext}(Y,Y')\operatorname{sgn}(F(Y'))|F(Y')|^{1/2},
\end{equation}
and write $\|\cdot\|_{\rm HS}$ for the Hilbert--Schmidt norm.

\begin{proof}[Proof of Proposition~\ref{prop:L2-centered}]
\textit{Hilbert--Schmidt bounds.} By Lemma~\ref{lem:L2-square},
\[
 \|\widetilde{\mathcal K}_F\|_{\rm HS}^2=\iint|F(Y)||K_{\rm ext}(Y,Y')|^2|F(Y')|\diff Y\diff Y'=\langle|F|,\mathcal Q|F|\rangle\le C\|F\|_2^2.
\]
For the difference, $\widetilde{\mathcal K}_F-\widetilde{\mathcal K}_G$ is the sum of the operators with kernels
\[
 \begin{gathered}
 \bigl(|F|^{1/2}-|G|^{1/2}\bigr)(Y)K_{\rm ext}(Y,Y')\bigl(\operatorname{sgn}(F)|F|^{1/2}\bigr)(Y'),\\
 |G(Y)|^{1/2}K_{\rm ext}(Y,Y')\bigl(\operatorname{sgn}(F)|F|^{1/2}-\operatorname{sgn}(G)|G|^{1/2}\bigr)(Y').
 \end{gathered}
\]
The pointwise inequalities $\bigl||F|^{1/2}-|G|^{1/2}\bigr|^2\le|F-G|$ and $\bigl|\operatorname{sgn}(F)|F|^{1/2}-\operatorname{sgn}(G)|G|^{1/2}\bigr|^2\le2|F-G|$ and the same computation bound their squared Hilbert--Schmidt norms by $C\|F-G\|_2\|F\|_2$ and $C\|G\|_2\|F-G\|_2$, which proves the corresponding bound for stationary test functions.

\textit{The determinant.} Let $F$ be real, bounded and compactly supported, and $V(y,u)=F(y,y^2-u)$. The map $\Psi(s,z)=(s,s^2-z)$ is a measure-preserving bijection of $\R^2$ with $V\circ\Psi=F$. At $(t,x,r)=(1,0,0)$ the function $\xi_y(w)$ of Proposition~\ref{lem:bounded-det} is $y^2-w$, which equals $z$ at $(y,w)=\Psi(s,z)$, so \eqref{eq:bounded-kernel} and \eqref{eq:bounded-kernel-minus} give $\widehat K^\gamma(\Psi(Y),\Psi(Y'))=K_{\rm ext}(Y,Y')$. Hence the unitary map $f\mapsto f\circ\Psi$ of $L^2(\R^2)$ conjugates the operator $\mathsf B_V$ of Proposition~\ref{lem:bounded-det} at $(1,0,0)$ to $\widetilde{\mathcal K}_F$, which leaves $\det_2$ unchanged. Since $\mathsf B_{-pV}=-p\mathsf B_V$ for real $p$ and $I_F=I_V^{1,0}(0)$, that proposition applied to $-pV$ gives
\[
 \E e^{pI_F}=Z_{-pV}(1,0,0)=e^{p\E I_F}\det{}_2(I+p\widetilde{\mathcal K}_F),
\]
the determinant identity for stationary test functions. The bound $|\det_2(I+A)|\le e^{\|A\|_{\rm HS}^2/2}$ for Hilbert--Schmidt $A$ \cite[Chapter~9]{Simon} and the Hilbert--Schmidt estimate above give the exponential bound. In particular $\E e^{p|X_F|}<\infty$ for every $p$, and since $\E X_F=0$, letting $p\to0$ in $p^{-2}(\E e^{pX_F}-1)\le p^{-2}(e^{Cp^2\|F\|_2^2}-1)$ gives $\E X_F^2\le2C\|F\|_2^2$.

\textit{Extension to $L^2$.} By linearity and the variance bound, $\E|X_F-X_G|^2\le C\|F-G\|_2^2$ for bounded compactly supported $F,G$. These functions are dense in $L^2(\R^2)$, so $F\mapsto X_F$ extends uniquely to a continuous linear map into $L^2(\mathbb P)$. Let $F_n\to F$ in $L^2$, with $F_n$ bounded and compactly supported. Then $X_{F_n}\to X_F$ in probability, and $\E e^{2pX_{F_n}}\le e^{4Cp^2\sup_n\|F_n\|_2^2}$, so the variables $e^{pX_{F_n}}$ are uniformly integrable and $\E e^{pX_{F_n}}\to\E e^{pX_F}$. By the Hilbert--Schmidt estimate above and the continuity of $\det_2$ in the Hilbert--Schmidt norm, $\det_2(I+p\widetilde{\mathcal K}_{F_n})\to\det_2(I+p\widetilde{\mathcal K}_F)$. Thus the stationary determinant identity passes to the limit. Taking $F=\widetilde V_{t,x,r}$, for which $\widetilde{\mathcal K}_F$ is the operator $\widetilde{\mathcal K}_V$ of Section~\ref{sec:L2}, and using $\|\widetilde V_{t,x,r}\|_2^2=t^{1/3}\|V\|_2^2$ gives \eqref{eq:L2-HS} and \eqref{eq:L2-det2} as stated in Proposition~\ref{prop:L2-centered}.
\end{proof}

For real $V\in L^2(\R^2)$, $t>0$ and $x,r\in\R$, let
\begin{equation}\label{eq:L2-rescaled}
 \widetilde V_{t,x,r}(s,z)=t^{2/3}V\bigl(x+t^{2/3}s,r+t^{1/3}(s^2-z)\bigr),\qquad\tau_V^{\rm ren}(t,x,r)=\E e^{-X_{\widetilde V_{t,x,r}}}.
\end{equation}
For $V\in C_c^\infty$, \eqref{eq:occ-scaled} gives $I_V^{t,x}(r)=I_{\widetilde V_{t,x,r}}$, so $X_{\widetilde V_{t,x,r}}=I_V^{t,x}(r)-\E I_V^{t,x}(r)$.

\begin{proof}[Proof of Lemma~\ref{lem:L2-tau}]
(i) The map $(s,z)\mapsto\bigl(x+t^{2/3}s,r+t^{1/3}(s^2-z)\bigr)$ is a bijection of $\R^2$ whose Jacobian determinant has absolute value $t$, so $\|\widetilde V_{t,x,r}\|_2^2=t^{4/3}t^{-1}\|V\|_2^2$. Proposition~\ref{prop:L2-centered} with $p=-1$ gives the determinant formula and the upper bound in \eqref{eq:L2-tau-bounds}. Since the centered occupation has mean zero, Jensen's inequality gives $\tau_V^{\rm ren}\ge1$. The map $(t,x,r)\mapsto \widetilde V_{t,x,r}$ is continuous from $(0,\infty)\times\R^2$ into $L^2(\R^2)$: for continuous compactly supported $V$ by dominated convergence, and for general $V$ by approximation, since $\|\widetilde V_{t,x,r}-\widetilde W_{t,x,r}\|_2=t^{1/6}\|V-W\|_2$. By \eqref{eq:L2-HS} and the continuity of $\det_2$ in the Hilbert--Schmidt norm, $\tau_V^{\rm ren}$ is continuous.

(ii) Since $\widetilde{V_n}_{t,x,r}-\widetilde V_{t,x,r}=\widetilde{(V_n-V)}_{t,x,r}$, \eqref{eq:L2-HS} and \eqref{eq:L2-tau-bounds} give
\[
 \sup_{0<t\le T,\ x,r\in\R}\|\widetilde{\mathcal K}_{\widetilde{V_n}_{t,x,r}}-\widetilde{\mathcal K}_{\widetilde V_{t,x,r}}\|_{\rm HS}^2\le CT^{1/3}(\|V_n\|_2+\|V\|_2)\|V_n-V\|_2\longrightarrow0.
\]
On a bounded set of Hilbert--Schmidt operators $\det_2$ is Lipschitz, $|\det_2(I+A)-\det_2(I+A')|\le\|A-A'\|_{\rm HS}\exp\{\frac12(\|A\|_{\rm HS}+\|A'\|_{\rm HS}+1)^2\}$ \cite[Theorem~9.2]{Simon}, so $\tau_{V_n}^{\rm ren}\to\tau_V^{\rm ren}$ uniformly on $(0,T]\times\R^2$. Since both are at least $1$, $|\log\tau_{V_n}^{\rm ren}-\log\tau_V^{\rm ren}|\le|\tau_{V_n}^{\rm ren}-\tau_V^{\rm ren}|$.

(iii) Take expectations in $e^{-I_V^{t,x}(r)}=e^{-\E I_V^{t,x}(r)}e^{-X_{\widetilde V_{t,x,r}}}$.
\end{proof}

Let $U_t$ be the linear KP-II group, the solution operator of \eqref{KP} without its nonlinear term. With $\widehat f(k,\xi)=\int e^{-i(kx+\xi r)}f(x,r)\diff x\diff r$, it is the unitary group on $L^2(\R^2)$ given by
\begin{equation}\label{eq:L2-group}
 \widehat{U_tV}(k,\xi)=\exp\Bigl\{it\Bigl(\frac{\xi^3}{12}-\frac{k^2}{4\xi}\Bigr)\Bigr\}\widehat V(k,\xi),\qquad\xi\ne0.
\end{equation}

\begin{proof}[Proof of Lemma~\ref{lem:L2-mean}]
\textit{The mean as a convolution.} By \eqref{eq:occ-density} with $u=r+a$,
\[
 \E I_V^{t,x}(r)=\iint_{\R^2}\Theta_t(x-y,a)V(y,r+a)\diff y\diff a,\qquad\Theta_t(h,a)=t^{-1/3}\rho_{\rm Ai}\bigl(t^{-1/3}(h^2/t-a)\bigr).
\]

\textit{The transform of $\Theta_t$.} Let $q>0$ and $k\in\R$. After the substitution $a=h^2/t-t^{1/3}w$, the integral in $w$ is $\mu(qt^{1/3})$ of \eqref{app:mu}, and the Gaussian integral in $h$ gives
\begin{equation}\label{eq:L2-laplace}
 \iint_{\R^2}e^{-ikh-qa}\Theta_t(h,a)\diff h\diff a=\frac1{2q^2}\exp\Bigl\{t\Bigl(\frac{q^3}{12}-\frac{k^2}{4q}\Bigr)\Bigr\}.
\end{equation}
The left side converges absolutely for $\operatorname{Re}q>0$, and both sides are analytic there, so \eqref{eq:L2-laplace} holds for $\operatorname{Re}q>0$.

\textit{A weighted Fourier transform.} Fix $0<\lambda\le1$ and put $q=\lambda-i\xi$. For a bounded compactly supported $W$, the function $e^{\lambda r}\E I_W^{t,x}(r)=\iint \Theta_t(h,a)e^{-\lambda a}(e^{\lambda r}W)(x-h,r+a)\diff h\diff a$ is the convolution of $e^{\lambda r}W\in L^1$ with a finite measure of mass $e^{t\lambda^3/12}/(2\lambda^2)$, by \eqref{eq:L2-laplace} at $k=0$. So it is integrable, and by \eqref{eq:L2-laplace} its Fourier transform is $\frac1{2q^2}\exp\{t(q^3/12-k^2/4q)\}\widehat{e^{\lambda r}W}$. Differentiating under the integral, which is legitimate since $V$ is smooth with compact support and $\Theta_t$ is integrable on $\R\times[-A,A]$ for every $A$, gives $\partial_r^2\E I_V=\E I_{\partial_r^2V}$. Apply the transform with $W=\partial_r^2V$, and use $e^{\lambda r}\partial_r^2V=(\partial_r-\lambda)^2(e^{\lambda r}V)$, whose transform is $q^2\widehat{e^{\lambda r}V}$. The factors $2q^2$ cancel:
\[
 \widehat{e^{\lambda r}2\partial_r^2\E I_V}(k,\xi)=\exp\Bigl\{t\Bigl(\frac{q^3}{12}-\frac{k^2}{4q}\Bigr)\Bigr\}\widehat{e^{\lambda r}V}(k,\xi).
\]

\textit{The limit $\lambda\downarrow0$.} The real part of the exponent is $t\lambda^3/12-t\lambda\xi^2/4-t\lambda k^2/(4(\lambda^2+\xi^2))\le t/12$. Since $V$ is smooth with compact support, $|\widehat{e^{\lambda r}V}(k,\xi)|\le C_N(1+|k|+|\xi|)^{-N}$ for every $N$, uniformly in $0<\lambda\le1$. So the right side is integrable, Fourier inversion applies to the continuous integrable function $e^{\lambda r}2\partial_r^2\E I_V$, and dominated convergence lets $\lambda\downarrow0$. For $\xi\ne0$ the exponent tends to $it(\xi^3/12-k^2/(4\xi))$, so the inversion integrals tend to that of \eqref{eq:L2-group}, which is $U_tV(x,r)$, while the left sides tend to $2\partial_r^2\E I_V^{t,x}(r)$.
\end{proof}

\textit{The standard $L^2$ flow.} We use one input from the $L^2$ theory of KP-II, \cite[Theorems~1.1 and~1.6]{Hadac2008} with dispersion exponent $2$; fixed rescalings of the variables and of the field convert Hadac's normalization into \eqref{KP}. For real $V\in L^2(\R^2)$ there is a solution $\Phi_t(V)$ of \eqref{KP} on a time interval whose length is bounded below in terms of $\|V\|_2$. It is unique in the resolution spaces of that theorem, depends real-analytically on $V$ as an element of $C_tL^2$, and satisfies the Duhamel formula
\[
 \Phi_t(V)=U_tV+\frac14\int_0^tU_{t-t'}\partial_r\bigl(\Phi_{t'}(V)^2\bigr)\diff t'.
\]
The $L^2$ norm of $\Phi_t(V)$ is conserved, so the solution extends to all times. We call the global solution $\Phi(V)\in C([0,\infty);L^2(\R^2))$ the standard solution. Below, $H^m$ and $W^{m,\infty}$ are the Sobolev spaces on $\R^2$, and $D^\alpha=\partial_x^{\alpha_1}\partial_r^{\alpha_2}$.

\begin{lemma}\label{lem:L2-flow}
Let $T>0$.
\begin{enumerate}
\item[\textup{(i)}] The map $V\mapsto\Phi(V)$ from real $L^2(\R^2)$ to $C([0,T];L^2(\R^2))$ is real analytic near every point.
\item[\textup{(ii)}] If $V_n\to V$ in $H^m$ for every $m$, then $\Phi(V_n)\to\Phi(V)$ in $C([0,T];H^m)$ for every $m$.
\end{enumerate}
\end{lemma}

\begin{proof}
(i) By conservation of the $L^2$ norm, every initial datum in a small $L^2$ ball around $V$ keeps norm at most $\|V\|_2+1$ along its trajectory, so the local theory applies from every point of the trajectory with one existence time $\theta>0$. On $[0,T]$, $\Phi$ is therefore a composition of finitely many local analytic flow maps, and is analytic.

(ii) Let $T_af=f(\cdot-a)$ for $a\in\R^2$. Uniqueness in the resolution spaces gives $\Phi(T_aV)=T_a\Phi(V)$. If $V\in H^m$, then $a\mapsto T_aV$ is $m$ times continuously differentiable into $L^2$, so by (i) $a\mapsto T_a\Phi(V)$ is $m$ times continuously differentiable into $C_tL^2$. Its derivatives at $a=0$ are, up to sign, the distributional derivatives $D^\alpha\Phi(V)$, $|\alpha|\le m$, which therefore lie in $C_tL^2$; by the chain rule each is a finite sum of derivatives of the analytic map $\Phi$ at $V$, applied to derivatives of $V$ of order at most $m$. The same formula for $V_n$, and the continuity of each derivative of $\Phi$, give $D^\alpha\Phi(V_n)\to D^\alpha\Phi(V)$ in $C_tL^2$.
\end{proof}

\begin{proposition}\label{prop:L2-weighted}
Let $V\in C_c^\infty(\R^2;\R)$ and $u(t)=\Phi_t(V)$. For every $T>0$, $\lambda>0$ and multi-index $\alpha$,
\[
 \sup_{0\le t\le T,\ x,r\in\R}e^{\lambda r_+}|D^\alpha u(t,x,r)|<\infty.
\]
Moreover $u$ is a classical solution of \eqref{KP} on $[0,\infty)\times\R^2$, and $D^\alpha u(t)\to D^\alpha V$ uniformly as $t\downarrow0$ for every $\alpha$.
\end{proposition}

\begin{proof}
\textit{Zero-mass approximants.} The energy estimate below needs a primitive of the solution in $L^2$, so we first cancel the mass of $V$ by a correction far below its support. Let $M(x)=\int V(x,r)\diff r$, choose $\psi\in C_c^\infty((-2,-1))$ with $\int\psi=1$, and for $R\ge1$ let $\psi_R(r)=R^{-1}\psi(r/R)$, $V_R=V-M\psi_R$ and $u_R=\Phi(V_R)$. Then $\int V_R(x,r)\diff r=0$, and $\|\partial_x^j\partial_r^k(M\psi_R)\|_2=O(R^{-k-1/2})$, so $V_R\to V$ in every $H^m$. By Lemma~\ref{lem:L2-flow}(ii) and Sobolev embedding, for every $m$ there is $R_m\ge1$ such that, as $R\to\infty$,
\begin{equation}\label{eq:L2-uniform-smooth}
 u_R\to u\text{ in }C([0,T];H^m),\qquad\sup_{R\ge R_m}\sup_{0\le t\le T}\|u_R(t)\|_{W^{m,\infty}}<\infty.
\end{equation}
For fixed $R$, $q_R(0)=\partial_r^{-1}V_R$ is smooth with compact support, so it lies in every $H^m$, although its norms grow with $R$. Let
\[
 q_R(t)=U_tq_R(0)+\frac14\int_0^tU_{t-t'}\bigl(u_R(t')^2\bigr)\diff t'.
\]
Then $q_R\in C_tH^m$ for every $m$, $\partial_rq_R=u_R$ by the Duhamel formula for $u_R$, and
\begin{equation}\label{eq:L2-KP-q}
 \partial_tu_R=\tfrac12u_R\partial_ru_R-\tfrac1{12}\partial_r^3u_R-\tfrac14\partial_x^2q_R
\end{equation}
in every $H^m$; in particular $u_R\in C^1_tH^m$ for every $m$.

\textit{The weighted energy identity.} Let $w=w(r)$ be smooth and bounded with bounded derivatives, and let $v\in C^1_tH^m$ and $p\in C_tH^m$ for all $m$, with $\partial_rp=v$ and $\partial_tv=g-\frac1{12}\partial_r^3v-\frac14\partial_x^2p$. Integration by parts, with all integrals over $\R^2$, gives
\begin{equation}\label{eq:L2-energy}
 \frac{d}{dt}\int wv^2=-\frac14\int w'(\partial_rv)^2+\frac1{12}\int w'''v^2-\frac14\int w'(\partial_xp)^2+2\int wvg.
\end{equation}
Here the dispersive term is $-\frac16\int wv\partial_r^3v=-\frac14\int w'(\partial_rv)^2+\frac1{12}\int w'''v^2$, and the nonlocal term is $-\frac12\int wv\partial_x^2p=\frac12\int w\partial_x\partial_rp\partial_xp=-\frac14\int w'(\partial_xp)^2$.

\textit{Weights.} Fix $\lambda>0$. Choose $\rho\in C^\infty(\R)$ with $\rho=0$ on $(-\infty,0]$, $0\le\rho'\le1$ and $\rho'=1$ on $[1,\infty)$, and $\chi\in C^\infty(\R)$ nonincreasing with $\chi=1$ on $(-\infty,1]$ and $\chi=0$ on $[2,\infty)$. For $N\ge1$ let
\[
 \rho_N(r)=\int_0^r\rho'(a)\chi(a/N)\diff a,\qquad w_N=e^{2\lambda\rho_N}.
\]
Each $w_N$ is bounded with bounded derivatives, and $w_N$ increases with $N$ to $w=e^{2\lambda\rho}$. Moreover $w_N'\ge0$ and $|w_N^{(j)}|\le C_\lambda w_N$ for $j=1,2,3$, uniformly in $N$.

\textit{The energy estimate.} For $|\alpha|\le m$, apply \eqref{eq:L2-energy} with $w=w_N$, $v=D^\alpha u_R$, $p=D^\alpha q_R$ and $g=\frac12D^\alpha(u_R\partial_ru_R)$, which is legitimate by \eqref{eq:L2-KP-q}. The first and third terms of \eqref{eq:L2-energy} are nonpositive, since $w_N'\ge0$, and we drop them; the second is at most $C_\lambda\int w_Nv^2$. By the Leibniz rule,
\[
 D^\alpha(u_R\partial_ru_R)=u_R\partial_rv+\sum_{0<\beta\le\alpha}\binom\alpha\beta(D^\beta u_R)(D^{\alpha-\beta}\partial_ru_R).
\]
The first term contributes $\int w_Nu_Rv\partial_rv=-\frac12\int(w_N'u_R+w_N\partial_ru_R)v^2$ to $2\int w_Nvg$. Each remaining term is a factor $D^\beta u_R$, bounded in $L^\infty$, times a derivative $D^{\alpha-\beta}\partial_ru_R$ of order at most $|\alpha|$, and the Cauchy--Schwarz inequality with weight $w_N$ bounds its contribution. Summing over $|\alpha|\le m$, $E_{m,N,R}(t)=\sum_{|\alpha|\le m}\int w_N|D^\alpha u_R(t)|^2$ satisfies
\[
 E_{m,N,R}'(t)\le C_{m,\lambda}\bigl(1+\|u_R(t)\|_{W^{m+1,\infty}}\bigr)E_{m,N,R}(t).
\]
The correction $M\psi_R$ is supported in $\{r<0\}$, where $w_N=1$, so
\[
 E_{m,N,R}(0)\le2\sum_{|\alpha|\le m}\Bigl(\int w|D^\alpha V|^2+\|D^\alpha(M\psi_R)\|_2^2\Bigr)
\]
is bounded uniformly in $N$ and $R$. By Gronwall's inequality and \eqref{eq:L2-uniform-smooth}, $E_{m,N,R}(t)$ is bounded uniformly in $N$, $R\ge R_{m+1}$ and $t\in[0,T]$, and by monotone convergence in $N$ the same holds with $w$ in place of $w_N$.

\textit{Pointwise decay.} The derivatives of $e^{\lambda\rho}$ are bounded by $Ce^{\lambda\rho}$, so the bound just proved, with $m+2$ in place of $m$, controls $e^{\lambda\rho}D^\alpha u_R$ in $H^2$ for $|\alpha|\le m$, uniformly in $R$ and $t$, and $H^2(\R^2)\subset L^\infty$ gives $|D^\alpha u_R|\le Ce^{-\lambda\rho(r)}\le Ce^{-\lambda r_+}$, since $|\rho(r)-r_+|\le1$. Letting $R\to\infty$ in \eqref{eq:L2-uniform-smooth} proves the bound of the proposition.

\textit{The equation.} Since $q_R(t)\in H^2(\R^2)$, it tends to $0$ as $r\to+\infty$, and $\partial_rq_R=u_R$ with the decay just proved gives $q_R=-\int_r^\infty u_R(t,x,a)\diff a$; the same holds for $\partial_x^2q_R$ and $\partial_x^2u_R$. The decay of $\partial_x^2u_R$, uniform in $R$, and the local uniform convergence $\partial_x^2u_R\to\partial_x^2u$ give $\partial_x^2q_R\to-\int_r^\infty\partial_x^2u(t,x,a)\diff a$ locally uniformly. Passing to the limit in \eqref{eq:L2-KP-q} in distributions, $u$ satisfies \eqref{KP} with the upper primitive, and its right side is continuous on $[0,\infty)\times\R^2$, so $u$ is continuously differentiable in $t$ and the equation holds classically. The convergence of $D^\alpha u(t)$ at $t=0$ follows from $u\in C([0,T];H^m)$ for every $m$.
\end{proof}

\begin{proof}[Proof of Theorem~\ref{thm:L2}]
\textit{Compact smooth data.} Let $V\in C_c^\infty(\R^2;\R)$, $u(t)=\Phi_t(V)$ and $T>0$. By Proposition~\ref{prop:L2-weighted} with some $\lambda>0$, $u$ is a classical solution of \eqref{KP}, its derivatives through order three satisfy $|D^\alpha u|\le Ce^{-\lambda r_+}\le Ce^{-\lambda r}$ on $[0,T]\times\R^2$, its negative part is bounded, and $u(t)\to V$ uniformly. So $u$ satisfies the hypotheses of Theorem~\ref{thm:representation}(ii) with $B=0$ and $\nu=0$, and that theorem gives $u(t)=-2\partial_r^2\log Z_V(t)$ for $0<t<T$. By \eqref{eq:L2-factor} and Lemma~\ref{lem:L2-mean},
\[
 \Phi_t(V)=-2\partial_r^2\log Z_V=2\partial_r^2\E I_V^{t,x}(r)-2\partial_r^2\log\tau_V^{\rm ren}=U_tV-2\partial_r^2\log\tau_V^{\rm ren}
\]
for $0<t<T$, and so for every $t>0$.

\textit{General data.} Let $V\in L^2(\R^2)$ be real, and choose real $V_n\in C_c^\infty$ with $V_n\to V$ in $L^2$. Fix $t>0$. By Proposition~\ref{prop:L2-centered} and \eqref{eq:L2-det2} applied to $V_n-V$, the centered occupations $X_{V_n}^{t,x}(r)=I_{V_n}^{t,x}(r)-\E I_{V_n}^{t,x}(r)$ converge in $L^2(\mathbb P)$ to $X_V^{t,x}(r)$; this limit is the extension of the centered occupation in Theorem~\ref{thm:L2}. As $n\to\infty$, $\Phi_t(V_n)\to\Phi_t(V)$ in $L^2$ by Lemma~\ref{lem:L2-flow}(i), $U_tV_n\to U_tV$ in $L^2$ since $U_t$ is unitary, and $\log\tau_{V_n}^{\rm ren}(t)\to\log\tau_V^{\rm ren}(t)$ uniformly on $\R^2$ by Lemma~\ref{lem:L2-tau}(ii), so that $\partial_r^2\log\tau_{V_n}^{\rm ren}(t)\to\partial_r^2\log\tau_V^{\rm ren}(t)$ in $\mathcal D'(\R^2)$. Passing to the limit in the identity for $V_n$ gives
\[
 \Phi_t(V)=U_tV-2\partial_r^2\log\tau_V^{\rm ren}(t,\cdot,\cdot)\quad\text{in }\mathcal D'(\R^2),\qquad t>0.
\]
All three terms are locally integrable on $(0,\infty)\times\R^2$, so the identity also holds in $\mathcal D'((0,\infty)\times\R^2)$. Its left side is the standard solution $\Phi(V)\in C([0,\infty);L^2(\R^2))$, which tends to $V$ in $L^2$ as $t\downarrow0$, and Lemma~\ref{lem:L2-tau} gives the stated properties of $\tau_V^{\rm ren}$.
\end{proof}

\subsection{Negative focusing near the initial time}\label{app:focusing}

We build $V$ as a sum of small nonnegative bumps whose base points approach the origin from below. As a bump shrinks, its effect on the field at its base point tends to an integral of the bump against an Airy kernel whose negative part is not integrable (Lemmas~\ref{lem:focusing-response} and~\ref{lem:focusing-kernel}), and the later bumps are so small in $L^1$ that together they change the field at each earlier observation point by less than $1$.

For a potential $F$ we write $\partial_rF$ for its derivative in the second variable.

\begin{lemma}\label{lem:focusing-response}
Let $B,f\in C_c^\infty(\R^2)$ be nonnegative, let $(x_0,r_0)\in\R^2$, and for $h>0$ put
\[
 V_h(x,r)=B(x,r)+f\Bigl(\frac{x-x_0}{h^2},\frac{r-r_0}{h}\Bigr).
\]
Then
\[
 \lim_{h\downarrow0}\bigl[\phi_{V_h}(h^3,x_0,r_0)-\phi_B(h^3,x_0,r_0)\bigr]=2\E I^{1,0}_{\partial_r^2f}(0).
\]
\end{lemma}

\begin{proof}
\textit{Scaling.} For a potential $W$ put $W^h(y,u)=h^2W(x_0+h^2y,r_0+hu)$. Applying \eqref{eq:occ-scaled} at time $h^3$ gives, for every realization of the lines,
\begin{equation}\label{eq:focusing-scaling}
 I^{h^3,x_0}_W(r_0+h\rho)=I^{1,0}_{W^h}(\rho),\qquad\rho\in\R.
\end{equation}
Taking the expectation of $e^{-I}$ on both sides gives $Z_W(h^3,x_0,r_0+h\rho)=Z_{W^h}(1,0,\rho)$. Put $Q_h=B^h$, so that $W^h=Q_h+h^2f$ for $W=V_h$ and $W^h=Q_h$ for $W=B$. Applying $-2\partial_\rho^2\log$ at $\rho=0$ gives
\[
 \phi_{V_h}(h^3,x_0,r_0)=h^{-2}\phi_{Q_h+h^2f}(1,0,0),\qquad\phi_B(h^3,x_0,r_0)=h^{-2}\phi_{Q_h}(1,0,0).
\]
The lemma therefore concerns the fields of $Q_h$ and $Q_h+h^2f$ at $(t,x,r)=(1,0,0)$. From now on every occupation is taken there, and we drop the superscripts.

\textit{The field.} For $U\in C_c^\infty(\R^2)$ let $\E_U$ be the expectation tilted by $e^{-I_U}$, that is, $\E_UF=\E[Fe^{-I_U}]/Z_U(1,0,0)$, and let $\operatorname{Var}_U$ and $\operatorname{Cov}_U$ be the variance and covariance under it. Differentiating $\log Z_U=\log\E e^{-I_U}$ twice in $r$ under the expectation, with $\partial_rI_U=I_{\partial_rU}$ (both by Lemma~\ref{lem:diff-under-E}), gives
\begin{equation}\label{eq:focusing-field}
 \phi_U(1,0,0)=2\E_UI_{\partial_r^2U}-2\operatorname{Var}_U(I_{\partial_rU}).
\end{equation}

\textit{The difference of the two fields.} Applied to $U=Q_h$ and $U=Q_h+h^2f$, \eqref{eq:focusing-field} involves the occupations of $Q_h$, $f$ and their first two $r$-derivatives. We name them, with the powers of $h$ that keep their moments bounded:
\[
 A_h=I_{Q_h},\qquad D_h=h^{-1}I_{\partial_rQ_h},\qquad C_h=h^{-2}I_{\partial_r^2Q_h},\qquad P_k=I_{\partial_r^kf}\quad(k=0,1,2).
\]
The tilts of $Q_h$ and of $Q_h+h^2f$ have the densities
\[
 q_h=\frac{e^{-A_h}}{\E e^{-A_h}},\qquad q_h^+=\frac{e^{-A_h-h^2P_0}}{\E e^{-A_h-h^2P_0}},
\]
and we write $\E_h$, $\operatorname{Var}_h$, $\operatorname{Cov}_h$ and $\E_h^+$, $\operatorname{Var}_h^+$, $\operatorname{Cov}_h^+$ for the expectation, variance and covariance under them. For $U=Q_h+h^2f$ the derivative occupations are $I_{\partial_rU}=hD_h+h^2P_1$ and $I_{\partial_r^2U}=h^2(C_h+P_2)$, and for $U=Q_h$ they are $hD_h$ and $h^2C_h$. Subtracting the two instances of \eqref{eq:focusing-field}, dividing by $h^2$ and using the scaling step gives
\begin{equation}\label{eq:focusing-difference}
 \begin{aligned}
 \phi_{V_h}(h^3,x_0,r_0)-\phi_B(h^3,x_0,r_0)
 &=2\E_h^+P_2+2(\E_h^+C_h-\E_hC_h)-2(\operatorname{Var}_h^+D_h-\operatorname{Var}_hD_h)\\
 &\quad-4h\operatorname{Cov}_h^+(D_h,P_1)-2h^2\operatorname{Var}_h^+P_1.
 \end{aligned}
\end{equation}
It remains to show that the first term on the right tends to $2\E P_2$ and the others tend to $0$.

\textit{Moments.} Since $\partial_r^kQ_h=h^k(\partial_r^kB)^h$, the identity \eqref{eq:focusing-scaling} at $\rho=0$ gives
\[
 A_h=I^{h^3,x_0}_B(r_0),\qquad D_h=I^{h^3,x_0}_{\partial_rB}(r_0),\qquad C_h=I^{h^3,x_0}_{\partial_r^2B}(r_0).
\]
By Lemma~\ref{lem:occ-moments}(i), every moment of $A_h$, $D_h$ and $C_h$ is bounded uniformly in $0<h\le1$, and $P_0$, $P_1$, $P_2$, which do not depend on $h$, have finite moments of every order. By Proposition~\ref{prop:lln}, $A_h\to\frac12\partial_r^{-2}B(x_0,r_0)$ in $L^2$ as $h\downarrow0$.

\textit{Changing the tilt.} Since $A_h,P_0\ge0$, the numerators of $q_h$ and $q_h^+$ are at most $1$, and by Jensen's inequality both denominators are at least $c=\exp\{-\sup_{0<h\le1}\E A_h-\E P_0\}>0$. Hence $0\le q_h,q_h^+\le c^{-1}$. Using $0\le1-e^{-h^2P_0}\le h^2P_0$ in the numerator and in the denominator gives $|q_h^+-q_h|\le c^{-1}h^2P_0+c^{-2}h^2\E P_0$, and the Cauchy--Schwarz inequality then gives, for every random variable $F$ with finite second moment,
\begin{equation}\label{eq:focusing-tilt}
 |\E_h^+F-\E_hF|\le Ch^2\|F\|_2,
\end{equation}
with $C$ independent of $h$ and $F$. For the first term of \eqref{eq:focusing-difference} we also need the limit of $q_h^+$. Its numerator is at most $1$ and tends to the constant $\exp\{-\frac12\partial_r^{-2}B(x_0,r_0)\}$ in probability, so its denominator tends to the same constant. Hence $q_h^+\to1$ in probability and, being bounded, in $L^2$, and $\E_h^+P_2=\E[q_h^+P_2]\to\E P_2$.

\textit{Conclusion.} In \eqref{eq:focusing-difference}, the first term tends to $2\E P_2$, as just shown. By \eqref{eq:focusing-tilt} with $F=C_h$, the second term is $O(h^2)$. So is the third, because
\[
 \operatorname{Var}_h^+D_h-\operatorname{Var}_hD_h=(\E_h^+D_h^2-\E_hD_h^2)-(\E_h^+D_h-\E_hD_h)(\E_h^+D_h+\E_hD_h)
\]
and the fourth moment of $D_h$ is bounded. The covariance and the variance in the last two terms are bounded, because the densities are bounded and the moments are bounded uniformly in $h$; so these terms are $O(h)$. The right side of \eqref{eq:focusing-difference} therefore tends to $2\E P_2=2\E I^{1,0}_{\partial_r^2f}(0)$.
\end{proof}

\begin{lemma}\label{lem:focusing-kernel}
For every $a>0$ and $N>0$ there is $f\in C_c^\infty(\R^2)$, with $0\le f\le a$ and $f(y,u)=0$ unless $|y|<1$ and $u>1$, such that
\[
 2\E I^{1,0}_{\partial_r^2f}(0)\le-N.
\]
\end{lemma}

\begin{proof}
\textit{The kernel.} For $F\in C_c^\infty(\R^2)$, \eqref{eq:occ-density} at $(t,x,r)=(1,0,0)$ gives
\[
 \E I^{1,0}_F(0)=\iint F(y,u)\rho_{\rm Ai}(y^2-u)\diff u\diff y.
\]
From $\rho_{\rm Ai}=(\Ai')^2-z\Ai^2$ and $\Ai''=z\Ai$ we get $\rho_{\rm Ai}'=-\Ai^2$ and $\rho_{\rm Ai}''=-2\Ai\Ai'$. Two integrations by parts in $u$ therefore give
\begin{equation}\label{eq:focusing-kernel}
 2\E I^{1,0}_{\partial_r^2F}(0)=\iint F(y,u)k(u-y^2)\diff u\diff y,\qquad k(w)=-4\Ai(-w)\Ai'(-w).
\end{equation}

\textit{Its negative part is not integrable.} The asymptotics of $\Ai$ and $\Ai'$ at $-\infty$ \cite[9.7.9--9.7.10]{DLMF} give
\[
 k(w)=\frac2\pi\cos\Bigl(\frac43w^{3/2}\Bigr)+O(w^{-3/2}),\qquad w\to\infty,
\]
so that $k_-(w)\ge\frac2\pi\bigl(\cos\frac43w^{3/2}\bigr)_--Cw^{-3/2}$ for $w\ge1$. The substitution $v=\frac43w^{3/2}$, with $\diff w=\frac12(\frac34v)^{-1/3}\diff v$, turns $\int_1^\infty(\cos\frac43w^{3/2})_-\diff w$ into the integral of $(\cos v)_-$, which is periodic with positive mean, against a decreasing weight that is not integrable; so it diverges, and $\int_1^\infty k_-(w)\diff w=\infty$.

\textit{The bump.} Fix $\chi\in C_c^\infty((-1,1))$ with $0\le\chi\le1$ and $\int\chi>0$, and choose $w_0>1$ with $\int_1^{w_0}k_-\ge1+N/(a\int\chi)$. Inner regularity on the open set $O=\{w\in(1,w_0):k(w)<0\}$ and a smooth cutoff give $g\in C_c^\infty(O)$ with $0\le g\le1$ and $\int gk_-\ge\int_1^{w_0}k_--1$. Since $g$ vanishes where $k\ge0$,
\[
 \int g(w)k(w)\diff w=-\int g(w)k_-(w)\diff w\le-\frac N{a\int\chi}.
\]
Put $f(y,u)=a\chi(y)g(u-y^2)$. Then $f\in C_c^\infty(\R^2)$, $0\le f\le a$, and $f(y,u)=0$ unless $|y|<1$ and $u>1+y^2$. By \eqref{eq:focusing-kernel} and the substitution $w=u-y^2$,
\[
 2\E I^{1,0}_{\partial_r^2f}(0)=a\int\chi(y)\diff y\int g(w)k(w)\diff w\le-N.
\]
\end{proof}

\begin{proof}[Proof of Remark~\ref{rem:focusing}]
\textit{The bumps.} Let $p_n=(0,-2^{-n})$ and $D=[-1,1]\times[-1,0]$. Lemma~\ref{lem:focusing-kernel} with $a=2^{-n}$ and $N=n+2$ gives $f_n$; fix $b_n>1$ such that $f_n(y,u)=0$ for $u\ge b_n$. For $h>0$ let
\[
 F_{n,h}(x,r)=f_n\Bigl(\frac x{h^2},\frac{r+2^{-n}}h\Bigr),
\]
the bump of Lemma~\ref{lem:focusing-response} with $f=f_n$ and $(x_0,r_0)=p_n$. Its $L^1$ norm is $h^3\|f_n\|_1$, and it vanishes unless $|x|<h^2$ and $-2^{-n}+h<r<-2^{-n}+hb_n$. If $h<2^{-n-3}/b_n$, it therefore vanishes outside the band
\[
 S_n=(-1,1)\times(-2^{-n},-7\cdot2^{-n-3})\subset D.
\]
The bands are pairwise disjoint, and $p_n$ lies above the closure of $S_j$ for every $j<n$.

\textit{Continuity in the data.} Let $\mathcal W$ be the set of measurable $W$ with $0\le W\le1$ that vanish outside $D$. Every $W\in\mathcal W$ satisfies $|W(x,r)|\le e^{-r_+}$, the hypothesis of Theorem~\ref{thm:bounded-data} with $M=\delta=1$. By Proposition~\ref{prop:bounded-stability}, if $W_k\in\mathcal W$ and $W_k\to W_0$ in $L^1$, then $\log Z_{W_k}\to\log Z_{W_0}$ in $C^\infty_{\rm loc}$ at positive time, so $\phi_{W_k}\to\phi_{W_0}$ pointwise. Consequently, for every $W_0\in\mathcal W$ and every $(t,x,r)$ with $t>0$ there is $\eta>0$ such that
\begin{equation}\label{eq:focusing-continuity}
 W\in\mathcal W,\quad\|W-W_0\|_1<\eta\quad\Longrightarrow\quad|\phi_W(t,x,r)-\phi_{W_0}(t,x,r)|<1.
\end{equation}

\textit{The inductive choice.} We choose $h_1,h_2,\ldots$ and $\eta_1,\eta_2,\ldots$ in turn. Suppose that $h_j$ and $\eta_j$ have been chosen for $j<n$, and put $t_j=h_j^3$ and $B_{n-1}=\sum_{j<n}F_{j,h_j}$, with $B_0=0$. Then $B_{n-1}$ is nonnegative, smooth and compactly supported, and $B_{n-1}(p_n)=0$ because $p_n$ lies above the closures of $S_1,\ldots,S_{n-1}$. By Theorem~\ref{thm:compact}, $\phi_{B_{n-1}}(h^3,p_n)\to B_{n-1}(p_n)=0$ as $h\downarrow0$, so Lemma~\ref{lem:focusing-response} with $B=B_{n-1}$, $f=f_n$ and $(x_0,r_0)=p_n$ gives
\[
 \lim_{h\downarrow0}\phi_{B_{n-1}+F_{n,h}}(h^3,p_n)=2\E I^{1,0}_{\partial_r^2f_n}(0)\le-(n+2).
\]
Choose $h_n<2^{-n-3}/b_n$ so small that, with $t_n=h_n^3$ and $B_n=B_{n-1}+F_{n,h_n}$,
\[
 \phi_{B_n}(t_n,p_n)\le-(n+1),\qquad t_n<\min(2^{-n},t_{n-1}),\qquad\|F_{n,h_n}\|_1<2^{-n}\eta_j\quad(j<n),
\]
where the conditions involving $t_{n-1}$ and $\eta_j$ are void for $n=1$. The summands of $B_n$ have disjoint supports in $D$ and are at most $\frac12$, so $B_n\in\mathcal W$. Let $\eta_n$ be given by \eqref{eq:focusing-continuity} with $W_0=B_n$ and $(t,x,r)=(t_n,p_n)$.

\textit{The potential.} Put $V=\sum_{n\ge1}F_{n,h_n}$; as for $B_n$, $V\in\mathcal W$. The support of $F_{n,h_n}$ lies within distance $h_n^2+2^{-n}$ of the origin and in the closure of $S_n$, which does not contain the origin. So near each point other than the origin $V$ is a finite sum of smooth summands, and for each $n$ it is at most $2^{-n}$ near the origin; thus $V$ is smooth away from the origin and continuous at the origin, where it vanishes. For each $n$, $V-B_n=\sum_{j>n}F_{j,h_j}$, so
\[
 \|V-B_n\|_1<\sum_{j>n}2^{-j}\eta_n<\eta_n,
\]
and \eqref{eq:focusing-continuity} gives
\[
 \phi_V(t_n,p_n)<\phi_{B_n}(t_n,p_n)+1\le-n.
\]
Since $t_n\downarrow0$ and $p_n\to(0,0)$, this proves the remark, with $(x_n,r_n)=p_n$.
\end{proof}

\section{Nonintersecting Brownian bridges and Toda integrability}
\label{app:section3}

This appendix contains the proofs for Section~\ref{sec:dyson-toda-kp}, and the limit \eqref{eq:bridge-limit} used in Section~\ref{sec:soliton}.

\subsection{The bridge limit}
\label{app:bridge-limit}\label{app:edge-scaled}

This subsection proves the limit \eqref{eq:db-target} and Proposition~\ref{prop:bridge-limit}. Let $\widetilde B_1>\dots>\widetilde B_N$ be nonintersecting Brownian bridges with diffusion coefficient $2$ on $[-T,T]$, all starting and ending at $0$, and let $t=(T^3/N)^{1/2}$. Define $G^N_1>\dots>G^N_N$ by
\begin{equation}\label{eq:db-edge-scaled}
 \widetilde B_i(y)=2\sqrt{NT}+t^{1/3}G_i^N\bigl(t^{-2/3}y\bigr),\qquad|y|\le T.
\end{equation}
Then, for each $k$,
\begin{equation}\label{eq:db-edge-limit}
 \bigl(G_1^N(\eta),\dots,G_k^N(\eta)\bigr)\Longrightarrow\bigl(\bar{\mathcal A}_1(\eta)-\eta^2,\dots,\bar{\mathcal A}_k(\eta)-\eta^2\bigr)
\end{equation}
for the topology of uniform convergence on compact sets of $\eta$. This follows from Theorem~3.1 of Corwin and Hammond \cite{CH} by Brownian scaling. Let $B_i$ be the bridges of Section~\ref{sec:dyson-brownian-scaling}. By \eqref{eq:db-endpoints},
\begin{equation}\label{eq:db-reflected}
 B_i(y)\stackrel d=2\sqrt{NT}+r+\frac{x^2}t-\frac{2xy}t-\widetilde B_i(y)
\end{equation}
jointly in $i$ and $y$. By \eqref{eq:db-edge-scaled} and \eqref{eq:db-edge-limit}, the right side converges to $r+(y-x)^2/t-t^{1/3}\bar{\mathcal A}_i(t^{-2/3}y)$, which is $r-\mathcal A_i^{t,x}(y)$ by \eqref{eq:markov-airy-marginal}. This proves \eqref{eq:db-target}.

Next we bound the number of bridges $G_i^N$ above a given height, uniformly in $N$. For $c\in\R$ and $|\eta|\le N^{1/3}$, let $\mathcal N_c(\eta)$ be the number of $i$ with $G_i^N(\eta)\ge c$.

\begin{lemma}\label{lem:bridge-counts}
For $c\in\R$, $|\eta|\le\frac12N^{1/3}$ and $\theta\ge0$,
\[
 \E e^{\theta\mathcal N_c(\eta)}\le\exp\bigl\{(e^\theta-1)C_0e^{-\min(c,2c)}\bigr\},\qquad C_0=(4\pi)^{-1/2}e^{1/12}.
\]
\end{lemma}

\begin{proof}
By Brownian scaling, the law of $G^N$ does not depend on $T$, so we take $T=1$. Then \eqref{eq:db-edge-scaled} reads $G_i^N(\eta)=N^{1/6}\bigl(\widetilde B_i(N^{-1/3}\eta)-2\sqrt N\bigr)$, and at a fixed time $y\in(-1,1)$ the heights $\widetilde B_i(y)$ have the law of the eigenvalues of $\sqrt{1-y^2}M$, where $M$ is an $N\times N$ GUE matrix with density proportional to $e^{-\operatorname{tr}M^2/2}$.

We first show that, for every $c\in\R$,
\begin{equation}\label{eq:bridge-mean-count}
 \E\mathcal N_c(0)\le C_0e^{-c}.
\end{equation}
The heights $G_i^N(0)$ form a determinantal point process whose kernel $P_N$ is the orthogonal projection onto the eigenfunctions of $H_N=-\partial_z^2+z+z^2/(4N^{2/3})$ with negative eigenvalues. Since these eigenvalues are negative, $P_N(z,z)\le e^{-H_N}(z,z)$. By the Feynman--Kac formula,
\[
 e^{-H_N}(z,z)=p_1(0)\,\E\exp\Bigl\{-\int_0^1\Bigl(\beta(\sigma)+\frac{\beta(\sigma)^2}{4N^{2/3}}\Bigr)\diff\sigma\Bigr\}\le p_1(0)\,\E\exp\Bigl\{-\int_0^1\beta(\sigma)\diff\sigma\Bigr\}=C_0e^{-z},
\]
where $\beta$ is a Brownian bridge with diffusion coefficient $2$ from $z$ to $z$ on $[0,1]$, and $\int_0^1\beta(\sigma)\diff\sigma$ is Gaussian with mean $z$ and variance $1/6$. Hence $\E\mathcal N_c(0)=\int_c^\infty P_N(z,z)\diff z\le C_0e^{-c}$.

Next let $|\eta|\le\frac12N^{1/3}$, and let $\gamma=1-N^{-2/3}\eta^2$, so that $\gamma\ge\frac34$. The heights $\widetilde B_i(N^{-1/3}\eta)$ have the law of $\sqrt\gamma\,\widetilde B_i(0)$, so the heights $G_i^N(\eta)$ have the law of $\sqrt\gamma\,G_i^N(0)-2N^{2/3}(1-\sqrt\gamma)$. Hence $\mathcal N_c(\eta)$ has the law of $\mathcal N_{c'}(0)$, where
\[
 c'=\frac{c+2N^{2/3}(1-\sqrt\gamma)}{\sqrt\gamma}\ge\frac c{\sqrt\gamma}\ge\min(c,2c).
\]

Finally, as in the proof of Lemma~\ref{lem:counts}(i), $\E e^{\theta\mathcal N_{c'}(0)}\le\exp\{(e^\theta-1)\E\mathcal N_{c'}(0)\}$. By \eqref{eq:bridge-mean-count} and since $c'\ge\min(c,2c)$, this proves the lemma.
\end{proof}

Next we show that the occupation of the bridges $G_i^N$ converges to that of the Airy lines. The limit \eqref{eq:db-edge-limit} gives this for the top $k$ bridges, and Lemma~\ref{lem:bridge-counts} controls the others uniformly in $N$.

For a bounded measurable function $W$ on $\R^2$ with compact support put
\[
 J_N(W)=\sum_{i=1}^N\int_{-N^{1/3}}^{N^{1/3}}W\bigl(\eta,G_i^N(\eta)\bigr)\diff\eta,\qquad J(W)=\sum_{i\ge1}\int_\R W\bigl(\eta,\bar{\mathcal A}_i(\eta)-\eta^2\bigr)\diff\eta.
\]

\begin{lemma}\label{lem:occupation-limit}
Let $W_0,\dots,W_m$ be continuous with compact support, and let $W_{j,N}\to W_j$ uniformly as $N\to\infty$, with all supports in a fixed compact set. Then, for every polynomial $P$ in $m$ variables,
\[
 \E\Bigl[P\bigl(J_N(W_{1,N}),\dots,J_N(W_{m,N})\bigr)e^{-J_N(W_{0,N})}\Bigr]\longrightarrow\E\Bigl[P\bigl(J(W_1),\dots,J(W_m)\bigr)e^{-J(W_0)}\Bigr],
\]
and the limit is finite.
\end{lemma}

\begin{proof}
Let all supports lie in $[-L,L]\times[c,\infty)$, and let $M$ bound all $|W_{j,N}|$ and $|W_j|$. Then $|J_N(W_{j,N})|\le MJ_N^*$, where $J_N^*=\int_{-L}^L\mathcal N_c(\eta)\diff\eta$. Jensen's inequality for the average over $[-L,L]$ and Lemma~\ref{lem:bridge-counts} give, for $p\ge0$ and $N\ge(2L)^3$,
\[
 \E e^{pJ_N^*}\le\frac1{2L}\int_{-L}^L\E e^{2pL\mathcal N_c(\eta)}\diff\eta\le\exp\bigl\{\bigl(e^{2pL}-1\bigr)C_0e^{-\min(c,2c)}\bigr\}.
\]
The same bound, with a different constant, holds for the Airy lines by Lemma~\ref{lem:counts}(i), with $J^*=\int_{-L}^LN_c(\eta)\diff\eta$ and $N_c(\eta)$ as in that lemma. So the limit in the lemma is finite, and the variables $P(\dots)e^{-J_N(W_{0,N})}$ are bounded in $L^2$ uniformly in $N$. It therefore suffices to prove that $\bigl(J_N(W_{0,N}),\dots,J_N(W_{m,N})\bigr)$ converges in law to $\bigl(J(W_0),\dots,J(W_m)\bigr)$.

First, $|J_N(W_{j,N})-J_N(W_j)|\le\|W_{j,N}-W_j\|_\infty J_N^*$, which tends to $0$ in $L^1$, so we may replace $W_{j,N}$ by $W_j$. Next fix $k\ge1$ and split $J_N(W_j)$ into the sum over $i\le k$ and the sum over $i>k$. The bridges are ordered, so at each $\eta$ at most $(\mathcal N_c(\eta)-k)_+$ of the bridges with $i>k$ meet the support of $W_j$. Hence
\[
 \E\Bigl|\sum_{i>k}\int W_j\bigl(\eta,G_i^N(\eta)\bigr)\diff\eta\Bigr|\le M\int_{-L}^L\E\bigl(\mathcal N_c(\eta)-k\bigr)_+\diff\eta\le\frac Mk\int_{-L}^L\E\mathcal N_c(\eta)^2\diff\eta,
\]
which is at most $C/k$ uniformly in $N$, by Lemma~\ref{lem:bridge-counts}; the same bound holds for the Airy lines. The sums over $i\le k$ are continuous functions of $(G_1^N,\dots,G_k^N)$ restricted to $[-L,L]$, in the uniform topology, so by \eqref{eq:db-edge-limit} they converge jointly in law to the corresponding sums for the Airy lines. Letting $N\to\infty$ and then $k\to\infty$ gives the convergence in law \cite[Theorem~3.2]{Bil}.
\end{proof}

\begin{proof}[Proof of Proposition~\ref{prop:bridge-limit}]
Let $V\in C_c^\infty(\R^2)$ vanish outside $[-Y,Y]\times[-U,U]$, and let $\mathcal C\subset(0,\infty)\times\R^2$ be compact. For $t>0$ and $x,r\in\R$ define
\[
 W_{t,x,r}(\eta,z)=t^{2/3}V\Bigl(t^{2/3}\eta,\ r+\frac{x^2}t-2xt^{-1/3}\eta-t^{1/3}z\Bigr).
\]
Then, with $T$, $u$ and $v$ given by \eqref{eq:db-endpoints},
\[
 Z_N^V(u,v)=\E e^{-J_N(W_{t,x,r})},\qquad Z_V(t,x,r)=\E e^{-J(W_{t,x,r})}.
\]
Indeed, by \eqref{eq:db-reflected} and \eqref{eq:db-edge-scaled}, the heights $B_i(y)$ have the law of $r+x^2/t-2xy/t-t^{1/3}G_i^N(t^{-2/3}y)$, jointly in $i$ and $y$, so the substitution $y=t^{2/3}\eta$ in the definition of $Z_N^V$ gives the first identity. The second follows in the same way from \eqref{eq:airy-occupation} and \eqref{eq:markov-airy-marginal}.

Next we show that $Z_V$ is finite, positive and smooth. For $(t,x,r)\in\mathcal C$, the function $W_{t,x,r}$ and its derivatives of every order in $t$, $x$ and $r$ are continuous in all variables and vanish outside a fixed compact set in $[-L,L]\times[c,\infty)$. A derivative of $e^{-J(W_{t,x,r})}$ in $t$, $x$ and $r$ is $e^{-J(W_{t,x,r})}$ times a polynomial in the occupations of derivatives of $W_{t,x,r}$, and, as in the proof of Lemma~\ref{lem:occupation-limit}, these occupations are bounded by a constant times $J^*=\int_{-L}^LN_c(\eta)\diff\eta$. Since $J^*$ has exponential moments of every order, dominated convergence allows differentiation under the expectation, and $Z_V$ is finite and smooth on $(0,\infty)\times\R^2$. Jensen's inequality gives $Z_V\ge e^{-\E J(|W_{t,x,r}|)}$, which is bounded below by a positive constant on $\mathcal C$. In the same way $Z_N^V$ is smooth in $x$ and $r$, its derivatives are expectations of the same polynomials in occupations $J_N$, and $Z_N^V\ge e^{-\E J_N(|W_{t,x,r}|)}$, which is bounded below on $\mathcal C$ uniformly in large $N$, by Lemma~\ref{lem:bridge-counts}.

It remains to prove the convergence. At fixed $N$ and $t$, derivatives in $x$ and $r$ act only on $W_{t,x,r}$. So the derivatives $\partial_x^j\partial_r^kZ_N^V(u,v)$ and $\partial_x^j\partial_r^kZ_V(t,x,r)$ are the two sides of Lemma~\ref{lem:occupation-limit}, with the same polynomial $P$ and with test functions that are derivatives of $W_{t,x,r}$. Let $(t_N,x_N,r_N)\in\mathcal C$ converge to $(t,x,r)$. The test functions at $(t_N,x_N,r_N)$ converge uniformly to those at $(t,x,r)$, with supports in a fixed compact set, so Lemma~\ref{lem:occupation-limit} gives $\partial_x^j\partial_r^kZ_N^V\to\partial_x^j\partial_r^kZ_V(t,x,r)$ along the sequence. Since the sequence is arbitrary and the derivatives of $Z_V$ are continuous, the convergence is uniform on $\mathcal C$. The positive lower bounds on $Z_N^V$ and $Z_V$ give the same convergence for the derivatives of $\log Z_N^V$ and $\log Z_V$.
\end{proof}

\subsection{The bridge limit for potentials with one-sided decay}
\label{app:s8-bridge}

The next proposition is the limit \eqref{eq:bridge-limit} used in the proof of Theorem~\ref{thm:onestep}.

\begin{proposition}\label{prop:s8-bridge-limit}
Let $V$ be bounded and continuous on $\R^2$, with $|V(y,u)|\le Ce^{c|y|-\delta u}$ for some $C,c,\delta>0$. Then, with $a$, $b$ and $R$ given by \eqref{eq:s8-edge}, $Z_N(a,b;R;V)\to Z_V(t,x,r)$ as $N\to\infty$.
\end{proposition}

By \eqref{eq:s8-edge}, we may take $X_i(y)=R-\widetilde B_i(y-x)$, where $\widetilde B_i$ are the bridges of Appendix~\ref{app:edge-scaled}. For a function $F$ on the plane let
\[
 I_{F,N}=\sum_{i=1}^N\int_a^bF\bigl(y,X_i(y)\bigr)\diff y,\qquad
 I_F=\sum_{i\ge1}\int_\R F\bigl(y,r-\mathcal A_i^{t,x}(y)\bigr)\diff y,
\]
so that $Z_N(a,b;R;V)=\E e^{-I_{V,N}}$ and $Z_V(t,x,r)=\E e^{-I_V}$.

We first prove the limit for test functions with compact support: for continuous $F_0,\dots,F_k$ with compact support and every polynomial $P$,
\begin{equation}\label{eq:s8-compact}
 \E\bigl[P(I_{F_1,N},\dots,I_{F_k,N})e^{-I_{F_0,N}}\bigr]\longrightarrow\E\bigl[P(I_{F_1},\dots,I_{F_k})e^{-I_{F_0}}\bigr].
\end{equation}
Indeed, by \eqref{eq:db-edge-scaled}, $X_i(x+t^{2/3}\eta)=r-t^{1/3}G_i^N(\eta)$, and by \eqref{eq:markov-airy-marginal} and the stationarity of the Airy line ensemble, the curves $r-\mathcal A_i^{t,x}(x+t^{2/3}\eta)$ have the law of $r+t^{1/3}\eta^2-t^{1/3}\bar{\mathcal A}_i(\eta)$, jointly in $i$ and $\eta$. So $I_{F,N}=J_N(W)$ and $I_F$ has the law of $J(W)$, with $W(\eta,z)=t^{2/3}F(x+t^{2/3}\eta,r-t^{1/3}z)$, and \eqref{eq:s8-compact} is Lemma~\ref{lem:occupation-limit}.

For general $V$ we need a bound on the bridges that is uniform in $N$.

\begin{lemma}\label{lem:s8-hermite}
For every $\delta>0$ there is a constant $C$, depending on $\delta$ and $t$, such that
\[
 \E\sum_{i=1}^Ne^{\delta(r-X_i(y))}\le Ce^{-\delta(y-x)^2/t}
\]
for all $N\ge1$ and $y\in[a,b]$.
\end{lemma}

\begin{proof}
By \eqref{eq:db-edge-scaled}, $r-X_i(y)=t^{1/3}G_i^N(\eta)$ with $\eta=t^{-2/3}(y-x)$. So, with $\theta=\delta t^{1/3}$, we have to show that $\E\sum_ie^{\theta G_i^N(\eta)}\le Ce^{-\theta\eta^2}$ for all $N\ge1$ and $|\eta|\le N^{1/3}$. We use the notation of the proof of Lemma~\ref{lem:bridge-counts}. We first show that, for $\theta/2\le\theta'\le\theta$,
\begin{equation}\label{eq:s8-centre}
 \E\sum_{i=1}^Ne^{\theta'G_i^N(0)}\le C.
\end{equation}
The left side is $\int_\R e^{\theta'z}P_N(z,z)\diff z$, and the proof of \eqref{eq:bridge-mean-count}, with $e^{-\lambda H_N}$ in place of $e^{-H_N}$, gives $P_N(z,z)\le(4\pi\lambda)^{-1/2}e^{\lambda^3/12-\lambda z}$ for every $\lambda>0$. With $\lambda=\theta+1$ this is at most $Ce^{-(\theta+1)z}$ for $z\ge0$, and with $\lambda=(1+|z|)^{-1}$ it is at most $C(1+|z|)^{1/2}$ for $z<0$, which proves \eqref{eq:s8-centre}.

Next let $\gamma=1-N^{-2/3}\eta^2$. As in that proof, the heights $G_i^N(\eta)$ have the law of $\sqrt\gamma\,G_i^N(0)-2N^{2/3}(1-\sqrt\gamma)$, and $2N^{2/3}(1-\sqrt\gamma)\ge N^{2/3}(1-\gamma)=\eta^2$. If $\gamma\ge\frac14$, then $\theta/2\le\theta\sqrt\gamma\le\theta$, so by \eqref{eq:s8-centre}
\[
 \E\sum_{i=1}^Ne^{\theta G_i^N(\eta)}=e^{-2\theta N^{2/3}(1-\sqrt\gamma)}\,\E\sum_{i=1}^Ne^{\theta\sqrt\gamma\,G_i^N(0)}\le Ce^{-\theta\eta^2}.
\]
If $\gamma<\frac14$, the left side is $e^{-2\theta N^{2/3}}\E\operatorname{tr}e^{\theta N^{1/6}\sqrt\gamma M}$. The law of $M$ is symmetric, so $\E\operatorname{tr}e^{cM}=\E\operatorname{tr}\cosh(cM)$, which increases with $c\ge0$. Since $\sqrt\gamma<\frac12$,
\[
 \E\sum_{i=1}^Ne^{\theta G_i^N(\eta)}\le e^{-2\theta N^{2/3}}\E\operatorname{tr}e^{\theta N^{1/6}M/2}=e^{-\theta N^{2/3}}\,\E\sum_{i=1}^Ne^{\theta G_i^N(0)/2}\le Ce^{-\theta N^{2/3}}\le Ce^{-\theta\eta^2},
\]
by \eqref{eq:s8-centre} and since $\eta^2\le N^{2/3}$.
\end{proof}

Next we show that, for $p>0$ and continuous $F$ with $0\le F(y,u)\le\min(C,Ce^{c|y|-\delta u})$,
\begin{equation}\label{eq:s8-expmoments}
 \sup_N\E e^{pI_{F,N}}<\infty.
\end{equation}
At a fixed time $y$ the heights $X_i(y)$ form a determinantal point process, so, as in the proof of Lemma~\ref{lem:counts}(i), for every bounded function $g\ge0$,
\begin{equation}\label{eq:s8-dpp}
 \E e^{\sum_ig(X_i(y))}\le\exp\Bigl\{\E\sum_i\bigl(e^{g(X_i(y))}-1\bigr)\Bigr\}\le\exp\Bigl\{e^{\sup g}\,\E\sum_ig\bigl(X_i(y)\bigr)\Bigr\}.
\end{equation}
Let $O_N(y)=\sum_iF(y,X_i(y))$ for $y\in[a,b]$ and $O_N(y)=0$ otherwise, and let $w$ be the probability density on $\R$ proportional to $(1+y^2)^{-3/4}$. By Jensen's inequality, $e^{pI_{F,N}}\le\int_\R w(y)e^{pO_N(y)/w(y)}\diff y$. By Lemma~\ref{lem:s8-hermite}, $\E O_N(y)\le Ce^{c|y|-\delta(y-x)^2/t}$. So \eqref{eq:s8-dpp} with $g=pF(y,\cdot)/w(y)$ gives
\[
 \E e^{pO_N(y)/w(y)}\le\exp\Bigl\{\frac C{w(y)}\exp\Bigl[c|y|-\frac{\delta(y-x)^2}t+\frac{pC}{w(y)}\Bigr]\Bigr\},
\]
with $C$ depending on $p$. Since $1/w(y)$ grows like $|y|^{3/2}$, the right side is bounded in $y$, and it does not depend on $N$. This proves \eqref{eq:s8-expmoments}.

It remains to remove the compact support. Let $\chi_k$, $k\ge1$, be continuous functions on $\R^2$ with $0\le\chi_k\le1$, increasing in $k$, equal to $1$ on $[-k,k]^2$ and to $0$ outside $[-2k,2k]^2$. We show that
\begin{equation}\label{eq:s8-tail}
 \sup_N\E I_{(1-\chi_k)|V|,N}\longrightarrow0\qquad(k\to\infty).
\end{equation}
Since $1-\chi_k$ vanishes on $[-k,k]^2$, and since $|V|\le Ce^{c|y|-\delta u}$ and $|V|\le C$,
\[
 (1-\chi_k)|V|(y,u)\le Ce^{-\delta k/2}e^{c|y|-\delta u/2}+Ce^{-k}e^{-u}+C\mathbf 1_{|y|>k}e^{c|y|-\delta u}:
\]
the first term bounds the left side where $u>k$, the second where $u<-k$, and the third elsewhere. By Lemma~\ref{lem:s8-hermite} with $\delta/2$, $1$ and $\delta$ in place of $\delta$, for all $N$ and $y\in[a,b]$,
\[
 \E\sum_{i=1}^N(1-\chi_k)|V|\bigl(y,X_i(y)\bigr)\le C\bigl(e^{-\delta k/2}+e^{-k}+\mathbf 1_{|y|>k}\bigr)e^{c|y|-\delta'(y-x)^2/t},\qquad\delta'=\min(\delta/2,1),
\]
and integrating over $y$ gives \eqref{eq:s8-tail}.

The bounds \eqref{eq:s8-expmoments} and \eqref{eq:s8-tail} also hold for the Airy lines, with $I_F$ in place of $I_{F,N}$ and without the supremum. Indeed, for $F$ as in \eqref{eq:s8-expmoments}, the limit \eqref{eq:s8-compact} gives $\E e^{pI_{\chi_jF}}=\lim_N\E e^{pI_{\chi_jF,N}}\le\sup_N\E e^{pI_{F,N}}$, and monotone convergence in $j$ gives $\E e^{pI_F}\le\sup_N\E e^{pI_{F,N}}$. In the same way, $\E I_{(1-\chi_k)|V|}\le\sup_N\E I_{(1-\chi_k)|V|,N}$.

\begin{proof}[Proof of Proposition~\ref{prop:s8-bridge-limit}]
Since $-I_{V,N}$ and $-I_{\chi_kV,N}$ are at most $I_{|V|,N}$ and differ by at most $I_{(1-\chi_k)|V|,N}$,
\[
 \bigl|e^{-I_{V,N}}-e^{-I_{\chi_kV,N}}\bigr|\le e^{I_{|V|,N}}\min\bigl(I_{(1-\chi_k)|V|,N},1\bigr).
\]
The square of the minimum is at most $I_{(1-\chi_k)|V|,N}$. So, by the Cauchy--Schwarz inequality, \eqref{eq:s8-expmoments} with $F=|V|$ and $p=2$, and \eqref{eq:s8-tail},
\[
 \sup_N\E\bigl|e^{-I_{V,N}}-e^{-I_{\chi_kV,N}}\bigr|\le\sup_N\bigl(\E e^{2I_{|V|,N}}\bigr)^{1/2}\bigl(\E I_{(1-\chi_k)|V|,N}\bigr)^{1/2}\longrightarrow0\qquad(k\to\infty),
\]
and in the same way $\E|e^{-I_V}-e^{-I_{\chi_kV}}|\to0$. For each $k$, $\E e^{-I_{\chi_kV,N}}\to\E e^{-I_{\chi_kV}}$ as $N\to\infty$, by \eqref{eq:s8-compact}. Hence $\E e^{-I_{V,N}}\to\E e^{-I_V}$.
\end{proof}

The next lemma shows that Proposition~\ref{prop:s8-bridge-limit} applies to $V_f$, and to the potentials $V_n$ of Section~\ref{sec:s8-several} whenever $W_n$ is a sum of exponentials with nonnegative coefficients, as for totally nonnegative data.

\begin{lemma}\label{lem:s8-soliton-potentials}
Let $\tau(y,u)=\sum_\alpha b_\alpha e^{q_\alpha u-d_\alpha y}$ be a finite sum with $b_\alpha\ge0$, not all zero, and let $V=-2\partial_u^2\log\tau$. Then $V$ is smooth and bounded, and $|V(y,u)|\le Ce^{c|y|-\delta u}$ for some $C,c,\delta>0$.
\end{lemma}

\begin{proof}
Under the probabilities proportional to the terms of $\tau$, $\partial_u^2\log\tau$ is the variance of $q_\alpha$. So $V$ is smooth and $-\frac12(q_{\max}-q_{\min})^2\le V\le0$, where $q_{\max}$ and $q_{\min}$ are the largest and smallest $q_\alpha$ with $b_\alpha>0$. The variance is at most the mean of $(q_\alpha-q_{\max})^2$, so
\[
 |V(y,u)|\le\frac{2(q_{\max}-q_{\min})^2}{\tau(y,u)}\sum_{q_\alpha<q_{\max}}b_\alpha e^{q_\alpha u-d_\alpha y}.
\]
If all $q_\alpha$ with $b_\alpha>0$ are equal, then $V=0$ and the conclusion is immediate. Otherwise, let $\delta>0$ be the smallest of the gaps $q_{\max}-q_\alpha$ among terms with $b_\alpha>0$ and $q_\alpha<q_{\max}$. Since $\tau(y,u)\ge b_{\alpha_*}e^{q_{\max}u-d_{\alpha_*}y}$ for some $\alpha_*$ with $q_{\alpha_*}=q_{\max}$ and $b_{\alpha_*}>0$, the right side is at most $Ce^{c|y|-\delta u}$ for $u\ge0$. For $u<0$ the bound holds because $V$ is bounded.
\end{proof}

\subsection{The Toda equation}
\label{app:toda-kp-proofs}

\begin{proof}[Proof of Proposition~\ref{prop:toda}]
We first prove \eqref{eq:db-ratio}. Since $V$ vanishes for $|y|\ge Y$,
\[
 P_{-T,T;V}(u,v)=\iint_{\R^2}p_{T-Y}(w-u)\,P_{-Y,Y;V}(w,w')\,p_{T-Y}(v-w')\diff w\diff w',
\]
which is smooth in $u$ and $v$, so the determinants $\tau_N^V$ are defined. Consider $N$ independent Brownian motions $X_1,\dots,X_N$ with diffusion coefficient $2$, started at heights $z_1<\dots<z_N$ at time $-T$. By the Karlin--McGregor formula \cite{KM}, $\det[P_{-T,T;V}(z_i,z_j')]_{i,j=1}^N$ is the density at $(z_1',\dots,z_N')$, with $z_1'<\dots<z_N'$, of their positions at time $T$, on the event that the paths do not intersect, weighted by $\exp\{-\sum_i\int_{-T}^TV(y,X_i(y))\diff y\}$. Dividing it by the same determinant with $V=0$ gives the expectation of this weight under nonintersecting bridges from the $z_i$ to the $z_j'$. Let the $z_i$ coalesce at $u$ and the $z_j'$ at $v$. After division by $\prod_{i<j}(z_j-z_i)(z_j'-z_i')$, both determinants tend to the corresponding $\tau_N$ times $\prod_{j<N}(j!)^{-2}$, which gives the ratio in \eqref{eq:db-ratio}. It is positive, since $Z_N^V$ is the expectation of a positive random variable.

Next we prove \eqref{eq:db-toda}. Differentiating $\tau_N^V$ in $u$ differentiates one row at a time, and the derivative of row $i<N-1$ equals row $i+1$, so only the derivative of the last row survives. Likewise, in $v$ only the derivative of the last column survives. Let $M$ be the $(N+1)\times(N+1)$ matrix $[\partial_{u}^i\partial_{v}^jP_{-T,T;V}]_{i,j=0}^N$, and let $M^{(i,j)}$ be $M$ without row $i$ and column $j$. Then
\[
 \begin{gathered}
 \tau_N^V=\det M^{(N,N)},\qquad\partial_{u}\partial_{v}\tau_N^V=\det M^{(N-1,N-1)},\\
 \partial_{u}\tau_N^V=\det M^{(N-1,N)},\qquad\partial_{v}\tau_N^V=\det M^{(N,N-1)},
 \end{gathered}
\]
$\tau_{N+1}^V=\det M$, and $\tau_{N-1}^V$ is the determinant of $M$ without rows $N-1$, $N$ and columns $N-1$, $N$. The Desnanot--Jacobi identity for $M$ and its last two rows and columns is \eqref{eq:toda-desnanot}, and hence
\begin{equation}\label{eq:toda-log-tau}
 \partial_u\partial_v\log\tau_N^V=\frac{\tau_{N+1}^V\tau_{N-1}^V}{(\tau_N^V)^2}.
\end{equation}
For $V=0$ a Gaussian computation gives $\tau_N^0=p_{2T}(v-u)^N\prod_{j=0}^{N-1}j!\,(4T)^{-j}$, so both sides of \eqref{eq:toda-log-tau} for $V=0$ equal $N/(4T)$. Subtracting \eqref{eq:toda-log-tau} for $V=0$ from \eqref{eq:toda-log-tau} for $V$, and using \eqref{eq:db-ratio}, gives \eqref{eq:db-toda}.
\end{proof}

\subsection{The KP limit of the Toda equation}\label{app:weak-toda-limit}

This subsection proves Theorem~\ref{thm:toda-kp}. Throughout, $h=1/T$, and $s$, $X$, $R$, $\mathsf S_h^\pm$ and $\delta_h$ are as in Section~\ref{sec:toda-kp-limit}. Let $F$ be a smooth function on $(0,\infty)\times\R^2$, and let $f(s,X,R)=F(t,x,r)$, where $t=s^{-1/2}$, $x=tX$ and $r=R-tX^2$. By the chain rule, $F$ satisfies \eqref{eq:toda-kp-log} if and only if
\begin{equation}\label{eq:db-cylindrical}
 \mathcal Lf=-\frac1{4s}\partial_X^2f-\frac1{2s^2}\bigl(\partial_R^2f\bigr)^2,\qquad\text{where}\quad\mathcal L=-\frac2{\sqrt s}\partial_s\partial_R+\frac1{2s^{3/2}}\partial_R+\frac1{12s^2}\partial_R^4.
\end{equation}

We now derive \eqref{eq:db-cylindrical} from the Toda equation \eqref{eq:db-scaled-toda} in the limit $h\to0$. The formal computation, which treats the solutions of \eqref{eq:db-scaled-toda} as smooth in $s$, appears in \cite[Section~2.5.1]{PNG-Toda} for the polynuclear growth model.

\begin{proposition}\label{prop:weak-toda}
For each small $h>0$, let $f_h$ be a real function on $h^3\mathbb N\times\R^2$, twice continuously differentiable in $(X,R)$, that satisfies \eqref{eq:db-scaled-toda} at every lattice point with $s\ge2h^3$. Let $f$ be smooth on $(0,\infty)\times\R^2$, and suppose that
\[
 f_h\to f,\qquad\partial_R^2f_h\to\partial_R^2f,\qquad\partial_X^2f_h\to\partial_X^2f
\]
uniformly on compact subsets of $(0,\infty)\times\R^2$ intersected with the lattice. Then $f$ satisfies \eqref{eq:db-cylindrical} on $(0,\infty)\times\R^2$.
\end{proposition}

\begin{proof}
We first show that, uniformly on compact subsets of the lattice,
\begin{equation}\label{eq:db-Lh-f}
 \delta_hf_h=\frac{h^2}s\partial_R^2f_h-h^4\Bigl(\frac1{4s}\partial_X^2f_h+\frac1{2s^2}\bigl(\partial_R^2f_h\bigr)^2\Bigr)+O(h^6).
\end{equation}
Indeed, \eqref{eq:db-scaled-toda} reads $\delta_hf_h=\log\bigl(1+\frac{h^2}s(\partial_R^2f_h-\frac{h^2}4\partial_X^2f_h)\bigr)$, and by the hypotheses $\partial_R^2f_h$ and $\partial_X^2f_h$ are bounded on compact subsets of the lattice, uniformly in $h$. So the expansion $\log(1+z)=z-z^2/2+O(z^3)$ gives \eqref{eq:db-Lh-f}. For a test function $\psi\in C_c^\infty\bigl((0,\infty)\times\R^2\bigr)$, Taylor's formula gives, uniformly in $(s,X,R)$,
\begin{equation}\label{eq:db-test-expansion}
 \delta_h\psi=\frac{h^2}s\partial_R^2\psi+h^4\mathcal L\psi+O(h^6).
\end{equation}
Here $\mathsf S_h^\pm\psi(s,X,R)=\psi(s\pm h^3,X,R+\rho_\pm)$, where
\[
 \rho_\pm=\frac{2\sqrt s-2\sqrt{s\pm h^3}}{h^2}=\mp\frac h{\sqrt s}+\frac{h^4}{4s^{3/2}}+O(h^7).
\]
The terms with $\pm h^3\partial_s\psi$ cancel and the odd powers of $\rho_\pm$ cancel to leading order, so
\[
 \delta_h\psi=(\rho_++\rho_-)\partial_R\psi+\tfrac12(\rho_+^2+\rho_-^2)\partial_R^2\psi+h^3(\rho_+-\rho_-)\partial_s\partial_R\psi+\tfrac1{24}(\rho_+^4+\rho_-^4)\partial_R^4\psi+O(h^6),
\]
which is \eqref{eq:db-test-expansion}. The supports of $\delta_h\psi$ lie in a fixed compact subset of $(0,\infty)\times\R^2$.

Next we show that
\begin{equation}\label{eq:db-adjoint}
 \langle\delta_hf_h,\psi\rangle_h=\langle f_h,\delta_h\psi\rangle_h,\qquad\text{where}\quad\langle g,\psi\rangle_h=h^3\sum_{N\ge1}\int_{\R^2}g\,\psi\diff X\diff R
\end{equation}
with the integrand evaluated at $s=Nh^3$, and $h$ is so small that $\psi$ vanishes for $s\le h^3$. Indeed, in the variables $(N,u,v)$ of Section~\ref{sec:toda-kp-limit} the operators $\mathsf S_h^\pm$ are the shifts $N\mapsto N\pm1$, and $\diff X\diff R=\frac h4\diff u\diff v$ for every $N$, so $\langle\mathsf S_h^\pm f_h,\psi\rangle_h=\langle f_h,\mathsf S_h^\mp\psi\rangle_h$. Integration by parts in $R$ gives the same identity with $\frac{h^2}s\partial_R^2$ in place of $\delta_h$. So, by \eqref{eq:db-Lh-f} and \eqref{eq:db-test-expansion},
\[
 \langle f_h,\mathcal L\psi\rangle_h=-\Bigl\langle\frac1{4s}\partial_X^2f_h+\frac1{2s^2}\bigl(\partial_R^2f_h\bigr)^2,\psi\Bigr\rangle_h+O(h^2),
\]
and, by the hypotheses, letting $h\to0$ gives
\[
 \int f\,\mathcal L\psi\diff s\diff X\diff R=-\int\Bigl(\frac1{4s}\partial_X^2f+\frac1{2s^2}\bigl(\partial_R^2f\bigr)^2\Bigr)\psi\diff s\diff X\diff R.
\]
The operator $\mathcal L$ is formally self-adjoint: its first two terms are $-\partial_s(s^{-1/2}\partial_R)-s^{-1/2}\partial_s\partial_R$, the sum of an operator and its adjoint, and $\partial_R^4$ is self-adjoint. Since $f$ is smooth, integration by parts gives $\int f\,\mathcal L\psi=\int(\mathcal Lf)\psi$. As $\psi$ is arbitrary, $f$ satisfies \eqref{eq:db-cylindrical}.
\end{proof}

The identity $\partial_r^{-1}\phi_V=-2\partial_rF$ in Theorem~\ref{thm:toda-kp} follows from the next lemma.

\begin{lemma}\label{lem:airy-tail}
Let $V\in C_c^\infty(\R^2)$ be real, and fix $t>0$ and $x\in\R$. Then $Z_V(t,x,r)\to1$ and $\partial_r^kZ_V(t,x,r)\to0$ for $k\ge1$ as $r\to+\infty$. In particular $\partial_r\log Z_V(t,x,r)\to0$.
\end{lemma}

\begin{proof}
Let $V$ vanish outside $[-Y,Y]\times(-\infty,U]$, and let $|\partial_u^jV|\le M_j$. By Lemma~\ref{lem:diff-under-E}, $\partial_r^kZ_V=\E\partial_r^ke^{-I_V}$, which is the expectation of $e^{-I_V}$ times a polynomial without constant term in the occupations $I_{\partial_u^jV}^{t,x}(r)$, $1\le j\le k$. By Lemma~\ref{lem:counts}(ii), with $c=t^{-1/3}(r-U)$,
\[
 \E\exp\bigl\{p\,I_{|\partial_u^jV|}^{t,x}(r)\bigr\}\le\exp\bigl\{\bigl(e^{2pM_jY}-1\bigr)\E N_c(0)\bigr\},
\]
and $\E N_c(0)\to0$ as $c\to\infty$. So, as $r\to+\infty$, the occupations $I_{\partial_u^jV}^{t,x}(r)$, $0\le j\le k$, tend to $0$ in every $L^p$, with exponential moments bounded uniformly in $r\ge0$. Hence $Z_V\to1$ and $\partial_r^kZ_V\to0$ for $k\ge1$.
\end{proof}

\begin{proof}[Proof of Theorem~\ref{thm:toda-kp}]
Let $f_h$ and $f$ be as in Section~\ref{sec:toda-kp-limit}. By Proposition~\ref{prop:bridge-limit}, $Z_V$ is finite, positive and smooth, so $f$ is smooth. The convergence in Proposition~\ref{prop:bridge-limit} is uniform for $t$ in compact sets, and at fixed $s$ the derivatives in $X$ and $R$ are combinations of derivatives in $x$ and $r$ with coefficients smooth in $(s,X,R)$. So $f_h$ and all its derivatives in $X$ and $R$ converge to those of $f$, uniformly on compact subsets of $(0,\infty)\times\R^2$ intersected with the lattice. By Proposition~\ref{prop:toda}, $f_h$ satisfies \eqref{eq:db-scaled-toda}. So Proposition~\ref{prop:weak-toda} gives \eqref{eq:db-cylindrical}, which is equivalent to \eqref{eq:toda-kp-log}.

By Lemma~\ref{lem:airy-tail}, $\partial_rF\to0$ as $r\to+\infty$, so $\partial_r^{-1}\phi_V=-2\partial_rF$. Then
\[
 \partial_t\phi_V-\frac12\phi_V\partial_r\phi_V+\frac1{12}\partial_r^3\phi_V+\frac14\partial_x^2\partial_r^{-1}\phi_V=-2\partial_r\Bigl(\partial_t\partial_rF+\frac14\partial_x^2F+\frac1{12}\partial_r^4F+\frac12\bigl(\partial_r^2F\bigr)^2\Bigr),
\]
which vanishes by \eqref{eq:toda-kp-log}.
\end{proof}

\section{Finite-band reduction and heat scattering}
\label{app:heat-scattering}

We prove Theorem~\ref{thm:lt-heat} by factoring the Green
operator into time-ordered and finite-band parts. Fix $q\in\R$,
$\kappa>0$, put $\widetilde V=W_q$, and choose a compact rectangle
$E\subset(a,b)\times\R$ containing its support. On $L^2(E)$,
write $G_{i\kappa}=\mathscr C_\kappa-\mathscr B_\kappa$,
where the kernel of $\mathscr C_\kappa$ vanishes unless the output $y$ exceeds the input $y$, and let
$M_{\widetilde V}$ denote multiplication by $\widetilde V$.

\textit{Volterra inverse and finite-band factorization.}
The variance-$2$ heat semigroup and the ordered time simplex imply
\begin{equation}\label{eq:appendix-Volterra}
 \|(M_{\widetilde V}\mathscr C_\kappa)^n\|
 \le e^{\kappa^2(b-a)}
      \frac{((b-a)\|\widetilde V\|_\infty)^n}{n!},
 \qquad n\ge0 .
\end{equation}
Thus $I+\lambda M_{\widetilde V}\mathscr C_\kappa$ is invertible
for every complex $\lambda$. Put
$S_\lambda=(I+\lambda M_{\widetilde V}\mathscr C_\kappa)^{-1}
\lambda M_{\widetilde V}$.

Factor $\mathscr B_\kappa=UV_0$, where the two Hilbert--Schmidt
maps between $L^2(-\kappa,\kappa)$ and $L^2(E)$ are
\begin{equation}\label{eq:appendix-UV}
\begin{aligned}
 (Uf)(y,u)&=\frac1{\sqrt{2\pi}}\int_{-\kappa}^{\kappa}
 e^{i(m-\kappa)u+(\kappa^2-m^2)y}f(m)\,dm,\\
 (V_0g)(m)&=\frac1{\sqrt{2\pi}}\int_E
 e^{-i(m-\kappa)u-(\kappa^2-m^2)y}g(y,u)\,dy\,du .
\end{aligned}
\end{equation}
The fundamental factorization is
\begin{equation}\label{eq:appendix-heat-factor}
 I+\lambda M_{\widetilde V}G_{i\kappa}
 =(I+\lambda M_{\widetilde V}\mathscr C_\kappa)
  (I-S_\lambda UV_0).
\end{equation}
In particular, $F_\lambda=V_0S_\lambda U$ is trace class.

The ordered Duhamel expansion of
\eqref{eq:tf-data} is
\begin{equation}\label{eq:tf-series}
\begin{aligned}
 T_V(k,\ell)=\frac1{2\pi}\sum_{n\ge1}(-1)^{n+1}
 \int_{a<y_1<\cdots<y_n<b}\!\diff\boldsymbol y
 \int_{\R^n}\!\diff\boldsymbol u\;&
 e^{-y_1k^2+y_n\ell^2+iku_1-i\ell u_n}\\[-2pt]
 &\times\prod_{j=1}^nV(y_j,u_j)
 \prod_{j=1}^{n-1}p_{y_{j+1}-y_j}(u_{j+1}-u_j).
\end{aligned}
\end{equation}
The ordered-simplex bound gives absolute, locally uniform
convergence on $\mathbb C^2$. Expanding the Volterra inverse
and comparing endpoint phases with \eqref{eq:tf-series} gives
\begin{equation}\label{eq:appendix-Ftranspose}
 F_\lambda(m,\ell)=\mathsf T_{q,\kappa;\lambda V}(\ell,m).
\end{equation}
The $e^{\pm i\kappa u}$ phases telescope along the heat
kernels, so this holds for every complex $\lambda$.

\textit{Invertibility.}
The factorization \eqref{eq:appendix-heat-factor} and the
Woodbury--Sylvester identity
\[
 (I-AB)^{-1}=I+A(I-BA)^{-1}B
\]
reduce invertibility of
$I+\lambda M_{\widetilde V}G_{i\kappa}$ to that of
$I-F_\lambda$, equivalently $I-\mathsf T_{q,\kappa;\lambda V}$.
The regular extension of the scattering coefficient in
\cite{ChenKP} does not by itself exclude homogeneous
Fredholm solutions; compact support allows us to do so.

\begin{lemma}[No exceptional finite-band points]
\label{lem:finite-band-invertible}
For every real $W\in C_c^\infty(\R^2)$, $\kappa>0$, and
$\lambda\in\R$, the finite-band operator $I-F_\lambda$ defined
above with $\widetilde V=W$ is invertible on
$L^2(-\kappa,\kappa)$.  More precisely,
\begin{equation}\label{eq:finite-band-positive-det}
 \det(I-F_\lambda)\ge
 \exp\left\{-\frac{\lambda\kappa}{\pi}
       \int_{\R^2}W(y,u)\diff y\diff u\right\}>0.
\end{equation}
\end{lemma}
\begin{proof}
It suffices to take $\lambda=1$, since $\lambda W$ is again a real
compactly supported smooth potential. Choose $a<b$ strictly outside
the $y$-support of $W$, let $H=-\partial_u^2$ on $L^2(\R)$, and set
$P_\kappa=\mathbf1_{\{H<\kappa^2\}}$.
Write $P_{a,b;W}$ for the Feynman--Kac propagator and
$D_W=e^{-(b-a)H}-P_{a,b;W}$. To keep the exponential multipliers
bounded, put $f_c(s)=e^{cs}\mathbf1_{\{0\le s<\kappa^2\}}$, $c\in\R$.
Then
\[
 K_{\kappa,W}=f_{-a}(H)D_Wf_b(H)
 =P_\kappa e^{-aH}D_We^{bH}P_\kappa,
\]
where the second expression uses bounded band multipliers.
Under the Fourier convention
$\widehat f(k)=(2\pi)^{-1/2}\int e^{-iku}f(u)\diff u$,
the frequency kernel of $K_{\kappa,W}$ is $T_W(-k,-\ell)$
by \eqref{eq:tf-data}, while \eqref{eq:appendix-Ftranspose}
gives $F_1(m,\ell)=T_W(\ell,m)$.
Reflection of both frequencies and transposition preserve the
Fredholm determinant, so
\begin{equation}\label{eq:finite-band-det-identification}
 \det(I-F_1)=\det(I-K_{\kappa,W}).
\end{equation}
The latter operator is trace class by the Duhamel factorization below.

Let $H_L=-\partial_u^2$ on $(-L,L)$ with Dirichlet boundary
conditions, normalized eigenfunctions $\varphi_{j,L}$, and
eigenvalues $\epsilon_{j,L}=(j\pi/(2L))^2$.
Write $N_L=\#\{j:\epsilon_{j,L}<\kappa^2\}$,
$Q_L=\mathbf1_{\{H_L<\kappa^2\}}$, and
$E_L=\sum_{j=1}^{N_L}\epsilon_{j,L}$; take $L$ large enough that
$N_L\ge1$. With $P_{a,b;W}^L$ the Dirichlet Feynman--Kac
propagator, set $D_W^L=e^{-(b-a)H_L}-P_{a,b;W}^L$,
$K_L=f_{-a}(H_L)D_W^Lf_b(H_L)$, and $\Delta_L=\det(I-K_L)$.
The finite-rank operator $K_L$ vanishes on $Q_L^\perp$, and
$I-K_L=Q_Le^{-aH_L}P_{a,b;W}^Le^{bH_L}Q_L$ on the range of $Q_L$.
Taking the determinant in the eigenfunction basis gives
\begin{equation}\label{eq:finite-band-finite-det}
 \Delta_L=e^{(b-a)E_L}
 \det[\langle\varphi_{i,L},P_{a,b;W}^L\varphi_{j,L}\rangle]_{i,j=1}^{N_L}.
\end{equation}

Consider the chamber
$\mathcal C_L=\{-L<u_1<\cdots<u_{N_L}<L\}$ and the Slater
determinant $h_L(\boldsymbol u)=\det[\varphi_{i,L}(u_j)]_{i,j=1}^{N_L}$.
Putting $\theta_j=\pi(u_j+L)/(2L)$, the identities
$\varphi_{i,L}(u)=L^{-1/2}\sin(i\theta)$ and
$\sin(i\theta)=\sin\theta\,U_{i-1}(\cos\theta)$ yield
\[
 h_L(\boldsymbol u)=C_{N_L,L}
 \prod_{j=1}^{N_L}\sin\theta_j
 \prod_{p<q}(\cos\theta_q-\cos\theta_p),\qquad C_{N_L,L}\ne0.
\]
Thus $h_L$ has constant nonzero sign on $\mathcal C_L$; choose
its sign positive. Orthonormality and antisymmetry give
$\int_{\mathcal C_L}h_L^2\diff\boldsymbol u=1$.
Since $h_L$ is a positive eigenfunction of the killed
$N_L$-particle Laplacian $-\sum_j\partial_{u_j}^2$ with
eigenvalue $E_L$, it is its ground state.

If $p_t^L$ denotes the Dirichlet heat kernel on $(-L,L)$,
the killed nonintersecting transition density is
$q_t^0(\boldsymbol u,\boldsymbol v)
=\det[p_t^L(u_i,v_j)]_{i,j=1}^{N_L}$. Its ground-state
transform is
\[
 q_t^h(\boldsymbol u,\boldsymbol v)
 =e^{tE_L}\frac{h_L(\boldsymbol v)}{h_L(\boldsymbol u)}
 q_t^0(\boldsymbol u,\boldsymbol v).
\]
This is a conservative Markov kernel with generator
$\sum_{j=1}^{N_L}\partial_{u_j}^2+
2\sum_{j=1}^{N_L}(\partial_{u_j}\log h_L)\partial_{u_j}$.
The Brownian variance is $2$, and the process is reversible
with invariant probability density $h_L^2$.

For the time-dependent potential $\sum_jW(y,u_j)$,
antisymmetrization gives the killed-chamber kernel
$q_{a,b}^W(\boldsymbol u,\boldsymbol v)
=\det[P_{a,b;W}^L(u_i,v_j)]_{i,j=1}^{N_L}$.
Indeed, both sides solve the same $N_L$-particle parabolic
equation with the same initial data and Dirichlet conditions
at collisions and at $\pm L$. Uniqueness applies also to real
sign-changing, time-dependent $W$.
Andr\'eief's identity applied to \eqref{eq:finite-band-finite-det},
followed by the ground-state transform, yields
\[
\begin{aligned}
 \Delta_L
 &=e^{(b-a)E_L}\int_{\mathcal C_L^2}
 h_L(\boldsymbol u)q_{a,b}^W(\boldsymbol u,\boldsymbol v)
 h_L(\boldsymbol v)\diff\boldsymbol u\diff\boldsymbol v\\
 &=\mathbb E_L^h\exp\left\{
 -\sum_{j=1}^{N_L}\int_a^b W(y,X_j(y))\diff y
 \right\},
\end{aligned}
\]
where $\mathbb E_L^h$ is expectation for the stationary
ground-state-transformed process. In particular $\Delta_L>0$.
Its one-particle intensity is
$\rho_L(u)=\sum_{j=1}^{N_L}|\varphi_{j,L}(u)|^2$,
so Jensen's inequality gives
\begin{equation}\label{eq:finite-band-Jensen}
 \Delta_L\ge\exp\left\{-\int_a^b\int_{-L}^L
 W(y,u)\rho_L(u)\diff u\diff y\right\}.
\end{equation}
From the explicit eigenfunctions,
\[
 \rho_L(u)=\frac{N_L}{2L}
 -\frac1{2L}\sum_{j=1}^{N_L}\cos\left(\frac{j\pi(u+L)}L\right).
\]
On every fixed compact interval the cosine sum is uniformly
bounded, while $N_L/(2L)\to\kappa/\pi$; hence
$\rho_L(u)\to\kappa/\pi$ uniformly on compact subsets of $\R$.

It remains to pass to the whole-line determinant.
Extend Dirichlet kernels by zero to $L^2(\R)$ and put
$W_y(u)=W(y,u)$. Duhamel's formula gives
\[
\begin{aligned}
 D_W&=\int_a^b e^{-(y-a)H}M_{W_y}P_{y,b;W}\diff y,\\
 D_W^L&=\int_a^b e^{-(y-a)H_L}M_{W_y}P_{y,b;W}^L\diff y.
\end{aligned}
\]
Factor $M_{W_y}=M_{\sqrt{|W_y|}}
M_{\operatorname{sgn}(W_y)\sqrt{|W_y|}}$ and set
\[
\begin{aligned}
 A(y)&=e^{-(y-a)H}M_{\sqrt{|W_y|}},&
 B(y)&=M_{\operatorname{sgn}(W_y)\sqrt{|W_y|}}P_{y,b;W},\\
 A_L(y)&=e^{-(y-a)H_L}M_{\sqrt{|W_y|}},&
 B_L(y)&=M_{\operatorname{sgn}(W_y)\sqrt{|W_y|}}P_{y,b;W}^L.
\end{aligned}
\]
Then $D_W=\int_a^b A(y)B(y)\diff y$ and
$D_W^L=\int_a^b A_L(y)B_L(y)\diff y$.
Since the variance-$2$ heat kernel satisfies
$\int_\R p_t(u)^2\diff u=(8\pi t)^{-1/2}$, Gaussian and
Feynman--Kac bounds give, uniformly in $L$,
\[
\begin{aligned}
 \|A_L(y)\|_{\mathrm{HS}}^2,\ \|A(y)\|_{\mathrm{HS}}^2
 &\le\frac{\|W_y\|_1}{\sqrt{8\pi(y-a)}},\\
 \|B_L(y)\|_{\mathrm{HS}}^2,\ \|B(y)\|_{\mathrm{HS}}^2
 &\le\frac{e^{2(b-y)\|W_-\|_\infty}\|W_y\|_1}
 {\sqrt{8\pi(b-y)}}.
\end{aligned}
\]
For each $a<y<b$, the Dirichlet heat kernels converge
pointwise to the whole-line kernels. The same holds for
the Feynman--Kac kernels, since a Brownian bridge stays
inside $(-L,L)$ with probability tending to $1$.
Gaussian domination therefore gives
$\|A_L(y)-A(y)\|_{\mathrm{HS}}\to0$ and
$\|B_L(y)-B(y)\|_{\mathrm{HS}}\to0$.
The trace-ideal product estimate
\[
\begin{aligned}
 \|A_LB_L-AB\|_1
 &\le\|A_L-A\|_{\mathrm{HS}}\|B_L\|_{\mathrm{HS}}\\
 &\quad+\|A\|_{\mathrm{HS}}\|B_L-B\|_{\mathrm{HS}}
\end{aligned}
\]
has the integrable majorant
$\frac{2e^{(b-y)\|W_-\|_\infty}}{\sqrt{8\pi}}\,
\frac{\|W_y\|_1}{[(y-a)(b-y)]^{1/4}}$.
Dominated convergence for the trace-class-valued integrals yields
\begin{equation}\label{eq:finite-band-defect-trace-limit}
 \|D_W^L-D_W\|_1\longrightarrow0.
\end{equation}

Finally, Dirichlet exhaustion gives strong-resolvent convergence
$H_L\to H$ under the natural embeddings. This can be formalized
by adjoining the constant $\kappa^2+1$ to $H_L$ on
$L^2((-L,L)^c)$. Since $H$ has no eigenvalues at the discontinuities
$0,\kappa^2$ of $f_c$, the spectral calculus gives
$f_c(H_L)\to f_c(H)$ strongly, with
$\sup_L\|f_c(H_L)\|\le e^{\max(c,0)\kappa^2}$.
These multipliers are self-adjoint, so their adjoints also
converge strongly. Write $R_L=f_{-a}(H_L)$,
$S_L=f_b(H_L)$, $R=f_{-a}(H)$, and $S=f_b(H)$.
Then
\[
 K_L-K_{\kappa,W}=R_L(D_W^L-D_W)S_L
 +(R_L-R)D_WS_L+RD_W(S_L-S).
\]
The first term tends to zero in trace norm by
\eqref{eq:finite-band-defect-trace-limit}; the other two do so
because multiplication of a fixed trace-class operator by
uniformly bounded strongly convergent operators (or their
adjoints on the right) preserves trace-norm convergence.
Thus $\|K_L-K_{\kappa,W}\|_1\to0$, and continuity of the
Fredholm determinant gives
\[
 \det(I-F_1)=\det(I-K_{\kappa,W})=\lim_{L\to\infty}\Delta_L.
\]
The compact support of $W$ now allows passage to the limit in
\eqref{eq:finite-band-Jensen}, proving
\eqref{eq:finite-band-positive-det} for $\lambda=1$.
Replacing $W$ by $\lambda W$ gives the bound for every real
coupling. Since $F_\lambda$ is trace class, positivity of
its determinant implies invertibility of $I-F_\lambda$.
\end{proof}

\textit{Endpoint amplitude and normalization.}
The lemma makes both resolvents below well-defined for every
real coupling. Solving through the two factors gives
\[
 h_\lambda=(I-S_\lambda UV_0)^{-1}S_\lambda1,\qquad
 h_\lambda=(I+\lambda M_{\widetilde V}G_{i\kappa})^{-1}
 \lambda\widetilde V .
\]
At $m=\kappa$ the $U$-kernel equals $1/\sqrt{2\pi}$,
so $S_\lambda1=\sqrt{2\pi}S_\lambda U(\cdot,\kappa)$.
Consequently,
$V_0h_\lambda=\sqrt{2\pi}(I-F_\lambda)^{-1}
F_\lambda(\cdot,\kappa)$. At $m=-\kappa$ the $V_0$-kernel is
$e^{2i\kappa u}/\sqrt{2\pi}$. Evaluating the preceding identity
there gives
\begin{equation}\label{eq:appendix-heat-coeff}
 \int_E e^{2i\kappa u}h_\lambda(y,u)\,dy\,du
 =2\pi[F_\lambda(I-F_\lambda)^{-1}](-\kappa,\kappa)
 =\mathcal A_q(\kappa;\lambda V).
\end{equation}
The smooth frequency kernels of $U$ and $V_0$ justify
endpoint evaluation. The transposition in
\eqref{eq:appendix-Ftranspose} fixes the order
$(\kappa,-\kappa)$ in \eqref{eq:tf-amp}; enlarging $E$
does not change the coefficient.

Finally, for $z=q+i\kappa$, shear $u=r+2qy$.
The normalized direct heat equation of \cite[(4.1)]{ChenKP}
becomes
\[
 (\partial_y-\partial_u^2-2i\kappa\partial_u)\chi
 =-W_q\chi,\qquad W_q(y,u)=V(y,u-2qy),
\]
whose Green kernel is \eqref{eq:tf-Green}; compare
\cite[Lemma~4.1]{ChenKP}. The corresponding initial KP field
in the conventions of \cite{ChenKP} is $u_0=-V$, so its heat
potential is $-u_0=V$. Its spectral frequency is
$\zeta(z)=(-2\kappa,-4q\kappa)$ in $(r,y)$ coordinates.
At $\lambda=1$, let $h=h_1$. Since $h\in L^2(E)$ and is supported in
the compact rectangle $E$, the Green representation
$\chi=1-G_{i\kappa}*h$ gives a bounded solution tending to $1$ at
spatial infinity. Uniqueness follows from the invertibility established
above.
Consequently the coefficient \cite[(4.12)]{ChenKP} is
\[
\begin{aligned}
 \mathsf{a}(q+i\kappa)
 &= -\frac{i}{2\pi}\int_{\R^2}
 u_0(r,y)\chi(z;r,y)
 e^{-i(r,y)\cdot\zeta(z)}\,dr\,dy\\
 &=\frac{i}{2\pi}\int_{\R^2}
 W_q(y,u)\chi(y,u)e^{2i\kappa u}\,dy\,du
 =\frac{i}{2\pi}\mathcal A_q(\kappa;V),
\end{aligned}
\]
by \eqref{eq:appendix-heat-coeff} at $\lambda=1$.

\section{Darboux transformations and Brownian interlacing}
\label{app:soliton-proof}

A positive solution of a backward heat equation with potential
defines a Darboux transformation of that potential. At the level
of nonintersecting paths, the transformation is implemented by an
interlacing link that changes the particle number by one. We
develop these links for a general base potential, compose them,
and examine their confluent endpoints. The resulting Gaussian
bordering law produces a finite-dimensional formula for Wronskian
potentials. Its large-$N$ asymptotics then yield the soliton
Airy expectation, independently of the KP equation.

The method has two separate limiting steps. The limit in which the bridge endpoints coincide must be taken
with the Vandermonde normalizations intact. Later, the lower and upper free endpoint
groups decouple only asymptotically; their finite-$N$ coupling
must be retained until the Laplace estimate is proved. The
Airy edge convergence itself is established in
Appendix~\ref{app:s8-bridge}.

\subsection{The one-particle Darboux transform}
\label{app:onepath}

The Darboux identity, tilted transition density, and
one-particle formula are \eqref{eq:darboux},
\eqref{eq:tilted-kernel}, and Lemma~\ref{lem:onepath-Q}
in the main text. We now prove the lemma.

\begin{proof}
The first claim is the Feynman--Kac formula for $\psi$, which solves the
backward equation for $Q$ and is at most $Ce^{C|u|}$.
Let $h$ be a smooth function on $\R$ with compact support.
Integrating \eqref{eq:Q-onepath} against $h(v)$, we see that
\eqref{eq:Q-onepath} is equivalent to
\begin{equation}\label{eq:Q-onepath-test}
 \int_\R P_{y,b;Q^\psi}(u,v)h(v)\diff v
 =-\psi(y,u)\partial_u\int_\R\frac{h(v)}{\psi(b,v)}
 \int_{-\infty}^{v}p^{\psi,Q}_{y,b}(u,w)\diff w\diff v
\end{equation}
for every such $h$. By the Feynman--Kac formula, the left side of
\eqref{eq:Q-onepath-test} is the only solution of the backward
equation for $Q^\psi$ that has terminal value $h$ at $y=b$ and
absolute value at most a constant multiple of $e^{C|u|}$. We show
that the right side of \eqref{eq:Q-onepath-test} is a solution with
these two properties.

First, we write the right side of \eqref{eq:Q-onepath-test} in terms
of a solution of the backward equation for $Q$. By
\eqref{eq:tilted-kernel} and Fubini's theorem, the double integral in
\eqref{eq:Q-onepath-test} is equal to $H(y,u)/\psi(y,u)$, where
\[
 H(y,u)=\int_\R P_{y,b;Q}(u,w)\psi(b,w)
 \int_w^\infty\frac{h(v)}{\psi(b,v)}\diff v\diff w.
\]
The right side of \eqref{eq:Q-onepath-test} is therefore
$-\psi\partial_u(H/\psi)=-(\partial_u-q)H$, where $q=\partial_u\log\psi$.
Next, $H$ solves the backward equation for $Q$, because the kernel
$P_{y,b;Q}(u,w)$ solves this equation in the variables $(y,u)$. By
the Darboux identity \eqref{eq:darboux}, with $L_Q=\partial_y+\partial_u^2-Q$,
$L_{Q^\psi}(\partial_u-q)H=(\partial_u-q)L_QH=0$, so the right side
of \eqref{eq:Q-onepath-test} solves the backward equation for
$Q^\psi$. At $y=b$ we have
$H(b,u)=\psi(b,u)\int_u^\infty h(v)/\psi(b,v)\diff v$,
so the right side of \eqref{eq:Q-onepath-test} is equal to $h(u)$
there.
Finally, the tilted kernel is the transition density of the diffusion
$dY_s=\sqrt2\,dB_s+2q(s,Y_s)\,ds$. Put
$g(w)=\int_w^\infty h(v)/\psi(b,v)\,dv$. Then
$H(y,u)/\psi(y,u)=\mathbb E_{y,u}g(Y_b)$.
Since $\partial_uq=(Q-Q^\psi)/2$ is bounded, the drift is globally
Lipschitz in $u$. Differentiating the diffusion flow in its starting
point, and using $\|g'\|_\infty\le\|h/\psi(b,\cdot)\|_\infty$, gives
\[
 \bigl|\partial_u(H/\psi)(y,u)\bigr|
 \le \|h/\psi(b,\cdot)\|_\infty
 e^{2\|\partial_uq\|_\infty(b-y)}.
\]
Thus the right side of \eqref{eq:Q-onepath-test} is at most a constant
multiple of $\psi(y,u)\le Ce^{C|u|}$ in absolute value.
\end{proof}

We now prove the one-more-path identity \eqref{eq:one-more-path} for
the potentials $Q$ and $Q^\psi$ of this subsection, for every $N\ge1$
and all increasing $u\in\R^{N+1}$ and $v\in\R^N$. We follow the proof
of Assiotis, O'Connell and Warren~\cite[Proposition~2.16]{AOW}, with
Lemma~\ref{lem:onepath-Q} in place of their one-path
duality~\cite[Lemma~2.2]{AOW}, whose proof uses time-homogeneity.

\begin{proof}[Proof of \eqref{eq:one-more-path}]
For $j=1,\dots,N$ let
$F_j(w)=\int_{-\infty}^{v_j}p^{\psi,Q}_{a,b}(w,z)\diff z$.
By Lemma~\ref{lem:onepath-Q} with $y=a$, the entries of the
determinant on the left of \eqref{eq:one-more-path} satisfy
$P_{a,b;Q^\psi}(w_i,v_j)/\psi(a,w_i)=-F_j'(w_i)/\psi(b,v_j)$.
The region $u\prec w$ is the product of the intervals
$(u_i,u_{i+1})$, and a determinant is linear in each row. Integrating
row $i$ over $w_i\in(u_i,u_{i+1})$ therefore shows that the left side
of \eqref{eq:one-more-path} is equal to
\begin{equation}\label{eq:one-more-path-rows}
 \frac1{\prod_{j=1}^N\psi(b,v_j)}
 \det[F_j(u_i)-F_j(u_{i+1})]_{i,j\le N}.
\end{equation}
It remains to write the determinant in \eqref{eq:one-more-path-rows}
as an integral of a determinant of size $N+1$. Let $\mathcal M$ be the
$(N+1)\times(N+1)$ matrix whose $i$th row is
$(1,F_1(u_i),\dots,F_N(u_i))$. Subtracting from each of the first $N$
rows of $\mathcal M$ the row below it, and expanding along the first
column, gives
\begin{equation}\label{eq:one-more-path-sign}
 \det[F_j(u_i)-F_j(u_{i+1})]_{i,j\le N}=(-1)^N\det\mathcal M.
\end{equation}
We now apply column operations to $\mathcal M$ in three steps:
subtract the last column from the first; then, for $j=N,N-1,\dots,2$,
subtract the column of $F_{j-1}$ from the column of $F_j$; then move
the first column to the end. The first two steps leave the determinant
unchanged, and the third multiplies it by $(-1)^N$, which cancels the
factor $(-1)^N$ in \eqref{eq:one-more-path-sign}. Since
$p^{\psi,Q}_{a,b}(u_i,\cdot)$ is a probability density, row $i$ of the
resulting matrix consists of the probabilities that it assigns to the
intervals $(-\infty,v_1)$, $(v_1,v_2)$, \dots, $(v_N,\infty)$. The
region $z\prec v$ is the product of these intervals, so linearity in
each column shows that the determinant of the resulting matrix is
$\int_{z\prec v}\det[p^{\psi,Q}_{a,b}(u_i,z_j)]_{i,j\le N+1}\diff z$.
By \eqref{eq:one-more-path-rows} and \eqref{eq:one-more-path-sign}, the
left side of \eqref{eq:one-more-path} is this integral divided by
$\prod_{j=1}^N\psi(b,v_j)$, and substituting \eqref{eq:tilted-kernel}
for the entries gives the right side of \eqref{eq:one-more-path}.
\end{proof}

\subsection{Interlacing links and their composition}
\label{app:interlacing}

This subsection proves Proposition~\ref{prop:n-endpoints}.

\textit{Normalized links.}
First we show that the one-more-path identity
\eqref{eq:one-more-path} can be written as an intertwining relation,
\eqref{eq:link-intertwining}. Applied $n$ times, it gives the identity
\eqref{eq:composite-intertwining} between $N$ paths in the potential
$V_n$ and $N+n$ free paths.

For $m\ge1$, let $\mathbb W_m=\{u_1<\cdots<u_m\}$, and
write $\Delta_m(u)=\prod_{i<j}(u_j-u_i)$ for the Vandermonde
determinant. For a positive function
$\psi$, define a kernel from $\mathbb W_{m+1}$ to $\mathbb W_m$ by
\begin{equation}\label{eq:normalized-link}
 \Lambda^{\psi}_{m,y}(u,\diff w)
 =m!\,\mathbf1_{\{u\prec w\}}
 \frac{\Delta_m(w)}{\Delta_{m+1}(u)}
 \frac{\prod_{i=1}^{m+1}\psi(y,u_i)}
 {\prod_{j=1}^{m}\psi(y,w_j)}\diff w.
\end{equation}
For $\psi=1$, this is a probability kernel, since
$\int_{u\prec w}\Delta_m(w)\diff w=\Delta_{m+1}(u)/m!$.
These are precisely the conditional eigenvalue kernels of consecutive
principal minors of a unitarily invariant Hermitian matrix.

For a potential $Q$, let
$K^Q_m(u,\diff v)=\det[P_{a,b;Q}(u_i,v_j)]_{i,j\le m}\diff v$ be its
Karlin--McGregor determinant, and use the Vandermonde-transformed kernel
\[
 \mathsf P_m^Q(a,u;b,\diff v)
 =\frac{\Delta_m(v)}{\Delta_m(u)}K^Q_m(u,\diff v).
\]
For $Q=0$, $\mathsf P_m^0$ is the Dyson Brownian transition kernel;
for nonzero $Q$, it is the corresponding positive Feynman--Kac
kernel. Multiplying
\eqref{eq:one-more-path} by the Vandermonde factors gives the exact
operator identity
\begin{equation}\label{eq:link-intertwining}
 \Lambda^{\psi}_{m,a}\mathsf P_m^{Q^\psi}
 =\mathsf P_{m+1}^{Q}\Lambda^{\psi}_{m,b}.
\end{equation}

\begin{lemma}[Positive Wronskians and coefficient minors]
\label{lem:positive-prefix-minors}
Let $p_1<\cdots<p_m$ and let $f_i,W_k$ be as in
\eqref{eq:Wronskians}, with real coefficients. Suppose that
$W_j(x,r)>0$ for every $(x,r)\in\R^2$ and $1\le j\le n$.
Then every ordered $k\times k$ minor of the first $k$ rows of
$(c_{ij})$ is nonnegative, for $1\le k\le n$, and at least
one such minor is positive for each $k$. In particular, each
$W_k$ is a positive finite exponential sum. On bounded $x$-intervals,
the potentials $V_k=-2\partial_r^2\log W_k$ are bounded,
$\partial_r\log\psi_k$ and its $r$-derivative are bounded, and
$\psi_k=W_k/W_{k-1}$ has at most exponential growth in $|r|$.
\end{lemma}

\begin{proof}
Fix $x$. Positivity of all leading Wronskians implies that for
$z_1<\cdots<z_k$,
\begin{equation}\label{eq:positive-collocation}
 \det[f_i(x,z_b)]_{i,b=1}^k>0.
\end{equation}
Here is a direct induction. The assertion for $k=1$ is $f_1>0$.
For $k\ge2$ put $h_i=\partial_r(f_{i+1}/f_1)$,
$1\le i<k$, with $x$ fixed. The elementary Wronskian identity
\[
 W(h_1,\ldots,h_j)
 =\frac{W_{j+1}(x,\cdot)}{f_1(x,\cdot)^{j+1}}>0,
 \qquad 1\le j<k,
\]
allows induction for the $h_i$. Divide column $b$ of the
collocation determinant by $f_1(x,z_b)$ and subtract each
preceding column from the next, working from right to left.
Expanding along the first row gives
\[
 \frac{\det[f_i(x,z_b)]_{i,b=1}^k}
 {\prod_{b=1}^k f_1(x,z_b)}
 =
 \int_{z_1}^{z_2}\!\!\cdots
 \int_{z_{k-1}}^{z_k}
 \det[h_i(u_b)]_{i,b=1}^{k-1}\,
 du_{k-1}\cdots du_1>0.
\]
The integration variables are strictly ordered, so positivity
follows by the inductive hypothesis.

Now fix an ordered $k$-subset $I=\{i_1<\cdots<i_k\}$ and
take $x=L$, $z_b=2p_{i_b}L$, where $L\to\infty$.
Dividing column $b$ of the positive determinant
\eqref{eq:positive-collocation} by $e^{p_{i_b}^2L}$,
its entry in row $a$ becomes
\[
 \sum_{j=1}^m c_{aj}e^{-L(p_j-p_{i_b})^2}
 \longrightarrow c_{a,i_b}.
\]
Thus $\det[c_{a,i_b}]_{a,b=1}^k\ge0$.
By Cauchy--Binet,
\[
 W_k(x,r)=
 \sum_{\substack{I\subset\{1,\ldots,m\}\\|I|=k}}
 \det[c_{a,i_b}]_{a,b=1}^k\,\Delta(p_I)\,
 \exp\left\{\sum_{i\in I}(p_i r-p_i^2x)\right\}.
\]
At least one coefficient is positive, since $W_k>0$.
Its logarithmic first derivative in $r$ is a weighted average
of the finitely many slopes $\sum_{i\in I}p_i$,
and its second derivative is their variance. Both are
uniformly bounded. The same holds for the differences
of these derivatives defining $\log\psi_k$.
Finally, $\psi_k(x,0)$ is bounded for $x$ in a fixed
compact interval, and the bound on $\partial_r\log\psi_k$
gives $\psi_k(x,r)\le Ce^{C|r|}$ there.
\end{proof}

Applying \eqref{eq:link-intertwining} $n$ times, along the chain
$\psi_k=W_k/W_{k-1}$, $V_k=-2\partial_u^2\log W_k$ of
\eqref{eq:psi-chain}, gives
\begin{equation}\label{eq:composite-intertwining}
 \mathcal L_{n,a}^{(N)}\mathsf P_N^{V_n}
 =\mathsf P_{N+n}^{0}\mathcal L_{n,b}^{(N)},
\end{equation}
where $\mathcal L_{n,y}^{(N)}$ is the kernel from $\mathbb W_{N+n}$ to
$\mathbb W_N$
\begin{equation}\label{eq:composite-link}
 \mathcal L_{n,y}^{(N)}
 =\Lambda^{\psi_1}_{N+n-1,y}
 \Lambda^{\psi_2}_{N+n-2,y}\cdots
 \Lambda^{\psi_n}_{N,y}.
\end{equation}
Indeed, \eqref{eq:link-intertwining} applies with $Q=V_{k-1}$,
$\psi=\psi_k$ and $m=N+n-k$, for $k=n,n-1,\dots,1$, and then
$Q^\psi=V_k$. The assumptions of Appendix~\ref{app:onepath} hold for
these pairs by Lemma~\ref{lem:positive-prefix-minors}: the potentials
are bounded and the positive functions $\psi_k$ have the required
growth bound.

For a family $(f_1,\ldots,f_n)$ with positive intermediate Wronskians,
consider a triangular interlacing array with top row $z\in\mathbb W_n$:
rows $z^{(k)}\in\mathbb W_{n-k}$, $1\le k\le n-1$, with
$z\prec z^{(1)}\prec\cdots\prec z^{(n-1)}$.
The set of these arrays is a subset of $\R^{n(n-1)/2}$ of volume
$\Delta_n(z)/C_n$.
Give the top row the weight $\prod_j\psi_1(z_j)$, and give the row
$z^{(k)}$ the weight
$\prod_j\tfrac{\psi_{k+1}}{\psi_k}(z^{(k)}_j)$.
Let $\mathcal I_n(z)$ denote the integral of the product of these
weights over the rows $z^{(1)},\ldots,z^{(n-1)}$. For instance,
\begin{equation}\label{eq:weighted-vol-2}
 \mathcal I_2(z_1,z_2)
 =f_1(z_1)f_1(z_2)
 \int_{z_1}^{z_2}\frac{\psi_2(w)}{\psi_1(w)}\diff w
 =\det[f_i(z_j)]_{i,j\le2}.
\end{equation}

\begin{lemma}[Divided-Wronskian volume identity]\label{lem:triangle}
For every $n\ge1$,
\begin{equation}\label{eq:triangle-identity}
 \mathcal I_n(z)=\det[f_i(z_j)]_{i,j=1}^n.
\end{equation}
Consequently $\mathcal I_n(z)$ divided by the volume of the set of
arrays is the divided Wronskian $\mathcal W(z)$ of
\eqref{eq:divided-W}.
\end{lemma}
\begin{proof}
The case $n=1$ is $f_1(z_1)$, and $n=2$ is
\eqref{eq:weighted-vol-2}. For the induction step, define
$g_i=(f_{i+1}/f_1)'$ for $1\le i\le n-1$. We will show that
\begin{equation}\label{eq:triangle-step}
 \mathcal I_n(z)=\left(\prod_{j=1}^nf_1(z_j)\right)
 \int_{z\prec w}\det[g_i(w_j)]_{i,j\le n-1}\diff w.
\end{equation}
The lemma then follows by integrating. Indeed,
after factoring $f_1(z_j)$ from each column of the determinant of
$f_i(z_j)$, subtract successive columns. The fundamental theorem of
calculus and determinant multilinearity turn it into the right side
of \eqref{eq:triangle-step}.

It remains to prove \eqref{eq:triangle-step}. The factor
$\prod_jf_1(z_j)$ is the weight of the top row, since \mbox{$\psi_1=f_1$}. By
the induction hypothesis, it suffices to show that the weights of the
other rows are those of an array for the $g_i$ with top row
$z^{(1)}$. This follows from the identity
\mbox{$\operatorname{Wr}(g_1,\ldots,g_k)=W_{k+1}/f_1^{k+1}$}, which is proved
either by column operations on a Wronskian matrix or by induction
using \mbox{$\operatorname{Wr}(f_1,f_i)=f_1^2(f_i/f_1)'$}. Indeed, the
Wronskians of the $g_i$ are then positive, and the $k$th divided by
the $(k-1)$st is $\tfrac{\psi_{k+1}}{\psi_1}$. So in that array
$z^{(1)}$ has the weight $\tfrac{\psi_2}{\psi_1}$, and $z^{(k)}$ for
$k\ge2$ has the weight $\tfrac{\psi_{k+1}}{\psi_k}$, since $\psi_1$
cancels.

The unweighted volume follows from the same recursion with
$f_i(u)=u^{i-1}/(i-1)!$, for which the weighted integrand is constant,
and equals $\Delta(z)/C_n$.
\end{proof}

\textit{Coinciding endpoints and Gaussian bordering.}
First we show that
\begin{equation}\label{eq:confluence-start}
 \mathcal L_{n,a}^{(N)}(R^{N+n},\diff w)
 =\left(\prod_{k=1}^n\psi_k(a,R)\right)\delta_{R^N}(\diff w)
 =W_n(a,R)\delta_{R^N}(\diff w),
\end{equation}
where $R^m=(R,\ldots,R)\in\R^m$ and the equality is understood as a
continuous extension from $\mathbb W_{N+n}$. To see this, note that by
\eqref{eq:composite-link} and the indicator in
\eqref{eq:normalized-link}, every variable integrated in
$\mathcal L_{n,a}^{(N)}(u,\diff w)$ has its points in
$[u_1,u_{N+n}]$, and this interval shrinks to $R$ as $u$ tends to
$R^{N+n}$. The $k$th kernel in \eqref{eq:composite-link} therefore
tends to a point mass, and by \eqref{eq:normalized-link} its total
mass tends to $\psi_k(a,R)$.

Let $\mathcal L_{n,b}^{0,(N)}$ be the kernel
\eqref{eq:composite-link} with every $\psi_k$ replaced by $1$. For
$v\in\mathbb W_N$, let $z(v)$ denote the eigenvalues of the random
matrix
\begin{equation}\label{eq:matrix-rep}
 \begin{pmatrix}\diag(v)&\xi\\\xi^*&RI_n+G\end{pmatrix},
\end{equation}
with $\xi$ and $G$ as in Proposition~\ref{prop:n-endpoints}. Next we
show that, for almost every $v$,
\begin{equation}\label{eq:fixed-block}
 W_n(a,R)\,
 \frac{\mathsf P_N^{V_n}(a,R^N;b,\diff v)}
 {\mathsf P_N^{0}(a,R^N;b,\diff v)}
 =\E_{\xi,G}\biggl[
 \frac{\mathcal L_{n,b}^{(N)}(z(v),\diff v)}
 {\mathcal L_{n,b}^{0,(N)}(z(v),\diff v)}\biggr].
\end{equation}
Indeed, by \eqref{eq:confluence-start}, the identity
\eqref{eq:composite-intertwining} at $u=R^{N+n}$ and the same identity
with every $\psi_k=1$ are
\begin{align}
 W_n(a,R)\,\mathsf P_N^{V_n}(a,R^N;b,\diff v)
 &=\int\mathsf P_{N+n}^{0}(a,R^{N+n};b,\diff z)\,
 \mathcal L_{n,b}^{(N)}(z,\diff v),\label{eq:G9-at-R}\\
 \mathsf P_N^{0}(a,R^N;b,\diff v)
 &=\int\mathsf P_{N+n}^{0}(a,R^{N+n};b,\diff z)\,
 \mathcal L_{n,b}^{0,(N)}(z,\diff v).\label{eq:G9-at-R-free}
\end{align}
Let $H$ be an $(N+n)\times(N+n)$ GUE matrix normalized as $G$. The
eigenvalues $z$ of $RI_{N+n}+H$ have law
$\mathsf P_{N+n}^{0}(a,R^{N+n};b,\diff z)$ and, by the GUE-corners
law~\cite{Baryshnikov}, given $z$ the eigenvalues $v$ of its
$N\times N$ leading principal block have law
$\mathcal L_{n,b}^{0,(N)}(z,\diff v)$. So the right-hand side of
\eqref{eq:G9-at-R-free} is the law of $v$, and the quotient of
\eqref{eq:G9-at-R} by \eqref{eq:G9-at-R-free} is the conditional
expectation, given $v$, of the ratio of
$\mathcal L_{n,b}^{(N)}(z,\diff v)$ to
$\mathcal L_{n,b}^{0,(N)}(z,\diff v)$. Given $v$, the eigenvalues $z$
have the law of $z(v)$: the leading block is unitarily conjugate to
$\diag(v)$, and the conjugation does not change the law of the
independent off-diagonal and bottom-right Gaussian blocks. This proves
\eqref{eq:fixed-block}.

Finally we let $v\to R^N$ in \eqref{eq:fixed-block}, through points at
which it holds. The left-hand side tends to $W_n(a,R)Z_N(a,b;R;V_n)$ by
the Karlin--McGregor formula~\cite{KM}, since bridges with common
endpoints are limits of bridges with distinct endpoints
\cite[Definition~2.13]{CH}. On the right-hand side, the
matrix \eqref{eq:matrix-rep} tends to $RI_{N+n}+M$, with $M$ as in
\eqref{eq:border}, and we show that, almost surely,
\begin{equation}\label{eq:final-link-ratio}
 \frac{\mathcal L_{n,b}^{(N)}(z(v),\diff v)}
 {\mathcal L_{n,b}^{0,(N)}(z(v),\diff v)}
 \longrightarrow
 \frac{\mathcal W_b(R+\lambda^-)\mathcal W_b(R+\lambda^+)}{W_n(b,R)}.
\end{equation}
The limit may be taken inside the expectation, because for $v$ near
$R^N$ the ratio in \eqref{eq:final-link-ratio} is at most a constant
multiple of $e^{C\|M\|}$, where $\|M\|$ is the operator norm, and the
entries of $M$ are Gaussian. So \eqref{eq:fixed-block} and
\eqref{eq:final-link-ratio} give \eqref{eq:n-endpoints}.

It remains to prove \eqref{eq:final-link-ratio}. By
\eqref{eq:normalized-link} and \eqref{eq:composite-link}, the ratio in
\eqref{eq:final-link-ratio} is the integral of a product of weights
over the arrays $z(v)\prec z^{(1)}\prec\cdots\prec z^{(n-1)}\prec v$
with $z^{(k)}\in\mathbb W_{N+n-k}$, divided by the volume of the set of
these arrays; the Vandermonde factors cancel. The rows
$z(v),z^{(1)},\ldots,z^{(n-1)}$ have the weights of
Lemma~\ref{lem:triangle} at time $b$, and the row $v$ has the weight
$\prod_j\psi_n(b,v_j)^{-1}$.

Almost surely, $z_i(v)\to R+\lambda^-_i$ and
$z_{N+i}(v)\to R+\lambda^+_i$ for $1\le i\le n$, so $z_n(v)<v_1$ and
$v_N<z_{N+1}(v)$ when $v$ is close to $R^N$. For such $v$, the
interlacing conditions hold if and only if the first $n-k$ entries of
the rows $z^{(k)}$ form a triangular interlacing array with top row
$(z_1(v),\ldots,z_n(v))$, the last $n-k$ entries form one with top row
$(z_{N+1}(v),\ldots,z_{N+n}(v))$, and the remaining entries interlace,
which places them in $[v_1,v_N]$. Hence the ratio in
\eqref{eq:final-link-ratio} is a product of three ratios. By
Lemma~\ref{lem:triangle}, the ratios for the two triangular arrays are
the divided Wronskians $\mathcal W_b$ of their top rows, which converge
to $\mathcal W_b(R+\lambda^-)$ and $\mathcal W_b(R+\lambda^+)$. In the
third ratio every weight is evaluated in $[v_1,v_N]$, so it converges
to the product of the weights at $R$,
\begin{equation}\label{eq:pinned-factor}
 \psi_1(b,R)^{N-n}\prod_{k=1}^{n-1}
 \biggl(\frac{\psi_{k+1}(b,R)}{\psi_k(b,R)}\biggr)^{N-n+k}
 \psi_n(b,R)^{-N}
 =\prod_{k=1}^n\psi_k(b,R)^{-1}=W_n(b,R)^{-1}.
\end{equation}
This proves \eqref{eq:final-link-ratio}.

\subsection{Gaussian endpoint laws and Laplace asymptotics}
\label{app:asymptotics}

The bordering matrix \eqref{eq:border} has two groups of
nonzero eigenvalues. We derive their joint density, evaluate
the exponential phases within each group, and control the
interaction between the groups. The cubic term in the
one-dimensional Laplace expansion is the source of the
KP time evolution.

\textit{The joint density.}
The GUE-corners law, followed by confluence at the repeated zero
eigenvalues, gives the following density for the nonzero eigenvalues of
\eqref{eq:border}:
\begin{equation}\label{eq:eigen-density}
 \mathrm{const}_{N,n,s}\,
 e^{-\sum_{j=1}^{2n}\lambda_j^2/(4s)}
 \prod_{j=1}^{2n}|\lambda_j|^{N-n}
 \Delta(\lambda)\Delta(\lambda^-)\Delta(\lambda^+),
\end{equation}
where $s=b-a$ and $\lambda=(\lambda^-,\lambda^+)$ is ordered, with $n$
negative and $n$ positive entries. To see the repeated-zero factors, add the Gaussian rows and columns
one at a time. At the $k$th border the remaining zero block has
multiplicity $N-k$. In the one-border secular equation, the residues
at simple eigenvalues have exponential laws, while the residue at zero
has a Gamma law of shape $N-k$. The change of variables from residues
to the new roots contributes a Vandermonde ratio and the zero-residue
factor
\[
 \frac{\prod_{\lambda\in\lambda^{(k+1)}}|\lambda|^{N-k-1}}
 {\prod_{\lambda\in\lambda^{(k)}}|\lambda|^{N-k}}.
\]
These factors telescope over the successive borders. Integrating the
negative and positive interlacing triangles contributes their respective
Vandermonde volumes, yielding \eqref{eq:eigen-density} up to its
normalizing constant. For $n=1$, putting $\lambda=(-L,U)$ reduces
\eqref{eq:eigen-density} to \eqref{eq:LU-density}.

Write the positive \emph{magnitudes} of the negative eigenvalues in
increasing order as $u_1<\dots<u_n$, and the positive eigenvalues as
$v_1<\dots<v_n$. Then \eqref{eq:eigen-density} is proportional to
\begin{equation}\label{eq:two-groups-density}
 \Bigl[\Delta(u)^2\prod_{i=1}^n
 u_i^{N-n}e^{-u_i^2/(4s)}\Bigr]
 \Bigl[\Delta(v)^2\prod_{i=1}^n
 v_i^{N-n}e^{-v_i^2/(4s)}\Bigr]
 \prod_{i,j=1}^n(u_i+v_j).
\end{equation}
Thus the two groups are not exactly independent. In
\eqref{eq:H-limit} we show that they are asymptotically independent,
in the sense that the last factor can be replaced by a constant as
$N\to\infty$.

For fixed $n$ let $\mu_{N,n,q}$ be the probability measure on
$(0,\infty)$ with density
\[
 \rho(u)\propto u^{N-n}e^{-u^2/(4s)+qu},\qquad
 s=2T=2t^{2/3}N^{1/3},\quad A=\sqrt{2sN}=2\sqrt{NT}.
\]
Write $
\nu_q=\E_{\mu_{N,n,q}}u$, and $\sigma_q^2$ for its variance.

\begin{lemma}[Scalar Laplace estimate]\label{lem:scalar-laplace}
As $N\to\infty$,
\begin{equation}\label{eq:one-dimensional-mgf}
 \log\E_{\mu_{N,n,0}}e^{qu}
 =Aq+\frac s2q^2+\frac t3q^3+o(1)
\end{equation}
and
\begin{equation}\label{eq:nuq-asymptotic}
 \nu_q=A+sq+tq^2+o(1),
\end{equation}
uniformly for $q$ in compact sets.
\end{lemma}

For $n=1$, Lemma~\ref{lem:scalar-laplace} is precisely the estimate
\eqref{eq:scalar-laplace} used in the main text.

\begin{proof}
We first prove \eqref{eq:nuq-asymptotic}. Using
$\int_0^\infty\partial_u\bigl(u\rho(u)\bigr)\diff u=0$, we obtain
\begin{equation}\label{eq:nuq-quadratic}
 \nu_q^2-2sq\nu_q=2s(N-n+1)-\sigma_q^2.
\end{equation}
Solving \eqref{eq:nuq-quadratic} for $\nu_q$, with $A^2=2sN$, gives
$\nu_q=sq+\sqrt{A^2+s^2q^2+O(s)}$, where the term $O(s)$ is
$2s(1-n)-\sigma_q^2$. Indeed, the logarithmic density has curvature at
least $1/(2s)$; integration by parts, or the one-dimensional covariance
bound for a strongly log-concave measure, gives $0\le\sigma_q^2\le2s$.
Expanding the square root, with $s^2/A=2t$, gives
\eqref{eq:nuq-asymptotic} with the error $O(s/A)=o(1)$, uniformly for
$q$ in compact sets. Finally,
$\nu_q=\partial_q\log\E_{\mu_{N,n,0}}e^{qu}$, so integrating
\eqref{eq:nuq-asymptotic} in $q$ from $0$ gives
\eqref{eq:one-dimensional-mgf}.
\end{proof}

Also, a direct Laplace expansion
of the density around its mode gives, with convergence of every fixed
moment,
\begin{equation}\label{eq:CLT-mode}
 \frac{u-A}{\sqrt s}\ \Longrightarrow\ \mathcal N(0,1)
 \quad\text{under }\mu_{N,n,0},
 \qquad
 \frac{u-A}{s}\ \longrightarrow\ q
 \quad\text{under }\mu_{N,n,q}.
\end{equation}
The second assertion also follows from the variance bound and
\eqref{eq:nuq-asymptotic}; the first is an elementary quadratic
Laplace approximation with Gaussian domination from log concavity.

Let $\mu_{N,n}$ be the probability law on
$0<u_1<\dots<u_n$ with density proportional to
\begin{equation}\label{eq:n-ensemble}
 \Delta(u)^2\prod_{j=1}^n u_j^{N-n}e^{-u_j^2/(4s)}.
\end{equation}
For strictly increasing $q_1<\cdots<q_n$ put
\begin{equation}\label{eq:exp-divided}
 \mathcal E_q(u)=C_n\frac{\det[e^{q_i u_j}]_{i,j\le n}}{\Delta(u)}.
\end{equation}
This function is positive.

\begin{lemma}[Many-endpoint Laplace estimate]\label{lem:multi-laplace}
For fixed $n$ and $q_1<\cdots<q_n$,
\begin{equation}\label{eq:multi-laplace}
 \displaystyle
 \E_{\mu_{N,n}}\mathcal E_q(u)
 =\Delta(q)\exp\left\{A\sum_iq_i
 +\frac s2\sum_iq_i^2+\frac t3\sum_iq_i^3\right\}
 (1+o(1)).
\end{equation}
The estimate is uniform as $q$ ranges in compact subsets with fixed
distinctness gaps.
\end{lemma}
\begin{proof}
Andr\'eief's determinant identity gives the \emph{exact} formula
\begin{equation}\label{eq:andreief}
 \E_{\mu_{N,n}}\mathcal E_q(u)
 =C_n\frac{\det[m_{j-1}(q_i)]_{i,j\le n}}
 {\det[m_{i+j-2}(0)]_{i,j\le n}},
\end{equation}
where $m_k(q)=\int_0^\infty u^{N-n+k}e^{-u^2/(4s)+qu}\diff u$.
First we show that the numerator of \eqref{eq:andreief} satisfies
\[
 \det[m_{j-1}(q_i)]
 =s^{n(n-1)/2}\left(\prod_i m_0(q_i)\right)
 \bigl(\Delta(q)+o(1)\bigr).
\]
Indeed, factor $m_0(q_i)$ from row $i$, so that the $(i,j)$ entry is
$\E_{\mu_{N,n,q_i}}u^{j-1}$. Subtract $A$ from $u$, and scale the
$j$th column by $s^{-(j-1)}$. By the second limit in
\eqref{eq:CLT-mode}, $(u-A)/s\to q_i$ in every fixed moment under
$\mu_{N,n,q_i}$, so the remaining determinant tends to
$\det[q_i^{j-1}]=\Delta(q)$.
Next, the denominator of \eqref{eq:andreief} satisfies
\[
 \det[m_{i+j-2}(0)]
 =m_0(0)^n s^{n(n-1)/2}(C_n+o(1)).
\]
Indeed, perform the same centering, but scale by $\sqrt s$. By the
first limit in \eqref{eq:CLT-mode}, the normalized Hankel determinant
tends to the Gaussian moment Hankel determinant, which is
$\prod_{j=0}^{n-1}j!=C_n$ (the product of the squared norms of the
monic Hermite polynomials).
The $C_n$ factors cancel in \eqref{eq:andreief}. Apply
\eqref{eq:one-dimensional-mgf} to each
$m_0(q_i)/m_0(0)=\E_{\mu_{N,n,0}}e^{q_iu}$ to obtain
\eqref{eq:multi-laplace}.
\end{proof}

\textit{The interaction between the two groups.}
Let $H(u,v)=\prod_{i,j}(u_i+v_j)$ be the last factor in
\eqref{eq:two-groups-density}, and let $\E$ denote expectation under
the product of two copies of $\mu_{N,n}$, one for $u$ and one for $v$.
By \eqref{eq:two-groups-density}, the joint law of the two groups is
this product law weighted by $H$, so under the joint law the
expectation of $\mathcal E_q(u)\mathcal E_{q'}(v)$ is
$\E[H\mathcal E_q(u)\mathcal E_{q'}(v)]/\E H$. We show that, for fixed
$q$ and $q'$ as in Lemma~\ref{lem:multi-laplace},
\begin{equation}\label{eq:H-limit}
 \frac{\E\bigl[H\mathcal E_q(u)\mathcal E_{q'}(v)\bigr]}{\E H}
 =\E\mathcal E_q(u)\E\mathcal E_{q'}(v)(1+o(1)).
\end{equation}
Indeed, \eqref{eq:H-limit} follows from the two limits
\begin{equation}\label{eq:H-two-limits}
 \frac{\E\bigl[H\mathcal E_q(u)\mathcal E_{q'}(v)\bigr]}
 {(2A)^{n^2}\E\mathcal E_q(u)\E\mathcal E_{q'}(v)}
 \longrightarrow1,
 \qquad
 \frac{\E H}{(2A)^{n^2}}\longrightarrow1.
\end{equation}
To justify the first limit in \eqref{eq:H-two-limits}, observe that
$\mathcal E_q$ has the Harish-Chandra--Itzykson--Zuber representation
\[
 \mathcal E_q(u)=\Delta(q)\int_{U(n)}
 \exp\{\operatorname{tr}(\diag(q)U\diag(u)U^*)\}\diff U,
\]
where Haar measure has total mass one. Since every $u_j$ is positive,
the trace lies between $q_1\sum_ju_j$ and $q_n\sum_ju_j$, so
\[
 \Delta(q)e^{q_1\sum_ju_j}\le\mathcal E_q(u)
 \le\Delta(q)e^{q_n\sum_ju_j}.
\]
In the variables $x_j=u_j/A$, the density \eqref{eq:n-ensemble} is
proportional to
\[
 \Delta(x)^2\prod_jx_j^{-n}e^{N\sum_j(\log x_j-x_j^2/2)},
\]
so $|u_j-A|\ge\varepsilon A$ has probability at most
$e^{-c_\varepsilon N}$, whereas the factors $e^{q_1\sum_ju_j}$ and
$e^{q_n\sum_ju_j}$ are $e^{q_1A\sum_jx_j}$ and $e^{q_nA\sum_jx_j}$,
with $A=O(N^{2/3})$.
It follows that, under $\mu_{N,n}$ reweighted by
$\mathcal E_q$, each $u_i/A\to1$ in probability, and the moments of
$u_i/A$ are bounded uniformly in $N$. The same holds for the
$v$-group. Since $H$ is a
polynomial of fixed degree $n^2$, uniform integrability gives
the first limit in \eqref{eq:H-two-limits}. The same concentration
argument, applied to the product law without the weights
$\mathcal E_q$ and $\mathcal E_{q'}$, gives the second limit.

\subsection{Totally nonnegative Grassmannian solitons}

The preceding interlacing and asymptotic statements apply to
positive Darboux chains. We next explain why a totally nonnegative
Grassmannian matrix admits such a chain, and then perform the
Airy edge evaluation.

\textit{Positive Darboux chains.}\label{app:positive}
The Darboux argument assumes positivity of each intermediate
Wronskian $W_k$, not merely that of $W_n$. We show that, after a
constant change of row basis, all maximal minors of the first $k$ rows
of $c$ are nonnegative for every $k\le n$. It follows that
$W_k(x,r)>0$ for all $x,r$. Indeed, by the Cauchy--Binet formula,
\[
 W_k(x,r)=\sum_{|I|=k}\det(c_{ij})_{i\le k,j\in I}
 \Delta(p_I)e^{\sum_{j\in I}(p_jr-p_j^2x)},
\]
where every term is nonnegative and at least one is positive, since
the first $k$ rows have rank $k$. The change of row basis multiplies
$W_n$ and every $c_I$ by the same constant, so it changes neither
$-2\partial_r^2\log W_n$ nor the right side of \eqref{eq:general-tau}.

For the change of row basis, first put $c$ in reduced row echelon
form, with pivot columns $j_1<\cdots<j_n$. The maximal minors of the
reduced matrix are still nonnegative, since the minor on the pivot
columns is $1$. Moreover, for every $k\le n$, all maximal minors of
the last $k$ rows of the reduced matrix are nonnegative. Indeed, the
last $k$ rows vanish to the left of column $j_{n-k+1}$, so their minor
on a column set $J$ is zero unless
$J\subset\{j_{n-k+1},\ldots,m\}$. For such $J$, the minor equals the
maximal minor of the reduced matrix on the columns
$\{j_1,\ldots,j_{n-k}\}\cup J$, because the columns
$j_1,\ldots,j_{n-k}$ of the reduced matrix are the first $n-k$ unit
vectors. Next, reverse the order of the rows and multiply the $i$th
row by $(-1)^{i-1}$. For the first $k$ rows, the row-reversal sign
$(-1)^{k(k-1)/2}$ is cancelled by the factor
$\prod_{i\le k}(-1)^{i-1}$, so all their maximal minors are again
nonnegative.

\textit{Airy edge evaluation.}
We now prove Theorem~\ref{thm:general}. After the change of row basis,
$W_k>0$ for every $k\le n$, so Proposition~\ref{prop:n-endpoints}
applies. Under the scaling \eqref{eq:s8-edge},
$Z_N(a,b;R;V_n)\to Z_{V_n}(t,x,r)$ by \eqref{eq:bridge-limit}. By
\eqref{eq:n-endpoints}, the formula \eqref{eq:general-tau} is therefore
equivalent to
\begin{equation}\label{eq:edge-goal}
 \frac{\E[\mathcal W_b(R+\lambda^-)\mathcal W_b(R+\lambda^+)]}
 {W_n(a,R)W_n(b,R)}
 \longrightarrow\frac{\tau_n(t,x,r)}{\alpha_{I_*}e^{\theta_{I_*}}}.
\end{equation}
The limit \eqref{eq:edge-goal} follows from the two estimates
\begin{equation}\label{eq:edge-numerator}
\begin{split}
 &\E[\mathcal W_b(R+\lambda^-)\mathcal W_b(R+\lambda^+)]\\
 &\qquad=\tau_n(t,x,r)\alpha_{I_*}
 \exp\left\{(2A+r)\sum_{j\in I_*}p_j-x\sum_{j\in I_*}p_j^2
 +\frac t3\sum_{j\in I_*}p_j^3\right\}(1+o(1))
\end{split}
\end{equation}
and
\begin{equation}\label{eq:edge-denominator}
 W_n(a,R)W_n(b,R)
 =\alpha_{I_*}^2
 \exp\left\{2(A+r)\sum_{j\in I_*}p_j
 -2x\sum_{j\in I_*}p_j^2\right\}(1+o(1)).
\end{equation}

To prove \eqref{eq:edge-numerator}, expand $\mathcal W_b$
using Cauchy--Binet:
\begin{equation}\label{eq:W-expand}
 \mathcal W_b(z)
 =\sum_{|I|=n}c_Ie^{-b\sum_{j\in I}p_j^2}
 C_n\frac{\det[e^{p_i z_j}]_{i\in I,\ j\le n}}{\Delta(z)}.
\end{equation}
Recall from Section~\ref{app:asymptotics} that
$\lambda^-=(-u_n,\ldots,-u_1)$ and $\lambda^+=(v_1,\ldots,v_n)$. At
$z=R+\lambda^-$, reversing the order of the rows and of the columns of
the determinant shows that the term $I$ of \eqref{eq:W-expand} is
$c_Ie^{\sum_{j\in I}(Rp_j-bp_j^2)}\mathcal E_q(u)$, where $q$ consists
of the phases $-p_j$, $j\in I$, in increasing order. By
Lemma~\ref{lem:multi-laplace}, the expectation of this term under
$\mu_{N,n}$ is
\begin{equation}\label{eq:edge-lower}
 \alpha_I\exp\left\{(R-A)\sum_{j\in I}p_j
 -\Bigl(b-\frac s2\Bigr)\sum_{j\in I}p_j^2
 -\frac t3\sum_{j\in I}p_j^3\right\}(1+o(1))
 =\alpha_Ie^{\theta_I}(1+o(1)),
\end{equation}
since $R-A=r$ and $b-s/2=x$. At $z=R+\lambda^+$, the term $J$ of
\eqref{eq:W-expand} is
$c_Je^{\sum_{j\in J}(Rp_j-bp_j^2)}\mathcal E_{q'}(v)$, where $q'$
consists of the phases $p_j$, $j\in J$, in increasing order. By
Lemma~\ref{lem:multi-laplace}, the expectation of this term under
$\mu_{N,n}$ is
\begin{equation}\label{eq:edge-upper}
 \alpha_J\exp\left\{(2A+r)\sum_{j\in J}p_j
 -x\sum_{j\in J}p_j^2
 +\frac t3\sum_{j\in J}p_j^3\right\}(1+o(1)),
\end{equation}
since $R+A=2A+r$ and $b-s/2=x$. By \eqref{eq:H-limit}, the expectation
of the product of these two terms under the joint law
\eqref{eq:two-groups-density} is the product of \eqref{eq:edge-lower}
and \eqref{eq:edge-upper}, up to a factor $1+o(1)$. All these
products are nonnegative, so summing over $I$ and $J$ shows that
$\E[\mathcal W_b(R+\lambda^-)\mathcal W_b(R+\lambda^+)]$ is the sum
of \eqref{eq:edge-lower} over $I$ times the sum of
\eqref{eq:edge-upper} over $J$, up to a factor $1+o(1)$. The first sum
is
$\tau_n(t,x,r)(1+o(1))$. In the second sum the term $J=I_*$ dominates
as $A\to\infty$, since $I_*$ is the unique maximizer of
$\sum_{j\in J}p_j$ among the $J$ with $\alpha_J>0$. This proves
\eqref{eq:edge-numerator}.

To prove \eqref{eq:edge-denominator}, use the Cauchy--Binet formula
$W_n(a,R)=\sum_{|I|=n}\alpha_Ie^{\sum_{j\in I}(Rp_j-ap_j^2)}$ and the
same formula with $b$ in place of $a$. Since $R=A+r$, $a=x-s/2$,
$b=x+s/2$ and $s/A\to0$, the term $I=I_*$ dominates both sums, and
\eqref{eq:edge-denominator} follows. This completes the proof of
Theorem~\ref{thm:general}.

\section{Proofs for Section~\ref{sec:kpz-fixed-point}}\label{app:kpz}

\subsection{Regularity and the initial-data limit}\label{app:kpz-wall}

\begin{lemma}\label{lem:kpz-smoothing}
Let $\phi$ be a distributional solution of \eqref{KP} on $(0,T)\times\R^2$, and let $J'\subset(0,T)$ be an interval. Suppose that there are $B<\infty$, $0\le\beta<2$, $\delta>0$ and, for each $\mu>0$, a constant $C_\mu<\infty$ such that, almost everywhere on $J'\times\R^2$,
\[
 |\phi(t,x,r)|\le C_\mu e^{B(1+|x|)^\beta+\mu r_--\delta r_+}.
\]
Then $\phi$ agrees almost everywhere with a function that is $C^\infty$ on $(\operatorname{int}J')\times\R^2$, and it satisfies the derivative bounds of Lemma~\ref{lem:representation-smoothing}.
\end{lemma}

\begin{proof}
The bound on $|\phi|$ assumed in Lemma~\ref{lem:representation-smoothing} enters its proof only through the first step, which shows that for $a>4B$ and $0<\lambda<\delta$ the functions $W_{a,\lambda}\phi$ and $W_{a/2,\lambda/2}\phi$ are bounded in $L^2(\R^2)\cap L^\infty(\R^2)$, uniformly for $t\in J'$, with the margin $e^{(2B-a/2)\langle x\rangle^\kappa}$ in $x$. Under the present bound with $\mu=\lambda/4$, both functions are at most $C_\mu e^{(2B-a/2)\langle x\rangle^\kappa}e^{-\lambda r_-/4-(\delta-\lambda)r_+}$, so this step holds, and the rest of the proof is unchanged.
\end{proof}

If $\phi$ is a classical solution of KP that satisfies $(\aleph 1)$, then the hypothesis of Lemma~\ref{lem:kpz-smoothing} holds on every compact time interval, with $B=0$. So the derivatives of $\phi$ satisfy the first bound in \eqref{eq:representation-bounds} with $\beta=1$, and the second bound holds since $\phi\ge0$.

\begin{lemma}\label{lem:kpz-countable}
There are finite sets $S_1\subset S_2\subset\cdots$ on which $h>-\infty$, with union $S$, such that $\sup_{y\in\R}(h(y)+g(y))=\sup_{y\in S}(h(y)+g(y))$ for every continuous $g\colon\R\to\R$.
\end{lemma}

\begin{proof}
For $m\ge1$ and $k\in\mathbb Z$, if $h$ is not identically $-\infty$ on $[k2^{-m},(k+1)2^{-m}]$, choose a point of this interval at which $h$ attains its maximum there; it exists because $h$ is upper semicontinuous. Let $S_n$ be the set of points chosen in intervals inside $[-n,n]$ with $m\le n$. If $h(y)>-\infty$, the points $y_m$ chosen in the intervals of length $2^{-m}$ that contain $y$ satisfy $h(y_m)\ge h(y)$ and $y_m\to y$, so $\liminf_m(h(y_m)+g(y_m))\ge h(y)+g(y)$.
\end{proof}

\begin{proof}[Proof of Lemma~\ref{lem:kpz-wall}]
We realize the curves for all base points $x$ from one stationary ensemble $\bar{\mathcal A}$, as in \eqref{eq:markov-airy-marginal}, and write $I_s=I_{\phi(s)}^{t,x}(r)$. Stationarity, the unit-interval upper-tail bound for the Airy$_2$ process (from the GOE maximum identity and the soft-edge right-tail bound \cite{LedouxRider2010}), and Borel--Cantelli imply $\limsup_{|y|\to\infty}\bar{\mathcal A}_1(y)/y^2\le0$ almost surely. Since $h$ is bounded above and $\mathcal A_1^{t,x}$ contains a negative parabola, the upper semicontinuous function $h+\mathcal A_1^{t,x}$ tends to $-\infty$ at $\pm\infty$, so its supremum is attained, and by Lemma~\ref{lem:kpz-countable} it is a supremum over a countable set. So the events $\{r-\mathcal A_1^{t,x}(y)\ge h(y)\text{ for all }y\}$ and $\mathcal S=\{r-\mathcal A_1^{t,x}(y)>h(y)\text{ for all }y\}$ are $\{\sup(h+\mathcal A_1^{t,x})\le r\}$ and $\{\sup(h+\mathcal A_1^{t,x})<r\}$. By \eqref{eq:kpz-wall}, and since $F_h$ is continuous \cite[Theorem~1.4]{QR}, both have probability $F_h(t,x,r)$. It therefore suffices to show that $\limsup_{s\downarrow0}Z_{\phi(s)}(t,x,r)$ is at most the probability of the first event and $\liminf_{s\downarrow0}Z_{\phi(s)}(t,x,r)\ge\mathbb P(\mathcal S)$.

\textit{Curves that cross the graph.} Suppose that $r-\mathcal A_1^{t,x}(y_0)<b<h(y_0)$ for some $y_0$ and $b$. By continuity, $r-\mathcal A_1^{t,x}<b$ on an interval $(y_0-\varepsilon,y_0+\varepsilon)$, and there $(\aleph 4)$ gives $\phi(s,y,r-\mathcal A_1^{t,x}(y))\ge\phi(s,y,b)$. So
\[
 I_s\ge\int_{y_0-\varepsilon}^{y_0+\varepsilon}\phi(s,y,b)\diff y\longrightarrow\infty
\]
by $(\aleph 3)$. Hence $e^{-I_s}\to0$ outside the first event, and bounded convergence gives the bound on the $\limsup$.

\textit{Curves above the graph.} For $L,M>0$, split $I_s=I_s^{\rm low}+I_s^{\rm high}+I_s^{\rm far}$ according to where the curves are: over $|y-x|\le L$ at heights at most $M$, over $|y-x|\le L$ at heights above $M$, and over $|y-x|>L$. We show below that
\begin{equation}\label{eq:kpz-three-parts}
 I_s^{\rm low}\to0\ \text{ on }\mathcal S,\qquad\sup_{0<s\le s_0}\E I_s^{\rm high}\xrightarrow[M\to\infty]{}0,\qquad\E\bigl[1-e^{-I_s^{\rm far}}\bigr]\le2Ce^{-\delta(r+L^2/t)}.
\end{equation}
Since $e^{-a-b-c}\ge e^{-a}-b-(1-e^{-c})$ for $a,b,c\ge0$,
\[
 Z_{\phi(s)}(t,x,r)\ge\E\bigl[\mathbf 1_{\mathcal S}e^{-I_s^{\rm low}}\bigr]-\E I_s^{\rm high}-2Ce^{-\delta(r+L^2/t)},
\]
and letting $s\downarrow0$, then $M\to\infty$, then $L\to\infty$ gives the bound on the $\liminf$.

On $\mathcal S$, the curves over $|y-x|\le L$ at heights at most $M$ lie in the compact set $\{(y,u):|y-x|\le L,\ r-\mathcal A_1^{t,x}(y)\le u\le M\}\subset\mathrm{epi}_S(h)$, on which $\phi(s)\to0$ uniformly by $(\aleph 2)$, and their total length $\sum_i\int_{x-L}^{x+L}\mathbf 1\{r-\mathcal A_i^{t,x}(y)\le M\}\diff y$ has finite mean by \eqref{eq:occ-density}. This gives the first claim in \eqref{eq:kpz-three-parts}.

For the third claim, let $y\ge x+L$. Then
\[
 \mathcal A_i^{t,x+L}(y)-\mathcal A_i^{t,x}(y)=\frac{(y-x)^2-(y-x-L)^2}t\ge\frac{L^2}t,
\]
so there the curves from $(x,r)$ lie above those from $(x+L,r+L^2/t)$. By $(\aleph 4)$, the part of $I_s^{\rm far}$ from $y\ge x+L$ is therefore at most
\[
 I_{\phi(s)}^{t,x+L}(r+L^2/t),
\]
and $Z_{\phi(s)}(t,x+L,r+L^2/t)\ge1-Ce^{-\delta(r+L^2/t)}$ by \eqref{eq:kpz-uniform-tail}. The part from $y\le x-L$ is bounded in the same way, with $x-L$ in place of $x+L$, and $1-e^{-b-c}\le(1-e^{-b})+(1-e^{-c})$.

For the second claim, fix $L$, and choose $C_L$ with $\mathbb P\{\max_{|y-x|\le L}(-\mathcal A_1^{t,x}(y))\le C_L\}\ge\frac12$. On this event the lowest curve from height $M-C_L$ stays at or below $M$ over $[x-L,x+L]$, so by $(\aleph 4)$ its occupation is at least $m_s(M)=\int_{x-L}^{x+L}\phi(s,y,M)\diff y$. By \eqref{eq:kpz-uniform-tail},
\[
 \tfrac12\bigl(1-e^{-m_s(M)}\bigr)\le1-Z_{\phi(s)}(t,x,M-C_L)\le Ce^{-\delta(M-C_L)},
\]
so $m_s(M)\le4Ce^{\delta C_L}e^{-\delta M}$ for $M\ge M_0$ and $0<s\le s_0$. For $|y-x|\le L$ and $u\ge0$, the density $\varrho_{t,x,r}(y,u)$ in \eqref{eq:occ-density} is at most $C'(1+\sqrt u)$, by \eqref{eq:airy-density-bound}. Hence, for $M\ge M_0$,
\[
 \E I_s^{\rm high}=\int_M^\infty\int_{x-L}^{x+L}\phi(s,y,u)\varrho_{t,x,r}(y,u)\diff y\diff u\le C'\int_M^\infty(1+\sqrt u)m_s(u)\diff u,
\]
which tends to $0$ as $M\to\infty$, uniformly in $s$.
\end{proof}

\subsection{The curvature inequality}\label{app:kpz-curvature}

\begin{proof}[Proof of Lemma~\ref{lem:kpz-gaussian}]
\textit{The tilt.} Put $\mathbf e=(1,\ldots,1)\in\R^n$ and $\Lambda=\{a:c_a<\infty\}$. In the variables $y_i=\theta\mathbf e-X_i$ the conditions become $y_{1a}<y_{2a}<\dots<y_{Na}$ for all $a$ and $y_{1a}\ge-c_a$ for $a\in\Lambda$. Let $\mathcal Y\subset\R^{Nn}$ be the set they define. Expanding the Gaussian density,
\[
 \log P(\theta)=\log Z(\theta)-\frac N2(\mathbf e^{\mathsf T}\Gamma\mathbf e)\theta^2+\theta\sum_{i=1}^N\mathbf e^{\mathsf T}\Gamma\mu_i+\text{const},\qquad Z(\theta)=\int_{\mathcal Y}e^{-\mathcal V_0(y)+\theta f\cdot y}\diff y,
\]
where $\mathcal V_0(y)=\frac12\sum_i(y_i+\mu_i)^{\mathsf T}\Gamma(y_i+\mu_i)$ and $f_{ia}=(\Gamma\mathbf e)_a\ge0$. So it suffices to show that $(\log Z)'$, the mean of $f\cdot y$ under the probability proportional to $\mathbf 1_{\mathcal Y}e^{-\mathcal V_0+\theta f\cdot y}$, is convex in $\theta$.

\textit{Soft constraints.} For $\gamma>0$ let
\[
 \mathcal V(y)=\mathcal V_0(y)+\sum_{a\in\Lambda}e^{-\gamma(y_{1a}+c_a)}+\sum_{a=1}^n\sum_{i=1}^{N-1}e^{-\gamma(y_{i+1,a}-y_{ia})},\qquad Z_\gamma(\theta)=\int_{\R^{Nn}}e^{-\mathcal V(y)+\theta f\cdot y}\diff y.
\]
As $\gamma\to\infty$, $e^{\mathcal V_0-\mathcal V}\to\mathbf 1_{\mathcal Y}$ off the boundary of $\mathcal Y$, and $e^{\mathcal V_0-\mathcal V}\le1$. By dominated convergence, $Z_\gamma$ and its derivatives in $\theta$ converge to those of $Z$, locally uniformly, and $Z>0$. So it suffices to show that $(\log Z_\gamma)'$ is convex.

\textit{Two properties of $\mathcal V$.} The terms added to $\mathcal V_0$ are convex functions of linear forms, so the Hessian $H(y)=\nabla^2\mathcal V(y)$ satisfies $H\ge\lambda I$, where $\lambda>0$ is the least eigenvalue of $\Gamma$. Let
\[
 \mathcal K=\bigl\{w\in\R^{Nn}:0\le w_{1a}\le w_{2a}\le\dots\le w_{Na}\text{ for }1\le a\le n\bigr\},
\]
which contains $f$. We claim:
\begin{enumerate}
\item[(a)] for every continuous path $u\mapsto y(u)$, the solution operators of $\dot w=-H(y(u))w$ and their transposes map $\mathcal K$ into $\mathcal K$;
\item[(b)] $\nabla^3\mathcal V(y)[w,z,z]\le0$ for $w\in\mathcal K$ and $z\in\R^{Nn}$.
\end{enumerate}
Fix $y$ and a site $a$, and write $\varpi_a=\gamma^2e^{-\gamma(y_{1a}+c_a)}$ if $a\in\Lambda$ and $\varpi_a=0$ otherwise, $\omega_i=\gamma^2e^{-\gamma(y_{i+1,a}-y_{ia})}$ for $1\le i<N$, and $\omega_0=\omega_N=0$. Then
\[
 (Hw)_{ia}=(\Gamma w_i)_a+\delta_{i1}\varpi_aw_{1a}+\omega_{i-1}(w_{ia}-w_{i-1,a})+\omega_i(w_{ia}-w_{i+1,a}).
\]
In the coordinates $v_{1a}=w_{1a}$ and $v_{ia}=w_{ia}-w_{i-1,a}$, $\mathcal K$ is the orthant $v\ge0$. For (a), it suffices that on each face of $\mathcal K$ the vector $-Hw$ points into $\mathcal K$: then $-H$ has nonnegative off-diagonal entries in these coordinates, so its solution operators have nonnegative entries. The transpose of a solution operator is the solution operator for the same Hessians in reversed time, so it is covered as well. On the face $w_{1a}=0$,
\[
 (Hw)_{1a}=\sum_{b\ne a}\Gamma_{ab}w_{1b}-\omega_1w_{2a}\le0.
\]
On the face $w_{ia}=w_{i+1,a}$, the terms with $\Gamma_{aa}$ and $\omega_i$ cancel, and
\begin{align*}
 (Hw)_{ia}-(Hw)_{i+1,a}={}&\sum_{b\ne a}(-\Gamma_{ab})(w_{i+1,b}-w_{ib})+\delta_{i1}\varpi_aw_{1a}\\
 &+\omega_{i-1}(w_{ia}-w_{i-1,a})+\omega_{i+1}(w_{i+2,a}-w_{i+1,a}),
\end{align*}
where each term is nonnegative for $w\in\mathcal K$. For (b), the third derivative of $e^{-\gamma\xi}$ is $-\gamma$ times its second derivative, so
\[
 \nabla^3\mathcal V(y)[w,z,z]=-\gamma\sum_{a\in\Lambda}\varpi_aw_{1a}z_{1a}^2-\gamma\sum_{a=1}^n\sum_{i=1}^{N-1}\omega_i(w_{i+1,a}-w_{ia})(z_{i+1,a}-z_{ia})^2\le0,
\]
with $\omega_i$ taken at the site $a$ of each term.

\textit{The Langevin dynamics.} Let $B$ be a Brownian motion in $\R^{Nn}$, and for $\theta\in\R$ let $Y_\theta$ solve
\[
 \diff Y_\theta(u)=\bigl(-\nabla\mathcal V(Y_\theta(u))+\theta f\bigr)\diff u+\sqrt2\diff B(u),\qquad Y_\theta(0)=0.
\]
Since $\mathcal V$ is smooth and $\langle\nabla\mathcal V(y)-\nabla\mathcal V(y'),y-y'\rangle\ge\lambda|y-y'|^2$, the difference $Y_\theta-\sqrt2B$ solves an ordinary differential equation whose vector field is continuous in $u$ and smooth in the other variables, and whose solutions are global. So $Y_\theta(u)$ is smooth in $\theta$. Its derivatives $z=\partial_\theta Y_\theta$ and $\zeta=\partial_\theta^2Y_\theta$ satisfy
\[
 \dot z=-Hz+f,\qquad\dot\zeta=-H\zeta-\nabla^3\mathcal V[z,z],\qquad z(0)=\zeta(0)=0,
\]
with $H$ and $\nabla^3\mathcal V$ evaluated at $Y_\theta(u)$. With $U(u,v)$ the solution operator of $\dot w=-Hw$ from time $v$ to time $u$,
\[
 \partial_\theta^2\bigl(f\cdot Y_\theta(u)\bigr)=f\cdot\zeta(u)=-\int_0^u\nabla^3\mathcal V(Y_\theta(v))\bigl[U(u,v)^{\mathsf T}f,z(v),z(v)\bigr]\diff v\ge0,
\]
by (a), (b) and $f\in\mathcal K$. So $f\cdot Y_\theta(u)$ is convex in $\theta$ for every $u$, and so is its mean. The probability $\pi_\theta$ with density proportional to $e^{-\mathcal V+\theta f\cdot y}$ is invariant for $Y_\theta$, and a stationary solution $Y'_\theta$ driven by the same $B$ satisfies $|Y_\theta(u)-Y'_\theta(u)|\le e^{-\lambda u}|Y'_\theta(0)|$. Since $\pi_\theta$ has Gaussian tails, $\E f\cdot Y_\theta(u)\to\int f\cdot y\diff\pi_\theta(y)=(\log Z_\gamma)'(\theta)$ as $u\to\infty$, and a pointwise limit of convex functions is convex.
\end{proof}

\begin{proof}[Proof of Theorem~\ref{thm:kpz-curvature}]
\textit{Passing to limits.} Let $P_k$ and $P$ be positive nondecreasing functions on $\R$ with $P_k\to P$ pointwise, and suppose that $(\log P_k)'''\ge0$ in the sense of distributions. For smooth $\chi\ge0$ supported in $[\alpha,\beta]$, $P_k\ge P_k(\alpha)\ge P(\alpha)/2$ on $[\alpha,\beta]$ for large $k$, so dominated convergence passes $-\int\log P_k\,\chi'''\ge0$ to the limit, and $(\log P)'''\ge0$ in the sense of distributions. An increasing affine change of the variable preserves this.

\textit{The Airy$_2$ process.} Let $\eta_1<\dots<\eta_m$ and $c_1,\dots,c_m\in\R$, and let
\[
 P(\theta)=\mathbb P\{\bar{\mathcal A}_1(\eta_j)\le\theta+c_j,\ 1\le j\le m\}.
\]
We show $(\log P)'''\ge0$ in the sense of distributions, by passing to the limit in three steps from Lemma~\ref{lem:kpz-gaussian}. In each step the limiting distribution functions are continuous, so the convergence holds at every $\theta$, and all functions involved are positive and nondecreasing.

Let $B_1,\dots,B_N$ be independent Brownian bridges with diffusion coefficient $2$ on $[-T,T]$, from $a_1>\dots>a_N$ to $b_1>\dots>b_N$, and let $-T<u_1<\dots<u_n<T$ be a grid containing $\eta_1,\dots,\eta_m$. With $u_0=-T$, $u_{n+1}=T$ and $d_a=u_a-u_{a-1}$, the vectors $(B_i(u_1),\dots,B_i(u_n))$ are independent Gaussian vectors with the common precision matrix
\[
 \Gamma_{aa}=\frac1{2d_a}+\frac1{2d_{a+1}},\qquad\Gamma_{a,a+1}=\Gamma_{a+1,a}=-\frac1{2d_{a+1}},
\]
and $\Gamma_{ab}=0$ for $|a-b|\ge2$. For $n\ge2$ its row sums are $1/(2d_1)$, $0,\dots,0$, $1/(2d_{n+1})$; for $n=1$ the single row sum is $1/(2d_1)+1/(2d_2)$. Thus Lemma~\ref{lem:kpz-gaussian} applies to the conditional probability that $B_1(\eta_j)\le\theta+c_j'$ for all $j$, given that $B_1>\dots>B_N$ on the grid.

First, refine the grid to a dense set. The order on the grid then decreases to $\{B_1\ge\dots\ge B_N\}$. This event differs from the event that the bridges do not intersect by a null set, since the difference of two of the bridges is a bridge with positive endpoints whose minimum has a continuous distribution, and the bridges do not intersect with positive probability. So the conditional probabilities converge to those of nonintersecting bridges. Second, let all $a_i$ and $b_i$ tend to $0$. The bridges converge in law to the nonintersecting bridges of Appendix~\ref{app:edge-scaled}. Third, take $T=N^{1/3}$, so that $t=1$ in \eqref{eq:db-edge-scaled}, and $c_j'=2\sqrt{NT}+c_j-\eta_j^2$. Then $B_1(\eta_j)\le\theta+c_j'$ is the event $G_1^N(\eta_j)\le\theta+c_j-\eta_j^2$ in the edge scaling \eqref{eq:db-edge-scaled}, and as $N\to\infty$ its probability tends to $P(\theta)$, by \eqref{eq:db-edge-limit}.

\textit{The KPZ fixed point.} If $h$ is finite only at finitely many spatial points $y_1,\dots,y_m$, \eqref{eq:kpz-finite-sources} writes $F_h(t,x,r)$ as $P(t^{-1/3}r)$, with $\eta_j=t^{-2/3}y_j$ and $c_j=t^{-1/3}((y_j-x)^2/t-h(y_j))$. For general $h$, let $S_n$ be as in Lemma~\ref{lem:kpz-countable}, and let $h_n=h$ on $S_n$ and $h_n=-\infty$ elsewhere. By \eqref{eq:kpz-wall} and Lemma~\ref{lem:kpz-countable}, applied with $g=\mathcal A_1^{t,x}$, $F_{h_n}(t,x,r)\downarrow F_h(t,x,r)$ for every $r$, and $F_h>0$ by \eqref{eq:kpz-TW}. So $(\log F_h(t,x,\cdot))'''\ge0$ in the sense of distributions, and pointwise since $F_h$ is smooth.
\end{proof}

\subsection{Identifying the parabolic lower-tail limit}\label{app:kpz-quadratic-tail}
Here we give the technical argument behind Theorem~\ref{thm:kpz-quadratic-tail}.
We use \(\bar{\mathcal A}_1\) for the stationary Airy$_2$ process in the
normalization of Section~\ref{sec:kpz-fixed-point}. The inputs are the KP
equation, the curvature inequality of Theorem~\ref{thm:kpz-curvature},
Johansson's flat-data identity, positive association of the Airy$_2$ process,
and nonsharp soft-edge tail estimates. In particular, no sharp
Tracy--Widom left-tail asymptotic is used.

\begin{lemma}[Exact parabolic reduction]\label{lem:kpz-quad-reduction}
For $b\ge1$, define
\[
G_b(s)=\mathbb P\!\left\{\sup_{y\in\R}\bigl(\bar{\mathcal A}_1(y)-by^2\bigr)\le s\right\},
\qquad g_b(s)=-2\partial_s^2\log G_b(s).
\]
Then with $b=1+at$ and $B(t)=a/(1+at)$,
\begin{equation}\label{eq:kpz-quad-exact}
F_a(t,x,r)=G_b\!\left(t^{-1/3}(r+B(t)x^2)\right),
\end{equation}
and for $b>1$,
\begin{equation}\label{eq:kpz-quad-g}
(b-1)g_{b,b}-\frac s3 g_{b,s}-\frac12g_bg_{b,s}
 +\frac1{12}g_{b,sss}
 +\left(\frac{b-1}{2b}-\frac23\right)g_b=0.
\end{equation}
At $b=1$ the stationary equation obtained by deleting the $g_{b,b}$ term holds, derived directly from flat-data KP and self-similarity.
\end{lemma}
\begin{proof}
The variational formula gives
$F_a(t,x,r)=\mathbb P\{\sup_y[t^{1/3}\bar{\mathcal A}_1(t^{-2/3}y)-ay^2-(x-y)^2/t]\le r\}$.
Complete the square:
\[
-ay^2-\frac{(x-y)^2}{t}
=-\frac{ax^2}{1+at}-\frac{1+at}{t}
 \left(y-\frac{x}{1+at}\right)^2.
\]
Stationarity of $\bar{\mathcal A}_1$ and the change of variables $y=t^{2/3}u$ give \eqref{eq:kpz-quad-exact}. Put $f_a(t,r)=\phi_a(t,0,r)$ and $q=r+B(t)x^2$. Then $\phi_a(t,x,r)=f_a(t,q)$, and
\[
\partial_x^2\partial_r^{-1}\phi_a=2Bf_a+4B^2x^2\partial_qf_a,
\quad B'=-B^2.
\]
Thus \eqref{KP} reduces to
$f_{a,t}-\tfrac12 f_af_{a,r}+\tfrac1{12}f_{a,rrr}+\tfrac B2 f_a=0$.
Insert $f_a(t,r)=t^{-2/3}g_{1+at}(t^{-1/3}r)$ and multiply by $t^{5/3}$. This gives \eqref{eq:kpz-quad-g}. If $a=0$, $b=1$ throughout and the flat-data equation directly gives its $b=1$ version.
\end{proof}

\begin{lemma}[Coarse tails and derivative bounds]\label{lem:kpz-quad-coarse}
There exist constants $c,C,S_0>0$ such that, for every $b\ge1$ and $S\ge S_0$,
\begin{align}\label{eq:kpz-quad-cubic}
 c S^3\le -\log G_b(-S)\le C S^3,\\
\label{eq:kpz-quad-linear}
 c S\le g_b(-S)\le C S.
\end{align}
For every fixed $B<\infty$ there is a $C_B$ such that
\begin{equation}\label{eq:kpz-quad-right}
\sup_{1\le b\le B}\left(g_b(0)+\int_0^\infty g_b(s)^2\,ds\right)\le C_B.
\end{equation}
Consequently, for $R\ge1$, $b\in[1,B]$ and $v_R(b,z)=R^{-1}g_b(Rz)$,
\[
0\le v_R(b,z)\le C_B(1+(-z)_+),\qquad
\int_0^\infty v_R(b,z)^2\,dz\le C_B R^{-3},
\]
and, for $R\ge S_0$ and $z\le-1$, $v_R(b,z)\ge c(-z)$.
\end{lemma}
\begin{proof}
Since $b\ge1$, increasing $b$ makes the event defining $G_b(s)$ more likely, while that event requires $\bar{\mathcal A}_1(0)\le s$ (the parabola vanishes at $y=0$). Hence
\[
G_1(s)=F_{\rm GOE}(4^{1/3}s)\le G_b(s)\le F_{\rm GUE}(s).
\]
The finite-$n$ $\beta$-Hermite bounds of Ledoux--Rider~\cite[Theorems~1 and~4]{LedouxRider2010}, for $\beta=1,2$, give two-sided estimates $e^{-CS^3}\le F_\beta(-S)\le e^{-cS^3}$ for $S$ large, by taking the soft-edge limit. This yields \eqref{eq:kpz-quad-cubic}.

Write $H_b=-\log G_b$. By the curvature inequality, $g_b\ge0$ and $g_b'\le0$; the right-tail decay and integration twice give
\begin{equation}\label{eq:kpz-quad-repr}
H_b(s)=\frac12\int_s^\infty(u-s)g_b(u)\,du.
\end{equation}
Since $H_b(0)\le H_1(0)$ and $H_b(-1)\le H_1(-1)$, monotonicity in $s$ gives $g_b(0)\le 4H_1(-1)$ and
$\int_0^\infty g_b(u)\,du\le2H_1(-1)$.
Thus $\int_0^\infty g_b(u)^2\,du\le8H_1(-1)^2$. Moreover, for $S\ge1$,
\[
H_b(-S)\le\tfrac14S^2g_b(-S)+H_b(0)+\tfrac S2\int_0^\infty g_b,
\]
where the negative-interval integral is bounded using $g_b'\le0$. Combining with $H_b(-S)\ge cS^3$ gives the lower bound in \eqref{eq:kpz-quad-linear}. Conversely
$H_b(-2S)\ge\tfrac14 S^2g_b(-S)$, and \eqref{eq:kpz-quad-cubic} gives the upper bound. The remaining assertions follow immediately by rescaling.
\end{proof}

\begin{lemma}[Flat constant from KP]\label{lem:kpz-quad-flat}
For $v_R(1,z)=R^{-1}g_1(Rz)$,
\begin{equation}\label{eq:kpz-quad-flat}
v_R(1,z)\longrightarrow2(-z)_+
\end{equation}
locally uniformly, and $-\log G_1(-S)=S^3/6+o(S^3)$.
\end{lemma}
\begin{proof}
Set $\varepsilon=R^{-3}$. The stationary instance of \eqref{eq:kpz-quad-g} is
\[
-\tfrac z3 v_R'-\tfrac12v_Rv_R'-\tfrac23v_R
 +\tfrac\varepsilon{12}v_R'''=0.
\]
By Lemma~\ref{lem:kpz-quad-coarse}, the $v_R$ are locally bounded, nonnegative and nonincreasing. Helly compactness gives a subsequential limit $v$ in $L^1_{\rm loc}$ and pointwise at continuity points. It vanishes on $z>0$ and obeys $c(-z)\le v(z)\le C(-z)$ for $z<0$. Distributionally the $\varepsilon v_R'''$ term tends to zero, and the nonlinear term passes by $L^1_{\rm loc}$ convergence of the bounded monotone functions. Hence
\[
-\tfrac z3v'-\tfrac14(v^2)'-\tfrac23v=0.
\]
The integration constant is zero because $v=0$ on $z>0$. Consequently,
\begin{equation}\label{eq:kpz-quad-identity}
4zv(z)+3v(z)^2=4\int_z^\infty v(w)\,dw
\end{equation}
for the one-sided traces of the monotone function. For $z<0$, the right side is strictly positive, so $v(z)>-4z/3$. At any downward jump at $z<0$, continuity of the right side of \eqref{eq:kpz-quad-identity} would require the sum of the two positive traces to equal $-4z/3$, impossible since each is strictly greater than $-4z/3$. At $z=0$ the identity gives the left trace $v(0-)=0$. Thus $v$ is continuous.

On $z<0$ the implicit relation \eqref{eq:kpz-quad-identity} is differentiable (the coefficient $4z+6v$ is nonzero by $v>-4z/3$), and the equation becomes $(2z+3v)v'=-4v$. Inverting gives
$d z/dv+z/(2v)=-3/4$, hence $z=-v/2+C/\sqrt v$. The boundary $v(z)\downarrow0$ as $z\uparrow0$ forces $C=0$. Thus $v(z)=-2z$ for $z<0$, and every subsequence has the same limit. Monotonicity and continuity of the limit improve this to local uniform convergence.

Finally use \eqref{eq:kpz-quad-repr} with $s=-S$, divide by $S^3$, and rescale $u=Sz$. On $z\in[-1,0]$, dominated convergence and \eqref{eq:kpz-quad-flat} apply; the integral over $u\ge0$ is $O(S)$ by Lemma~\ref{lem:kpz-quad-coarse}. It follows that
$S^{-3}H_1(-S)\to\frac12\int_{-1}^0(1+z)2(-z)\,dz=1/6$.
\end{proof}

\begin{lemma}[Positive association and the entrance bound]\label{lem:kpz-quad-entrance}
For $0<\delta<1/4$, uniformly over $1\le b\le1+\delta$,
\begin{equation}\label{eq:kpz-quad-entrance-prob}
G_1(-S)\le G_b(-S)\le
\frac{G_1(-(1-4\delta)S)}{1-Ce^{-cS^{3/2}}},
\qquad S\ge S_\delta,
\end{equation}
with universal constants $c,C>0$ after increasing $S_\delta$ if necessary. For any smooth positive weight $\psi$ with $\psi(z)\asymp(1+|z|)^{-p}$ for $z\ll0$, $p>3$, and $\psi=1$ for $z\ge-1$,
\begin{equation}\label{eq:kpz-quad-entrance-E}
\lim_{\delta\downarrow0}\limsup_{R\to\infty}
\sup_{1\le b\le1+\delta}
\int_\R\psi(z)\left|v_R(b,z)-(1+b^{-1/2})(-z)_+\right|^2\,dz=0.
\end{equation}
\end{lemma}
\begin{proof}
For fixed times $y_1,\ldots,y_m$, the Airy$_2$ vector, up to
deterministic shifts and positive scaling, is a distributional limit
of geometric last-passage times to several endpoints in a common
i.i.d. environment \cite{Johansson03}. Each last-passage time is
coordinatewise nondecreasing in the independent weights. The
Harris--FKG inequality therefore gives positive association of the
finite-dimensional vectors, for decreasing as well as increasing
events. The property passes to the Airy$_2$ finite-dimensional
limits. Refining finite grids and using continuity of the sample
paths proves positive association for events requiring the process
to remain below a specified function on a closed set, including
unbounded closed sets.

Fix $L=2\sqrt S$, $q=1-4\delta$, and let $E_{\rm in}$ and $E_{\rm out}$ be the events that $\bar{\mathcal A}_1(y)\le-qS+y^2$ for $|y|\le L$ and $|y|\ge L$, respectively. For $b\in[1,1+\delta]$ and $|y|\le L$,
$-S+by^2\le-qS+y^2$. Hence $G_b(-S)\le\mathbb P(E_{\rm in})$. Positive association yields
$G_1(-qS)=\mathbb P(E_{\rm in}\cap E_{\rm out})\ge\mathbb P(E_{\rm in})\mathbb P(E_{\rm out})$.

The exterior event has probability $\mathbb P(E_{\rm out})\ge1-Ce^{-cS^{3/2}}$. Indeed split $\{|y|\ge L\}$ into unit intervals. By stationarity and the GOE identity, for $u>2$,
\[
\mathbb P\{\sup_{y\in[j,j+1]}\bar{\mathcal A}_1(y)>u\}
=\mathbb P\{\sup_{v\in[0,1]}\bar{\mathcal A}_1(v)>u\}
\le1-G_1(u-1)\le Ce^{-cu^{3/2}},
\]
where the last inequality follows from the soft-edge right-tail bound in \cite[Theorem~1]{LedouxRider2010}. On an exterior interval at distance $j\ge L-1$, the function $y^2-qS$ is at least $j^2-qS$, which is of order $S+j^2$. Summing the stretched-exponential bounds proves the claim, and hence \eqref{eq:kpz-quad-entrance-prob}.

By Lemma~\ref{lem:kpz-quad-flat}, \eqref{eq:kpz-quad-entrance-prob} implies, uniformly for $1\le b\le1+\delta$,
\[
\frac{q^3}{6}+o(1)
\le \frac{H_b(-S)}{S^3}
\le\frac16+o(1)\qquad(S\to\infty).
\]
Since $H_b''=g_b/2$ is nonincreasing, the centered second-difference bounds for a function with decreasing second derivative show, for fixed $x>0$ and $\eta_0\in(0,1/3)$,
\[
\limsup_{R\to\infty}\sup_{1\le b\le1+\delta}
\left|v_R(b,-x)-2x\right|
\le Cx\left(\eta_0+\frac{1-q^3}{\eta_0^2}\right).
\]
Take $\eta_0=(1-q^3)^{1/3}$ when $\delta$ is small. For $x\le0$, where $\psi=1$, the corresponding part of the integral in \eqref{eq:kpz-quad-entrance-E} is $\int_0^\infty v_R(b,z)^2\,dz\le CR^{-3}$, by Lemma~\ref{lem:kpz-quad-coarse}. For $x>0$, the uniform linear bounds from Lemma~\ref{lem:kpz-quad-coarse} and $p>3$ permit dominated convergence in the weighted $L^2$ norm, and $1+b^{-1/2}=2+O(\delta)$. This proves \eqref{eq:kpz-quad-entrance-E}.
\end{proof}

\begin{lemma}[Weighted stability]\label{lem:kpz-quad-stability}
For every $B>1$ and $p>3$, choose a smooth nondecreasing $\psi>0$ such that $\psi=1$ on $[-1,\infty)$ and $\psi(z)\asymp(1+|z|)^{-p}$ for $z\le-2$, with $|z|\psi'\le C\psi$ and $|\psi'''|\le C\psi$. Set
\[
v_*(b,z)=K(b)(-z)_+ ,\qquad K(b)=1+b^{-1/2}.
\]
Then
\begin{equation}\label{eq:kpz-quad-L2selection}
\lim_{R\to\infty}\sup_{1\le b\le B}
\int_\R\psi(z)|v_R(b,z)-v_*(b,z)|^2\,dz=0.
\end{equation}
\end{lemma}
\begin{proof}
Set $\varepsilon=R^{-3}$ and $A=b-1$, $\kappa(b)=(b-1)/(2b)-2/3$. The exact rescaled equation is
\begin{equation}\label{eq:kpz-quad-v}
A\partial_bv_R-\frac z3\partial_zv_R
-\frac12v_R\partial_zv_R+\kappa v_R
+\frac\varepsilon{12}\partial_z^3v_R=0.
\end{equation}
On $z<0$, the inviscid profile $v_*=-Kz$ satisfies the equation because
\begin{equation}\label{eq:kpz-quad-K}
(b-1)K'=\left(1-\frac{b-1}{2b}\right)K-\frac12K^2.
\end{equation}
It is also an inviscid weak solution across $z=0$.

\emph{Smooth comparison profile.} Let $\rho\ge0$ be a symmetric $C_c^\infty([-1,1])$ mollifier of integral one. Write $\eta=R^{-1}=\varepsilon^{1/3}$ and $W_\eta(b,z)=K(b)(\rho_\eta*(-\cdot)_+)(z)$. For $z\le-\eta$ this equals $-Kz$, and for $z\ge\eta$ it is zero. In $|z|\le\eta$ its size is $O(\eta)$ and $|(W_\eta)_z|\le C$. Its residual
\[
\mathcal R_\eta=A(W_\eta)_b-\frac z3(W_\eta)_z-
\frac12W_\eta(W_\eta)_z+\kappa W_\eta+
\frac\varepsilon{12}(W_\eta)_{zzz}
\]
is supported in $[-\eta,\eta]$. The inviscid part has size $O_B(\eta)$ there, while
$\|(W_\eta)_{zzz}\|_2\le C_B\eta^{-3/2}$. Hence
\begin{equation}\label{eq:kpz-quad-residual}
\|\mathcal R_\eta\|_{L^2(\R)}\le C_B(\eta^{3/2}+\varepsilon\eta^{-3/2})
\le C_B R^{-3/2}.
\end{equation}
Lemma~\ref{lem:kpz-quad-coarse} also gives $|v_R(b,z)|\le C_B\eta$ when $|z|\le\eta$. Therefore, with $e=v_R-W_\eta$,
\begin{equation}\label{eq:kpz-quad-corner}
\int_{-\eta}^\eta e^2\,dz\le C_B\eta^3=C_B\varepsilon,
\qquad
\left|\int\psi e\mathcal R_\eta\,dz\right|\le C_B\varepsilon.
\end{equation}
This is the estimate that makes smoothing the corner harmless.

\emph{Energy identity.} The following identity is exact for any smooth positive weight $w$ independent of $b$, provided the boundary terms vanish:
\begin{align}\label{eq:kpz-quad-energy}
\frac A2\frac d{db}\int w e^2
+\frac\varepsilon8\int w'e_z^2
={}&\int w\left(-\frac16-\kappa+\frac14(W_\eta)_z\right)e^2
\\&-\frac1{12}\int w'(2z+W_\eta+2v_R)e^2
+\frac\varepsilon{24}\int w'''e^2
-\int w e\mathcal R_\eta.\notag
\end{align}
To derive the identity, subtract the equation for $W_\eta$ from \eqref{eq:kpz-quad-v} and note that $v_R(v_R)_z-W_\eta(W_\eta)_z=W_\eta e_z+(W_\eta)_ze+ee_z$. Multiply by $we$ and integrate. The transport and cubic terms are integrated once, and the dispersive term uses the exact integration-by-parts formula
\[
-\frac\varepsilon{12}\int we e_{zzz}= -\frac\varepsilon8\int w'e_z^2+\frac\varepsilon{24}\int w'''e^2.
\]
To justify it without assumptions on the left-tail derivatives, choose $w=\psi_M$, where $\psi_M=0$ on $(-\infty,-2M]$, $\psi_M=\psi$ on $[-M,\infty)$, and $\psi_M$ is smooth and nondecreasing. One may take $\psi_M(z)=\psi(z)\theta(2+z/M)$ for a smooth nondecreasing cutoff $\theta$ that is zero at arguments $\le0$ and one at arguments $\ge1$. For each compact interval $b\in[b_0,B]\subset(1,\infty)$, take the fixed initial profile $h(y)=-y^2$ and time $t=b-1$. By Lemma~\ref{lem:kpz-smoothing}, its first three height derivatives decay exponentially at $+\infty$, uniformly in $b$, so integration by parts at $+\infty$ is justified; all terms at $-\infty$ vanish because $\psi_M$ does. The bounds below are uniform in $M$, and then $M\to\infty$.

\emph{Signs.} On $z\le-\eta$, $(W_\eta)_z=-K$, and
\[
-\frac16-\kappa+\frac14(W_\eta)_z
=\frac1{2b}-\frac{1+b^{-1/2}}4\le0\qquad(b\ge1).
\]
On $z\ge\eta$, $W_\eta=0$, and the potentially positive first term of \eqref{eq:kpz-quad-energy} is bounded by $C_B\int_0^\infty v_R^2\le C_B\varepsilon$. On the corner, \eqref{eq:kpz-quad-corner} yields the same bound. Thus the first term on the right of \eqref{eq:kpz-quad-energy} is at most $C_B\varepsilon$. The residual term is bounded by \eqref{eq:kpz-quad-corner}, and $|\psi'''|\le C\psi$ bounds the third term by $C\varepsilon E_R$, where $E_R(b)=\int\psi e^2$.

By Lemma~\ref{lem:kpz-quad-coarse}, $v_R(b,z)\ge c(-z)$ on $z\le-1$, for $R$ large. Choose $\delta_0>0$ so small that $2-K(b)\le2c$ for $1\le b\le1+\delta_0$. Then, on the support of $\psi_M'\subset(-\infty,-1]$, where $W_\eta=v_*$,
\[
2z+W_\eta+2v_R\ge[-(2-K)+2c](-z)\ge0.
\]
Because $\psi_M'\ge0$, the nonlinear flux term of \eqref{eq:kpz-quad-energy} is nonpositive on this interval. The uniform bound $E_R\le C_B$, which follows from Lemma~\ref{lem:kpz-quad-coarse} and $p>3$, then gives, for any fixed $1<b_0\le b\le1+\delta_0$,
\begin{equation}\label{eq:kpz-quad-near}
E_R(b)\le E_R(b_0)+C_B\varepsilon\log\frac{b-1}{b_0-1}.
\end{equation}
Here one first proves the same inequality with the cutoff $\psi_M$ and then lets $M\to\infty$; the extra $\psi_M'''$ contributions in the transition region are $O_B(\varepsilon M^{-p})$, and the energy converges by dominated convergence.

On $1+\delta_0\le b\le B$, the sign of the flux term is unnecessary: use $|2z+W_\eta+2v_R|\le C_B(1+|z|)$ for $z\le-1$ and $|z|\psi'\le C\psi$. This bounds the flux term by $C_B E_R$. The cutoff transition $[-2M,-M]$ contributes at most $C_B M^{3-p}$, since $|e|\le C_B(1+|z|)$ and $\psi_M'=O(M^{-p-1})$ there; this vanishes as $M\to\infty$. As $b-1\ge\delta_0$, Gronwall gives
\begin{equation}\label{eq:kpz-quad-far}
\sup_{1+\delta_0\le b\le B}E_R(b)
\le C_B\left(E_R(1+\delta_0)+\varepsilon\right).
\end{equation}

Finally, Lemma~\ref{lem:kpz-quad-entrance} says that
$\lim_{b_0\downarrow1}\limsup_{R\to\infty}E_R(b_0)=0$,
since $\|W_\eta-v_*\|_{L^2(\psi)}^2=O_B(\eta^3)=O_B(\varepsilon)$. First let $R\to\infty$ in \eqref{eq:kpz-quad-near}, then $b_0\downarrow1$; this gives $E_R(1+\delta_0)\to0$. Equation \eqref{eq:kpz-quad-far} gives convergence on $[1+\delta_0,B]$. The entrance estimate and \eqref{eq:kpz-quad-near}, applied with arbitrarily small $b_0-1$, also give uniform convergence on $[1,1+\delta_0]$. We conclude \eqref{eq:kpz-quad-L2selection}.
\end{proof}

\begin{proof}[Proof of Theorem~\ref{thm:kpz-quadratic-tail}]
By Lemma~\ref{lem:kpz-quad-stability}, $v_R(b,\cdot)\to K(b)(-\cdot)_+$ in weighted $L^2$, uniformly for $b$ in a compact interval. Because $v_R(b,\cdot)$ is nonincreasing and the limiting profile is uniformly Lipschitz, this implies locally uniform convergence in $(b,z)$.

For $b\ge1$ and fixed $z\le0$, use \eqref{eq:kpz-quad-repr} and rescale $u=Rw$:
\[
R^{-3}H_b(Rz)=\frac12\int_z^\infty(w-z)v_R(b,w)\,dw
\longrightarrow \frac{K(b)}{12}(-z)^3.
\]
The part $w\ge0$ vanishes uniformly by \eqref{eq:kpz-quad-right} and \eqref{eq:kpz-quad-repr}; the part $w\in[z,0]$ converges by the local uniform result. For $z\ge0$, monotonicity gives $0\le H_b(Rz)\le H_b(0)\le H_1(0)$ and the same zero limit. Thus, locally uniformly in $(b,z)$,
\[
-R^{-3}\log G_b(Rz)\longrightarrow\frac{1+b^{-1/2}}{12}(-z)_+^3.
\]
Insert the exact identity \eqref{eq:kpz-quad-exact}, with $x=\sqrt R X$ and $r=Rz$, to obtain \eqref{eq:kpz-lower-tail-quadratic}. Finally, the curvature inequality makes $\partial_r^2\log F_a$ monotone in $r$. The standard second-difference squeeze, applied to the locally uniform cubic asymptotic, gives the locally uniform limit \eqref{eq:kpz-lower-tail-phi}.
\end{proof}

\subsection{Approximate solutions}\label{app:kpz-transport}

\begin{proof}[Proof of Theorem~\ref{thm:kpz-stability}]
Fix $0<t<T$ and $x,r\in\R$, and let $0<s\le s_0$, with $s_0<\min\{t_*,T-t\}$. The transport identity \eqref{eq:kpz-transport} on $(0,T)$, with $\tau=s+t$ and $s_2=\tau$, gives
\begin{equation}\label{eq:kpz-defect}
 \exp\{-\tfrac12\partial_r^{-2}\phi(s+t,x,r)\}-Z_{\phi(s)}(t,x,r)=-\int_s^{s+t}\E\Bigl[e^{-I_{\phi(\sigma)}}I_{R_\phi(\sigma)}\Bigr]\diff\sigma,
\end{equation}
with the occupations at $(s+t-\sigma,x,r)$. Since $\phi\ge0$, the integrand is at most $\nu(\sigma)\E I_{e^{-\delta u}}^{s+t-\sigma,x}(r)$ in absolute value, and \eqref{eq:occ-density} with \eqref{app:mu} gives
\begin{equation}\label{eq:kpz-exp-occupation}
 \E I_{e^{-\delta u}}^{t',x}(r)=\frac{e^{\delta^3t'/12-\delta r}}{2\delta^2},\qquad t'>0.
\end{equation}
So the right side of \eqref{eq:kpz-defect} is at most $\frac1{2\delta^2}e^{\delta^3(s+t)/12-\delta r}\int_s^{s+t}\nu$ in absolute value. Together with the bound on $\phi$ on $[t,t+s_0]$, this gives \eqref{eq:kpz-uniform-tail}, as in the proof of uniqueness, so Lemma~\ref{lem:kpz-wall} applies, and $Z_{\phi(s)}(t,x,r)\to F_h(t,x,r)$. Letting $s\downarrow0$ in \eqref{eq:kpz-defect} gives \eqref{eq:kpz-stability}.
\end{proof}

\subsection{The thinned laws}\label{app:kpz-thinning}

\begin{proof}[Proof of Proposition~\ref{prop:kpz-thinned}]
Write $\log F^{(p)}=g_p(\sigma)$. By \cite[(1.15)]{TWAiry},
\[
 g_p(\sigma)=-\int_\sigma^\infty(z-\sigma)q_p(z)^2\diff z,\qquad q_p''(z)=zq_p(z)+2q_p(z)^3,\qquad q_p(z)\sim\sqrt p\Ai(z)\ \ (z\to+\infty).
\]
Since $\det(I-pK_{\rm Ai})_{L^2(\sigma,\infty)}$ is positive and real analytic in $\sigma$, $q_p^2=-g_p''$ is real analytic on $\R$. By \cite[(1.18) and Remark~1.6]{Bothner}, as $w\to+\infty$,
\[
 q_p(-w)=d_pw^{-1/4}\cos\varphi_p(w)+O(w^{-7/10}),\qquad d_p^2=-\frac{\log(1-p)}\pi,
\]
where $\varphi_p(w)=\frac23w^{3/2}-\frac34d_p^2\log w+c_p$ for a constant $c_p$, and therefore
\begin{equation}\label{eq:kpz-AS}
 q_p(-w)^2=d_p^2w^{-1/2}\cos^2\varphi_p(w)+O(w^{-19/20}).
\end{equation}
In particular $q_p$ is bounded, and it changes sign infinitely often on $(-\infty,0)$.

\textit{The equation.} With $\ell=\log F^{(p)}$, KP in the form \eqref{eq:hirota-physical}, divided by $24(F^{(p)})^2$, reads $\ell_{tr}+\frac14\ell_{xx}+\frac1{12}\ell_{rrrr}+\frac12\ell_{rr}^2=0$. Substituting $\ell=g_p(\sigma)$, its left side is
\[
 \frac{t^{-4/3}}{12}\bigl(g_p''''+6(g_p'')^2-4\sigma g_p''+2g_p'\bigr)(\sigma).
\]
Here $g_p''=-q_p^2$, so the Painlev\'e equation gives $g_p''''=-2(q_p')^2-2\sigma q_p^2-4q_p^4$, and $g_p'=(q_p')^2-\sigma q_p^2-q_p^4$, since both sides tend to $0$ at $+\infty$ and have derivative $-q_p^2$. So the bracket vanishes. Differentiating in $r$ and using $\ell_r=t^{-1/3}g_p'(\sigma)=-\frac12\partial_r^{-1}\phi^{(p)}$ gives \eqref{KP} for
\[
 \phi^{(p)}=-2\ell_{rr}=2t^{-2/3}q_p(\sigma)^2.
\]

\textit{Monotonicity, $(\aleph 1)$ and $(\aleph 2)$.} The zeros of $q_p$ are isolated, and $q_p$ has some, so $q_p^2$ vanishes at isolated points and is positive between them; since $\sigma$ is an increasing affine function of $r$, $\phi^{(p)}(t,x,\cdot)$ is not monotone. Since $q_p$ is bounded and decays like $\Ai$ at $+\infty$, $\phi^{(p)}$ is bounded on compact time intervals and decays faster than every exponential as $r\to\infty$, so $(\aleph 1)$ holds. A compact subset of $\mathrm{epi}_S(h)$ lies in $\{|x-y_0|\ge\rho\}\cup\{r\ge h_0+\rho\}$ for some $\rho>0$, and on it $\sigma\ge t^{-1/3}(r-h_0)+t^{-4/3}\rho^2$ or $\sigma\ge t^{-1/3}\rho$. So $\sigma\to+\infty$ uniformly there, and $F^{(p)}\to1$ and $\phi^{(p)}\to0$ uniformly. At $x=y_0$ and $r<h_0$, $\sigma\to-\infty$, and
\[
 -g_p(\sigma)\ge\int_{\sigma/2}^0(z-\sigma)q_p(z)^2\diff z\ge\frac{|\sigma|}2\int_{\sigma/2}^0q_p(z)^2\diff z\longrightarrow\infty
\]
by \eqref{eq:kpz-AS}; so $F^{(p)}(t,y_0,r)\to0$.

\textit{The integral below the tip.} Let $b<h_0$ and $v=t^{-1/3}(h_0-b)$. On each side of $y_0$, substitute $z=\sigma=-v+t^{-4/3}(x-y_0)^2$:
\[
 \int_\R\phi^{(p)}(t,x,b)\diff x=2\int_{-v}^\infty\frac{q_p(z)^2}{\sqrt{z+v}}\diff z.
\]
For $|x-y_0|\ge\varepsilon$, $z\ge t^{-4/3}\varepsilon^2-v\to+\infty$, so the part of the integral outside $(y_0-\varepsilon,y_0+\varepsilon)$ tends to $0$. The part with $z\ge0$ is at most $2v^{-1/2}\int_0^\infty q_p^2\to0$. Fix $w_0\ge1$. The part with $-w_0\le z\le0$ is $O(v^{-1/2})$, and by \eqref{eq:kpz-AS} the error term contributes $O(v^{-9/20})$ to the part with $z=-w$, $w_0\le w\le v$. In the main term substitute $w=vu$:
\[
 2d_p^2\int_{w_0}^v\frac{\cos^2\varphi_p(w)}{\sqrt{w(v-w)}}\diff w=2d_p^2\int_{w_0/v}^1\frac{\cos^2\varphi_p(vu)}{\sqrt{u(1-u)}}\diff u.
\]
Write $\cos^2=\frac12+\frac12\cos2\varphi_p$. The constant part tends to $d_p^2\int_0^1\frac{\diff u}{\sqrt{u(1-u)}}=\pi d_p^2=-\log(1-p)$. On $[\kappa,1-\kappa]$ the $u$-derivative of $\varphi_p(vu)$ is at least $cv^{3/2}\sqrt\kappa$ for large $v$, so an integration by parts makes the oscillating part tend to $0$ there, and the end intervals contribute $O(\sqrt\kappa)$. Letting $v\to\infty$ and then $\kappa\downarrow0$ gives the limit $-\log(1-p)$.

\textit{The limit in $\mathcal D'(\R^2)$.} For $b>h_0$ the same substitution, now with $z\ge t^{-1/3}(b-h_0)$, shows that $\int_\R\phi^{(p)}(t,x,b)\diff x\to0$. By \eqref{eq:kpz-AS}, $q_p(-w)^2\le C\min\{1,w^{-1/2}\}$, so the integral $2\int_{-v}^\infty q_p(z)^2(z+v)^{-1/2}\diff z$ is bounded uniformly in $v\in\R$. For $\eta\in C_c^\infty(\R^2)$ and $b\ne h_0$, $\int\phi^{(p)}(t,x,b)\eta(x,b)\diff x$ therefore tends to $-\log(1-p)\eta(y_0,b)\mathbf 1_{\{b<h_0\}}$, boundedly in $b$, and integrating in $b$ gives $\phi^{(p)}(t)\to V_p$ in $\mathcal D'(\R^2)$.
\end{proof}

\section*{AI disclosure and author notes}
\addcontentsline{toc}{section}{AI disclosure and author notes}
Corollary~\ref{cor:kpz-restart} was discovered by ChatGPT when the authors asked it to prove Theorem~\ref{thm:kpz-curvature}. The authors recognized that this observation could lead to a very general representation for KP.

OpenAI's ChatGPT and Anthropic's Claude were used extensively in preparing the manuscript. The assistance included proposing and developing proof arguments, expanding and simplifying proofs, checking calculations and identifying possible gaps, locating references, organizing the manuscript, and drafting and revising the exposition and \LaTeX{} source. Python scripts and \LaTeX{} tools were used to assemble and typeset the document.

Responsibility for the correctness and adequacy of the mathematical arguments, and for the accuracy and completeness of the references and attribution, remains entirely with the authors.

A forthcoming article by the authors will address applications of the representation to the zero dispersion limit of the KP-II equation.
The authors welcome comments about this first version of the manuscript.

\phantomsection

\bigskip
{\footnotesize
(J.~Quastel) \textsc{Department of Mathematics, University of Toronto, 40 St.~George Street, Toronto, Ontario M5S 2E4, Canada}\par\nopagebreak
\textit{Email address}: \href{mailto:jeremy.quastel@utoronto.ca}{\texttt{jeremy.quastel@utoronto.ca}}\par
\bigskip
(C.~A.~Rodriguez) \textsc{Department of Mathematics, University of Michigan, 2074 East Hall, 530 Church Street, Ann Arbor, MI 48109-1043, USA}\par\nopagebreak
\textit{Email address}: \href{mailto:rodrigca@umich.edu}{\texttt{rodrigca@umich.edu}}\par}

\end{document}